%% file: arxiv.tex
\pdfoutput=1
\documentclass[10pt]{aart}

\usepackage{arxiv-setting}
\usepackage{booktabs}
\usepackage{placeins}
\usepackage{needspace}
\usepackage{tabularx}

\cslet{blx@noerroretextools}\empty
\usepackage{autonum}

\renewcommand{\epsilon}{\varepsilon}

\def\stj{\mathfrak{m}} 
\def\res{{\mG}}
\def\eres{{\emG}}
\def\resmax{{\mathrm{K}}}
\def\reserr{{\mathfrak{E}}} 

\newcommand{\bulk}{\bm{\mathfrak{B}}_{\kappa}} 
\newcommand{\specdomain}{\bm{\mathfrak{S}}_{\kappa}} 
\def\scA{{\rmH}}
\def\hatH{{\widehat{\scA}^{(\vtau)}}}

\def\ratio{{\gamma}} 
\def\qr{{\vartheta}} 
\def\signal{{\mathfrak{f}}} 
\def\ledge{\lambda_{\mathrm{L}}} 
\def\redge{\lambda_{\mathrm{R}}} 
\newcommand{\gapR}{\gap_{\mathrm R}}
\newcommand{\gapL}{\gap_{\mathrm L}}
\newcommand{\invR}[1]{\Lambda_{#1}^{-1,\mathrm R}}
\newcommand{\invL}{\Lambda_{\qr}^{-1,\mathrm L}}
\newcommand{\rateR}{\mathcal I_{\mathrm R}}
\newcommand{\rateL}{\mathcal I_{\mathrm L}}
\newcommand{\rhoR}{\rho_{\mathrm R}}
\newcommand{\rhoL}{\rho_{\mathrm L}}
\newcommand{\bcurve}[1]{\mathcal B_{\mathrm{#1}}}
\newcommand{\bcurveL}{\bcurve{L}}
\newcommand{\bcurveR}{\bcurve{R}}
\newcommand{\bcurveRi}[1]{\mathcal B_{\mathrm R,#1}}
\newcommand{\bcrit}[1]{b_{\mathrm{#1}}^\star}
\newcommand{\bcritL}{\bcrit{L}}
\newcommand{\bcritR}{\bcrit{R}}
\newcommand{\bcritRi}[1]{b_{\mathrm R,#1}^\star}
\newcommand{\bcritF}[1]{b_{#1}^\star}
\newcommand{\bcritIso}{\bcrit{iso}}
\newcommand{\bcritFork}{\bcrit{fork}}
\def\benrate{\mathrm{h}} 

\def\error{\mathfrak{\xi}}
\def\radius{\mathrm{r}_{\star}}

\def\fka{\mathfrak{a}}
\def\fkb{\mathfrak{b}}

\def\Vtyp{\mathcal{T}} 
\def\Vonehigh#1{{\gV_{1}^{(\geq #1)}}}
\def\Vtwolow#1{{\gV_{2}^{(\leq #1)}}}
\def\Vtwohigh#1{{\gV_{2}^{(\geq #1)}}}

\def\hatvtau{{\widehat{\rvv}^{(\vtau)}}}

\DeclareMathAlphabet{\mathcalligra}{T1}{calligra}{m}{n}
\DeclareMathAlphabet{\mathpzc}{OT1}{pzc}{m}{it}

\makeatletter
\define@key{todonotes}{ID}[]{%
    \setkeys{todonotes}{author=ID, inline, color=magenta}}%
\makeatother
\makeatletter
\define@key{todonotes}{HW}[]{%
    \setkeys{todonotes}{author=HW, inline, color=orange}}%
\makeatother

\newcommand{\toromanlower}[1]{\ifcase#1 \or i\or ii\or iii\or iv\or v\or vi\or vii\or viii\or ix\or x\fi}

\title{Singular vectors of sparse random rectangular matrices:\\
Transition between delocalized and semilocalized phases}

\author{
Ioana Dumitriu%
\thanks{\mbox{\scriptsize Department of Mathematics, University of California San Diego, La Jolla, CA 92093, USA; \texttt{idumitriu@ucsd.edu}.}}
\and
Hai-Xiao Wang\orcidlink{0000-0003-2730-1439}%
\thanks{\mbox{\scriptsize Department of Applied Mathematics, University of Washington, Seattle, WA 98195, USA; \texttt{haixwang@uw.edu}.}}
}

\date{\today}

\begin{document}

\maketitle

\begin{abstract}
  Consider the bipartite Erd\H{o}s--R\'enyi graph $\mathds{G}(n,m,d/\sqrt{mn})$ at the critical sparsity $d = b\log(N)$ for some $b > 0$ with $N = n + m$, and constant aspect ratio $\gamma \coloneqq n/m \geq 1$. We show that the spectrum of its centered and rescaled bi-adjacency matrix, whose limiting law is the Mar\v{c}enko--Pastur (MP) distribution, is split into two phases: in the bulk, singular vectors are completely delocalized, while near the edges, singular vectors are semilocalized, with their mass concentrating on small disjoint balls, centered at vertices with extreme degrees. As $b$ decreases, we identify sharp thresholds for the emergence of outlier singular values. For explicit $b_{\rm R}^{\star}$ and $b_{\rm L}^{\star}$ as functions of $\gamma$, we prove that when $b > b_{\rm R}^{\star}$, there is no outlier outside the bulk; when $b_{\rm L}^{\star}< b < b_{\rm R}^{\star}$, only right outliers emerge; and when $b < b_{\rm L}^{\star}$, outliers are present on both sides, all with high probability. This removes an extra technical condition from our previous work \cite{dumitriu2025singular}, providing a complete characterization of the parameter space. In sparser regimes, in terms of left outliers, we provide some nontrivial generating mechanisms, demonstrating that local graph structure plays a key role in understanding extreme singular values. 

  Our proofs rely on two technical inputs: (1) an entrywise local law in the MP bulk, where the deterministic equivalent of the resolvent depends on both variance profiles and normalized vertex degrees, different from the Gaussian setting; (2) spectral radius upper bounds for nonbacktracking matrices constructed from regularized random graphs, recently developed in \cite{wang2026regularized}. Our results reveal fruitful phenomena for sparse random rectangular matrices at criticality, which are beyond the universality class of their Gaussian counterparts.

\end{abstract}


\newpage

\tableofcontents

\newpage
\input{intro}
\subsection*{Acknowledgments}
Part of this research was performed while H.X.W. was visiting the Institute for Pure and Applied Mathematics (IPAM) in Los Angeles, California, in August 2025, supported by NSF award DMS-1925919, and while H.X.W. was attending the workshop ``Recent Advances in Random Networks'' at the Institute for Mathematical and Statistical Innovation (IMSI) in Chicago, Illinois, in January 2026, supported by NSF award DMS-2425650. Part of this research was completed while the authors were attending FoCM 2026 in Vienna, Austria. The authors acknowledge the support of NSF award DMS-2609677. H.X.W. thanks Mark Rudelson for pointing out the reference \cite{basak2021sharp} and for helpful discussions on the invertibility of sparse random matrices during H.X.W.'s May 2026 visit to the National Institute for Theory and Mathematics in Biology (NITMB) in Chicago, Illinois, for the workshop ``Recent Developments in Discrete Probability''. Both authors would like to thank Petar Nizic-Nikolac for helpful discussions on \cite{bandeira2026matrix} during their visit to Simons Institute for the Theory of Computing in October 2026, during the workshop ``Probabilistic Phenomena on Graphs, Groups, and Manifolds''.

\input{delocalized}
\input{semilocalized}
\input{transition}

\printbibliography

\newpage
\appendix
\phantomsection
\addcontentsline{toc}{section}{Appendices}

\input{appendix}
\input{app-delocalized}
\input{app-semilocalized}
\input{app-transition}
\end{document}

%% file: intro.tex

\section{Introduction}\label{sec:intro}

Sparse random matrices arise broadly in many fields, such as numerical linear algebra \cite{martinsson2020randomized}, community detection \cite{abbe2018community}, and modern attention mechanisms \cite{liu2025deepseek}. In the past several decades, considerable effort has been devoted to understanding the spectral properties of sparse random matrices, with particular focus on determining the precise locations of extreme eigenvalues.
One of the most fundamental and extensively studied models is the Erd\H{o}s--R\'enyi random graph $\gG\sim\mathds{G}(N, d/N)$, in which each edge is included independently with probability $d/N$. This simple but mathematically rich model has attracted significant attention since its introduction in the 1950s \cite{renyi1959random}.
It has been proved that when $d/\log(N)\to\infty$, the extreme eigenvalues of the centered and rescaled adjacency matrix approach the semicircle edges \cite{furedi1981eigenvalues, vu2007spectral, benaych2020spectral}; in particular, $\|\mA-\E\mA\| = (2 + o(1))\sqrt{d}$ where $\mA$ is the adjacency matrix of $\gG$. In addition, all the eigenvectors of $\mA$ are delocalized \cite{he2019local}. When $d/\log(N)\to 0$, $\|\mA-\E\mA\|$ fails to concentrate at the $\sqrt{d}$ scale \cite{feige2005spectral, lei2015consistency, le2017concentration} due to the fact that $\|\mA\|$ is determined by the maximum degree of $\gG$ \cite{krivelevich2003largest, benaych2019largest}, which is much larger than $\|\E\mA\| = d$. As a consequence, the leading eigenvector of $\mA$ is localized \cite{alt2023localized}. In the crossover regime $d = b \log(N)$ for some constant $b > 0$, Alt, Ducatez, and Knowles \cite{alt2021extremal}, as well as Tikhomirov and Youssef \cite{tikhomirov2021outliers}, identified the critical threshold $b_{*} = 1/\ln(4/e)$: spectral outliers induced by high-degree vertices emerge if and only if $b < b_{*}$.

Another topic of interest concerns the geometric structure of the eigenvectors of $\mA$, with particular emphasis on their \emph{spatial localization}. For an $\ell^{2}$-normalized vector $\vw = [\evw_{x}]_{x \in [N]} \in \R^N$, we say that it is
\begin{enumerate}[label=(\roman*)]
    \item \emph{delocalized} if $\|\vw\|_{\infty}^{2} = N^{-1+o(1)}$,
    \item \emph{semilocalized} if $\|\vw\|_{\infty}^{2} \geq N^{-\rho}$ for some constant $\rho \in (0,1)$,
    \item \emph{localized} if $\|\vw\|_{\infty}^{2} \asymp 1$.
\end{enumerate}
Note that $\vw = [\evw_{x}]_{x \in [N]} \in \R^N$ naturally induces a probability measure $x \mapsto \evw_{x}^{2}$ on the set of vertices $[N]$. Informally speaking, $\vw$ is delocalized if its mass is approximately uniformly distributed over all vertices, whereas it is localized when most of its mass is confined to a small number of vertices. More recently, Alt, Ducatez, and Knowles \cite{alt2021delocalization, alt2022completely, alt2023poisson, alt2023localized} identified the \emph{mobility edge} for localization-delocalization transition: edge eigenvectors are localized while bulk eigenvectors remain delocalized.

In this paper, we focus on the bipartite Erd\H{o}s--R\'enyi graph $\gG\sim\mathds{G}(n, m, d/\sqrt{mn})$, where each edge, containing exactly one vertex in $\gV_{1} = [n]$ and another one in $\gV_{2} = [m]$, is included independently with probability $d/\sqrt{mn}$.  We study the spectral properties and singular vector behaviors of its \emph{bi-adjacency} matrix $\widetilde{\rmA}$ at criticality $d \asymp \log(N)$. Without loss of generality, we assume $n \geq m$ and denote the \emph{aspect ratio} by $\ratio \coloneqq n/m \geq 1$, with $\qr\coloneqq\ratio^{1/4}$. Previously, Dumitriu and Zhu \cite{dumitriu2024extreme} have shown that when $d \gg \log(N)$, the largest and smallest singular values of $(\widetilde{\rmA} - \E\widetilde{\rmA})/\sqrt{d}$ converge to the right and left edges of the \emph{Mar\v{c}enko--Pastur (MP)} distribution, respectively. At the critical sparsity $d = b \log(N)$ for a constant $b \equiv b(\ratio)$, Dumitriu, Wang, Wang and Zhu \cite{dumitriu2025singular} established the critical threshold for the emergence of outlier singular values: high-degree vertices in either class produce right outliers, while low-degree vertices in the smaller class produce left outliers when above the connectivity threshold, and an additional left-edge lower bound. Therefore, two main problems remain unaddressed.
\begin{enumerate}
    \item Detailed behaviors of the singular vectors of $\widetilde{\rmA}$ at the critical sparsity $d \asymp \log(N)$.
    \item The mechanism for generating the left outliers, especially the regime below the connectivity threshold and the additional left-edge restriction.
\end{enumerate}

The difficulties arise from two aspects. First, it is now a standard procedure to establish the local law before proving delocalization. When controlling the resolvent of the Hermitian block of $\widetilde{\rmA}$, its deterministic equivalent, a diagonally dominant matrix, is degree-dependent at criticality due to the emergence of vertices with extreme degrees. Therefore, the local law established for classical random matrices, for example \cite{knowles2013isotropic,knowles2017anisotropic} cannot be applied directly. 
Second, the extra lower bound on sparsity imposed in \cite{dumitriu2025singular} effectively prevents the presence of very high-degree vertices in the larger vertex class. However, as we now know, such vertices are not responsible for producing left outliers, thus this restriction is unnatural and redundant. 
At a technical level, this additional assumption remains necessary in \cite{dumitriu2025singular} since it ensures the positive semidefiniteness of certain quadratic forms, which is crucial when establishing the lower bounds on the least singular values. 
An alternative approach is to first prune the high-degree vertices that are not going to produce left outliers, and then proceed with the proof strategy established in \cite{dumitriu2025singular}, utilizing the Ihara--Bass formula applied to the associated nonbacktracking matrices.
However, the challenge with this alternative approach is that, after pruning, the entries of the resulting active submatrix become dependent. 
As far as we know, in the existing literature, such dependencies have generally not been handled in a detailed manner.

In this paper, we address those challenges by combining the proof strategy in \cite{alt2021delocalization, dumitriu2025singular}, and the recent development by the second author \cite{wang2026regularized}, which allows us to handle the dependencies introduced by the pruning process. Our contributions can be summarized as follows.

\begin{enumerate}

    \item \emph{Bulk delocalization.}
    We prove that all the bulk singular vectors are delocalized (Theorem~\ref{thm:eigenvector_delocalization}), which remains valid even at the subcritical sparsity regime when $d \gg \sqrt{\log N} \log\log(N)$.
    Along the proof, we establish an entrywise local law (Theorem~\ref{thm:bipartite-sparse-local-law}), where the deterministic equivalent of the resolvent is a diagonally dominant matrix depending on both vertex classes and normalized degrees.

    \item \emph{Edge semilocalization.} For outlier singular vectors at both right and left edges, we establish lower bounds $\|\vw\|_{\infty}^{2} \geq N^{-\rho}$ for some function $\rho(\lambda, b) \in (0,1)$ in Theorems~\ref{thm:right_vector_semilocalization} and~\ref{thm:left_vector_semilocalization}, respectively. Indeed, they are localized on a small number of vertices, with their mass concentrated on a small number of disjoint balls, decaying exponentially with the distance from the center.

    \item \emph{Transition boundaries.}
    We establish sharp transition boundaries for the emergence of outlier singular values at right and left edges in Theorems~\ref{thm:right-edge-transition}, as well as the boundaries for the delocalization-localization transition in Theorem~\ref{thm:left-edge-transition}. When $b$ is below certain threshold, zero columns emerge with high probability. We identify two other mechanisms, small components and star-tuning forks on the giant component, inducing positive singular values of order $d^{-1/2}$ with localized singular vectors. See Figure~\ref{fig:edge-singular-vector-transition-two-aspects} for transition boundaries for the emergence of outlier singular values, and Figure~\ref{fig:combined-aspect-ratio-thresholds} for the delocalization-localization transition.
\end{enumerate}

\subsection{Bipartite Erd\H{o}s--R\'enyi graph}\label{subsec:intro-model}
For the bipartite Erd\H{o}s--R\'enyi graph $\gG=(\gV,\gE)\sim\mathds{G}(n,m,p)$, we denote its vertex classes by $\gV=\gV_1\sqcup\gV_2$ with $|\gV_1|=n$, $|\gV_2|=m$ and $N=n+m$. Let $\widetilde{\rmA}\in\{0,1\}^{n\times m}$ denote its bi-adjacency matrix where $\widetilde{\rmA}_{ij} \sim_{\textnormal{i.i.d.}} \Ber(p)$. Throughout the paper, we assume proportional dimensions.
\begin{assumption}[Proportional regime]
\label{ass:proportional_regime}
There is some constant $\Gamma > 0$ such that $1 \leq  \ratio \leq \Gamma$, and
      \begin{align}
     n &\, = \lfloor \ratio m \rfloor,
     &
     \qr &\, \coloneqq\ratio^{1/4}.
      \label{eqn:proportional-block-regime}
     \end{align}
\end{assumption}
We further define the sparsity parameter by $d = d_N \coloneqq p\sqrt{mn}$, which is the geometric mean of the expected degrees in the two vertex classes:
\begin{align}
    d_1\coloneqq mp = d\sqrt{m/n} = d\qr^{-2}(1+o(1)),
    \qquad
    d_2\coloneqq np = d\sqrt{n/m} = d\qr^2(1+o(1)).
\end{align}
We denote the centered and rescaled bi-adjacency matrix by $\mX\coloneqq (\widetilde{\rmA}-\E\widetilde{\rmA})/\sqrt d$. The adjacency matrix of $\gG$ and its centered and rescaled version are
\begin{align}
    \rmA
    &\coloneqq
    \begin{bmatrix}
        \bzero & \widetilde{\rmA} \\
        \widetilde{\rmA}^{\sT} & \bzero
    \end{bmatrix},
    &
    \scA
    &\coloneqq
    \begin{bmatrix}
        \bzero&\mX\\
        \mX^{\sT}&\bzero
    \end{bmatrix} = \frac{1}{\sqrt{d}}(\mA - \E \mA).
    \label{eqn:scAdj}
\end{align}
Note that for $s>0$ and unit vectors $\rvu\in\R^n$ and $\rvv\in\R^m$, if they form a singular-vector pair of $\mX$, then they also form a singular-vector pair of $\scA$, and vice versa. Formally, we have
\begin{align}
    \mX\rvv=s\rvu, \quad \mX^{\sT}\rvu=s\rvv \quad \iff \quad \scA\frac1{\sqrt2} \begin{bmatrix}\rvu\\ \pm\rvv\end{bmatrix} = \pm s\frac1{\sqrt2} \begin{bmatrix}\rvu\\ \pm\rvv\end{bmatrix}. \label{eqn:singular-eigenvector-correspondence}
\end{align}
Therefore, nonzero eigenvalues of $\scA$ come in pairs $\pm \lambda$, with $\lambda$ being a positive singular value of $\mX$. We denote the $m$ largest singular values of $\mX$, including zeros, by $\lambda_1(\scA) \geq \lambda_2(\scA) \geq \cdots \geq \lambda_m(\scA) \geq 0$. The last $n-m$ eigenvalues of $\scA$ are zero, called \emph{structural zeros}. Furthermore, their kernel spaces are related by
\begin{align}
    \dim(\ker\scA)=\dim(\ker\mX) + \dim(\ker\mX^{\sT}) = n - m + 2\dim(\ker\mX).
    \label{eqn:linearization-zero-multiplicity}
\end{align}
Consequently, each right-kernel direction $\vv\in\ker\mX$ contributes two zero eigenvalues of $\scA$ in addition to the $n-m$ structural zeros.

Consider a rectangular random matrix $\widetilde{\mX}\in\R^{n\times m}$ with i.i.d. entries with zero mean and variance $\sigma^2>0$. For the matrix $\widetilde{\mY} \coloneqq \widetilde{\mX}^{\sT}\widetilde{\mX}/n \in \R^{m \times m}$, let $\mu_{m}(\widetilde{\mY}) \coloneqq m^{-1}\sum_{j=1}^{m} \delta_{\lambda_{j}(\widetilde{\mY})}$ denote the  \emph{Empirical Spectral Distribution} (ESD) of the matrix $\widetilde{\mY}$, which converges weakly to the Mar\v{c}enko--Pastur (MP) law \cite{marchenko1967distribution} below
\begin{align}
    \mu_{\rm{MP}}(dx)
    =\frac{\ratio}{2\pi\sigma^2 x}
      \sqrt{(\lambda_{+}-x)(x-\lambda_{-})}
      \indi{x\in[\lambda_{-},\lambda_{+}]}\,dx,
    \label{eqn:MP_law_original}
\end{align}
where $\lambda_{\pm}=\sigma^2(1\pm\ratio^{-1/2})^2$ are the right and left squared-singular-value edges.

We consider another scaling $\mX \coloneqq \widetilde{\mX}/[\sigma(mn)^{1/4}]$. By applying an elementary change of variables in \eqref{eqn:MP_law_original}, we obtain the limiting density for the positive singular values as
\begin{align}
    f(s)&\coloneqq
    \frac{\qr^2}{\pi s}\sqrt{(s^2-\ledge^2)(\redge^2-s^2)}
    \indi{s\in[\ledge,\redge]},
    &
    \ledge&\coloneqq\qr-\qr^{-1},
    &
    \redge&\coloneqq\qr+\qr^{-1}.
    \label{eqn:MP-singular-value-edges}
\end{align}
Due to the correspondence between the singular values of $\mX$ and the eigenvalues of $\scA$ in \eqref{eqn:singular-eigenvector-correspondence}, the limiting spectral distribution of $\mH$ is given by the symmetrized singular-value law, accompanied by a Dirac mass at zero arising from the $n-m$ surplus rows:
\begin{align}
    \mu_{\rm SMP}(dx)&\coloneqq
    \frac{\ratio-1}{\ratio+1}\delta_0(dx)
    +\frac{1}{\ratio+1}f(|x|)\,dx.
    \label{eqn:symmetric_MP_measure}
\end{align}
Let $\mu_{N}(\scA)\coloneqq N^{-1}\sum_{j=1}^{N} \delta_{\lambda_{j}(\scA)}$ denote the ESD of $\scA$. In the sparse regime $p\to0$ and $d\to\infty$, it was shown in \cite{nadutkina2024marchenko} that $\mu_{N}(\scA)$ converges weakly to $\mu_{\rm SMP}$ as $N\to\infty$. We further characterize the Stieltjes transform of $\mu_{\rm SMP}$ by the following lemma, with the proof deferred to Appendix~\ref{sec:proof-of-lemma-stieltjes-smp-biregular-tree}.
\begin{lemma}\label{lem:stieltjes-smp}
For $z\in\C_+$, the Stieltjes transform
$\stj(z)\coloneqq\int_\R(x-z)^{-1}\mu_{\rm SMP}(dx)$ of the measure in \eqref{eqn:symmetric_MP_measure} is
\begin{align}
    \stj(z)=\frac{\sqrt\ratio}{(1+\ratio)z}
    \left(-z^2+\sqrt{(z^2-\ledge^2)(z^2-\redge^2)}\right).
    \label{eqn:stieltjes-smp-explicit}
\end{align}
Here the square root denotes the analytic branch on $\C_+$ that is asymptotic to $z^2$ as $|z|\to\infty$.
\end{lemma}

\subsection{Bulk delocalization}\label{subsec:intro-bulk}
For any fixed
$0<\kappa<\qr^{-1}$, define the
\emph{$\kappa$-bulk} region by
\begin{align}
    \bulk \coloneqq [-\redge + \kappa, - \ledge - \kappa] \cup [\ledge + \kappa, \redge - \kappa]. \label{eqn:k-bulk-region},
\end{align}
where $\redge$ and $\ledge$ are the right and left MP edges of the support in \eqref{eqn:MP-singular-value-edges}, respectively. For a fixed constant $L \geq 1$, we define the
\emph{bulk-restricted} spectral region by
\begin{align}
    \specdomain \coloneqq \{E+\ii\eta\in\C: E\in \bulk,\ N^{-1+\kappa}\leq \eta\leq L\}\,. \label{eqn:spectral_region}
\end{align}
\begin{assumption}[Subcritical sparsity regime]\label{ass:subcritical_regime}
    For some fixed $B>0$, assume that
    \begin{align}
        \gC_0\sqrt{\log N}\log\log N
        \leq d\leq B\log N,
        \label{eqn:subcritical_regime}
    \end{align}
    For each result with failure exponent $\nu>0$, the constant $\gC_0>0$ is chosen sufficiently large depending on $B$, $\qr$, $\nu$, and the other fixed parameters, before $N\to\infty$.
\end{assumption}
\begin{theorem}[Delocalized phase]\label{thm:eigenvector_delocalization}
    For any fixed $\nu>0$ and $0<\kappa<\qr^{-1}$, there exists some constant $\gC$ such that the following holds. Suppose Assumption~\ref{ass:proportional_regime} and Assumption~\ref{ass:subcritical_regime} with some sufficiently large $\gC_0$ depending on $B$, $\qr$, $\nu$.  For every normalized eigenvector $\rvw$ of $\scA$ with corresponding eigenvalue $\lambda \in \bulk$, there exist some $N_0\in\N$ such that, for every $N\geq N_0$, the following holds with probability at least $1-\gC N^{-\nu}$
    \begin{align}
        \|\rvw\|_{\infty}^{2} \leq N^{-1 + \kappa}.
    \end{align}
\end{theorem}

\subsection{Edge semilocalization}\label{subsec:intro-edge}
For a vertex $x\in\gV$, let $\rD_x\coloneqq\sum_{y\in\gV}\emA_{xy}$ denote the degree of $x$ and $\alpha_x\coloneqq\rD_x/d$ denote the normalized degree. Define two \emph{location maps}:
\begin{align}
    \Lambda_{\qr}(t)
    &\coloneqq
    \sqrt{t+\qr^{-2}+\frac{1}{t-\qr^2}},
    &
    \Lambda_{\qr^{-1}}(t)
    &\coloneqq
    \sqrt{t+\qr^2+\frac{1}{t-\qr^{-2}}}.
    \label{eqn:Lambda-maps}
\end{align}
When $\qr>1$, the map $\Lambda_{\qr}$ is increasing from
$(0,\qr^2-1]$ onto $(0,\ledge]$; it is also increasing from
$[\qr^2+1,\infty)$ onto $[\redge,\infty)$. The map
$\Lambda_{\qr^{-1}}$ is increasing from
$[\qr^{-2}+1,\infty)$ onto $[\redge,\infty)$.
We define the \emph{error scale} by
\begin{align}
    \error_d\coloneqq\frac{\sqrt{\log N}}{d}\log d.
    \label{eqn:error}
\end{align}

\subsubsection{Right edge results}
Define the vertices that induce right outliers by
\begin{subequations}\label{eqn:right-outlier-profile-sets}
\begin{align}
    \gR_1
    &\coloneqq
    \{x\in\gV_1:\alpha_x>\qr^{-2}+1,\
      \Lambda_{\qr^{-1}}(\alpha_x)
      \geq\redge+\tfrac18\error_d^{1/4}\},
      \label{eqn:right-outlier-profile-set-V1}\\
    \gR_2
    &\coloneqq
    \{x\in\gV_2:\alpha_x>\qr^2+1,\
      \Lambda_{\qr}(\alpha_x)
      \geq\redge+\tfrac18\error_d^{1/4}\}.
      \label{eqn:right-outlier-profile-set-V2}
\end{align}
\end{subequations}
Denote $\Lambda^{\sR_{1}}\coloneqq\{\Lambda_{\qr^{-1}}(\alpha_x):x\in\gR_1\}$ and $\Lambda^{\sR_{2}}\coloneqq\{\Lambda_{\qr}(\alpha_x):x\in\gR_2\}$. We arrange the multiset
    $
        \Lambda^{\sR}\coloneqq\Lambda^{\sR_{1}}\cup\Lambda^{\sR_{2}}
    $
    in non-increasing order as $\Lambda^{\sR}_{1}\geq\cdots\geq\Lambda^{\sR}_{|\sR_{1}|+|\sR_{2}|}$.
The location theorem matches these profiles to distinct ordered singular-value positions. Near the right edge, the profile values can interleave with complementary singular values.

\begin{theorem}[Locations of right outliers]\label{thm:right_edge_behavior}
    Under Assumptions~\ref{ass:proportional_regime} and~\ref{ass:subcritical_regime}, with $\gC_0$ sufficiently large depending on $(B,\qr,\nu)$, for every $\nu>0$ there is a constant $\gC_{B,\nu}>0$ such that the following holds with probability at least $1-\gC_{B,\nu}N^{-\nu}$.
\begin{enumerate}[label=(\roman*)]
    \item If $\gR_1\cup\gR_2\neq\emptyset$, there are indices
    $1\leq j_1<\cdots<j_{|\gR_1|+|\gR_2|}\leq m$ such that
    \begin{align}
        |\lambda_{j_k}(\scA)-\Lambda^{\sR}_{k}|
        \leq
        \gC_{B,\nu}
        \big(\error_d+ \error_d^{3/4}/\log(d) \big),
        \qquad 1\leq k\leq|\gR_1|+|\gR_2|.
        \label{eqn:right-outlier-locations}
    \end{align}
    \item If $|\gR_1|+|\gR_2|<m$, we have
    \begin{align}
        \lambda_{|\gR_1|+|\gR_2|+1}(\scA)
        \leq\redge+\error_d^{1/4}.
        \label{eqn:right-no-additional-outliers}
    \end{align}
\end{enumerate}
\end{theorem}
\begin{remark}
    The subcritical extension of
\cite[Theorem~1.9]{dumitriu2025singular} uses quantitative pruning
and the regularized nonbacktracking bound in
Section~\ref{sec:semilocalized_phase}. The component reduction enters
the left-edge argument. The error $\error_d$ dominates when
$d\asymp\sqrt{\log N}\log\log N$, while
$\error_d^{3/4}/\log d$ dominates when $d\asymp\log N$.
\end{remark}

For $0<\lambda\leq\ledge$ or $\lambda\geq\redge$, define
\begin{align}
    \Delta_{\Lambda}(\lambda)
    \coloneqq
    \sqrt{(\lambda^2-\redge^2)(\lambda^2-\ledge^2)}.
    \label{eqn:profile-gap-function}
\end{align}
For $\lambda\geq\redge$ and $s\in\{\qr^{-1},\qr\}$, the inverse branch of the right location map is defined as
\begin{align}
    \invR{s}(\lambda)
    &\coloneqq
    \frac12\left(
      \lambda^2+ s^2 - s^{-2}+\Delta_{\Lambda}(\lambda)
    \right),
    \label{eqn:Lambda-inverse-right-edge}
\end{align}
For $s>0$ and $\alpha\geq0$, define the \emph{transition rate function}
by
\begin{align}
    I_s(\alpha)
    &\coloneqq
    \indi{\alpha>0}
      \bigl(\alpha\log(s^{-2}\alpha)-\alpha+s^2\bigr)
      +\indi{\alpha=0}s^2.
    \label{eqn:transition-rate-function}
\end{align}
The degree rate function for the right edge is defined for $t\geq\redge$,
\begin{align}
    \rateR(t;\qr)
    &\coloneqq
    \min_{s\in\{\qr^{-1},\qr\}}
      I_s\bigl(\invR{s}(t)\bigr),
    \label{eqn:right-degree-rate}
\end{align}
Note that $\rateR(t;\qr)$ is strictly increasing on $[\redge,\infty)$ since the right inverse branches in \eqref{eqn:Lambda-inverse-right-edge} are increasing and $I_s'(\alpha)=\log(\alpha/s^2)>0$ for $\alpha\geq s^2+1$.

At critical sparsity, the right semilocalization exponent has an explicit degree-rate formula.
\begin{assumption}[Critical sparsity]\label{ass:critical_sparsity}
    There is some constant $b>0$ such that $d=b\log N$.
\end{assumption}
\begin{theorem}[Right singular vector semilocalization]
\label{thm:right_vector_semilocalization}
Under Assumptions~\ref{ass:proportional_regime}
and~\ref{ass:subcritical_regime}, with $\gC_0$ sufficiently large
depending on $(B,\qr,\nu)$ for the prescribed failure exponent,
each assertion below holds
simultaneously for all normalized eigenpairs $(\lambda,\rvw)$ of
$\scA$ in its stated spectral range and all $2\leq p\leq\infty$.
\begin{enumerate}[label=(\roman*)]
    \item \emph{Critical sparsity.} Suppose in addition that Assumption~\ref{ass:critical_sparsity} holds. For every $\epsilon,\nu>0$ there is $\gC_{B,\epsilon,\nu}>0$ such that if $\lambda\geq\redge+\gC_{B,\epsilon,\nu}\error_d^{1/4}$, the following holds with probability at least $1-\gC_{B,\epsilon,\nu}N^{-\nu}$
    \begin{align}
        \|\rvw\|_p^2
        \geq
        N^{(2/p-1)\rhoR(\lambda,b)-\epsilon},
        \label{eqn:right-lp-semilocalization-critical}
    \end{align}
    where the transition exponent function $\rhoR(\lambda,b)$ for $\lambda\geq\redge$ is defined by
    \begin{align}
        \rhoR(\lambda,b)
        &\coloneqq
          \max\bigl\{1-b\rateR(\lambda;\qr), \quad 0\bigr\}.
        \label{eqn:right-transition-count-exponent}
    \end{align}
    In particular, $\|\vw\|_{\infty}^{2}\geq N^{-\rhoR(\lambda,b) - \epsilon}$ as $p\to\infty$.

    \item \emph{Subcritical sparsity.} Suppose $d/\log(N) \to 0$.
    For every $\nu>0$, choose constants $K_{B,\nu},\gC_{B,\nu}>0$
    sufficiently large. For $\lambda\geq\redge+\gC_{B,\nu}\error_d^{1/4}$, define the \emph{resonance width}
    \begin{align}
        \delta_\lambda
        \coloneqq K_{B,\nu}\error_d
        \left(1+\frac1{(\lambda-\redge)\log d}\right).
        \label{eqn:right-intro-resonance-width}
    \end{align}
    Then $\delta_\lambda \leq (\lambda-\redge)/2$, and the following holds with probability at least $1-\gC_{B,\nu}N^{-\nu}$:
    \begin{align}
        \|\rvw\|_p^2
        \geq \frac12 N^{2/p-1}
        \exp\left[
          \left(1-\frac2p\right)
          \left\{
            \min\bigl\{
              d\rateR(\lambda-\delta_\lambda;\qr),
              \log N
            \bigr\}
            -\gC_{B,\nu}
              \bigl(\sqrt{d\log d}+\log\log N\bigr)
          \right\}
        \right].
        \label{eqn:right-lp-semilocalization-sublog}
    \end{align}
    In particular, $\|\rvw\|_{\infty}^{2}\geq (2N)^{-1} \exp(d\rateR(\redge;\qr))$ as $p\to\infty$.
\end{enumerate}
\end{theorem}

\subsubsection{Left edge results}\label{subsec:intro-zero-modes}
At the left edge, isolated vertices and finite configurations can produce
rank deficiency and positive singular values approaching zero. An
$n\times m$ matrix with $n\geq m$ has column rank deficiency when its
rank is less than $m$, equivalently when one of its $m$ singular values
is zero. This is distinct from the $n-m$ structural zero eigenvalues
of its Hermitian linearization. An isolated vertex in $\gV_2$ gives a
zero column of the raw matrix $\widetilde{\rmA}$; centering changes the
corresponding kernel, as described below.

The isolation threshold in $\gV_2$ is given by the following lemma,
proved in Appendix~\ref{subsec:bipartite-vertex-counts}.

\begin{lemma}
\label{lem:bipartite-isolated-column-threshold}
Suppose Assumptions~\ref{ass:proportional_regime}
and~\ref{ass:critical_sparsity} hold with fixed $\qr$. Let $\sI_2\coloneqq\{y\in\gV_2:\rD_y=0\}$ denote the set of isolated vertices in $\gV_2$. Define the threshold function $\bcritIso(\qr)\coloneqq\qr^{-2}$. The following statements hold.
\begin{enumerate}[label=(\roman*)]
    \item If $b>\bcritIso(\qr)$, then $|\sI_2| = 0$ with probability at least $1 - ((1+\qr^4)^{-1}+o(1))N^{1-b\qr^2}$.

    \item If $b=\bcritIso(\qr)$, then $|\sI_2| \xrightarrow{\mathrm d} \Pois(({1+\qr^4})^{-1})$.

    \item If $b<\bcritIso(\qr)$, then for every $0<\varepsilon<(1+\qr^4)^{-1}$ and every $\nu>0$, there is a constant $C_{b,\qr,\varepsilon,\nu}>0$ such that, for all sufficiently large $N$,
    \begin{align}
        \P\left(
            |N^{-1+b\qr^2}|\sI_2|
         - (1+\qr^4)^{-1}| \leq \varepsilon\right)
        \geq1-C_{b,\qr,\varepsilon,\nu}N^{-\nu}.
        \label{eqn:bipartite-nonstructural-zero-count}
    \end{align}
\end{enumerate}
\end{lemma}

\paragraph{Full rank regime.}
For $\qr>1$ and $b>\bcritIso(\qr)$,
Theorem~\ref{thm:left_edge_behavior} identifies locations of left outliers induced by low-degree vertices in $\gV_2$. In its proof, the least singular value is bounded away from zero with polynomially small failure probability.

\begin{assumption}\label{ass:no-isolated-vertices-in-V2}
    Under Assumptions~\ref{ass:proportional_regime} and ~\ref{ass:critical_sparsity}, we assume $b>\qr^{-2} \eqqcolon \bcritIso(\qr)$.
\end{assumption}

Define the vertices that induce left outliers by
\begin{align}
    \gL_2
    \coloneqq
    \{x\in\gV_2:0<\alpha_x<\qr^2-1,\
      \Lambda_{\qr}(\alpha_x)
      \leq\ledge-\tfrac18\error_d^{1/4}\}.
    \label{eqn:left-outlier-profile-set}
\end{align}
Arrange the multiset $\Lambda^{\gL} \coloneqq \{\Lambda_{\qr}(\alpha_x):x\in\gL_2\}$ in non-increasing order as $\Lambda^{\gL}_1\geq\cdots\geq \Lambda^{\gL}_{|\gL_2|}$.

\begin{theorem}[Locations of left outliers]
\label{thm:left_edge_behavior}
Suppose $\ratio > 1$ and Assumption~\ref{ass:no-isolated-vertices-in-V2} holds. There are constants $\nu=\nu(b,\qr)>0$ and $\gC=\gC(b,\qr)>0$ such that the following hold with probability at least $1-\gC N^{-\nu}$.
\begin{enumerate}[label=(\roman*)]
    \item If $\gL_2\neq\emptyset$, there are indices
    $1\leq j_1<\cdots<j_{|\gL_2|}\leq m$ such that
    \begin{align}
        |\lambda_{j_k}(\scA)-\Lambda^{\gL}_k|
        \leq\gC \error_d^{3/4}/\log(d),
        \qquad 1\leq k\leq |\gL_2|,
        \label{eqn:left-outlier-locations}
    \end{align}
    \item For all $1\leq j \leq m - |\gL_2|$, we have
    \begin{align}
        \lambda_{j}(\scA)
        \geq\ledge-\error_d^{1/4}.
        \label{eqn:left-no-additional-outliers}
    \end{align}
\end{enumerate}
\end{theorem}
\begin{remark}
In our previous work, a technical condition $b > b_{\mathrm I}$ \cite[Assumption 1.6]{dumitriu2025singular} was imposed, where $b_{\mathrm I}\coloneqq \qr^2/\benrate(\qr^2-1)$ with $\benrate(t)\coloneqq(1+t)\log(1+t)-t$. Here, Theorem~\ref{thm:left_edge_behavior} removes the technical condition by incorporating the recent development \cite{wang2026regularized} for spectral radius estimates of non-backtracking matrices constructed from regularized sparse random matrices. See Figure~\ref{fig:combined-aspect-ratio-thresholds} for an illustration.
\end{remark}

\paragraph{Semilocalization of left singular vectors.}
For $0<\lambda\leq\ledge$, define the inverse branch of the left location map by
\begin{align}
    \invL(\lambda)
    &\coloneqq
    \frac12\left(
      \lambda^2+\qr^2-\qr^{-2}-\Delta_{\Lambda}(\lambda)
    \right).
    \label{eqn:Lambda-inverse-lower}
\end{align}
The left degree rate function is
\begin{align}
    \rateL(t;\qr)
    \coloneqq I_{\qr}\bigl(\invL(t)\bigr).
    \label{eqn:left-degree-rate}
\end{align}
Under Assumption~\ref{ass:critical_sparsity}, the corresponding profile-count exponent function for $0<\lambda<\ledge$ is
\begin{align}
    \rhoL(\lambda,b)
    \coloneqq
    \max\left\{1-b\rateL(\lambda;\qr),\,0\right\}.
    \label{eqn:left-transition-count-exponent}
\end{align}

\begin{theorem}[Left singular vector semilocalization]
\label{thm:left_vector_semilocalization}
Suppose $\ratio>1$ and Assumption~\ref{ass:no-isolated-vertices-in-V2} holds. Fix a compact interval $\sK\subset(0,\ledge)$ of positive length. There is a constant
$\nu=\nu(b,\qr,\sK)>0$ such that, for every $\epsilon>0$, there
is a constant $\gC_\epsilon=\gC_\epsilon(b,\qr,\sK)>0$ for which the
following holds with probability at least $1-\gC_\epsilon N^{-\nu}$.
For every normalized eigenpair $(\lambda,\rvw)$ of $\scA$ with
$\lambda\in\sK$ and all $2\leq p\leq\infty$, we have
\begin{align}
    \|\rvw\|_p^2
    \geq
    N^{(2/p-1)\rhoL(\lambda,b)-\epsilon}.
    \label{eqn:left-lp-semilocalization}
\end{align}
In particular, $\|\rvw\|_\infty^2\geq N^{-\rhoL(\lambda,b)-\epsilon}$ as $p\to\infty$.
\end{theorem}

\paragraph{Kernel space induced by isolated vertices in $\gV_2$.}
When $b<\bcritIso(\qr)$,
\eqref{eqn:bipartite-nonstructural-zero-count} gives polynomially many
isolated vertices in $\gV_2$ with overwhelming probability. Their
columns in $\mX$ are identical and nonzero, so kernel vectors supported
on these coordinates must have zero coordinate sum. The following
lemma characterizes this subspace; its proof is in
Appendix~\ref{subsec:zero-and-near-zero-modes}.

\begin{lemma}
\label{lem:bipartite-isolated-column-zero-modes}
Let $\sK_{\rm iso}\coloneqq \{\vv\in\R^{m}: \operatorname{supp}\vv\subset\sI_2, \mX\vv=0\}$ be the kernel subspace induced by the set of isolated vertices $\sI_2$ in Lemma~\ref{lem:bipartite-isolated-column-threshold}. Then
\begin{align}
    \sK_{\rm iso}
    &=\bigl\{\vv\in\R^m:
      \operatorname{supp}\vv\subset\sI_2,\quad
      \sum\nolimits_{y\in\sI_2} \evv_y=0\bigr\},
    \qquad
    \dim(\sK_{\rm iso}) = (|\sI_2|-1)_{+}.
    \label{eqn:bipartite-isolated-V2-kernel}
\end{align}
In particular, $\dim(\ker\mX) \geq \dim(\sK_{\rm iso})$.
\end{lemma}
\begin{remark}
    Suppose $|\sI_2|\geq2$. Every unit $\vv\in\sK_{\rm iso}$ satisfies $\|\vv\|_\infty^2\geq |\sI_2|^{-1}$. For distinct $i,j\in\sI_2$, the vector $(\ones_i-\ones_j)/\sqrt2$ belongs to this subspace and has squared infinity norm $1/2$. A balanced sign vector on the largest even subset of $\sI_2$, normalized to have norm one, has squared infinity norm at most $2/|\sI_2|$. Thus, below the threshold, \eqref{eqn:bipartite-nonstructural-zero-count} gives the existence of both localized unit kernel vectors and diffuse ones with $\|\vv\|_\infty^2 \asymp N^{b\qr^2-1}$ with overwhelming probability.
\end{remark}

\paragraph{Singular values induced by small components.}
After centering, columns of $\mX$ corresponding to isolated vertices in $\gV_2$ are identical and nonzero. Consequently, $\mX$ becomes rank deficient when $|\sI_2|\geq2$; see Corollary~\ref{cor:least-singular-value-regimes}. Beside that, small components are bounded trees with probability tending to one, inducing positive singular values with scale $d^{-1/2}$; see Corollary~\ref{cor:components-size}.

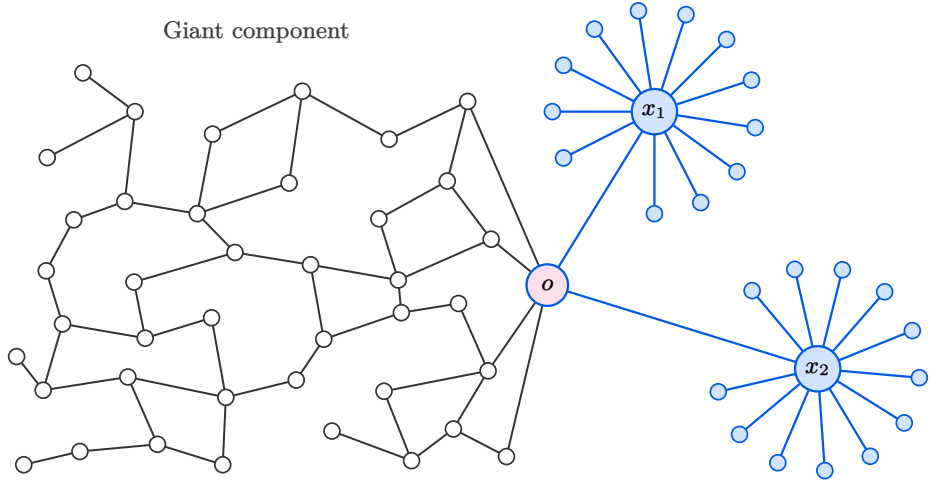
\begin{figure}[!htbp]
    \centering
    \definecolor{forkSemilocalized}{HTML}{D1E3FF}
    \definecolor{forkSemilocalizedBorder}{HTML}{005AE5}
    \begin{tikzpicture}[
        x=1.35cm,y=1.35cm,
        giant edge/.style={draw=black!80,line width=0.75pt},
        fork edge/.style={draw=forkSemilocalizedBorder,line width=0.85pt},
        giant vertex/.style={circle,draw=black!80,fill=white,
            minimum size=2.1mm,inner sep=0pt,line width=0.65pt},
        fork leaf/.style={circle,draw=forkSemilocalizedBorder,fill=forkSemilocalized,
            minimum size=2.1mm,inner sep=0pt,line width=0.65pt},
        fork base/.style={circle,draw=forkSemilocalizedBorder,fill=magenta!15,
            minimum size=5.5mm,inner sep=0pt,line width=0.85pt,font=\small},
        fork hub/.style={circle,draw=forkSemilocalizedBorder,fill=forkSemilocalized,
            minimum size=6mm,inner sep=0pt,line width=0.85pt,font=\small}
    ]
        \coordinate (base-position) at (0,0);
        \coordinate (hub-one-position) at (1.05,1.70);
        \coordinate (hub-two-position) at (2.65,-0.82);

        \foreach \i/\x/\y in {
            1/-4.90/1.25,2/-4.04/1.70,3/-3.28/1.48,
            4/-2.40/1.90,5/-1.55/1.43,6/-0.78/1.80,
            7/-0.98/1.02,8/-1.65/0.65,9/-2.53/1.00,
            10/-3.43/0.70,11/-4.14/0.82,12/-4.64/0.64,
            13/-4.91/0.14,14/-4.05/0.03,15/-3.06/0.32,
            16/-2.32/0.20,17/-1.46/0.05,18/-0.55/0.45,
            19/-0.87/-0.18,20/-1.43/-0.27,21/-2.19/-0.53,
            22/-3.29/-0.32,23/-3.94/-0.52,24/-4.75/-0.38,
            25/-4.94/-1.03,26/-4.11/-0.90,27/-3.15/-1.10,
            28/-2.46/-0.93,29/-1.60/-1.04,30/-0.58/-0.84,
            31/-0.92/-1.41,32/-1.33/-1.72,33/-2.11/-1.43,
            34/-3.18/-1.76,35/-3.82/-1.56,36/-4.58/-1.63,
            37/-4.55/2.08,38/-0.40/-1.68,
            39/-5.13/-1.76,40/-5.20/-0.70}{
            \coordinate (giant-\i) at (\x,\y);
        }
        \foreach \a/\b in {
            1/2,3/4,4/5,5/6,12/11,11/10,10/9,8/7,
            14/15,15/16,16/17,17/18,24/23,23/22,21/20,20/19,
            25/26,26/27,27/28,29/30,36/35,35/34,33/32,32/31,
            2/11,3/10,4/9,6/7,12/13,10/15,8/17,7/18,
            13/24,14/23,16/21,17/20,24/25,22/27,21/28,19/30,
            26/35,27/34,29/32,30/31,2/37,31/38,36/39,25/40}{
            \draw[giant edge] (giant-\a)--(giant-\b);
        }
        \foreach \i in {6,18,30,38}{
            \draw[giant edge] (base-position)--(giant-\i);
        }

        \draw[fork edge] (base-position)--(hub-one-position);
        \draw[fork edge] (base-position)--(hub-two-position);
        \foreach \i in {1,...,12}{
            \pgfmathsetmacro{\forkAngleOne}{270+27*(\i-1)}
            \pgfmathsetmacro{\forkAngleTwo}{193+27*(\i-1)}
            \coordinate (leaf-one-\i) at
                ($(hub-one-position)+(\forkAngleOne:1)$);
            \coordinate (leaf-two-\i) at
                ($(hub-two-position)+(\forkAngleTwo:1)$);
            \draw[fork edge] (hub-one-position)--(leaf-one-\i);
            \draw[fork edge] (hub-two-position)--(leaf-two-\i);
        }

        \foreach \i in {1,...,40}{
            \node[giant vertex] at (giant-\i) {};
        }
        \foreach \i in {1,...,12}{
            \node[fork leaf] at (leaf-one-\i) {};
            \node[fork leaf] at (leaf-two-\i) {};
        }
        \node[fork base] at (base-position) {$o$};
        \node[fork hub] at (hub-one-position) {$x_1$};
        \node[fork hub] at (hub-two-position) {$x_2$};
        \node[font=\small,text=black!75] at (-2.85,2.48)
            {Giant component};
    \end{tikzpicture}
    \caption{An $F$-fork with $F=12$ leaves per hub, attached to the giant through its base $o$. The base and leaves lie in one vertex class, and the hubs $x_1,x_2$ in the other.}
    \label{fig:bipartite-pendant-configurations}
\end{figure}

\paragraph{Kernel spaces from star-tuning forks.}
A \emph{star-tuning fork} has a base $o$ joined to the hubs $x_1,x_2$ of two disjoint stars, each with $F\in\N_0$ leaves in sets $R_1,R_2$. We call it an \emph{$F$-fork}; see Figure~\ref{fig:bipartite-pendant-configurations}. It is
\emph{rooted in $\gG$} if
\begin{align}
    \rD_{x_1}=\rD_{x_2}=F+1,
    \qquad \rD_y=1 \quad\text{for every }y\in R_1\cup R_2.
    \label{eqn:star-tuning-fork-degree-conditions}
\end{align}
As highlighted in \cite{alt2021delocalization}, $F$-forks in $\gG(N, d/N)$ induce vanishing eigenvalues when $d< \log(N)/2$. We give the full descriptions for the kernel space induced by $F$-forks in Lemma~\ref{lem:bipartite-star-tuning-fork-modes}, with its proof deferred to Section~\ref{subsec:bipartite-component-structure}.

\begin{lemma}[Singular values and kernel vectors of a star-tuning fork]
    \label{lem:bipartite-star-tuning-fork-modes}
    Let an $F$-fork be rooted in $\gG$.
    \begin{enumerate}[label=(\roman*),topsep=3pt,itemsep=3pt,parsep=0pt]
        \item \emph{$0$-forks.} A unit kernel vector is
        \begin{align}
            \rvw_0
            \coloneqq
            \frac{1}{\sqrt2}
            (\ones_{x_1}-\ones_{x_2})
            \in
            \begin{cases}
                \ker\mX,&o\in\gV_1,\\
                \ker\mX^{\sT},&o\in\gV_2.
            \end{cases}
            \label{eqn:star-tuning-fork-zero-mode}
        \end{align}
        In particular, when $o\in\gV_1$, at least one of the $m$ singular values of
        $\mX$ is zero, and $\rvw_0\notin\sK_{\rm iso}$ from
        Lemma~\ref{lem:bipartite-isolated-column-zero-modes}, since the
        hubs $x_{1}$ and $x_{2}$ are not isolated vertices. When $o\in\gV_2$ and $n>m$, this left-kernel vector may lie in the structural row kernel and does not by itself imply column rank deficiency of $\mX$. In the square case, either orientation forces a zero singular value.
    
        \item \emph{Positive singular values ($F\geq1$).}
    Define vectors $\vh_F$ and $\vll_F$, supported on the hubs and leaves, respectively,
        \begin{align}
            \vh_F&\coloneqq (\ones_{x_1}-\ones_{x_2})/\sqrt2,
            \qquad
            \vll_F\coloneqq(\ones_{R_1}-\ones_{R_2})/\sqrt{2F},
            \label{eqn:star-tuning-fork-singular-profiles}
        \end{align}
        with support sizes $2$ and $2F$ and norms $\|\vh_F\|_\infty=1/\sqrt2$ and $\|\vll_F\|_\infty=1/\sqrt{2F}$. The positive singular value $s_F\coloneqq\sqrt{F/d}$ of $\mX$ has the corresponding singular pair $(\vu_F, \vv_F)\in\R^n\times\R^m$, satisfying $\mX \vv_F = s_F \vu_F$ and $\mX^{\sT}\vu_F=s_F\vv_F$, where
        \begin{align}
            (\vu_F,\vv_F)&\coloneqq (\vll_F,\vh_F), \quad \text{when }o\in\gV_1,
            &
            (\vu_F,\vv_F)&\coloneqq (\vh_F,\vll_F), \quad \text{when }o\in\gV_2.
            \label{eqn:star-tuning-fork-singular-vectors}
        \end{align}
    
        \item \emph{Leaf-supported kernel spaces ($F\geq2$).} Suppose
        $o\in\gV_j$ for some $j\in\{1,2\}$. For $i\in\{1,2\}$, define
        \begin{align}
            \sK^{\rm leaf}_{i}&\coloneqq
            \bigl\{
              \rvv\in\R^{\gV_j}:
              \operatorname{supp}\rvv\subset R_i,\quad
              \sum\nolimits_{y\in R_i}\evv_y=0
            \bigr\},
            \label{eqn:star-tuning-fork-branch-zero-modes}
        \end{align}
        with $\dim(\sK^{\rm leaf}_{i})=F-1$. Furthermore, $\sK^{\rm leaf}_{1}$ and $\sK^{\rm leaf}_{2}$ are orthogonal to each other, and both are orthogonal to $\vll_F$ in item~(ii). Thus,
        \begin{align}
            \sK^{\rm leaf}_{1}\oplus\sK^{\rm leaf}_{2}&\subset \ker\mX^{\sT}, \quad \text{when }o\in\gV_1,
            &
            \sK^{\rm leaf}_{1}\oplus\sK^{\rm leaf}_{2}&\subset \ker\mX, \quad \text{when }o\in\gV_2.
            \label{eqn:star-tuning-fork-branch-kernels}
        \end{align}
         In particular, when $o\in\gV_2$, at least $2(F-1)$ of the $m$ singular values of $\mX$ are zero. When $o\in\gV_1$ and $n>m$, these left-kernel spaces do not alone imply column rank deficiency of $\mX$.
    \end{enumerate}
\end{lemma}

For fixed $F$, singular vectors constructed from $F$-forks are localized, with exact zero or vanishing singular values. 
Below, we establish thresholds for emergence of $F$-forks with base $o$ in either $\gV_{1}$ and $\gV_{2}$.

\begin{lemma}[Emergence of $F$-forks]
\label{lem:bipartite-star-tuning-fork-threshold}
Suppose that Assumptions~\ref{ass:proportional_regime} and \ref{ass:critical_sparsity} hold, and fix $F\in\N_0$. Let
\begin{align}
    \bcritF{F, 1}(\qr)&\coloneqq\frac{1}{2(\qr^2+F\qr^{-2})},
    \qquad
    \bcritF{F, 2}(\qr)\coloneqq\frac{1}{2(\qr^{-2}+F\qr^2)}.
    \label{eqn:star-tuning-fork-orientation-threshold-definitions}
\end{align}
denote the thresholds for the emergence of $F$-forks based in $\gV_1$ and $\gV_2$, respectively. There is a constant $C=C(b,\qr,F)>0$ such that, for all sufficiently large $N$,
\begin{enumerate}[label=(\alph*)]
    \item when $b>\bcritF{F, i}(\qr)$, an $F$-fork based in $\gV_i$ is absent with probability at least $1-CN^{1-b/\bcritF{F, i}(\qr)}(\log N)^{2F+2}$;
    \item when $b\leq\bcritF{F, i}(\qr)$, an $F$-fork based in $\gV_i$ emerges with probability at least $1-CN^{b/\bcritF{F, i}(\qr)-1}(\log N)^{-2F-2}$.
\end{enumerate}
\end{lemma}
The preceding probability bounds concern forks anywhere in the graph. For each fixed $F$ and orientation $i$, Corollary~\ref{cor:bipartite-star-tuning-fork-giant-emergence} additionally shows that such a fork occurs in the giant component with probability $1-o(1)$.
\begin{remark}[Rank deficiency induced by $F$-forks]
Two families that force column rank deficiency are
$0$-forks based in $\gV_1$ and $F$-forks with $F\geq2$ based in
$\gV_2$. Their first thresholds are
$\bcritF{0,1}(\qr)=1/(2\qr^2)$ and
$\bcritF{2,2}(\qr)=1/(2\qr^{-2}+4\qr^2)$, respectively.
Both lie below $\bcritIso(\qr)$, while the $0$-fork threshold is
above the positive $1$-fork threshold $\bcritFork(\qr)$ defined below.
\end{remark}

Among the $F$-forks, the first positive singular values are induced by $1$-forks, with both orientations having the same threshold,
\begin{align}
    \bcritFork(\qr)
    \coloneqq\bcritF{1, 1}(\qr)
    =1/(2(\qr^2+\qr^{-2})).
    \label{eqn:nonzero-star-tuning-fork-threshold}
\end{align}
At criticality $d \asymp \log(N)$, one can prove that $F \lesssim \log(N)$ (Lemma~\ref{lem:bipartite-tuning-forks}) for every fork on the giant, hence $s_F = \sqrt{F/d} = O(d^{-1/2})$ with overwhelming probability. Consequently, these fork-induced singular values eventually lie below every fixed compact interval in $(0,\ledge)$ used in Theorem~\ref{thm:left_vector_semilocalization}.

\subsection{Transition boundaries}\label{sec:localized-delocalized-transition}
\subsubsection{Phase diagrams}

\begin{figure}[p]
    \centering
    \definecolor{thresholdorange}{RGB}{205,125,0}
    \definecolor{thresholdteal}{RGB}{0,128,105}
    \definecolor{thresholdpurple}{RGB}{132,66,155}
    \definecolor{thresholdpeach}{RGB}{252,224,211}
    \tikzset{
        threshold axis/.style={-{Stealth[length=1.7mm]},thick,
            line cap=round,line join=round},
        threshold boundary/.style={line width=1pt,
            line cap=round,line join=round},
        threshold tick/.style={font=\footnotesize},
        threshold endpoint/.style={circle,inner sep=0pt,outer sep=0pt,
            minimum size=3.2pt,draw=white,line width=0.25pt},
        threshold number/.style={draw,circle,fill=white,
            line width=0.85pt,inner sep=1pt,minimum size=7.5mm,
            font=\fontsize{14}{15}\selectfont},
        threshold key number/.style={threshold number,
            minimum size=6mm,font=\fontsize{11}{12}\selectfont},
        threshold key text/.style={anchor=west,inner sep=0pt,
            font=\fontsize{10.5}{12}\selectfont},
        threshold key reference/.style={draw,circle,inner sep=0.8pt,
            line width=0.6pt,minimum size=4.8mm,font=\small},
        threshold legend auto/.style={fill=white,draw=black!30,
            line width=0.35pt,rounded corners=1pt,
            inner xsep=5pt,inner ysep=4pt,font=\small}
    }
    \resizebox{0.985\textwidth}{!}{%
    \begin{tikzpicture}[
        x=1cm,y=1cm,>=Stealth,
        declare function={
            combinedright(\x)=1/((\x*\x+1)*ln(1+1/(\x*\x))-1);
            combinedrone(\x)=1/((1+1/(\x*\x))*ln(1+\x*\x)-1);
            combinedleft(\x)=1/((\x*\x-1)*ln(1-1/(\x*\x))+1);
            combinediso(\x)=1/(\x*\x);
            combinedfork(\x)=1/(2*(\x*\x+1/(\x*\x)));
        }
    ]
        \useasboundingbox (0,0.6) rectangle (16.1,9.1);
        \node[threshold legend auto,anchor=north] at (8.05,9.0) {%
            \begin{tabular}{@{}c@{\hspace{3pt}}l@{\hspace{12pt}}c@{\hspace{3pt}}l@{\hspace{12pt}}c@{\hspace{3pt}}l@{\hspace{12pt}}c@{\hspace{3pt}}l@{\hspace{12pt}}c@{\hspace{3pt}}l@{}}
                \tikz[baseline=-0.55ex]\draw[threshold boundary,blue!65!black]
                    (0,0)--(4.2mm,0); & $\bcritR(\qr)=\bcritRi{2}(\qr)$ &
                \tikz[baseline=-0.55ex]\draw[threshold boundary,thresholdorange]
                    (0,0)--(4.2mm,0); & $\bcritRi{1}(\qr)$ &
                \tikz[baseline=-0.55ex]\draw[threshold boundary,red!85!black]
                    (0,0)--(4.2mm,0); & $\bcritL(\qr)$ &
                \tikz[baseline=-0.55ex]\draw[threshold boundary,thresholdteal]
                    (0,0)--(4.2mm,0); & $\bcritIso(\qr)$ &
                \tikz[baseline=-0.55ex]\draw[threshold boundary,thresholdpurple,dashed]
                    (0,0)--(4.2mm,0); & $\bcritFork(\qr)$
            \end{tabular}%
        };

        \begin{scope}[shift={(1.1,1.25)},x=15.6cm,y=1.7cm]
            \begin{scope}
                \clip (0,0) rectangle (0.5,4);
                \fill[black!5] (0,0) rectangle (0.5,4);
                \fill[blue!12] (0,0)
                    --plot[domain=1:1.5,samples=180,smooth]
                        ({\x-1},{combinedright(\x)})--(0.5,0)--cycle;
                \fill[thresholdorange!18] (0,0)
                    --plot[domain=1:1.5,samples=180,smooth]
                        ({\x-1},{combinedrone(\x)})--(0.5,0)--cycle;
                \fill[red!11] (0,0)--(0,1)
                    --plot[domain=1.001:1.5,samples=180,smooth]
                        ({\x-1},{combinedleft(\x)})--(0.5,0)--cycle;
                \fill[thresholdpeach] (0,0)--(0,1)
                    --plot[domain=1.001:1.5,samples=180,smooth]
                        ({\x-1},{min(combinedleft(\x),combinedrone(\x))})
                    --(0.5,0)--cycle;
                \fill[thresholdteal!15] (0,0)
                    --plot[domain=1:1.5,samples=180,smooth]
                        ({\x-1},{combinediso(\x)})--(0.5,0)--cycle;
                \fill[thresholdteal!15] (0,0)
                    --plot[domain=1:1.5,samples=180,smooth]
                        ({\x-1},{combinedfork(\x)})--(0.5,0)--cycle;
                \path[pattern=north east lines,pattern color=thresholdpurple] (0,0)
                    --plot[domain=1:1.5,samples=180,smooth]
                        ({\x-1},{combinedfork(\x)})--(0.5,0)--cycle;
                \draw[black!65,line width=0.85pt,dash pattern=on 3pt off 2pt]
                    (0,1.9993166721)--(0.1822169196,1.9993166721)
                    --(0.1822169196,0);
                \draw[threshold boundary,blue!65!black]
                    plot[domain=1:1.5,samples=180,smooth]
                        ({\x-1},{combinedright(\x)});
                \draw[threshold boundary,thresholdorange]
                    plot[domain=1:1.5,samples=180,smooth]
                        ({\x-1},{combinedrone(\x)});
                \draw[threshold boundary,red!85!black] (0,1)
                    --plot[domain=1.001:1.5,samples=180,smooth]
                        ({\x-1},{combinedleft(\x)});
                \draw[threshold boundary,thresholdteal]
                    plot[domain=1:1.5,samples=180,smooth]
                        ({\x-1},{combinediso(\x)});
                \draw[threshold boundary,thresholdpurple,dashed]
                    plot[domain=1:1.5,samples=180,smooth]
                        ({\x-1},{combinedfork(\x)});
            \end{scope}
            \draw[threshold axis] (0,0)--(0.52,0);
            \node[font=\small,anchor=west] at (0.52,-0.15) {$\qr$};
            \draw[threshold axis] (0,0)--(0,4.13);
            \foreach \x/\label in {0/1,0.1822169196/1.182,0.5/1.5}{
                \draw (\x,0)--++(0,-0.036);
                \node[threshold tick,anchor=north] at (\x,-0.065) {$\label$};
            }
            \foreach \y in {0,1,2}{
                \draw (0,\y)--++(-0.004,0);
                \node[threshold tick,anchor=east] at (-0.009,\y) {$\y$};
            }
            \draw (0,4)--++(-0.004,0);
            \node[font=\small,anchor=east] at (-0.009,4) {$b$};
            \draw (0,2.588699)--++(-0.004,0);
            \node[font=\scriptsize,anchor=east] at (-0.0125,2.588699)
                {$\dfrac{1}{\ln4-1}$};
            \node[threshold endpoint,fill=black] at (0,2.588699) {};
            \node[threshold endpoint,fill=red!85!black] at (0,1) {};
            \node[threshold endpoint,fill=black] at (0.1822169196,1.9993166721) {};

            \node[threshold number,draw=black!60,text=black!70]
                at (0.095,3.55) {1};
            \node[threshold number,draw=blue!65!black,text=blue!65!black]
                at (0.25524,3.08600) {2a};
            \node[threshold number,draw=thresholdorange,text=thresholdorange]
                at (0.042,1.70) {2b};
            \node[threshold number,draw=red!85!black,text=red!85!black]
                at (0.39363,2.36733) {3a};
            \node[threshold number,draw=thresholdorange,
                text=thresholdorange!80!black]
                at (0.3125,1.21) {3b};
            \node[threshold number,draw=thresholdteal,text=thresholdteal]
                at (0.05,0.57) {4};
            \node[threshold number,draw=thresholdpurple,text=thresholdpurple]
                (phase-five) at (-0.045,0.38) {5};
            \draw[thresholdpurple,line width=0.7pt]
                (phase-five.south east)--(0.02,0.10);
            \node[threshold endpoint,fill=thresholdpurple]
                at (0.02,0.10) {};
        \end{scope}

        \node[threshold key number,draw=black!60,text=black!70]
            at (9.6,7.55) {1};
        \node[threshold key text] at (10.1,7.55) {No outliers};
        \node[threshold key number,draw=blue!65!black,text=blue!65!black]
            at (9.6,6.65) {2a};
        \node[threshold key text] at (10.1,6.65)
            {Right outliers only ($\gV_2$)};
        \node[threshold key number,draw=thresholdorange,text=thresholdorange]
            at (9.6,5.75) {2b};
        \node[threshold key text] at (10.1,5.75)
            {Right outliers only ($\gV_1,\gV_2$)};
        \node[threshold key number,draw=red!85!black,text=red!85!black]
            at (9.6,4.85) {3a};
        \node[threshold key text] at (10.1,4.85)
            {Left ($\gV_2$) $+$ right outliers ($\gV_2$)};
        \node[threshold key number,draw=thresholdorange,
            text=thresholdorange!80!black] at (9.6,3.95) {3b};
        \node[threshold key text] at (10.1,3.95)
            {Left ($\gV_2$) $+$ right outliers ($\gV_1,\gV_2$)};
        \node[threshold key number,draw=thresholdteal,text=thresholdteal]
            at (9.6,3.05) {4};
        \node[threshold key text,text width=6cm,align=left] at (10.1,3.05)
            {Region \tikz[baseline=(reference-three-b.base)]{
                \node[threshold key reference,draw=thresholdorange,
                    text=thresholdorange!80!black] (reference-three-b) {3b};}
                $+$ zero columns};
        \node[threshold key number,draw=thresholdpurple,text=thresholdpurple]
            at (9.6,2.15) {5};
        \node[threshold key text,text width=6cm,align=left] at (10.1,2.15)
            {Region \tikz[baseline=(reference-four.base)]{
                \node[threshold key reference,draw=thresholdteal,
                    text=thresholdteal] (reference-four) {4};}
                $+$ $1$-fork singular pairs};
        \node[font=\small,inner sep=0pt,anchor=north west,
            text width=6cm,align=left] at (10.1,1.35)
            {Sources at fixed positive locations are listed in parentheses.};
    \end{tikzpicture}%
    }
    \caption{Selected transition thresholds for $1\leq\qr\leq1.5$
    and $0\leq b\leq4$. Dashed grey guides mark the crossing of
    $\bcritRi{1}$ and $\bcritL$ at
    $(\qr,b)\approx(1.1822,1.9993)$. The right thresholds meet at
    $(1,1/(\ln4-1))$, and the left-outlier and isolation thresholds
    meet at $(1,1)$. The dashed purple curve marks the positive
    $1$-fork threshold. At $\qr=1$, the left threshold is interpreted
    by continuity and there is no positive left gap.}
    \label{fig:combined-aspect-ratio-thresholds}

    \vspace{0.5em}

    \definecolor{phaseDelocalized}{HTML}{FFBDCB}
    \definecolor{phaseSemilocalized}{HTML}{D1E3FF}
    \definecolor{phaseDelocalizedBorder}{HTML}{E50031}
    \definecolor{phaseSemilocalizedBorder}{HTML}{005AE5}
    \colorlet{phaseLeft}{red!85!black}
    \colorlet{phaseRight}{blue!65!black}
    \definecolor{phaseOrange}{RGB}{205,125,0}
    \definecolor{phaseTeal}{RGB}{0,128,105}
    \definecolor{phasePurple}{RGB}{132,66,155}
    \definecolor{phasePeach}{RGB}{252,224,211}
    \tikzset{
        singularvalue axis/.style={-{Stealth[length=1.7mm]},thick},
        singularvalue boundary/.style={line width=1.6pt,line cap=round,line join=round},
        singularvalue text/.style={font=\scriptsize,inner sep=1.5pt},
        singularvalue number/.style={circle,draw=black!60,text=black!70,fill=white,
            font=\fontsize{9}{10}\selectfont,minimum size=5mm,
            inner sep=0.6pt,line width=0.65pt},
        singularvalue guide/.style={dash pattern=on 3pt off 2pt,
            black!65,line width=0.85pt},
        singularvalue point/.style={circle,inner sep=0pt,minimum size=2.8pt,
            draw=white,line width=0.2pt},
        singularvalue legend/.style={fill=white,draw=black!30,
            rounded corners=1pt,inner xsep=7pt,inner ysep=5pt,font=\small}
    }
    \resizebox{\textwidth}{!}{%
    \begin{tikzpicture}[x=1cm,y=1cm,>=Stealth,
        declare function={
            singularvaluex(\u,\q)=sqrt(\q*\q+1/(\q*\q)+\u+1/\u);
            singularvaluertwo(\u,\q)=1/((\q*\q+\u)*ln(1+\u/(\q*\q))-\u);
            singularvaluerone(\u,\q)=1/((1/(\q*\q)+\u)*ln(1+\u*\q*\q)-\u);
            singularvalueleftx(\a,\q)=sqrt(max(0,\a*(\a-\q*\q+1/(\q*\q))/(\a-\q*\q)));
            singularvalueleftb(\a,\q)=1/(\a*ln(\a/(\q*\q))-\a+\q*\q);
        }]
        \useasboundingbox (0,-0.65) rectangle (16.3,7.6);
        \node[singularvalue legend,anchor=north] at (8.15,0.18) {%
            \begin{tabular*}{\dimexpr16.3cm-14pt\relax}{@{\extracolsep{\fill}}llllll@{}}
                \tikz[baseline=-0.5ex]\filldraw[fill=phaseDelocalized,
                    draw=phaseDelocalizedBorder] (0,0) rectangle (3.5mm,1.7mm);
                    \hspace{3pt}Delocalized &
                \tikz[baseline=-0.5ex]\filldraw[fill=phaseSemilocalized,
                    draw=phaseSemilocalizedBorder] (0,0) rectangle (3.5mm,1.7mm);
                    \hspace{3pt}Semilocalized &
                \tikz[baseline=-0.5ex]\filldraw[fill=black!18,
                    draw=black!55] (0,0) rectangle (3.5mm,1.7mm);
                    \hspace{3pt}Uncovered &
                \tikz[baseline=-0.5ex]\draw[singularvalue boundary,phaseLeft]
                    (0,0)--(4mm,0); \hspace{3pt}$\bcurveL(\lambda;\qr)$ &
                \tikz[baseline=-0.5ex]\draw[singularvalue boundary,phaseOrange]
                    (0,0)--(4mm,0); \hspace{3pt}$\bcurveRi{1}(\lambda;\qr)$ &
                \tikz[baseline=-0.5ex]\draw[singularvalue boundary,phaseRight]
                    (0,0)--(4mm,0); \hspace{3pt}$\bcurveRi{2}(\lambda;\qr)$
            \end{tabular*}%
        };

        \def\drawSingularValuePanel#1#2#3#4#5{%
            \pgfmathsetmacro{\singularvalueQ}{#1}
            \pgfmathsetmacro{\singularvalueLeft}{\singularvalueQ-1/\singularvalueQ}
            \pgfmathsetmacro{\singularvalueRight}{\singularvalueQ+1/\singularvalueQ}
            \pgfmathsetmacro{\singularvalueAMax}{\singularvalueQ*\singularvalueQ-1}
            \pgfmathsetmacro{\singularvalueBL}{singularvalueleftb(\singularvalueAMax,\singularvalueQ)}
            \pgfmathsetmacro{\singularvalueBRone}{singularvaluerone(1,\singularvalueQ)}
            \pgfmathsetmacro{\singularvalueBRtwo}{singularvaluertwo(1,\singularvalueQ)}
            \pgfmathsetmacro{\singularvalueIso}{1/(\singularvalueQ*\singularvalueQ)}
            \pgfmathsetmacro{\singularvalueFork}{1/(2*(\singularvalueQ*\singularvalueQ+1/(\singularvalueQ*\singularvalueQ)))}
            \pgfmathsetmacro{\singularvalueMiddleHigh}{max(\singularvalueBL,\singularvalueBRone)}
            \pgfmathsetmacro{\singularvalueMiddleLow}{min(\singularvalueBL,\singularvalueBRone)}
            \pgfmathsetmacro{\singularvalueStripLeft}{(#3-#2)/1.9066667}
            \pgfmathsetmacro{\singularvalueStripRight}{\singularvalueStripLeft+0.60/1.9066667}
            \pgfmathsetmacro{\singularvalueStripCenter}{(\singularvalueStripLeft+\singularvalueStripRight)/2}
            \ifdim\singularvalueBL pt<\singularvalueBRone pt
                \colorlet{singularvalueMiddleColor}{phaseOrange}
                \def\singularvalueMiddleFill{phaseOrange!18}
            \else
                \colorlet{singularvalueMiddleColor}{phaseLeft}
                \def\singularvalueMiddleFill{red!11}
            \fi
            \begin{scope}[shift={(#2,1.0)},x=1.9066667cm,y=1.38cm]
                \begin{scope}
                    \clip (0,0) rectangle (3,4.2);
                    \fill[phaseSemilocalized] (0,0)--(0,\singularvalueIso)
                        --plot[domain=0.001:\singularvalueAMax,samples=240,smooth,variable=\a]
                            ({singularvalueleftx(\a,\singularvalueQ)},{singularvalueleftb(\a,\singularvalueQ)})
                        --(\singularvalueLeft,0)--cycle;
                    \fill[phaseSemilocalized] (\singularvalueRight,0)
                        --plot[domain=1:7,samples=280,smooth,variable=\u]
                            ({singularvaluex(\u,\singularvalueQ)},{singularvaluertwo(\u,\singularvalueQ)})
                        --(3.1,0)--cycle;
                    \fill[phaseDelocalized]
                        ({\singularvalueLeft+0.075},0) rectangle ({\singularvalueRight-0.075},4.2);
                    \fill[black!18]
                        ({\singularvalueLeft-0.035},0) rectangle ({\singularvalueLeft+0.075},4.2);
                    \fill[black!18]
                        ({\singularvalueRight-0.075},0) rectangle ({\singularvalueRight+0.035},4.2);
                    \draw[singularvalue guide] (0,\singularvalueBL)--(\singularvalueLeft,\singularvalueBL);
                    \draw[singularvalue guide] (0,\singularvalueBRone)--(\singularvalueRight,\singularvalueBRone);
                    \draw[singularvalue guide] (0,\singularvalueBRtwo)--(\singularvalueRight,\singularvalueBRtwo);
                    \draw[black!45,line width=0.35pt]
                        (\singularvalueLeft,0)--(\singularvalueLeft,4.2)
                        (\singularvalueRight,0)--(\singularvalueRight,4.2);
                    \draw[singularvalue boundary,phaseLeft] (0,\singularvalueIso)
                        --plot[domain=0.001:\singularvalueAMax,samples=240,smooth,variable=\a]
                            ({singularvalueleftx(\a,\singularvalueQ)},{singularvalueleftb(\a,\singularvalueQ)});
                    \draw[singularvalue boundary,phaseRight]
                        plot[domain=1:7,samples=280,smooth,variable=\u]
                            ({singularvaluex(\u,\singularvalueQ)},{singularvaluertwo(\u,\singularvalueQ)});
                    \draw[singularvalue boundary,phaseOrange]
                        plot[domain=1:7,samples=280,smooth,variable=\u]
                            ({singularvaluex(\u,\singularvalueQ)},{singularvaluerone(\u,\singularvalueQ)});
                \end{scope}
                \draw[singularvalue axis] (0,0)--(3.123,0);
                \draw[singularvalue axis] (0,0)--(0,4.36);
                \node[font=\small,anchor=base] at (0,4.48) {$b$};
                \node[font=\small,anchor=base] at (1.5,4.48) {(#5) $\qr=#1$};
                \node[font=\footnotesize,anchor=north east,inner sep=1.5pt]
                    at (-0.025,-0.055) {$0$};
                \foreach \x/\label in {\singularvalueLeft/\ledge,\singularvalueRight/\redge}{
                    \draw (\x,0)--++(0,-0.04);
                    \node[font=\footnotesize,anchor=north] at (\x,-0.07) {$\label$};
                }
                \node[font=\small,anchor=north] at (3,-0.07) {$\lambda$};
                \foreach \y/\label/\singularvalueTickColor in {4.2/4.2/black,
                    \singularvalueFork/{\bcritFork}/phasePurple,
                    \singularvalueIso/{\bcritIso}/phaseTeal,
                    \singularvalueBL/{\bcritL}/phaseLeft,
                    \singularvalueBRone/{\bcritRi{1}}/phaseOrange,
                    \singularvalueBRtwo/{\bcritRi{2}}/phaseRight}{
                    \draw[\singularvalueTickColor] (0,\y)--++(-0.018,0);
                    \node[font=\footnotesize,anchor=east,text=\singularvalueTickColor]
                        at (-0.06,\y) {$\label$};
                }
                \draw[phaseTeal,line width=1.6pt] (0,0)--(0,\singularvalueIso);
                \draw[phasePurple,line width=1.6pt] (0,0)--(0,\singularvalueFork);
                \node[singularvalue point,fill=phaseTeal,minimum size=4.5pt,line width=0.35pt]
                    at (0,\singularvalueIso) {};
                \node[singularvalue point,fill=phasePurple,minimum size=4.5pt,line width=0.35pt]
                    at (0,\singularvalueFork) {};
                \node[singularvalue point,fill=phaseLeft] at (\singularvalueLeft,\singularvalueBL) {};
                \node[singularvalue point,fill=phaseOrange] at (\singularvalueRight,\singularvalueBRone) {};
                \node[singularvalue point,fill=phaseRight]
                    at (\singularvalueRight,\singularvalueBRtwo) {};

                \fill[phaseTeal!15] (\singularvalueStripLeft,0)
                    rectangle (\singularvalueStripRight,\singularvalueFork);
                \path[pattern=north east lines,pattern color=phasePurple]
                    (\singularvalueStripLeft,0) rectangle (\singularvalueStripRight,\singularvalueFork);
                \fill[phaseTeal!15] (\singularvalueStripLeft,\singularvalueFork)
                    rectangle (\singularvalueStripRight,\singularvalueIso);
                \fill[phasePeach] (\singularvalueStripLeft,\singularvalueIso)
                    rectangle (\singularvalueStripRight,\singularvalueMiddleLow);
                \fill[\singularvalueMiddleFill] (\singularvalueStripLeft,\singularvalueMiddleLow)
                    rectangle (\singularvalueStripRight,\singularvalueMiddleHigh);
                \fill[blue!12] (\singularvalueStripLeft,\singularvalueMiddleHigh)
                    rectangle (\singularvalueStripRight,\singularvalueBRtwo);
                \fill[black!5] (\singularvalueStripLeft,\singularvalueBRtwo)
                    rectangle (\singularvalueStripRight,4.2);
                \draw[black!45,line width=0.4pt]
                    (\singularvalueStripLeft,0) rectangle (\singularvalueStripRight,4.2);
                \foreach \y in {\singularvalueFork,\singularvalueIso,\singularvalueBL,\singularvalueBRone,\singularvalueBRtwo}{
                    \draw[black!45,line width=0.4pt] (\singularvalueStripLeft,\y)--(\singularvalueStripRight,\y);
                }
                \draw[singularvalueMiddleColor,line width=0.7pt]
                    (\singularvalueStripLeft,\singularvalueMiddleLow)
                    rectangle (\singularvalueStripRight,\singularvalueMiddleHigh);
                \node[singularvalue number]
                    at (\singularvalueStripCenter,{(4.2+\singularvalueBRtwo)/2}) {1};
                \node[singularvalue number,draw=phaseRight,text=phaseRight]
                    at (\singularvalueStripCenter,{(\singularvalueBRtwo+\singularvalueMiddleHigh)/2}) {2a};
                \node[singularvalue number,draw=singularvalueMiddleColor,text=singularvalueMiddleColor,
                    line width=0.8pt]
                    at (\singularvalueStripCenter,{(\singularvalueMiddleHigh+\singularvalueMiddleLow)/2}) {#4};
                \node[singularvalue number,draw=phaseOrange,text=phaseOrange!80!black]
                    at (\singularvalueStripCenter,{(\singularvalueMiddleLow+\singularvalueIso)/2}) {3b};
                \node[singularvalue number,draw=phaseTeal,text=phaseTeal]
                    at (\singularvalueStripCenter,{(\singularvalueIso+\singularvalueFork)/2}) {4};
                \node[singularvalue number,draw=phasePurple,text=phasePurple] (singularvalueFive#5)
                    at (\singularvalueStripCenter,-0.25) {5};
                \draw[phasePurple,line width=0.6pt] (singularvalueFive#5.north)
                    --(\singularvalueStripCenter,{\singularvalueFork/2});
                \fill[phasePurple] (\singularvalueStripCenter,{\singularvalueFork/2})
                    circle[radius=0.65pt];
            \end{scope}%
        }
        \drawSingularValuePanel{1.10}{1.2}{7.50}{2b}{i}
        \drawSingularValuePanel{1.25}{10.15}{8.30}{3a}{ii}
    \end{tikzpicture}%
    }
    \caption{Phase diagrams in the $(\lambda,b)$ plane for $\qr=1.10$ (left) and $\qr=1.25$ (right). The central strips use the region key of Figure~\ref{fig:combined-aspect-ratio-thresholds}. The left semilocalization theorem applies for $b>\bcritIso$. Grey strips show the uncovered exterior $C\error_d^{1/4}$-margins and interior $\kappa$-windows, enlarged and not to scale. Green vertical marks indicate rank deficiency from isolated vertices in $\gV_2$ for $b<\bcritIso$. Purple marks indicate the emergence of the positive singular values $d^{-1/2}$ induced by $1$-fork when $b\leq\bcritFork$.}
    \label{fig:edge-singular-vector-transition-two-aspects}
\end{figure}

For fixed $\qr\geq1$, define the right transition curves using
the rate function in \eqref{eqn:transition-rate-function}:
\begin{align}
    \bcurveRi{1}(\lambda;\qr)
    &\coloneqq\frac{1}{I_{\qr^{-1}}\bigl(\invR{\qr^{-1}}(\lambda)\bigr)},
    &
    \bcurveRi{2}(\lambda;\qr)
    &\coloneqq\frac{1}{I_{\qr}\bigl(\invR{\qr}(\lambda)\bigr)}
    \label{eqn:right-class-transition-boundaries}
\end{align}
where $\lambda>\redge$. The overall curve is
$\bcurveR(\lambda;\qr)=\max_{i\in\{1,2\}}\bcurveRi{i}(\lambda;\qr)
=\bcurveRi{2}(\lambda;\qr)$; equivalently, it solves
$\rhoR(\lambda,b)=0$ for $b$.
The corresponding thresholds for the emergence of right outliers induced by high-degree vertices in $\gV_1$ and $\gV_2$ are obtained by sending $\lambda$ to $\redge$:
\begin{subequations}\label{eqn:edge-transition-thresholds}
    \begin{align}
        \bcritRi{1}(\qr)
        &\coloneqq\lim_{\lambda\downarrow\redge}\bcurveRi{1}(\lambda;\qr)=[(\qr^{-2}+1)\log(1+\qr^2)-1]^{-1},
        \label{eqn:right_threshold_1}\\
        \bcritRi{2}(\qr)
        &\coloneqq\lim_{\lambda\downarrow\redge}\bcurveRi{2}(\lambda;\qr)=[(\qr^2+1)\log(1+\qr^{-2})-1]^{-1}.
        \label{eqn:right_threshold_2}
    \end{align}
\end{subequations}
Consequently, for any $\qr \geq 1$, the threshold for the first emergence of right outliers is
\begin{align}
    \bcritR(\qr)
    =\max\{\bcritRi{1}(\qr),\bcritRi{2}(\qr)\}
    =\bcritRi{2}(\qr).
    \label{eqn:transition-threshold-identification}
\end{align}
For $\qr>1$, decreasing $b$ first produces right outliers from
$\gV_2$, then from $\gV_1$; at $\qr=1$, the thresholds coincide.
For $\qr>1$, the left transition curve solves $\rhoL(\lambda,b)=0$
for $b$:
\begin{align}
    \bcurveL(\lambda;\qr) \coloneqq 1/\sI_{\mathrm L}(\lambda;\qr),
    \label{eqn:left-class-transition-boundary}
\end{align}
where $0<\lambda<\ledge$. The threshold for the emergence of left outliers is obtained by sending $\lambda$ to $\ledge$:
\begin{align}
    \bcritL(\qr)
    &\coloneqq\lim_{\lambda\uparrow\ledge}\bcurveL(\lambda;\qr) =[(\qr^2-1)\log(1-\qr^{-2})+1]^{-1},
    \label{eqn:left_threshold}
\end{align}
which corresponds to the emergence of low degrees in $\gV_2$.

Figure~\ref{fig:combined-aspect-ratio-thresholds} shows five selected thresholds: $\bcritRi{1}$, $\bcritRi{2}$, $\bcritL$, the isolation threshold $\bcritIso$, and the positive $1$-fork threshold $\bcritFork$. For $b<\bcritIso$, isolated vertices in $\gV_2$ produce a rank deficient centered matrix; the conditional centering estimate and the additional tree and fork thresholds are given in Subsection~\ref{subsec:intro-zero-modes}. At $b=\bcritIso$, the isolated-vertex count has the Poisson limit in
Lemma~\ref{lem:bipartite-isolated-column-threshold}, and
$|\sI_2|\geq2$ guarantees rank deficiency. The figure's fifth region also contains positive $1$-fork singular values $d^{-1/2}$,
including inside the giant with probability $1-o(1)$ when
$b\leq\bcritFork$.

The panels in Figure~\ref{fig:edge-singular-vector-transition-two-aspects} display two orderings of $\bcritRi{1}$ and $\bcritL$. Region~2b requires $\bcritL<b<\bcritRi{1}$, whereas region~3a requires $\bcritRi{1}<b<\bcritL$. Only one of these intervals can be nonempty at a fixed $\qr$, and both collapse when the thresholds $\bcritRi{1}$ and $\bcritL$ coincide.

\subsubsection{Transition for outlier statistics}
For a self-adjoint matrix $\rmT$, a centre $\lambda\in\R$, and a width
$\delta>0$, define
\begin{align}
    \mathrm{N}_{\lambda,\delta}(\rmT)
    \coloneqq
    \#\bigl(\operatorname{Spec}(\rmT)
    \cap[\lambda-\delta,\lambda+\delta]\bigr),
    \label{eqn:transition-eigenvalue-count}
\end{align}
with eigenvalues counted according to multiplicity. For $0<\delta<\lambda$,
$\mathrm{N}_{\lambda,\delta}(\scA)$ also counts the singular values of
$\mX$ in the same interval. The positive-count asymptotics below
concern fixed positive locations strictly inside the corresponding
outlier intervals and deterministic widths $\delta=\delta_N>0$
satisfying
\begin{align}
    \delta_N&\longrightarrow0,
    &{\error_d}/{\delta_N}&\longrightarrow0.
    \label{eqn:transition-spectral-window}
\end{align}

At critical sparsity, the spectral boundaries are determined by the extreme normalized degrees, which can be expressed in terms of the Lambert $W$ function below.

\begin{definition}[Lambert $W$ function]
\label{def:lambert-w}
The Lambert $W$ function is the solution to the equation $W(z)\mathrm e^{W(z)}=z$. It has two real branches. The principal branch $W_0:[-\mathrm e^{-1},\infty)\to[-1,\infty)$ is increasing, and the lower branch $W_{-1}:[-\mathrm e^{-1},0)\to(-\infty,-1]$ is decreasing. Both branches equal $-1$ at $z=-\mathrm e^{-1}$, while $W_0(0)=0$ and $W_{-1}(z)\to-\infty$ as $z\uparrow0$.
\end{definition}

For the right edge boundary, define the upper degree root using $W_0$ for fixed $\qr\geq1$ and $b>0$,
\begin{align}
    a_{\mathrm R}(b;\qr)
    &\coloneqq
    \qr^2\exp\!\left[
      1+W_0\!\left(\frac{1/(b\qr^2)-1}{\mathrm e}\right)
    \right].
    \label{eqn:right-explicit-degree-root}
\end{align}
It is the unique solution of $I_{\qr}(a)=b^{-1}$ with $a>\qr^2$.
At critical sparsity, the maximum normalized degree in $\gV_2$ converges
in probability to this root. Define the right spectral boundary by
\begin{align}
    \redge(b;\qr)
    &\coloneqq
    \Lambda_{\qr}\!\left(
      \max\{a_{\mathrm R}(b;\qr),\qr^2+1\}
    \right).
    \label{eqn:global-right-transition-boundary}
\end{align}
For $b\geq\bcritR(\qr)$, this boundary equals the bulk edge $\redge$.
For $0<b<\bcritR(\qr)$, the root exceeds $\qr^2+1$, and
$\redge(b;\qr)>\redge$ is the upper endpoint of the positive right
outlier interval $(\redge,\redge(b;\qr))$. The right-outlier count
aggregates the contributions from both vertex classes.

\begin{theorem}[Right-edge transition]
\label{thm:right-edge-transition}
Under Assumptions~\ref{ass:proportional_regime}
and~\ref{ass:critical_sparsity}, fix $b>0$ with
$b\neq\bcritR(\qr)$.
\begin{enumerate}[label=(\roman*)]
    \item \emph{Spectral boundary.} For every $\epsilon>0$,
    there are constants $\gC_\epsilon,c_\epsilon>0$, depending on
    $b$ and $\qr$, such that, for all sufficiently large $N$,
    \begin{align}
        \P\left(
          \left|\lambda_1(\scA)-\redge(b;\qr)\right|\leq\epsilon
        \right)
        \geq1-\gC_\epsilon N^{-c_\epsilon}.
        \label{eqn:rightmost-eigenvalue-transition}
    \end{align}
    \item \emph{Outlier counts.} If $b>\bcritR(\qr)$, then for every fixed $\lambda>\redge$
    and every deterministic sequence $\delta_N\downarrow0$,
    there are constants $\gC_\lambda,c_\lambda>0$ such that,
    for all sufficiently large $N$,
    \begin{align}
        \P\left(\mathrm{N}_{\lambda,\delta}(\scA)=0\right)
        \geq1-\gC_\lambda N^{-c_\lambda}.
        \label{eqn:right-high-eigenvalue-absence}
    \end{align}

    If $b<\bcritR(\qr)$, fix
    $\lambda\in(\redge,\redge(b;\qr))$ and let $\delta_N>0$
    satisfy \eqref{eqn:transition-spectral-window}.
    For every $\nu>0$, there exist a constant
    $\gC_{\lambda,\nu}>0$ and a deterministic sequence
    $\eta_N\downarrow0$ such that, with probability at least
    $1-\gC_{\lambda,\nu}N^{-\nu}$ for all sufficiently large $N$,
    \begin{align}
        N^{\rhoR(\lambda,b)-\eta_N}
        \leq\mathrm{N}_{\lambda,\delta}(\scA)
        \leq N^{\rhoR(\lambda,b)+\eta_N}.
        \label{eqn:right-transition-eigenvalue-count}
    \end{align}
\end{enumerate}
\end{theorem}
We now turn attention to the left-edge transition. For fixed $\qr>1$,
define the lower degree root using the branch $W_{-1}$ from
Definition~\ref{def:lambert-w} by
\begin{align}
    a_{\mathrm L}(b;\qr)
    &\coloneqq
    \begin{cases}
      0,&0<b\leq\bcritIso(\qr),\\
      \displaystyle
      \qr^2\exp\!\left[
        1+W_{-1}\!\left(\frac{1/(b\qr^2)-1}{\mathrm e}\right)
      \right],&b>\bcritIso(\qr).
    \end{cases}
    \label{eqn:left-explicit-degree-root}
\end{align}
For $b>\bcritIso(\qr)$, it is the unique solution of
$I_{\qr}(a)=b^{-1}$ in $(0,\qr^2)$. For every $b>0$, the minimum
normalized degree in $\gV_2$ converges in probability to this root at
critical sparsity. Extending the lower location map
continuously by $\Lambda_{\qr}(0)=0$, define the left spectral boundary by
\begin{align}
    \ledge(b;\qr)
    &\coloneqq
    \Lambda_{\qr}\!\left(
      \min\{a_{\mathrm L}(b;\qr),\qr^2-1\}
    \right).
    \label{eqn:global-left-transition-boundary}
\end{align}
This boundary equals $\ledge$ for $b\geq\bcritL(\qr)$, lies in
$(0,\ledge)$ for $\bcritIso(\qr)<b<\bcritL(\qr)$, and equals zero
for $0<b\leq\bcritIso(\qr)$. The positive left-outlier interval is
$(\ledge(b;\qr),\ledge)$.
When the lower endpoint is zero, the theorem concerns each fixed
positive location in $(0,\ledge)$;
it does not locate the least positive singular value. Specific
near-zero families and their multiplicity bounds are given in
Subsection~\ref{subsec:intro-zero-modes}.

\begin{theorem}[Left-edge transition]
\label{thm:left-edge-transition}
Under Assumptions~\ref{ass:proportional_regime}
and~\ref{ass:critical_sparsity}, fix $\qr>1$ and $b>0$ with
$b\neq\bcritL(\qr)$.
\begin{enumerate}[label=(\roman*)]
    \item \emph{Spectral absence below the boundary.} For every compact set
    $K\subset(0,\ledge(b;\qr))$, there are constants
    $\gC_K,c_{\mathrm L,K}>0$ such that $K$ is free of
    eigenvalues with probability at least
    $1-\gC_KN^{-c_{\mathrm L,K}}$ for all sufficiently large $N$.
    This assertion is vacuous when $\ledge(b;\qr)=0$.

    \item \emph{Outlier counts.} If $b>\bcritL(\qr)$, then for every fixed
    $\lambda\in(0,\ledge)$ and every deterministic sequence
    $\delta_N\downarrow0$, there are constants
    $\gC_\lambda,c_{\mathrm L,\lambda}>0$ such that, for all
    sufficiently large $N$,
    \begin{align}
        \P\left(\mathrm{N}_{\lambda,\delta}(\scA)=0\right)
        \geq1-\gC_\lambda N^{-c_{\mathrm L,\lambda}}.
        \label{eqn:left-high-eigenvalue-absence}
    \end{align}

    If $0<b<\bcritL(\qr)$, fix
    $\lambda\in(\ledge(b;\qr),\ledge)$ and let $\delta_N>0$
    satisfy \eqref{eqn:transition-spectral-window}.
    There are constants $\gC_\lambda,c_{\mathrm L,\lambda}>0$
    and a deterministic sequence $\eta_N\downarrow0$ such that,
    with probability at least $1-\gC_\lambda N^{-c_{\mathrm L,\lambda}}$, the following holds for all sufficiently large $N$,
    \begin{align}
        N^{\rhoL(\lambda,b)-\eta_N}
        \leq\mathrm{N}_{\lambda,\delta}(\scA)
        \leq N^{\rhoL(\lambda,b)+\eta_N}.
        \label{eqn:left-transition-eigenvalue-count}
    \end{align}
\end{enumerate}
\end{theorem}

The explicit degree roots, the boundary identities, and the two
transition theorems are proved in
Section~\ref{sec:transition-from-delocalization-to-localization}.

\subsection{Numerical simulations}\label{sec:numerical-results}
\begin{figure}[h]
    \centering
    \includegraphics[width=\textwidth]{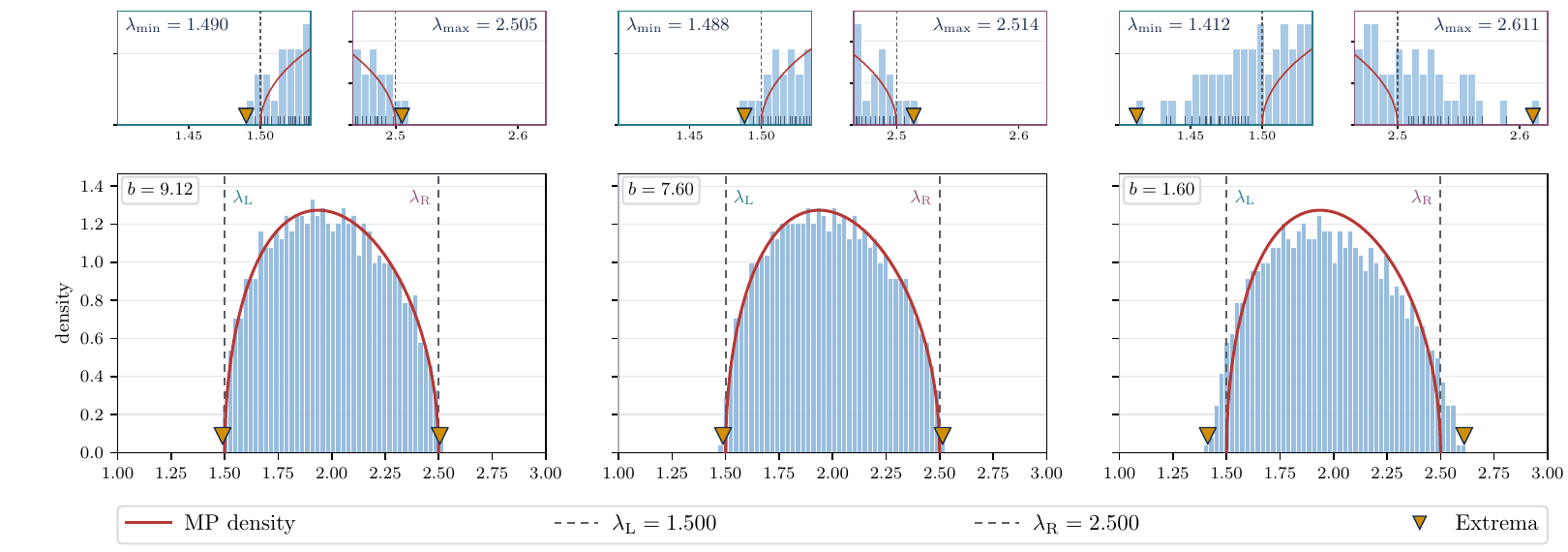}
    \caption{From left to right, we choose the sparsity parameter
    $b\in\{9.12,7.60,1.60\}$. The red curve is the MP density in \eqref{eqn:MP-singular-value-edges} with $\qr=2$, and dashed lines mark the MP edges $\ledge=3/2$ and $\redge=5/2$. Enlarged triangles mark the observed extrema. All axes are shared across the three panels, with $1\leq\lambda\leq3$.}
    \label{fig:bipartite-outliers-ratio16}
\end{figure}

\begin{figure}[h]
    \centering
    \includegraphics[width=\textwidth]{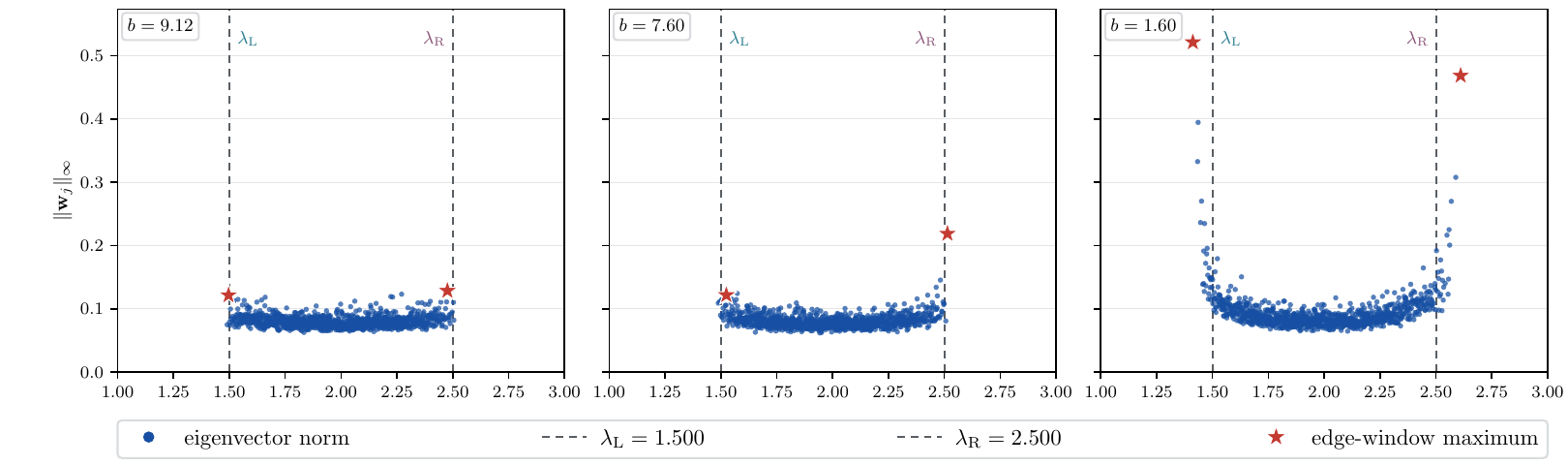}
    \caption{The pairs $(\lambda_j,\|\rvw_j\|_\infty)$ for the same realizations with $\ratio = 16$ as in Figure~\ref{fig:bipartite-outliers-ratio16}.  The red stars mark the maximum $\ell^\infty$ norm separately in the common left and right edge windows.}
    \label{fig:eigenvector-norm-transition-ratio16}
\end{figure}

We simulate $\ratio=16$, $m=1000$, and seed $777$. Here $\qr=2$,
the MP edges are $3/2$ and $5/2$, and the right- and left-outlier
thresholds are approximately $8.64$ and $7.30$. The choices
$b\in\{9.12,7.60,1.60\}$ sample the no-outlier, right-only, and
two-sided regimes. Figures~\ref{fig:bipartite-outliers-ratio16}
and~\ref{fig:eigenvector-norm-transition-ratio16} compare the
spectra and eigenvector infinity norms across these regimes.

\FloatBarrier
\subsection{Related literature}\label{sec:literature}

\paragraph{Classical random matrices}
Marchenko and Pastur \cite{marchenko1967distribution} identified the limiting empirical eigenvalue distribution of sample covariance matrices in the proportional regime. Yin, Bai, and Krishnaiah \cite{yin1988limit} and Bai and Yin \cite{bai1993limit} established the corresponding almost-sure limits for normalized largest and smallest singular values of matrices with centered i.i.d.\ entries and finite fourth moments. For the normalized smallest singular value, Tikhomirov \cite{tikhomirov2015limit} obtained almost-sure convergence for real centered i.i.d.\ entries of finite nonzero variance when the limiting aspect ratio is strictly greater than one.
Tracy and Widom \cite{tracy1994level,tracy1996orthogonal} identified
limiting largest-eigenvalue distributions for Gaussian ensembles,
and Pillai and Yin \cite{pillai2014universality} proved soft-edge
universality for real and complex sample covariance matrices with
independent entries and subexponential tails. Their lower-edge
result excludes limiting aspect ratio one.
At the hard edge, Tao and Vu \cite{tao2010smallest} established
universality for the smallest singular value of square matrices and
matrices with fixed dimension difference under sufficiently high
finite moment assumptions on independent entries.
These laws provide benchmarks for the degree-induced outliers studied
here; see \cite{yau2026randomwigner} for broader universality background.

\paragraph{Local laws and delocalization in classical random matrices.}
For Hermitian Wigner matrices, Erd\H{o}s, Schlein, and Yau established
short-scale semicircle laws and eigenvector delocalization
\cite{erdos2009shortscales}, followed by complete bulk delocalization
\cite{erdos2009localsemicircle} and bulk
local laws at optimal spectral resolution, Wegner estimates, and level
repulsion \cite{erdos2010wegner}.
Erd\H{o}s, Yau, and Yin \cite{erdos2012rigidity} obtained a strong
local semicircle law and eigenvalue rigidity uniformly through the
spectral edges for generalized Wigner matrices, with edge-universality
consequences under second-moment matching.
Knowles and Yin \cite{knowles2013isotropic} proved the isotropic
semicircle law, and Bloemendal and coauthors
\cite{bloemendal2014isotropic} established isotropic laws for sample
covariance and generalized Wigner matrices. These laws control
resolvent quadratic forms in deterministic directions.
Vu and Wang \cite{vu2015projections} sharpened entrywise
delocalization bounds for bulk eigenvectors of Bernoulli
Wigner matrices; see \cite{orourke2016survey} for a survey.
Knowles and Yin \cite{knowles2017anisotropic} developed anisotropic
local laws for covariance and deformed Wigner matrices.
He, Knowles, and Rosenthal \cite{he2018isotropic} derived isotropic
self-consistent equations for mean-field models with comparable
variances and arbitrary expectation.
Ajanki, Erd\H{o}s, and Kr\"uger \cite{ajanki2019stability} proved
bulk local laws, rigidity, and delocalization for matrices with
rapidly decaying correlations through stability analysis of the
matrix Dyson equation; Erd\H{o}s \cite{erdos2019dyson} surveys the
vector and matrix Dyson equation framework.
Benigni and Lopatto \cite{benigni2022optimal} obtained optimal
isotropic bounds of order $(\log N)/N$ for squared overlaps in
generalized Wigner matrices with subexponential entries.
Osman \cite{osman2026precise} proved precise delocalization tails
and Gumbel laws for maximal eigenvector entries in Wigner matrices.
For rectangular matrices, Alt, Erd\H{o}s, and Kr\"uger
\cite{alt2017gram} established optimal entrywise and averaged bulk
local laws for Gram matrices with independent centered entries and
general variance profiles.

\paragraph{Local law and delocalization in sparse random matrices.}
Tran, Vu, and Wang \cite{tran2013sparse} obtained eigenvalue-count
estimates for sparse random graphs and entrywise eigenvector bounds
for Erd\H{o}s--R\'enyi graphs.
Erd\H{o}s, Knowles, Yau, and Yin \cite{erdos2013spectral} established
averaged and entrywise local semicircle laws down to nearly
microscopic scales under logarithmic-power degree assumptions.
He, Knowles, and Marcozzi \cite{he2019local} reduced the expected-degree
requirement to $d\geq C\log N$, proving a local law and complete
entrywise delocalization. Dumitriu and Zhu
\cite{dumitriu2019sparse} treated nonconstant variance profiles in the
bulk when $d/\log N\to\infty$.
He, Wang, and Zhu \cite{he2024sparse} extended logarithmic-degree
local laws and delocalization to sparse sub-Gaussian weights, and
obtained corresponding results under larger polylogarithmic degree
requirements for sparse $\alpha$-subexponential weights with
$0<\alpha<2$.
At critical sparsity, the local law of Alt, Ducatez, and Knowles
\cite{alt2021delocalization} approximates diagonal resolvent entries
by degree-dependent functions in a bulk region separated from zero
and the spectral edges.
For sparse sample covariance matrices, Hwang, Lee, and Schnelli
\cite{hwang2019local} proved local laws up to the upper spectral edge
at polynomial expected degrees.
For aspect ratios bounded away from one, Tikhomirov and Timushev
\cite{tikhomirov2022local} reached polylogarithmic expected degrees
under finite $(4+\delta)$-moment and truncation assumptions, with
$\delta>0$. G\"otze, Tikhomirov, and Timushev
\cite{gotze2024locallaw} removed the truncation assumption for bulk
spectral parameters under the same degree and moment requirements.
Isotropic local laws and delocalization at polynomial degrees were
obtained by He, Huang, and Wang \cite{he2026extremal}.
For $C(\log N)^2\leq d\leq N^{1-\epsilon}$, with fixed $\epsilon>0$,
Henheik and Knowles \cite{henheik2026bernoulli} used the Bernoulli
flow to prove optimal isotropic bulk bounds of order $(\log N)/N$
for squared overlaps with deterministic unit vectors.

\paragraph{Universality and fluctuations for sparse random matrices.}
Bulk universality for sparse matrices was established at polynomial
expected degrees \cite{erdos2012spectralii,huang2015bulk} and,
using the Bernoulli flow, for
$(\log N)^{2+\epsilon}\leq d\leq N^{1-\epsilon}$ with fixed
$\epsilon>0$ \cite{henheik2026bernoulli}.
He \cite{he2020bulkfluctuations} identified collective Gaussian
fluctuations of bulk eigenvalues away from zero and the spectral edges.
At the edge, Lee and Schnelli \cite{lee2018local} proved Tracy--Widom
fluctuations after a deterministic correction to the spectral edge.
Huang, Landon, and Yau \cite{huang2020transition} identified the
crossover to Gaussian fluctuations at $p\asymp N^{-2/3}$, and
He and Knowles \cite{he2021extremefluctuations} proved joint Gaussian
fluctuations of all nontrivial eigenvalues away from zero for
$N^{\epsilon}\leq Np\leq N^{1/3-\epsilon}$.
Lee \cite{lee2024higher} constructed higher-order random edge
corrections, and Huang and Yau \cite{huang2026edge} established
Tracy--Widom fluctuations for centered sparse matrices after an
explicit random shift at polynomial degrees.
For sparse sample covariance matrices, Hwang, Lee, and Schnelli
\cite{hwang2019local} obtained upper-edge Tracy--Widom laws after
a deterministic edge correction when $p\geq N^{-2/3+\epsilon}$
for fixed $\epsilon>0$.
Lee and Yeo \cite{lee2025local} proved local laws, eigenvalue rigidity,
and asymptotic normality of extremal eigenvalues for a class of
sparse matrices with random or deterministic diagonal deformations
independent of the sparse part.

\paragraph{Degree-induced outliers, semilocalization, and mobility edges.}
Spectral-norm estimates for bounded symmetric random matrices were
established by F\"uredi and Koml\'os \cite{furedi1981eigenvalues}
and refined by Vu \cite{vu2007spectral}.
For Erd\H{o}s--R\'enyi graphs, Benaych-Georges, Bordenave, and Knowles
identified a crossover at logarithmic mean degree:
below this scale, the extreme eigenvalues of the centered adjacency
matrix are governed by the largest degrees \cite{benaych2019largest};
above it, the normalized nontrivial spectral edges converge to the
semicircle edges \cite{benaych2020spectral}.
At criticality, Alt, Ducatez, and Knowles established a correspondence
between high-degree vertices and nontrivial outliers
\cite{alt2021extremal}, and proved bulk delocalization away from zero
together with semilocalization around resonant vertices near the
spectral edges \cite{alt2021delocalization}.
They subsequently determined the thresholds for complete
delocalization \cite{alt2022completely}, proved Poisson edge statistics
and exponential localization \cite{alt2023poisson}, and established
a localized phase reaching the mobility edge in suitable degree
regimes \cite{alt2023localized}.
For inhomogeneous graphs with prescribed vertex weights,
Buc-d'Alch\'e and Knowles \cite{buc2026semilocalization} proved
semilocalization near the spectral edges and localization around
a single vertex for extreme eigenvalues under empirical moment
conditions on the weights.

\paragraph{Sparse rectangular matrices.}
For balanced bipartite Erd\H{o}s--R\'enyi graphs, Basak and Rudelson
\cite{basak2021sharp} identified zero rows and columns as the
asymptotic obstruction to invertibility.
G\"otze and Tikhomirov \cite{gotze2022largest} obtained
extreme-singular-value bounds for sparse rectangular matrices with
fixed limiting aspect ratio greater than one and expected degree
$\asymp(\log N)^\alpha$, $\alpha>1$, under moment and truncation
conditions. Dumitriu and Zhu
\cite{dumitriu2024extreme} derived two-sided bounds for inhomogeneous
rectangular matrices using the nonbacktracking operator; for
bipartite Erd\H{o}s--R\'enyi graphs, these exclude outliers
when the expected degree grows faster than logarithmically.
Nadutkina, Tikhomirov, and Timushev \cite{nadutkina2024marchenko}
proved a global symmetrized Marchenko--Pastur law for weighted
bipartite graphs.
At logarithmic degree, Dumitriu, Wang, Wang, and Zhu
\cite{dumitriu2025singular} characterized degree-induced outlier
locations and counts. Their left-edge results require connectivity
and an additional positive-definiteness condition
\cite[Assumptions~1.5--1.6]{dumitriu2025singular}, which is removed in the present work.

\paragraph{Random regular and biregular graphs.}
Random regular and biregular graphs provide comparisons in which
the prescribed degrees suppress the fluctuations responsible for
degree-induced outliers.
Early second-eigenvalue estimates \cite{friedman1989second} were
followed by Bordenave's new proof of Friedman's sharp fixed-degree
bound \cite{bordenave2020new} and spectral-gap estimates for growing
degrees \cite{cook2018size,tikhomirov2019spectral,sarid2023spectral}.
Bauerschmidt, Huang, Knowles, and Yau \cite{bauerschmidt2017bulk}
proved GOE bulk universality for
$N^\alpha\leq d\leq N^{2/3-\alpha}$, with fixed $\alpha>0$.
Bourgade and Huang \cite{bourgade2026loop} used loop equations to
extend bulk universality to $(\log N)^{24}\ll d\leq N^{1/2}$.
Edge universality was established for dense random regular graphs
by He \cite{he2024spectral} and for every fixed degree $d\geq3$
by Huang, McKenzie, and Yau \cite{huang2024ramanujan}.
For bipartite biregular graphs with growing degrees, Dumitriu and Johnson
\cite{dumitriu2016marvcenko} proved global and local symmetrized
Marchenko--Pastur laws, and Tran \cite{tran2020local} obtained a
local law for diverging degrees.
Spectral-gap estimates were established at fixed degrees by Brito,
Dumitriu, and Harris \cite{brito2022spectral} and in growing-degree
regimes by Zhu \cite{zhu2022second}.

\paragraph{Other models and transition results.}
For non-Hermitian matrices, Rudelson and Vershynin
\cite{rudelson2015delocalization} proved entrywise eigenvector
delocalization for independent centered
entries with uniform lower variance bounds and exponential tail control. For i.i.d.\ entries
with uniformly bounded moments of all orders, Cipolloni and Landon
\cite{cipolloni2025optimal} obtained the optimal isotropic scale
$\sqrt{(\log N)/N}$ for both left and right eigenvectors.
Local circular laws treat comparable variance profiles in the bulk
\cite{alt2018inhomogeneous} and at the spectral edge
\cite{alt2021spectralradius}, as well as decaying correlations
\cite{alt2021correlatedcircular}.
For complex i.i.d.\ matrices, extremal eigenvalues have universal
Gumbel and Poisson statistics \cite{cipolloni2023extremal}; edge
universality also holds at typical edge points of deterministic
deformations \cite{campbell2025spectraledge}.
In real elliptic Ginibre matrices, inverse participation ratios
depend on the eigenvalue's distance from the real axis
\cite{benigni2026delocalization}. Directed inhomogeneous graphs have
inhomogeneous spectral laws and finite-rank outlier fluctuations under
logarithmic degree-growth conditions \cite{hazra2026spectrum}.

For a class of one-dimensional band matrices, bulk universality was established
at $W\asymp N$ \cite{bourgade2017banduniversality} and, together with
delocalization and quantum unique ergodicity, at
$W\geq N^{3/4+\epsilon}$ \cite{bourgade2020random}, using
generalized-resolvent and fluctuation-averaging estimates
\cite{bourgade2019random,yang2021random}.
The flow method gave quantum diffusion and bulk delocalization for
Gaussian bands with $W\gg N^{8/11}$ \cite{dubova2026diffusion}.
Subsequent work reached $W\geq N^{1/2+\epsilon}$ for Gaussian block
bands \cite{yau2025delocalization} and then for more general real and
complex band models \cite{erdos2025zigzag}.
For block bands with regular entry densities, exponential
localization occurs on the sharp scale $W^2$ when $W^2\ll N$
\cite{drogin2025localization}.
In higher dimensions, renormalized resolvent estimates and $T$-expansions
gave delocalization of most bulk eigenvectors and quantum diffusion for Gaussian bands in
dimensions at least eight
\cite{yang2021delocalization,yang2022delocalization}; later work
obtained bulk delocalization, quantum unique ergodicity, and
universality at polynomial bandwidths in dimensions two and higher
\cite{dubova2025delocalization,dubova2025nonmeanfield}.
Under suitable coupling conditions, Gaussian block Schr\"odinger
operators exhibit bulk delocalization in dimensions at least seven
\cite{yang2025blockschrodinger} and in dimension two, with
one-dimensional lower bounds on bulk localization lengths
\cite{truong2025localization}.
Gaussian band-matrix delocalization in dimensions one and two has also been extended toward
and through the spectral edges under suitable bandwidth conditions
\cite{yang2025delocalization}; coupled Wigner blocks exhibit
transitions whose coupling threshold depends on the spectral
parameter \cite{fan2025localization}.

For the Sachdev--Ye--Kitaev model, Benigni and Cipolloni
\cite{benigni2026weak} proved mesoscopic local laws for even interaction
orders much smaller than the square root of the fermion count,
and bulk eigenvector delocalization for fixed interaction order.
He \cite{he2026spectral} determined the asymptotic spectral edge
for fixed even interaction order at least four.

Mobility edges have been established for L\'evy matrices with
tail index sufficiently close to zero or one from below
\cite{aggarwal2022mobility}, Anderson models with unbounded potentials
on Bethe lattices of sufficiently large degree
\cite{aggarwal2025bethe}, and Rosenzweig--Porter-type models
\cite{cipolloni2026rosenzweig}.
For rectangular matrices with i.i.d.\ power-law entries and limiting aspect
ratio greater than one, bottom right singular vectors undergo a
localization--delocalization transition as the tail exponent crosses
two \cite{bao2024phasebs}.

\paragraph{Matrix concentration inequalities for sum of random matrices}
Tropp established matrix concentration inequalities for sums of independent random self-adjoint matrices, providing noncommutative analogues of classical results through the matrix Laplace-transform method \cite{tropp2012user,tropp2015introduction}. Oliveira \cite{oliveira2009concentration} and Tropp \cite{tropp2011freedman} established matrix Freedman inequalities for matrix-valued martingales. Moreover, Tropp \cite{tropp2018second} refined matrix Khintchine inequality by using \emph{matrix alignment parameter} to reduce the dimensional dependence. More recently, Tropp obtained one-sided Gaussian comparisons for extreme eigenvalues \cite{tropp2025comparison,tropp2026extreme} and an exchangeable-counterpart proof of Brailovskaya and van Handel's universality theorems \cite{tropp2026exchangeable}.
Bandeira, Boedihardjo, and van Handel \cite{bandeira2021matrix} established nonasymptotic upper bounds for the spectra of general Gaussian random matrices, using asymptotic freeness \cite{Voi91} and strong convergence \cite{HT05,And13} to account for matrix noncommutativity. Corresponding lower bounds for the spectral edges were subsequently established in \cite{bandeira2024matrix}, and both were recently sharpened in \cite{bandeira2026matrix}. For sum of random matrices with non-Gaussian entries, \cite{brailovskaya2022universality} proved concentration inequalities by comparing them to sums of independent Gaussian matrices having the same means and covariances. Exponential concentration inequalities for sums of dependent random matrices have been established using an alternative exchangeable-pair method \cite{MJCFT14}, while \cite{SvW23} demonstrates Gaussian universality for sums arising from $\psi$-mixing Markov chains. See the recent reviews \cite{bordenave2026sparse, tropp2026applied, van2026strong} for more details.

\subsection{Notation}\label{sec:notations}
Throughout the paper, $[k]\coloneqq\{1,\ldots,k\}$ for
$k\in\N$, $|S|$ denotes the cardinality of a finite set $S$, and
$S^{\complement}$ denotes its complement. We use $\C_+ = \{z\in\C:\Im z>0\}$ for the upper half-plane. All asymptotic statements are understood as $N\to\infty$ unless another limit is specified. The symbols $c,C>0$ denote constants that may change from line to line, where a subscript records the parameters on which a constant may depend. We use $a\vee b\coloneqq\max\{a,b\}$ and $a\wedge b\coloneqq\min\{a,b\}$ for the maximum and minimum of two quantities $a$ and $b$, respectively. For nonnegative quantities $a$ and $b$, we write $a\lesssim_{\theta}b$, equivalently $a=O_\theta(b)$, if $a\leq C_{\theta}b$, $a\asymp_{\theta}b$ if both $a\lesssim_{\theta}b$ and $b\lesssim_{\theta}a$, and $a\ll b$ if $a/b\to0$.

For a vector $\vv$, let $\|\vv\|_p$ denote its $\ell^p$-norm. For a matrix $\mH$, $\|\mH\|$ denotes its operator norm, while $\mH^*$ and $\mH^{\sT}$ denote its adjoint and transpose, respectively.
For a finite set $S$, $\ones_S$ denotes its all-ones vector.  When $S\subset\gV$, we regard $\ones_S$ as a vector in $\R^{\gV}$ by extending it by zero outside $S$, and write $\ones_x\coloneqq\ones_{\{x\}}$ for a singleton.

We use $\indi{A}$ for the indicator of an event $A$, which is $1$ if $A$ occurs and $0$ otherwise. Denote probability and expectation by $\P$ and $\E$. We use $\rX\sim\Ber(p)$, $\rY\sim\Bin(k,p)$, and $\rZ\sim\Pois(\lambda)$ for Bernoulli, binomial, and Poisson random variables, respectively.

\begin{definition}[Overwhelming probability]
\label{def:overwhelming-probability}
An event $A=A_N$ holds with overwhelming probability if, for every $\nu>0$, there is a constant $\gC_\nu>0$ depending on $\nu$ such that $\P(A^{\complement})\leq\gC_\nu N^{-\nu}$ for all sufficiently large $N$.
\end{definition}

\subsection{Organization}
Section~\ref{sec:delocalized_phase} proves the bulk local law and
delocalization, with technical proofs in
Section~\ref{sec:delocalized-deferred-proofs}.
Section~\ref{sec:semilocalized_phase} establishes the edge profile
decomposition and semilocalization bounds; its deferred proofs appear in
Section~\ref{sec:semilocalized-deferred-proofs}.
Sections~\ref{sec:transition-from-delocalization-to-localization}
and~\ref{sec:transition-deferred-proofs} prove the transition theorems.
The appendices contain the tree calculations, the required random-graph
estimates, and the spectral comparison with the giant component.

%% file: delocalized.tex
\section{The Delocalized Phase}\label{sec:delocalized_phase}

We prove Theorem~\ref{thm:eigenvector_delocalization} using an
entrywise local law for sparse rectangular matrices. We first state
the local law and deduce delocalization, then develop its proof through
typical vertices, a controlled Schur equation, and a continuity
bootstrap. Auxiliary estimates are provided in Section~\ref{sec:delocalized-deferred-proofs}.

\subsection{Local law and bulk eigenvector delocalization}
\label{sec:local-law-and-delocalization}

\begin{definition}[Sparse rectangular noise]
\label{def:sparse-random-rectangular-matrix}
A random matrix $\mX\in\R^{n\times m}$ with independent entries has
sparsity $0<d<N$ if each entry $\emX_{ij}$, $i\in[n]$, $j\in[m]$,
satisfies $\E\emX_{ij}=0$ and
\begin{align}
    \max_{i,j}|\emX_{ij}|\leq Cd^{-1/2}\,\,\text{almost surely}, \quad
    \max_{i,j}\left|\sqrt{mn}\,\E|\emX_{ij}|^2-1\right| \leq\frac{Cd}{N},\quad
    \E|\emX_{ij}|^k\leq \frac{C}{Nd^{(k-2)/2}} (k\geq 2),
    \label{eqn:sparse-matrix-variance-defect}
\end{align}
where $C>0$ is a fixed constant independent of $n$, $m$, $d$, and $k$.
\end{definition}

We write $\evsigma_x=j$ for $x\in\gV_j$, $j=1,2$, and let $\mH$
be the linearization of $\mX$ as in \eqref{eqn:scAdj}. For a
deterministic signal $\signal\geq0$, define
\begin{align}
    \emM_{xy}\coloneqq\emH_{xy}
      +\frac{2\signal}{N}\indi{\evsigma_x\neq\evsigma_y},
    \qquad x,y\in[N].
    \label{eqn:two-spike-model}
\end{align}
The rescaled adjacency matrix $\mM=\rmA/\sqrt d$ is obtained by
taking $\signal=Np/(2\sqrt d)=O_{\qr}(\sqrt d)$ and
$\mX=(\widetilde{\mA}-\E\widetilde{\mA})/\sqrt d$.
For $z\in\C_+$, define the \emph{resolvent} of $\mM$ by
\begin{align}
    \res(z)\coloneqq(\mM-z)^{-1}.
    \label{eqn:resolvent}
\end{align}

Set $c_1\coloneqq\qr^{-2}$ and $c_2\coloneqq\qr^2$. Define the
class profiles as the solutions with positive imaginary part and
asymptotic behavior $\stj^{(j)}(z)\sim-z^{-1}$ of
\begin{align}
    \stj^{(j)}(z)=-\frac{1}{z+c_j\stj^{(3-j)}(z)}
    \label{eqn:stieltjes-m1-m2-definition}
\end{align}
for $j=1,2$. Solving this coupled system gives
\begin{subequations}
\begin{align}
    \stj^{(1)}(z)&=
    \frac{1}{2z\qr^2} \Big(-z^2-\qr^2+\qr^{-2}
      +\sqrt{(z^2-\redge^2)(z^2-\ledge^2)}\,\, \Big),
    \label{eqn:stieltjes-V1}\\
    \stj^{(2)}(z)&=
    \frac{\qr^{2}}{2z}\Big(-z^2+\qr^2-\qr^{-2}
      +\sqrt{(z^2-\redge^2)(z^2-\ledge^2)} \,\, \Big),
    \label{eqn:stieltjes-V2}
\end{align}
\end{subequations}
where the square root is the holomorphic branch in $\C_+$
asymptotic to $z^2$ at infinity. The scalar transform in
Lemma~\ref{lem:stieltjes-smp} is their class-weighted average
$\stj=(\ratio\stj^{(1)}+\stj^{(2)})/(1+\ratio)$, describing the limiting ESD.

Individual diagonal resolvent entries also depend on the realized
squared row norm
$\zeta_x\coloneqq\sum_{y\in[N]}|\emH_{xy}|^2$. The variance
normalization and Assumption~\ref{ass:proportional_regime} give
\begin{align}
    \E\zeta_x=c_j+O_{\qr}(d/N),\qquad x\in\gV_j.
    \label{eqn:zeta-class-means}
\end{align}
In this sparsity regime, these squared row norms need not be uniformly
close to their class means. We therefore retain their realized values
through the comparison functions
\begin{align}
    \stj^{(j)}_\alpha(z)\coloneqq
    -\frac{1}{z+\alpha\stj^{(3-j)}(z)},
    \qquad \alpha\geq0,\quad j=1,2.
    \label{eqn:stieltjes-alpha-inverse}
\end{align}
Thus $\stj^{(j)}_{c_j}=\stj^{(j)}$, and the local law compares
$\eres_{xx}$ with $\stj^{(\evsigma_x)}_{\zeta_x}$.

\begin{theorem}[Local law]\label{thm:bipartite-sparse-local-law}
Fix $0<\kappa<\qr^{-1}$ and $L\geq1$. Let $\mM$ be given by
\eqref{eqn:two-spike-model}, with $0\leq\signal\leq N^{\kappa/6}$ and
$\mX$ satisfying Definition~\ref{def:sparse-random-rectangular-matrix}.
Under Assumptions~\ref{ass:proportional_regime}
and~\ref{ass:subcritical_regime}, with overwhelming probability,
\begin{align}
    \sup_{z\in\specdomain}\max_{x,y\in[N]}
    \left|\eres_{xy}(z)
      -\indi{x=y}\stj_{\zeta_x}^{(\evsigma_x)}(z)\right|
    \lesssim\left(\frac{\log N}{d^2}\right)^{1/3}.
    \label{eqn:sparse-entrywise-local-law}
\end{align}
\end{theorem}

The comparison functions are uniformly bounded in the bulk by
Lemma~\ref{lem:stieltjes-bound}, even for atypical squared row norms.
The local law therefore yields bulk eigenvector delocalization.

\begin{proof}[Proof of Theorem~\ref{thm:eigenvector_delocalization}]
Apply Theorem~\ref{thm:bipartite-sparse-local-law} to $\mM=\mA/\sqrt{d}$ with $\kappa/2$ in place of $\kappa$ and $L=1$. For any normalized eigenpair
$(\lambda,\rvw)$ with $\lambda\in\bulk$, put
$\eta=N^{-1+\kappa/2}$ and $z=\lambda+\ii\eta$.
This point lies in the enlarged spectral domain. The local law,
Lemma~\ref{lem:stieltjes-bound}, and
$(\log N/d^2)^{1/3}=o(1)$ give
\begin{align}
    \frac{|\rvw(x)|^2}{\eta}
    \leq\Im\eres_{xx}(\lambda+\ii\eta)\leq C,
    \qquad x\in[N],
\end{align}
where the first inequality follows from the spectral theorem, also for
multiple eigenvalues. Hence
$\|\rvw\|_\infty^2\leq CN^{-1+\kappa/2}\leq N^{-1+\kappa}$ for large
$N$. The local-law event is uniform in $z$, so this holds simultaneously
for all such eigenvectors, with any prescribed polynomial failure bound.
\end{proof}

\subsection{Proof strategy for the local law}
\label{sec:outline-local-law}

The main task is to control the opposite-class background while
retaining the realized squared row norm in the diagonal comparison
\eqref{eqn:stieltjes-alpha-inverse}.

\paragraph{The Schur equation at local scales.}
Schur's complement suggests the approximate relation
\begin{align}
    \frac1{\eres_{xx}}
    \approx-z-\sum_y^{(x)}|\emH_{xy}|^2\eres_{yy},
\end{align}
where the approximation suppresses the off-diagonal terms in the minor
resolvent, the signal terms, and the error from removing the minor.
The first difficulty is to justify this reduction at spectral scales down to
$\Im z=N^{-1+\kappa}$. At these scales, the deterministic bound
$\|\res(z)\|\leq(\Im z)^{-1}$ is too large for the required error
estimates. In the adjacency-matrix case, local tree structure alone
also does not control the resolvent at local scales.
We work under the a priori weak local law $\reserr(z)\leq\epsilon$,
where $\reserr$ is the entrywise error defined in
\eqref{eqn:local-law-control-scale}. This is the event $\phi(z)=1$ in
\eqref{eqn:res-priori-bound}. It supplies bounded resolvent entries and
proximity to the desired comparison profiles. Minor comparison, Ward's
identity, and sparse large-deviation estimates then control the Schur
remainder and the deletion errors.

\paragraph{Instability of the neighborhood equation.}
The sparse neighborhood equation does not admit an entrywise stability
estimate uniform over sparse graphs. The regular-tree example in
\cite[Section~2.3 and Appendix~A.10]{alt2021delocalization} illustrates
this obstruction in the Erd\H{o}s--R\'enyi setting. A biregular benchmark
shows that weak proximity to the two class profiles and a small residual
do not yield a graph-uniform entrywise improvement for general
approximate solutions of the neighborhood equation. At logarithmic
sparsity, residuals of order $\varphi_{\fka}$ can coexist with errors
of order $\sqrt{\varphi_{\fka}}$, where $\varphi_{\fka}$ is the
control scale in \eqref{eqn:local-law-control-scale}; see
\eqref{eqn:neighborhood-nonlinear-instability} and its construction in
Subsection~\ref{sec:delocalized-neighborhood-instability-proof}.
This motivates the reduction to class averages.

\paragraph{Typical vertices and class averaging.}
Following the typical-vertex strategy of \cite{alt2021delocalization},
we replace the random weighted neighborhood sums by averages over the
opposite vertex class. A squared row norm close to its mean is
insufficient for this replacement: the resolvent-weighted fluctuation
must also be small. The parameters $\Phi_x$ and $\Psi_x$ in
\eqref{eqn:typicality-parameters} test these two fluctuations and define
the typical and atypical sets in
\eqref{eqn:typical-vertices}--\eqref{eqn:atypical-vertices}.

Under the bounded-resolvent event,
Proposition~\ref{prop:size-typical-vertices} bounds the proportion of
atypical vertices and their total weight in every row by
$O(\varphi_{\fka})$, with overwhelming probability.
The cardinality bound permits replacement of full class averages by
typical-class averages, while the row-weight bound controls the
atypical contribution to the weighted Schur equation.

On typical vertices, the Schur equation reduces to a two-class
averaging system that is uniformly stable in the bulk by
Lemma~\ref{lem:stability-self-consistent}. The a priori proximity
therefore improves to an $O(\varphi_{\fka})$ bound on the difference
between the typical diagonal entries and their class profiles.
Substitution into the original weighted Schur equation gives
the controlled equation in
Proposition~\ref{prop:controlled-schur-system}, retaining the full
squared row norm $\zeta_x$ even when $x$ itself is atypical.

\paragraph{Decoupling and the control scale.}
Under the bounded-resolvent event, concentration gives a
small probability of atypicality for a fixed vertex, but its tail is
not small enough to prove that every vertex is typical by a union bound.
Moreover, the typicality events are dependent through the minor
resolvents in $\Psi_x$.
For a small deterministic collection $\sU$ of vertices, we replace
these different minors by the common minor with $\sU$ deleted.
Conditional on that minor, the remaining row tests are independent.
Quantitative minor comparison controls the replacement and preserves
the atypicality test for $|\sU|\lesssim\varphi_{\fka}d$.
This gives joint bounds of the form
\begin{align}
    \P\bigl(\sU\subset\Vtyp_{\fka/2}^{\complement},\,
       \theta_1^{(\emptyset)}=1\bigr)
    \leq e^{-c\varphi_{\fka}^2d|\sU|}+CN^{-\nu};
\end{align}
see Lemma~\ref{lem:typical-vertices-probability}.
The largest usable collections yield an exponent of order
$\varphi_{\fka}^3d^2$. Requiring polynomially small failure
probabilities explains the balance
\begin{align}
    \varphi_{\fka}^3d^2=\fka^3\log N,
    \qquad
    \varphi_{\fka}=\fka\left(\frac{\log N}{d^2}\right)^{1/3}.
\end{align}
A conditional weight-level decomposition in the minor deleting the
tested row converts these counting estimates into a bound on atypical
weight in every row.

\paragraph{Closing the bootstrap.}
The weighted profile bound \eqref{eqn:weighted-stieltjes-bound}
compensates for the factor $1+\zeta_x$ in the Schur error. Together
with minor comparison, this yields the fixed-parameter improvement
in Proposition~\ref{prop:resolvent-strong-bound}. Starting at large
imaginary part, where the weak bound holds deterministically,
continuity propagates it to successive spectral parameters, and the
proposition restores the stronger estimate at each step. Iteration
reaches $\Im z=N^{-1+\kappa}$, and a spectral net makes the estimate
uniform throughout $\specdomain$.

\subsection{Typical vertices}
\label{sec:local-law-typical-vertices}

Throughout the proof, assume the hypotheses of
Theorem~\ref{thm:bipartite-sparse-local-law}. For $\setT\subset[N]$,
let $\mM^{(\setT)}$ be the minor obtained by deleting $\setT$, and write
\begin{align}
    \res^{(\setT)}(z)\coloneqq(\mM^{(\setT)}-z)^{-1},
    \qquad \sum_x^{(\setT)}\coloneqq\sum_{x\notin\setT}.
    \label{eqn:minor-notation}
\end{align}
The same convention applies to $\mH$. We omit braces for singleton
sets and use juxtaposition for unions, as in $\setT x$ and $\setS\setT$.
Fix a deterministic $z\in\specdomain$ and suppress it when convenient.
For a constant $\fka>0$, set
\begin{align}
    \varphi_{\fka}\coloneqq\fka\left(\frac{\log N}{d^2}\right)^{1/3},
    \qquad
    \reserr(z)\coloneqq\max_{x,y\in[N]}
    \left|\eres_{xy}(z)-\indi{x=y}\stj_{\zeta_x}^{(\evsigma_x)}(z)\right|.
    \label{eqn:local-law-control-scale}
\end{align}
For $0<\epsilon\leq1$, define
\begin{align}
    \phi(z)\coloneqq\indi{\reserr(z)\leq\epsilon}.
    \label{eqn:res-priori-bound}
\end{align}
The event $\phi(z)=1$ is the a priori weak local law used in the
controlled Schur equation and the one-step improvement.
Choose $\resmax\geq C_{\rm St}+1$, where $C_{\rm St}$ is from
Lemma~\ref{lem:stieltjes-bound}\textup{(1)}, and define
\begin{align}
    \theta_c^{(\setT)}\coloneqq
    \indi{\max_{x,y\notin\setT}|\eres_{xy}^{(\setT)}|\leq c\resmax},
    \qquad c\geq1.
    \label{eqn:indi-max-res-general}
\end{align}
By \eqref{eqn:indi-max-res-general} and the bound on the comparison profiles,
\begin{align}
    \phi(z)\leq\theta_1^{(\emptyset)}.
    \label{eqn:phi-less-than-theta}
\end{align}
All constants below may depend on the fixed model parameters. For each
prescribed failure exponent, $\fka$ is chosen sufficiently large and
then kept fixed as $N\to\infty$.

\paragraph{Typicality parameters.}
For deterministic $\setT\subset[N]$ and every $x\in[N]$, define
\begin{align}
    \Phi_x^{(\setT)}&\coloneqq\sum_y^{(\setT x)}
    (|\emH_{xy}|^2-\E|\emH_{xy}|^2),
    &
    \Psi_x^{(\setT)}&\coloneqq\sum_y^{(\setT x)}
    (|\emH_{xy}|^2-\E|\emH_{xy}|^2)\eres_{yy}^{(\setT x)}.
    \label{eqn:typicality-parameters}
\end{align}
Write $\Phi_x=\Phi_x^{(\emptyset)}$ and
$\Psi_x=\Psi_x^{(\emptyset)}$; in particular,
$\Phi_x=\zeta_x-\E\zeta_x$.
The typical vertices of a minor and their class restrictions are
\begin{align}
    \Vtyp_{\fka}^{(\setT)}&\coloneqq
    \{x\notin\setT:|\Phi_x^{(\setT)}|\vee|\Psi_x^{(\setT)}|
       \leq\varphi_{\fka}\},
    &
    \Vtyp_{\fka}&\coloneqq\Vtyp_{\fka}^{(\emptyset)},
    &
    \Vtyp_{j,\fka}&\coloneqq\Vtyp_{\fka}\cap\gV_j.
    \label{eqn:typical-vertices}
\end{align}
The corresponding atypical set is the relative complement
\begin{align}
    (\Vtyp_{\fka}^{(\setT)})^{\complement}\coloneqq
    \{x\notin\setT:|\Phi_x^{(\setT)}|\vee|\Psi_x^{(\setT)}|
       >\varphi_{\fka}\}.
    \label{eqn:atypical-vertices}
\end{align}
The parameters $\Psi_x^{(\setT)}$ and the typical and atypical sets
depend on $z$; we suppress this dependence in our notation.
Both parameters remain defined for $x\in\setT$; this is
needed when testing deleted rows against a common minor.
The condition on $\Phi_x$ controls the row norm, while the condition
on $\Psi_x$ controls the fluctuation of the weighted diagonal average.

\paragraph{Typical-class averages.}
When both typical classes are nonempty, define
\begin{align}
    \overline G_j\coloneqq\frac1{|\Vtyp_{j,\fka}|}
       \sum_{y\in\Vtyp_{j,\fka}}\eres_{yy},\qquad j=1,2.
    \label{eqn:typical-class-average}
\end{align}

\subsection{Main estimates and the controlled Schur equation}
\label{sec:local-law-resolvent-inputs}

We state the estimates used to derive the controlled Schur equation
at a fixed $z$. Their proofs are deferred to
Section~\ref{sec:delocalized-deferred-proofs}.

\begin{lemma}[Minor-resolvent comparison]
\label{lem:res-max-difference-bound}
Fix $\nu>0$. There are constants $c_*,C>0$, depending only on
$\kappa,\qr,L,\nu$, such that, for each deterministic $\setT\subset[N]$
satisfying
\begin{align}
    |\setT|\leq c_*\resmax^{-2}d,
    \label{eqn:setT-size}
\end{align}
with probability at least $1-CN^{-\nu}$,
\begin{align}
    \theta_1^{(\emptyset)}
    \max_{x,y\notin\setT}|\eres_{xy}^{(\setT)}|\leq2\resmax,
    \qquad
    \theta_1^{(\emptyset)}\max_{x\neq y}|\eres_{xy}|\leq Cd^{-1/2}.
    \label{eqn:resolvent-entry-upper-bound}
\end{align}
For every deterministic $u\notin\setT$, with the same probability bound,
\begin{align}
    \theta_1^{(\emptyset)}
    \max_{x,y\notin\setT u}
    |\eres_{xy}^{(\setT)}-\eres_{xy}^{(\setT u)}|\leq Cd^{-1}.
    \label{eqn:resolvent-schur-difference-upper-bound}
\end{align}
In particular, $\theta_1^{(\emptyset)}\leq\theta_2^{(\setT)}$
outside an event of probability at most $CN^{-\nu}$.
\end{lemma}

The next proposition gives the two typical-set properties used in the
proof. Sharper exponential estimates are established in
Subsection~\ref{sec:delocalized-atypical-set-proof}.

\begin{proposition}[Size and edge weight of the atypical set]
\label{prop:size-typical-vertices}
For each $\nu>0$, one can choose a fixed $\fka>0$ such that, on
$\{\theta_1^{(\emptyset)}=1\}$ and outside an event of probability
at most $CN^{-\nu}$,
\begin{align}
    \frac{|\Vtyp_{\fka}^{\complement}|}{N}
    &\lesssim\varphi_{\fka},
    \label{eqn:atypical-set-cardinality}\\
    \max_{x\in[N]}\sum_{y\in\Vtyp_{\fka}^{\complement}}|\emH_{xy}|^2
    &\lesssim\varphi_{\fka}.
    \label{eqn:atypical-set-weight}
\end{align}
In particular, $|\Vtyp_{j,\fka}|\geq|\gV_j|/2$ for $j=1,2$ and
large $N$.
\end{proposition}

\begin{lemma}[Schur-complement equations]
\label{lem:approximate-self-consistent}
For $x\in\gV_j$, the Schur complement can be written in the forms
\begin{align}
    \frac1{\eres_{xx}}
    &=-z-\sum_y^{(x)}|\emH_{xy}|^2\eres_{yy}^{(x)}+\rY_x,
    \label{eqn:res-inverse-expansion-raw}\\
    \frac1{\eres_{xx}}
    &=-z-c_j\overline G_{3-j}+\widehat\epsilon_x,
    \label{eqn:res-inverse-expansion-error-epsilon}
\end{align}
These identities define the residuals $\rY_x$ and
$\widehat\epsilon_x$; the second requires both typical classes to be
nonempty. The residual $\rY_x$ contains the off-diagonal minor terms and
the signal terms. The expanded formulas are given in
Subsection~\ref{sec:delocalized-schur-proof}.
\end{lemma}

For the row-dependent comparison, define
\begin{align}
    \epsilon_x\coloneqq\rY_x+\sum_y^{(x)}|\emH_{xy}|^2
       (\stj^{(3-\evsigma_x)}-\eres_{yy}^{(x)}).
    \label{eqn:error-epsilon-zeta-x}
\end{align}
Then \eqref{eqn:res-inverse-expansion-raw} is equivalently
\eqref{eqn:controlled-schur-diagonal}.

\begin{lemma}[Schur error bounds]
\label{lem:schur-complement-error}
With overwhelming probability, on $\{\theta_1^{(\emptyset)}=1\}$,
both typical classes are nonempty and
\begin{align}
    \max_{x\in[N]}|\rY_x|\lesssim d^{-1},
    \qquad
    \max_{x\in\Vtyp_{\fka}}|\widehat\epsilon_x|
    \lesssim\varphi_{\fka}.
    \label{eqn:schur-complement-error-bounds}
\end{align}
\end{lemma}

The following deterministic lemma gives bulk stability of the averaged
equation.

\begin{lemma}[Bulk stability]
\label{lem:stability-self-consistent}
Let $\sX_1,\sX_2$ be nonempty finite sets and
$\sX=\sX_1\sqcup\sX_2$. Suppose that $\Re z\in\bulk$, $\Im z>0$,
and $(\evg_x)_{x\in\sX}$ satisfies
\begin{align}
    \frac1{\evg_x}=-z-\frac{c_j}{|\sX_{3-j}|}
       \sum_{y\in\sX_{3-j}}\evg_y+\epsilon_x,
    \qquad x\in\sX_j,\quad j=1,2.
    \label{eqn:self-consistent-equation-stability}
\end{align}
There are $\epsilon_*,C>0$, depending only on $\kappa,\qr$, such that
\begin{align}
    \max_{j=1,2}\max_{x\in\sX_j}|\evg_x-\stj^{(j)}(z)|
    \leq\epsilon_*
    \quad\Longrightarrow\quad
    \max_{j=1,2}\max_{x\in\sX_j}|\evg_x-\stj^{(j)}(z)|
    \leq C\max_{x\in\sX}|\epsilon_x|.
    \label{eqn:stability-self-consistent}
\end{align}
\end{lemma}

\begin{proposition}[Controlled Schur equation]
\label{prop:controlled-schur-system}
There is $\epsilon>0$, depending only on $\kappa$ and $\qr$, such that,
with overwhelming probability, on $\{\phi(z)=1\}$,
\begin{align}
    \frac1{\eres_{xx}(z)}
    &=-z-\zeta_x\stj^{(3-\evsigma_x)}(z)+\epsilon_x,
    \qquad x\in[N],
    \label{eqn:controlled-schur-diagonal}\\
    \max_{x\in[N]}\frac{|\epsilon_x|}{1+\zeta_x}
    &\lesssim\varphi_{\fka}.
    \label{eqn:controlled-schur-error}
\end{align}
\end{proposition}

\begin{proof}[Proof of Proposition~\ref{prop:controlled-schur-system}]
On $\{\phi(z)=1\}$, we have $\theta_1^{(\emptyset)}=1$ by
\eqref{eqn:phi-less-than-theta}. It remains to bound the residual in
\eqref{eqn:error-epsilon-zeta-x}.

\paragraph{Typical vertices.}
For $y\in\Vtyp_{j,\fka}$, typicality and
\eqref{eqn:zeta-class-means} give
$\zeta_y=c_j+O(\varphi_{\fka}+d/N)=c_j+O(\varphi_{\fka})$.
The Lipschitz bound in Lemma~\ref{lem:stieltjes-bound}\textup{(2)} and
$\reserr(z)\leq\epsilon$ therefore imply
$|\eres_{yy}-\stj^{(j)}|\leq\epsilon+C\varphi_{\fka}$.
Choose $\epsilon$ small enough for Lemma~\ref{lem:stability-self-consistent}
and then take $N$ large. Applying that lemma to
\eqref{eqn:res-inverse-expansion-error-epsilon}, with
$\sX_j=\Vtyp_{j,\fka}$, yields
\begin{align}
    \Delta_{\rm typ}\coloneqq
    \max_{j=1,2}\max_{y\in\Vtyp_{j,\fka}}
       |\eres_{yy}-\stj^{(j)}|
    \lesssim\max_{y\in\Vtyp_{\fka}}|\widehat\epsilon_y|
    \lesssim\varphi_{\fka}.
    \label{eqn:typical-diagonal-control}
\end{align}

\paragraph{All vertices.}
For every $x$, nonzero $\emH_{xy}$ have
$y\in\gV_{3-\evsigma_x}$. Split \eqref{eqn:error-epsilon-zeta-x} into
typical and atypical neighbors. Using
\eqref{eqn:resolvent-entry-upper-bound},
\eqref{eqn:resolvent-schur-difference-upper-bound}, and the bounded
comparison profiles gives
\begin{align}
    |\epsilon_x|
    &\lesssim d^{-1}+\zeta_x(\Delta_{\rm typ}+d^{-1})
       +\sum_{y\in\Vtyp_{\fka}^{\complement}}|\emH_{xy}|^2\notag\\
    &\lesssim d^{-1}+\zeta_x\varphi_{\fka}
       +\varphi_{\fka}
    \lesssim(1+\zeta_x)\varphi_{\fka}.
\end{align}
The second line uses \eqref{eqn:typical-diagonal-control} and
\eqref{eqn:atypical-set-weight}, together with $d^{-1}=o(\varphi_{\fka})$.
This proves \eqref{eqn:controlled-schur-error}.
\end{proof}

\subsection{One-step improvement and continuity}
\label{sec:local-law-bootstrap}

Proposition~\ref{prop:controlled-schur-system} yields the following
improvement at a fixed spectral parameter.

\begin{proposition}[One-step resolvent improvement]
\label{prop:resolvent-strong-bound}
For $\epsilon$ as in Proposition~\ref{prop:controlled-schur-system} and
every deterministic $z\in\specdomain$, with overwhelming probability,
\begin{align}
    \phi(z)\reserr(z)\lesssim\varphi_{\fka}.
    \label{eqn:resolvent-strong-bound}
\end{align}
\end{proposition}

\begin{proof}
Work on $\{\phi(z)=1\}$. By \eqref{eqn:phi-less-than-theta} and
\eqref{eqn:resolvent-entry-upper-bound}, the off-diagonal entries are
$O(d^{-1/2})=O(\varphi_{\fka})$. For the diagonal entries, the explicit
profiles have imaginary parts bounded below on $\specdomain$.
Together with Lemma~\ref{lem:stieltjes-bound},
\eqref{eqn:stieltjes-alpha-inverse} therefore gives
\begin{align}
    \sup_{\alpha\geq0}\max_{j=1,2}
    (1+\alpha)|\stj_\alpha^{(j)}(z)|\lesssim1.
    \label{eqn:weighted-stieltjes-bound}
\end{align}
Indeed, the uniform profile bound handles $\alpha\leq1$, while
$|\stj_\alpha^{(j)}|\leq
(\alpha\Im\stj^{(3-j)})^{-1}$ handles $\alpha>1$.
Subtracting the inverse identities
\eqref{eqn:stieltjes-alpha-inverse} and
\eqref{eqn:controlled-schur-diagonal}, with $j=\evsigma_x$, yields
\begin{align}
    |\eres_{xx}-\stj_{\zeta_x}^{(j)}|
    =|\eres_{xx}\stj_{\zeta_x}^{(j)}\epsilon_x|
    \lesssim\frac{|\epsilon_x|}{1+\zeta_x}
    \lesssim\varphi_{\fka},
\end{align}
using $|\eres_{xx}|\leq\resmax$ and
\eqref{eqn:controlled-schur-error}. This proves the claim.
\end{proof}

\begin{proof}[Proof of Theorem~\ref{thm:bipartite-sparse-local-law}]
Let $\epsilon>0$ be the threshold in \eqref{eqn:res-priori-bound},
chosen as in Proposition~\ref{prop:resolvent-strong-bound}.
Enlarging $L$ if necessary, it is enough to prove the theorem under
the additional assumption $L\geq2\epsilon^{-1}+1$, since this only
enlarges the spectral domain.

\paragraph{Fixed spectral parameter.}
We first prove the estimate at a fixed deterministic point
$z_0\in\specdomain$. Set
\begin{align}
    J\coloneqq\min\{j\in\N_0:\Im z_0+jN^{-3}\geq2\epsilon^{-1}\},
    \qquad z_j\coloneqq z_0+\ii jN^{-3},\quad 0\leq j\leq J.
\end{align}
Then $J\leq C_{\kappa,L}N^3$, and all points $z_j$ belong to the
enlarged spectral domain. We claim that, with overwhelming probability,
\begin{align}
    \reserr(z_j)\lesssim\varphi_{\fka},\qquad 0\leq j\leq J.
    \label{eqn:local-law-lattice-bound}
\end{align}
The proof is a downward induction in $j$. At the top point $z_J$,
the resolvent bound and
\eqref{eqn:stieltjes-uniform-resolvent-bound} give
\begin{align}
    |\eres_{xy}(z_J)|\leq(\Im z_J)^{-1}\leq\epsilon/2,
    \qquad
    |\stj_{\zeta_x}^{(k)}(z_J)|\leq(\Im z_J)^{-1}\leq\epsilon/2,
    \quad k=1,2.
\end{align}
Hence $\reserr(z_J)\leq\epsilon$, so $\phi(z_J)=1$.
Applying Proposition~\ref{prop:resolvent-strong-bound} at $z_J$
gives \eqref{eqn:local-law-lattice-bound} for $j=J$.

Assume now that \eqref{eqn:local-law-lattice-bound} holds at $z_j$
for some $1\leq j\leq J$. By
Lemma~\ref{lem:resolvent-identities}\textup{(5)}, the functions
$\eres_{xy}(z)$ are Lipschitz on $\specdomain$ with constant at most
$N^2$, uniformly in $x,y$. Since $\zeta_x\geq0$ for every $x\in[N]$,
Lemma~\ref{lem:stieltjes-bound}\textup{(3)} gives the same Lipschitz
constant for $\stj_{\zeta_x}^{(k)}(z)$, uniformly in $x$ and $k$.
These bounds hold for all realized row norms $\zeta_x$.
Since $|z_j-z_{j-1}|=N^{-3}$, we obtain
\begin{align}
    \reserr(z_{j-1})
    \leq\reserr(z_j)+2N^2|z_j-z_{j-1}|
    =\reserr(z_j)+2N^{-1}
    \lesssim\varphi_{\fka}+N^{-1}.
\end{align}
By \eqref{eqn:subcritical_regime}, $\varphi_{\fka}=o(1)$, and the
last quantity is at most $\epsilon$ for sufficiently large $N$.
Thus $\phi(z_{j-1})=1$, and another application of
Proposition~\ref{prop:resolvent-strong-bound} yields
\eqref{eqn:local-law-lattice-bound} at $z_{j-1}$. This closes the
induction. For any target exponent $\nu>0$, invoke each input estimate
with failure exponent $\nu+4$. A union bound over the $O(N^3)$
mesh points then gives a failure probability $O(N^{-\nu})$ at $z_0$.

\paragraph{Uniformity in the spectral parameter.}
Let $\mathcal N\subset\specdomain$ be a deterministic $N^{-3}$-net
with $|\mathcal N|\leq C_{\kappa,L}N^6$. For any target exponent
$\nu>0$, apply the fixed-$z_0$ estimate with failure exponent $\nu+7$
at every $z_0\in\mathcal N$. The constant $\fka$ is chosen once,
sufficiently large for the resulting input exponent $\nu+11$.
A union bound gives, with probability at least $1-O(N^{-\nu})$,
\begin{align}
    \max_{z\in\mathcal N}\reserr(z)\lesssim\varphi_{\fka}.
\end{align}
For any $z\in\specdomain$, choose $z_0\in\mathcal N$ with
$|z-z_0|\leq N^{-3}$. The same uniform Lipschitz estimates give
\begin{align}
    \reserr(z)
    \leq\reserr(z_0)+2N^2|z-z_0|
    \leq\reserr(z_0)+2N^{-1}
    \lesssim\varphi_{\fka},
\end{align}
where $N^{-1}=o(\varphi_{\fka})$ by
\eqref{eqn:subcritical_regime}. Taking the supremum over
$z\in\specdomain$ and recalling \eqref{eqn:local-law-control-scale}
proves \eqref{eqn:sparse-entrywise-local-law}, since
$\varphi_{\fka}=\fka(\log N/d^2)^{1/3}$ and $\fka$ is fixed.
\end{proof}

%% file: semilocalized.tex
\section{The Semilocalized Phase}\label{sec:semilocalized_phase}
Radial profiles and complementary-spectrum bounds give outlier locations
in Theorems~\ref{thm:right_edge_behavior} and~\ref{thm:left_edge_behavior}.
Resonant projection, degree counts and support bounds give
Theorems~\ref{thm:right_vector_semilocalization} and
\ref{thm:left_vector_semilocalization}.

At the left edge, we restrict auxiliary operators to the giant component and extend them by zero; by
\eqref{eqn:giant-component-centered-comparison}, this costs $O_{b,\qr}(\error_d)$ in operator norm. This reduction is justified in Appendix~\ref{sec:giant-component-reduction-and-spectrum-near-zero} and used in the left input events below.

\subsection{Pruned graph and singular vector profiles}
\label{subsec:semilocalized-profile-projections}

\subsubsection{Existence of pruned graph}
For a triple $\vtau \coloneqq\bigl(\tau_1,\tau_2^{\mathrm{L}},\tau_2^{\mathrm{R}}\bigr)$, define the vertex set by
\begin{align}
  \gV^{(\vtau)}
&\coloneqq\Vonehigh{\tau_1}\cup\Vtwohigh{\tau_2^{\mathrm{R}}}
      \cup\Vtwolow{\tau_2^{\mathrm{L}}},
    \label{eqn:energy-dependent-pruning-set}
\end{align}
where the vertex sets are defined by
\begin{align}
    \Vonehigh{\tau_1}
    &\coloneqq\{x\in\gV_1:\alpha_x\geq\tau_1\},
    &
    \Vtwohigh{\tau_2^{\mathrm{R}}}
    &\coloneqq\{x\in\gV_2:\alpha_x\geq\tau_2^{\mathrm{R}}\},
    &
    \Vtwolow{\tau_2^{\mathrm{L}}}
    &\coloneqq\{x\in\gV_2:\alpha_x\leq\tau_2^{\mathrm{L}}\}.
\end{align}
Let $\gap>0$ denote the distance between the spectral edge and the target eigenvalue. We define the gap-dependent triple  $\vtau(\gap)\coloneqq\bigl(\tau_1(\gap),\tau_2^{\mathrm{L}}(\gap),\tau_2^{\mathrm{R}}(\gap)\bigr)$ by
\begin{align}
    \tau_1(\gap)&\coloneqq\qr^{-2}+\frac{1}{8\qr}s(\gap),
    \qquad
    \tau_2^{\mathrm{R}}(\gap)\coloneqq\qr^2+\frac{\qr}{8}s(\gap),
    \qquad
    \tau_2^{\mathrm{L}}(\gap)\coloneqq\qr^2-\frac{\qr}{8}s(\gap).
    \label{eqn:energy-dependent-thresholds}
\end{align}
where $s(\gap)\coloneqq\gap\wedge1$. We further define the radius by
\begin{align}
    \radius
    \coloneqq
    \left\lfloor(\qr^{-2}/4)\sqrt{d/\log(d)}\right\rfloor.
    \label{eqn:radius}
\end{align}

The degree displacements are $s/(8\qr)$ on side $1$ and
$\qr s/8$ on side $2$. Their local spectral buffers and the selected
profiles are recorded in Table~\ref{tab:semilocalized-buffers}.
\begin{table}[htbp]
\centering\small
\begin{tabularx}{\textwidth}{@{}lXXX@{}}
\toprule
Application & Regular local boundary & Target & Selected profiles\\
\midrule
Right window & $\redge+s/8+C\error_d$ & $\redge+s/2$
 & $\gR_1\cup\gR_2$, intersected with high-degree cutoffs\\
Right location & $\redge+s/8+C\error_d$ & $\redge+s$
 & All of $\gR_1\cup\gR_2$\\
Left window & Squared floor at $\ledge-s/8-C\error_d$
 & $\ledge-s/4$
 & $\inf\sK/2\leq\Lambda_x\leq\ledge-2\error_d^{1/4}$,
 intersected with low-degree cutoff\\
Global left & The same squared floor & $\ledge-s$, including zero
 & All of $\gL_2$, with retained root degree at least $a_0/2$\\
\bottomrule
\end{tabularx}
\caption{Deterministic buffers. The statements include their form and
comparison errors. The two location applications use
$s=3\error_d^{1/4}/4$; the left window has $\inf\sK>0$.}
\label{tab:semilocalized-buffers}
\end{table}

For fixed $a>0$, the inputs hold on one simultaneous event over
deterministic gaps with $s(\gap)\geq a\error_d^{1/4}$, including
$a=3/4$ for locations. Subsection~\ref{subsec:semilocalized-pruning-proofs}
constructs this event before any eigenvalue-dependent gap is selected.
For $\gap\geq1$,
the pruning is fixed with error $C\error_d$, and the spectral and
resolvent parameters are treated deterministically.

For every deterministic admissible cutoff, the one-vertex Chernoff
bound is
\begin{align}
 \max_x\P(x\in\gV^{(\vtau)})
 \leq 2e^{-c_{\qr}ds^2}.
 \label{eqn:pruning-single-vertex-tail}
\end{align}

\begin{lemma}[Quantitative pruning]
\label{lem:existence_pruned_graph}
Fix $a>0$ and $\nu>0$. Under
\eqref{eqn:subcritical_regime}, with $\gC_0$ sufficiently large depending
on $(B,\qr,\nu,a)$, put $s=s(\gap)$ and suppose
$s\geq a\error_d^{1/4}$. There is an event of probability at least $1-N^{-\nu}$ on which,
simultaneously for all these cutoff sets, a subgraph
$\gG^{(\vtau)}\subset\gG$ has the following properties.
\begin{enumerate}
\item Distinct exceptional vertices have retained distance at least
$4\radius+1$. Their retained $2\radius$-balls are disjoint.
\item The retained $4\radius$-ball at every exceptional vertex is a tree.
\item Every deleted edge is incident to an exceptional vertex. The
deleted graph $\gD^{(\vtau)}$ admits a decomposition
$\gD^{(\vtau)}=\gF^{(\vtau)}\cup\gU^{(\vtau)}$ with
$\gF^{(\vtau)}$ a forest and
\begin{align}
 \max_x\rD_x^{\gD^{(\vtau)}}
 &\leq C_{B,\qr,\nu,a}\frac{\log N}{ds^2},\nonumber\\
 \max_x\rD_x^{\gF^{(\vtau)}}
 &\leq C_{B,\qr,\nu,a}\frac{\log N}{ds^2},
 &\max_x\rD_x^{\gU^{(\vtau)}}&\leq C_{B,\qr,\nu,a}.
 \label{eqn:removed-graph-degrees}
\end{align}
\item For every exceptional $x$ and $j\leq2\radius$,
$\layer_j^{(\vtau)}(x)\subset\layer_j^{\gG}(x)$. For
$y\in\ball_{2\radius}^{(\vtau)}(x)\setminus\{x\}$ and
$j\leq2\radius$,
\begin{align}
 \layer_1^{(\vtau)}(y)\cap\layer_j^{(\vtau)}(x)
 =\layer_1^{\gG}(y)\cap\layer_j^{\gG}(x).
 \label{eqn:retained-branch-neighbours}
\end{align}
\item Put $m_x=\rD_x^{\gD^{(\vtau)}}$ and
$\ell_j(x)=|\layer_j^{\gG}(x)\setminus\layer_j^{(\vtau)}(x)|$.
For every exceptional root,
\begin{align}
 \ell_1(x)&=m_x,
 &\ell_j(x)&\leq C_{B,\qr,\nu,a}
       (dm_x+\log N)d^{j-2},\qquad 2\leq j\leq2\radius.
 \label{eqn:pruned-layer-loss}
\end{align}
The transfer to nonempty, comparable retained layers is asserted only
for roots with $\rD_x\geq c_0d$, as in
Lemma~\ref{lem:semilocalized-layer-concentration}.
\end{enumerate}
The same event includes the crude all-root ball bound
\begin{align}
 \max_x|\ball_{8\radius+2}^{\gG}(x)|
 \leq (8\radius+3)(1+C_{B,\qr,\nu}\log N)^{8\radius+2}
 =N^{o(1)}.
 \label{eqn:pruning-crude-ball-size}
\end{align}
\end{lemma}
The construction removes root branches of unpruned breadth-first trees
containing a nearby exceptional vertex or a nontree-edge endpoint.
Concentration of their total sizes gives \eqref{eqn:pruned-layer-loss}.

Let $\rmA^{(\vtau)}$ be the pruned adjacency matrix and let
$\chi^{(\vtau)}$ project onto the complement of
$\bigcup_{x\in\gV^{(\vtau)}}\ball_{2\radius}^{(\vtau)}(x)$.
To preserve the tree operator inside these balls, define
\begin{align}
    \scA^{(\vtau)}
    \coloneqq
    \frac1{\sqrt d}
    \left(\rmA^{(\vtau)}-
    \chi^{(\vtau)}(\E\rmA)\chi^{(\vtau)}\right) = \begin{bmatrix}
      \bzero & \mX^{(\vtau)} \\
      (\mX^{(\vtau)})^{\sT} & \bzero
    \end{bmatrix}.
    \label{eqn:def_H_tau}
\end{align}
The following norm estimate transfers the pruned construction to the full centered operator.

For $\gap>0$, define the profile-error scale by
\begin{align}
    \error_{\mathrm{prof}}(\gap)
    \coloneqq\frac{\sqrt{\log N}}{d\gap}
    =\frac{\error_d}{\gap\log d}.
    \label{eqn:profile-error-scale}
\end{align}

\begin{lemma}[Pruning input]
\label{lem:semilocalized-pruning-input}
Fix $a>0$ and $\nu>0$. Under \eqref{eqn:subcritical_regime}, choose
$\gC_0$ sufficiently large depending on $(B,\qr,\nu,a)$.
With probability at least $1-N^{-\nu}$,
simultaneously for every deterministic cutoff with
$s(\gap)\geq a\error_d^{1/4}$, the pruned graph in
Lemma~\ref{lem:existence_pruned_graph} satisfies
\begin{align}
    \|\scA-\scA^{(\vtau)}\|
    \leq C_{B,\qr,\nu,a}
      \bigl(\error_d+\error_{\mathrm{prof}}(\gap)\bigr).
    \label{eqn:semilocalized-pruning-input}
\end{align}
Under the standing left-edge assumptions, restriction to the giant
component and extension by zero satisfy the same estimate,
simultaneously for $s(\gap)\geq a\error_d^{1/4}$, with probability
at least $1-CN^{-\nu_0}$ and constants depending on $(b,\qr,a)$.
\end{lemma}

\subsubsection{Active centres and radial profiles}

For the left fixed-window comparison, fix a compact interval
$\sK\subset(0,\ledge)$ of positive length. Constants in this clause
may depend on $(b,\qr,\sK)$. The global left comparison instead uses
a degree floor $a_0>0$ and constants depending on $(b,\qr,a_0)$.
Right constants may depend on $(B,\qr,\nu)$; dependence on the
fixed aspect-ratio bound $\Gamma$ is suppressed in subscripts.

For a selected triple $\vtau$, the \emph{active centres} are the following subsets of the right-outlier
vertices from \eqref{eqn:right-outlier-profile-sets} and the left-branch
vertices in the prescribed positive spectral interval:
\begin{subequations}\label{eqn:profile-centres}
\begin{align}
    \gV_{1,\mathrm{R}}^\star(\vtau)
    &\coloneqq
    \gR_1\cap\Vonehigh{\tau_1},
    \label{eqn:V1-right-profile-centres}\\
    \gV_{2,\mathrm{R}}^\star(\vtau)
    &\coloneqq
    \gR_2\cap\Vtwohigh{\tau_2^{\mathrm{R}}},
    \label{eqn:V2-right-profile-centres}\\
    \gV_{2,\mathrm{L}}^\star(\vtau)
    &\coloneqq
    \{x\in\gV_2:0<\alpha_x<\qr^2-1,
      (\inf\sK)/2\leq\Lambda_{\qr}(\alpha_x)
      \leq\ledge-2\error_d^{1/4}\}
      \cap\Vtwolow{\tau_2^{\mathrm{L}}}.
    \label{eqn:V2-left-profile-centres}
\end{align}
\end{subequations}
For spectral parameters outside the positive MP interval, write
\begin{align}
    \gapR&\coloneqq\lambda-\redge,
    &&\lambda>\redge,
    \label{eqn:right-gap}\\
    \gapL&\coloneqq\ledge-\lambda,
    &&0<\lambda<\ledge.
    \label{eqn:left-gap}
\end{align}
The right and left window comparisons use $\vtau=\vtau(\gapR)$ and
$\vtau=\vtau(\gapL)$, respectively. The choice for outlier locations is
specified in \eqref{eqn:calibrated-comparison-scales}.
Suppress $\vtau$ from the active sets and group the centres by edge:
\begin{align}
    \gV_{\mathrm{R}}^\star
    &\coloneqq\gV_{1,\mathrm{R}}^\star\cup\gV_{2,\mathrm{R}}^\star,
    &\gV_{\mathrm{L}}^\star
    &\coloneqq\gV_{2,\mathrm{L}}^\star,
    &\gV^\star&\coloneqq\gV_{\mathrm{R}}^\star\cup\gV_{\mathrm{L}}^\star.
    \label{eqn:all-profile-centres}
\end{align}
For a centre $x\in\gV^{(\vtau)}$ on a positive outlier branch,
namely $x\in\gR_1\cup\gR_2$ or $x\in\gV_2$ with
$0<\alpha_x<\qr^2-1$, define the profile value
\begin{align}
    \Lambda_x
    \coloneqq
    \Lambda_{\qr^{-1}}(\alpha_x)\indi{x\in\gV_1}
    +\Lambda_{\qr}(\alpha_x)\indi{x\in\gV_2}.
    \label{eqn:profile-energy}
\end{align}

The degree--eigenvalue relation identifies which spectral candidates
are active. Write $s=s(\gap)$ and $u=|\alpha-\qr^2|>1$ on either
$\gV_2$ branch, with $\lambda_{\rm edge}$ its corresponding edge. Then
\begin{align}
    |\Lambda_\qr(\alpha)^2-\lambda_{\rm edge}^2|
    =\frac{(u-1)^2}{u},
    \qquad
    \lambda_{\rm edge}\in\{\ledge,\redge\}.
    \label{eqn:profile-energy-degree-separation}
\end{align}
The cutoff $u\geq\qr s/8$ follows from $u>1$ when $\qr s/8\leq1$.
Otherwise $\qr>8$, so $\ledge\geq\qr/2$. For the candidates used in
resonant projection, the profile value $E$ satisfies
$E\geq\redge+\gap/2$ on a right branch or
$E\leq\ledge-\gap/2$ on the left branch. Hence
\eqref{eqn:profile-energy-degree-separation} gives
$u\geq|E^2-\lambda_{\rm edge}^2|\geq\qr\gap/4\geq\qr s/8$.
On the $\gV_1$ branch, the displacement exceeds one and the cutoff
$s/(8\qr)$ is below one.

The ideal recurrence has coefficients
\begin{align}
    \omega_1&=\sqrt{\alpha_x},
    &\omega_{2j}&=\qr\indi{x\in\gV_1}+\qr^{-1}\indi{x\in\gV_2},
    &\omega_{2j+1}&=\qr^{-1}\indi{x\in\gV_1}+\qr\indi{x\in\gV_2},
    \qquad j\geq1,
    \label{eqn:ideal-radial-coefficients}
\end{align}
Let $(\ervu_j)_{j\geq0}$ be the square-summable solution with
$\ervu_0=1$ of
\begin{align}
    \Lambda_x\ervu_0&=\omega_1\ervu_1,
    &\Lambda_x\ervu_j&=\omega_j\ervu_{j-1}+\omega_{j+1}\ervu_{j+1},
    \qquad j\geq1.
    \label{eqn:ideal-radial-recurrence}
\end{align}
Suppressing the dependence on $x$, we recall the ideal profiles from
Lemma~\ref{lem:eigenZ1Z2} and
\eqref{eqn:Z1Eigenvec}--\eqref{eqn:Z2Eigenvec}, with $\ervu_0=1$.
The coefficient on the first layer is $\ervu_1=\Lambda_x/\sqrt{\alpha_x}$; the
geometric tail is described by
\[
 q_x=\begin{cases}
 (\alpha_x-\qr^{-2})^{-1},&x\in\gV_1,\\
 (\alpha_x-\qr^2)^{-1},&x\in\gV_2.
 \end{cases}
\]
Then, for $k\geq0$,
\begin{align}
 \ervu_{2k+1}&=\ervu_1q_x^k,&
 \ervu_{2k+2}&=\ervu_2q_x^k,\\
 \ervu_2&=\sqrt{\alpha_x}
   (\qr\indi{x\in\gV_1}+\qr^{-1}\indi{x\in\gV_2})q_x,&
 \Lambda_x^2&=\alpha_x\left[1+
 (\qr^2\indi{x\in\gV_1}+\qr^{-2}\indi{x\in\gV_2})q_x\right].
 \label{eqn:signed-profile-heuristic}
\end{align}
On the left branch $q_x<0$, so the profile alternates sign every
two layers; its decay is governed by $|q_x|<1$.
We choose the terminal coefficient from the alternating partial sum
to make the two signed profiles orthogonal. Its positivity is proved
in Subsection~\ref{subsec:semilocalized-profile-proofs}.

Let $J$ be the largest odd integer not exceeding $\radius$ and define
\begin{align}
    \widetilde a_j(x)&=\ervu_j, &&0\leq j<J,\notag\\
    \widetilde a_J(x)&=
    \operatorname{sgn}(\ervu_J)
    \left(\sum_{j=0}^{J-1}(-1)^j \ervu_j^2\right)^{1/2},
    &\widetilde a_j(x)&=0, &&j>J.
    \label{eqn:explicit-terminal-coefficient}
\end{align}
Normalize by
$a_j(x)=\widetilde a_j(x)/(\sum_l\widetilde a_l(x)^2)^{1/2}$.
Since $J$ is odd, the terminal adjustment gives
\begin{align}
    \sum_{j\textnormal{ even}}a_j(x)^2
    =\sum_{j\textnormal{ odd}}a_j(x)^2=\frac12.
    \label{eqn:profile-parity-balance}
\end{align}
For $\sigma\in\{+,-\}$, set
\begin{align}
    \hatvtau_\sigma(x)
    \coloneqq
    \sum_{j=0}^{\radius}
    \sigma^j a_j(x)
    \frac{\ones_{\layer_j^{(\vtau)}(x)}}
    {|\layer_j^{(\vtau)}(x)|^{1/2}},
    \label{eqn:pruned-localization-profile}
\end{align}
with zero contribution from an empty terminal layer. Each profile is
supported on $\ball_{\radius}^{(\vtau)}(x)$. Disjoint balls separate
different centres, and parity balance makes the two signs at each
centre orthogonal.

\begin{lemma}[Profile properties]
\label{lem:semilocalized-profile-properties}
The family
$\{\hatvtau_\sigma(x):x\in\gV^\star,\ \sigma=\pm\}$ is orthonormal.
For every $\nu>0$, uniformly for
$\gap\geq\gC_{B,\nu}\error_d^{1/4}$, with probability at least
$1-N^{-\nu}$,
for $x\in\gV_{\mathrm{R}}^\star$,
\begin{align}
    \|\bigl(\scA^{(\vtau)}-\sigma\Lambda_x\bigr)
    \hatvtau_\sigma(x)\|
    \leq C_{B,\qr,\nu}\error_d.
    \label{eqn:pruned-profile-residual}
\end{align}
At the left edge, the same bound, with constant $C_{b,\qr,\sK}$, holds for
$x\in\gV_{\mathrm{L}}^\star$, uniformly for
$\gap\geq\gC\error_d^{1/4}$, with probability at least
$1-CN^{-\nu_0}$.
The residual in \eqref{eqn:pruned-profile-residual} is supported in
$\ball_{\radius+1}^{(\vtau)}(x)$. Moreover, for
$0\leq r\leq\radius$,
\begin{align}
    \sum_{y\notin\ball_r(x)}|\hatvtau_+(x)_y|^2
    \leq{}& C_{\qr}\bigl(
    \indi{x\in\gV_1}(\alpha_x-\qr^{-2})^{-r}
+\indi{x\in\gV_2,\,\alpha_x>\qr^2}
       (\alpha_x-\qr^2)^{-r}
+\indi{x\in\gV_2,\,\alpha_x<\qr^2}
       (\qr^2-\alpha_x)^{-r}\bigr).
    \label{eqn:profile-tail-bound}
\end{align}
On the same events, uniformly over the active centres and gaps in the
stated ranges, every profile layer up to depth $\radius$ is nonempty and
\begin{align}
    |\operatorname{supp}\hatvtau_\sigma(x)|
    \leq\exp\!\left[
       C\bigl(\sqrt{d\log d}+\log\log N\bigr)\right]
    =N^{o(1)}.
    \label{eqn:profile-support-size}
\end{align}
The constant has the same parameter dependence as the corresponding
right- or left-edge estimate above.
\end{lemma}
The proof and the refined residual estimates are given in
Subsection~\ref{subsec:semilocalized-profile-proofs}.

\subsection{Block approximation}
\label{subsec:semilocalized-block-calibration}

\subsubsection{Orthogonal projections and block approximation}

For a selected centre set $\mathcal C$ whose signed profiles are
orthonormal, define
\begin{align}
    \Pi_{\mathcal C}^{(\vtau)}
    &\coloneqq\sum_{x\in\mathcal C}\sum_{\sigma=\pm}
      \hatvtau_\sigma(x)(\hatvtau_\sigma(x))^*,
    &Q_{\mathcal C}^{(\vtau)}
    &\coloneqq\id_N-\Pi_{\mathcal C}^{(\vtau)}.
    \label{eqn:profile-projections}
\end{align}
\begin{align}
    \hatH_{\mathcal C}
    &\coloneqq\sum_{x\in\mathcal C}\sum_{\sigma=\pm}
      \sigma\Lambda_x\hatvtau_\sigma(x)(\hatvtau_\sigma(x))^*
      +Q_{\mathcal C}^{(\vtau)}\scA^{(\vtau)}
        Q_{\mathcal C}^{(\vtau)}.
    \label{eqn:right-hatH}
\end{align}
The comparison operator acts by the ideal eigenvalues on the profiles
and by the compression of $\scA^{(\vtau)}$ to their complement. For the signed
pairs in \eqref{eqn:pruned-localization-profile}, parity balance makes
$Q_{\mathcal C}^{(\vtau)}$ block diagonal and $\hatH_{\mathcal C}$
bipartite. At gap $\gap$, abbreviate
$\Pi_{\mathrm E}^{(\vtau)}=\Pi_{\gV_{\mathrm E}^\star}^{(\vtau)}$,
$\overline\Pi_{\mathrm E}^{(\vtau)}=Q_{\gV_{\mathrm E}^\star}^{(\vtau)}$,
and $\hatH_{\mathrm E}(\gap)=\hatH_{\gV_{\mathrm E}^\star}$,
where $\vtau=\vtau(\gap)$ and $\mathrm E\in\{\mathrm L,\mathrm R\}$.

\begin{proposition}[Block approximation]
\label{prop:semilocalized-block-approximation}
For the construction in \eqref{eqn:profile-projections}, suppose the
selected profiles are orthonormal, their residual supports
are disjoint between distinct centres, and
\[
    \max_{x\in\mathcal C,\,\sigma=\pm}
    \|(\scA^{(\vtau)}-\sigma\Lambda_x)\hatvtau_\sigma(x)\|
    \leq\varepsilon_{\mathrm{res}},
\]
with maximum zero when $\mathcal C$ is empty. For an absolute constant
$C$, the pruned block error and compression defect satisfy
\begin{align}
    \|\scA^{(\vtau)}-\hatH_{\mathcal C}\|
    +\|\Pi_{\mathcal C}^{(\vtau)}\scA^{(\vtau)}
       Q_{\mathcal C}^{(\vtau)}\|
    &\leq C\varepsilon_{\mathrm{res}},
    \label{eqn:pruned-profile-block-error}\\
    \|\scA-\hatH_{\mathcal C}\|
    &\leq\|\scA-\scA^{(\vtau)}\|+C\varepsilon_{\mathrm{res}}.
    \label{eqn:full-profile-block-error}
\end{align}
\end{proposition}
The proof is given at the end of
Subsection~\ref{subsec:semilocalized-profile-proofs}.

\paragraph{Ordered profile transfer.}
For a bipartite comparison operator $\hatH_{\mathcal C}$, parity
balance gives $r=|\mathcal C|$ orthonormal column profiles, so $r\leq m$.
Its $r$ positive profile values occur, with multiplicity, at distinct
positions $j_1<\cdots<j_r\leq m$ in the nonincreasing list of all $m$
column singular values. Complementary values may interleave with them.
If the profile values in this order are $\Lambda_1\geq\cdots\geq\Lambda_r$,
singular-value Weyl gives
\begin{align}
 |\lambda_{j_k}(\scA)-\Lambda_k|
 \leq\|\scA-\hatH_{\mathcal C}\|,
 \qquad 1\leq k\leq r.
 \label{eqn:ordered-profile-transfer}
\end{align}

\subsection{Nonbacktracking and quadratic-form bounds}
\label{subsec:semilocalized-loewner-inputs}

Row regularization replaces high-degree rows by zero. The nonbacktracking
bounds of \cite{wang2026regularized} then give quadratic forms that
transfer to the pruned operator by the pruning norm estimate.

\subsubsection{Row regularization and nonbacktracking radius}

Fix a gap $\gap>0$ and write $s=s(\gap)$ and
$\vtau=\vtau(\gap)$, with thresholds as in
\eqref{eqn:energy-dependent-thresholds}. We suppress the dependence
on $\gap$ in the cutoff and retained set below.

\begin{definition}[Row regularization]
\label{def:semilocalized-row-regularization}
Define the integer row cutoff and the set of retained rows by
\begin{align}
    U_1&\coloneqq\lceil\tau_1d\rceil,
    &I_1&\coloneqq\{x\in\gV_1:\rD_x<U_1\}
      =\gV_1\setminus\Vonehigh{\tau_1}.
    \nonumber
\end{align}
Let $\mP^{(\tau_1)}$ be the coordinate projection onto $I_1$.
The \emph{row regularization} of $\mX$ is
\begin{align}
    \mP^{(\tau_1)}&\coloneqq
      \operatorname{diag}(\indi{x\in I_1}:x\in\gV_1),
    &\mX^{(\tau_1)}&\coloneqq\mP^{(\tau_1)}\mX,
    \nonumber\\
    \emX^{(\tau_1)}_{xy}&\coloneqq
      \indi{x\in I_1}\frac{\emA_{xy}-p}{\sqrt d},
    & &x\in\gV_1,\quad y\in\gV_2.
    \label{eqn:semilocalized-row-projection}
\end{align}
Applying the same projection to the pruned block in
\eqref{eqn:def_H_tau} gives its row regularization
\begin{align}
    \mZ\coloneqq\mP^{(\tau_1)}\rmX^{(\vtau)}.
    \label{eqn:semilocalized-pruned-row-regularization}
\end{align}
\end{definition}
Since degrees are integers, $\rD_x<U_1$ is equivalent to
$\rD_x<\tau_1d$. Retained entries of $\mX$ keep the original centering term
$p$, and both regularizations have size $n\times m$.
The pruning estimate \eqref{eqn:semilocalized-pruning-input} bounds their difference. For any $\mY\in\C^{n\times m}$, let
$\mH_{\mY}=\left[\begin{smallmatrix}0&\mY\\\mY^*&0\end{smallmatrix}\right]$
and define its nonbacktracking matrix on
$\vec{\gE}=(\gV_1\times\gV_2)\sqcup(\gV_2\times\gV_1)$ by
\begin{align}
    \mB(\mY)_{(u,v),(v',w)}
    \coloneqq\indi{v=v'}\indi{w\neq u}(\mH_{\mY})_{vw}.
    \label{eqn:semilocalized-nb-definition}
\end{align}
Write $\mB^{(\tau_1)}=\mB(\mX^{(\tau_1)})$.
A deterministic family of gaps $\mathfrak G_N$ is \emph{admissible} if
\begin{align}
    \inf_{\gap\in\mathfrak G_N}
      \frac{d\,s(\gap)^2}{\log d}\longrightarrow\infty.
    \label{eqn:semilocalized-admissible-gap-family}
\end{align}

\begin{proposition}[Nonbacktracking input,
    {\cite[nonbacktracking-radius proposition]{wang2026regularized}}]
\label{prop:semilocalized-nb-radius-input}
Under Assumptions~\ref{ass:proportional_regime} and
\ref{ass:subcritical_regime}, fix $B>0$, $\qr\geq1$, and $\nu>0$.
For every deterministic admissible family $\mathfrak G_N$, with
probability at least $1-N^{-\nu}$ for all sufficiently large $N$,
\begin{align}
    \sup_{\gap\in\mathfrak G_N}\rho(\mB^{(\tau_1(\gap))})
    \leq1+C_{B,\qr,\nu}d^{-1/2}.
    \label{eqn:semilocalized-nb-input}
\end{align}
\end{proposition}
The specialization and its hypothesis checks are given in
Subsection~\ref{subsec:semilocalized-loewner-proofs}.

\subsubsection{Quadratic forms for the regularized matrix}

Deleting all edges incident to $\gV_1\setminus I_1$ gives the row-regularized graph. Let $\mH^{(\tau_1)}\coloneqq\mH_{\mX^{(\tau_1)}}$ be the linearization of $\mX^{(\tau_1)}$. We further define
\begin{align}
    \emD^{(1,\tau_1)}_{xx}
    &\coloneqq\indi{x\in I_1}\rD_x,
    &\emD^{(2,\tau_1)}_{yy}
    &\coloneqq\sum_{x\in I_1}\emA_{xy},
    \nonumber\\
    \mD^{(a,\tau_1)}&\coloneqq
      \operatorname{diag}(\emD^{(a,\tau_1)}_{zz}:z\in\gV_a),
    &\mD^{(\tau_1)}&\coloneqq
      \operatorname{diag}(\mD^{(1,\tau_1)},\mD^{(2,\tau_1)}).
    \label{eqn:semilocalized-regularized-degrees-and-linearization}
\end{align}
Every retained row has normalized degree below
$\tau_1$, and regularization can only decrease column degrees.
For the left-edge application, suppose $\qr>1$ and
$0<\gap<\ledge$. Then
\begin{align}
    \tau_1<\frac{1+7\qr^{-2}}8<1.
    \label{eqn:good-row-margin}
\end{align}
The corresponding inverse row weight is
\begin{align}
    \mM^{(\tau_1)}
    \coloneqq\sum_{x\in I_1}
      \frac{1}{1-\rD_x/d}\,\ones_x\ones_x^*.
    \label{eqn:semilocalized-regularized-inverse}
\end{align}
Then $0\preceq\mM^{(\tau_1)}\preceq(1-\tau_1)^{-1}\id_n$.

\begin{proposition}[Regularized Loewner bounds]
\label{prop:semilocalized-regularized-loewner-bounds}
Fix $\nu>0$ and a deterministic admissible family of gaps.
Under the respective assumptions below, each assertion holds with
probability at least $1-N^{-\nu}$ for all sufficiently large $N$,
simultaneously over the family.
\begin{enumerate}
    \item Under Assumptions~\ref{ass:proportional_regime} and
    \ref{ass:subcritical_regime}, fix $B>0$ and $\qr\geq1$.
    \begin{align}
        \mH^{(\tau_1)}
        &\preceq\id_N+d^{-1}\mD^{(\tau_1)}+\mE_+^{(\tau_1)},
        \label{eqn:upper_bound_H}\\
        0\preceq\mE_+^{(\tau_1)}
        &\preceq C_{B,\qr,\nu}d^{-1/2}
           \bigl(\id_N+d^{-1}\mD^{(\tau_1)}\bigr).
        \label{eqn:upper-bound-H-loewner-error}
    \end{align}

    \item Under Assumptions~\ref{ass:proportional_regime} and
    \ref{ass:critical_sparsity}, fix $b>0$ and $\qr>1$,
    and suppose the family lies in $0<\gap<\ledge$.
    \begin{align}
        (\mX^{(\tau_1)})^*\mM^{(\tau_1)}\mX^{(\tau_1)}
        &\succeq d^{-1}\mD^{(2,\tau_1)}-\id_m-\mE^{(\tau_1)},
        \nonumber\\
        0\preceq\mE^{(\tau_1)}
        &\preceq C_{b,\qr,\nu}d^{-1/2}\id_m.
        \label{eqn:lower_bound_H}
    \end{align}
\end{enumerate}
\end{proposition}
The square-sum calculations and both form derivations are proved in
Subsection~\ref{subsec:semilocalized-loewner-proofs}.

\subsubsection{Column-degree loss under row regularization}
\label{subsec:semilocalized-row-regularization-transfer}

The lower form involves retained column degrees. To use it in the
left complementary-spectrum argument, we compare these with the
original column degrees. Write
$\mD^{(2)}=\operatorname{diag}(\rD_y:y\in\gV_2)$ for the original
column-degree matrix. The loss at a column is the number of its
neighbours in the discarded row set.

\begin{lemma}[Column-degree loss]
\label{lem:semilocalized-column-degree-loss}
Under Assumptions~\ref{ass:proportional_regime} and
\ref{ass:subcritical_regime}, fix $B>0$ and $\nu>0$.
For $s=s(\gap)$ and $\tau_1=\tau_1(\gap)$, set
\begin{align}
    \ell_y(s)&\coloneqq
      \sum_{x\in\gV_1}\emA_{xy}\indi{\rD_x\geq\tau_1d},
      \qquad y\in\gV_2,
    &\ell_s&\coloneqq d^{-1}+\frac{\log N}{d^2s^2}.
    \label{eqn:regularized-loss-scale}
\end{align}
For every deterministic admissible family
\eqref{eqn:semilocalized-admissible-gap-family}, with probability at
least $1-N^{-\nu}$ for all sufficiently large $N$, simultaneously over
its gaps,
\begin{align}
    \max_{y\in\gV_2}\ell_y(s)
    &\leq C_{B,\qr,\nu}
      \left(1+\frac{\log N}{ds^2}\right),
    \label{eqn:regularized-incidence}\\
    0\preceq d^{-1}(\mD^{(2)}-\mD^{(2,\tau_1)})
    &\preceq C_{B,\qr,\nu}\ell_s\id_m.
    \label{eqn:regularized-degree-loss}
\end{align}
\end{lemma}
Its proof is given in
Subsection~\ref{subsec:semilocalized-loewner-proofs}.

\subsubsection{Conditional complementary-spectrum form tests}

The form tests use approximate singular vectors of $\mZ$ and bounds
on their exceptional-coordinate mass. We record the scalar algebra
shared with the local forest tests.

\begin{lemma}[Scalar form-testing algebra]
\label{lem:semilocalized-scalar-form-algebra}
Put $B=\frac12(c\redge+\qr\tau_1+\qr^{-1}\tau_2)$ for real
$c,\tau_1,\tau_2$. If $\mu\geq0$, $0\leq\theta\leq1/2$,
$\epsilon\geq0$ and $\qr\mu(1-\theta)\leq\qr B+\epsilon$, then
\begin{align}
 \mu\leq\frac{B+\epsilon/\qr}{1-\theta}
       \leq B+C_{M,\qr}(\theta+\epsilon)
 \label{eqn:semilocalized-scalar-upper-conclusion}
\end{align}
whenever $|c|+|\tau_1|+|\tau_2|\leq M$. At $\theta=0$, the
exact bound is $\mu\leq B+\epsilon/\qr$.
If $1-\tau_1>0$ and
$b-\epsilon_-\leq L\leq\mu^2/(1-\tau_1)+\epsilon_+$, then
\begin{align}
 \mu^2\geq(1-\tau_1)(b-\epsilon_--\epsilon_+).
 \label{eqn:semilocalized-scalar-lower-conclusion}
\end{align}
\end{lemma}
\begin{proof}
Divide the upper inequality by $\qr(1-\theta)$, using
$\theta/(1-\theta)\leq2\theta$ and $|B|\leq C_{M,\qr}$.
Combine the lower inequalities and multiply by $1-\tau_1$.
\end{proof}

\begin{lemma}[Conditional complementary-spectrum form tests]
\label{lem:semilocalized-complement-form-tests}
Fix a gap $\gap$ at which the pruning construction is valid, and use
the row regularizations in Definition~\ref{def:semilocalized-row-regularization}.
Write $\beta_y=d^{-1}\emD^{(2,\tau_1)}_{yy}$ for
$y\in\gV_2$. Fix $M>0$ and suppose \eqref{eqn:upper_bound_H} and
\eqref{eqn:upper-bound-H-loewner-error} hold with the Loewner-error
constant at most $M$, and
\[
    \|\mZ-\mX^{(\tau_1)}\|\leq M\eta,
    \qquad d^{-1/2}\leq M\eta,
    \qquad 0<\eta\leq\eta_0,
\]
where $\eta_0>0$ is sufficiently small depending on $(M,\qr)$.
Then the following implications hold with $C=C_{M,\qr}$.
\begin{enumerate}
    \item \emph{Right edge.}
    Let $\mu\geq0$ and $w=(w^{(1)},w^{(2)})\in\C^n\oplus\C^m$
    satisfy $\|w^{(1)}\|^2=\|w^{(2)}\|^2=1/2$.
    If
    \begin{gather}
        \mZ w^{(2)}=\mu\mP^{(\tau_1)}w^{(1)}+r,
        \qquad\|r\|\leq M\eta,
        \qquad\sum_{x\in\gV^{(\vtau)}}(1+\alpha_x)|w_x|^2\leq M\eta,
    \end{gather}
    then the following inequality holds:
    \begin{align}
        \mu\leq\frac12\left(
          \redge+\qr\tau_1+\qr^{-1}\tau_2^{\mathrm R}\right)+C\eta.
        \label{eqn:right-regularized-quadratic-test}
    \end{align}

    \item \emph{Left edge.}
    Suppose also that $\qr>1$, $1-\tau_1\geq c_\qr>0$,
    $\|\mX^{(\tau_1)}\|\leq M$, and the lower form
    \eqref{eqn:lower_bound_H} holds with error at most $M\eta$.
    Assume the degree comparison
    $\beta_y\geq\alpha_y-M\eta$ for every $y\in\gV_2$.
    Let $0\leq\mu\leq\ledge$, $v\in\C^m$ with $\|v\|=1$,
    and $u\in\C^n$ with $\|u\|\leq1$. If
    \begin{gather}
        \mZ v=\mu\mP^{(\tau_1)}u+r,
        \qquad\|r\|\leq M\eta,
        \qquad
        \sum_{y\in\Vtwolow{\tau_2^{\mathrm L}}}
          (1+\alpha_y)|v_y|^2\leq M\eta,
    \end{gather}
    then the following inequality holds:
    \begin{align}
        \mu^2\geq(1-\tau_1)(\tau_2^{\mathrm L}-1)-C\eta.
        \label{eqn:left-complement-squared-edge}
    \end{align}
\end{enumerate}
\end{lemma}

\begin{proof}[Proof of Lemma~\ref{lem:semilocalized-complement-form-tests}]
\emph{Right edge.}
The discarded rows form $\Vonehigh{\tau_1}\subset\gV^{(\vtau)}$,
so the mass hypothesis gives
\[
    \|(\id_n-\mP^{(\tau_1)})w^{(1)}\|^2\leq M\eta,
    \qquad
    \|\mP^{(\tau_1)}w^{(1)}\|^2\geq\frac12-M\eta.
\]
Test the upper form at $h=(\qr w^{(1)},w^{(2)})$.
The first-block relation and pruning transfer yield
\begin{align}
    \langle h,\mH^{(\tau_1)}h\rangle
    &=2\qr\operatorname{Re}
      \langle w^{(1)},\mX^{(\tau_1)}w^{(2)}\rangle\\
    &=2\qr\mu\|\mP^{(\tau_1)}w^{(1)}\|^2
      +O_{M,\qr}(\eta)\\
    &\geq\qr\mu(1-2M\eta)-C\eta.
\end{align}
For the diagonal form, the retained row degrees satisfy
$\alpha_x<\tau_1$, while $\beta_y\leq\alpha_y$ for all columns.
Columns outside $\Vtwohigh{\tau_2^{\mathrm R}}$ have
$\alpha_y<\tau_2^{\mathrm R}$, and the weighted mass hypothesis
controls the contribution from this exceptional set. Hence
\begin{align}
    \langle h,(\id_N+d^{-1}\mD^{(\tau_1)})h\rangle
    &=\frac{\qr^2+1}{2}
      +\qr^2\sum_{x\in I_1}\alpha_x|w_x^{(1)}|^2
      +\sum_{y\in\gV_2}\beta_y|w_y^{(2)}|^2\\
    &\leq\frac{\qr^2+1+\qr^2\tau_1+\tau_2^{\mathrm R}}2+C\eta.
\end{align}
By \eqref{eqn:upper-bound-H-loewner-error} and
$d^{-1/2}\leq M\eta$, the Loewner error at $h$ is at most $C\eta$.
Combining the two form bounds gives
\[
    \qr\mu(1-2M\eta)
    \leq\frac{\qr^2+1+\qr^2\tau_1+\tau_2^{\mathrm R}}2+C\eta.
\]
Apply \eqref{eqn:semilocalized-scalar-upper-conclusion} with
$c=1$, $\theta=2M\eta$ and $\epsilon=C\eta$. Taking $\eta_0$
sufficiently small proves \eqref{eqn:right-regularized-quadratic-test}.

\emph{Left edge.}
The row margin gives $\|\mM^{(\tau_1)}\|\leq c_\qr^{-1}$.
Together with the pruning transfer and
$\|\mX^{(\tau_1)}\|\leq M$, this implies
\begin{align}
    &\left|\langle\mZ v,\mM^{(\tau_1)}\mZ v\rangle
      -\langle\mX^{(\tau_1)}v,
        \mM^{(\tau_1)}\mX^{(\tau_1)}v\rangle\right|
    \leq\|\mM^{(\tau_1)}\|\,
      \|\mZ-\mX^{(\tau_1)}\|\,
      \bigl(\|\mZ\|+\|\mX^{(\tau_1)}\|\bigr)
      \leq C\eta.
\end{align}
The lower form \eqref{eqn:lower_bound_H} therefore yields
\[
    \langle\mZ v,\mM^{(\tau_1)}\mZ v\rangle
    \geq\sum_{y\in\gV_2}(\beta_y-1)|v_y|^2-C\eta.
\]
On the other hand, the first-block relation, $\|u\|\leq1$, and
$\mu\leq\ledge$ give
\begin{align}
    \langle\mZ v,\mM^{(\tau_1)}\mZ v\rangle
    &=\mu^2\langle\mP^{(\tau_1)}u,
      \mM^{(\tau_1)}\mP^{(\tau_1)}u\rangle+O_{M,\qr}(\eta)
    \leq\frac{\mu^2}{1-\tau_1}+C\eta.
\end{align}
Outside $\Vtwolow{\tau_2^{\mathrm L}}$, the degree comparison gives
$\beta_y-1\geq\tau_2^{\mathrm L}-1-M\eta$, while
$\beta_y-1\geq-1$ everywhere. The low-column mass bound thus gives
\[
    \sum_{y\in\gV_2}(\beta_y-1)|v_y|^2
    \geq\tau_2^{\mathrm L}-1-C\eta.
\]
Apply \eqref{eqn:semilocalized-scalar-lower-conclusion} with
$b=\tau_2^{\mathrm L}-1$ and errors $C\eta$ to obtain
\eqref{eqn:left-complement-squared-edge}. Column normalization
suffices, including at $\mu=0$.
\end{proof}

\subsection{Complementary spectrum and comparison operators}
\label{subsec:semilocalized-comparison-inputs}

Window comparison uses the active profiles; location comparison selects
all outlier profiles at a deterministic gap. The mass estimates are
proved in Subsection~\ref{subsec:semilocalized-radial-mass-proof}, and
the comparison hypotheses are verified in
Subsection~\ref{subsec:semilocalized-comparison-proofs}.

\subsubsection{Comparison in spectral windows}
\label{subsec:semilocalized-window-comparison}

\begin{proposition}[Window comparison]
\label{prop:semilocalized-complement-spectrum}
The following estimates hold for the active centre sets in
\eqref{eqn:profile-centres}.
\begin{enumerate}
    \item For every $\nu>0$, uniformly for
    $\gap\geq\gC_{B,\nu}\error_d^{1/4}$, with probability at least
    $1-C_{B,\qr,\nu}N^{-\nu}$,
    \begin{align}
        \|\scA-\hatH_{\mathrm R}(\gap)\|
        &\leq C_{B,\qr,\nu}
          (\error_d+\error_{\mathrm{prof}}(\gap)),
        \label{eqn:full-block-approximation-right}\\
        \left\|\overline\Pi_{\mathrm R}^{(\vtau)}\scA^{(\vtau)}
          \overline\Pi_{\mathrm R}^{(\vtau)}\right\|
        &\leq\redge+\frac{s(\gap)}2
          +C_{B,\qr,\nu}(\error_d+\error_{\mathrm{prof}}(\gap)).
        \label{eqn:right-complement-spectrum}
    \end{align}
    \item At critical sparsity, for $\qr>1$ and any fixed $b>0$,
    uniformly for $0<\gap<\ledge$ satisfying
    $\gap\geq\gC\error_d^{1/4}$, with probability at least
    $1-CN^{-\nu_0}$,
    \begin{align}
        \|\scA-\hatH_{\mathrm L}(\gap)\|
        &\leq C_{b,\qr,\sK}
          (\error_d+\error_{\mathrm{prof}}(\gap)),
        \label{eqn:full-block-approximation-left}\\
        \operatorname{Spec}_+\!\left(
          \overline\Pi_{\mathrm L}^{(\vtau)}\scA^{(\vtau)}
          \overline\Pi_{\mathrm L}^{(\vtau)}\right)
        &\cap\left[\frac{3\inf\sK}{4},\ledge-\frac{s(\gap)}4
          -C_{b,\qr,\sK}(\error_d+\error_{\mathrm{prof}}(\gap))\right)
          =\varnothing.
        \label{eqn:left-complement-spectrum}
    \end{align}
\end{enumerate}
On these events, the pruned block error and compression defect in
\eqref{eqn:pruned-profile-block-error} are bounded by $C\error_d$,
with the parameter dependence of the corresponding clause.
\end{proposition}
The left statement concerns a positive spectral window and permits
zero and smaller positive complementary eigenvalues.

\subsubsection{Comparison at the location scale}
\label{subsec:semilocalized-location-comparison}

For the eigenvalue locations, prune at
\begin{align}
    \gap&=\frac34\error_d^{1/4},
&\vtau&=\vtau\!\left(\frac34\error_d^{1/4}\right),
&
\error_d+\error_{\mathrm{prof}}\!\left(\frac34\error_d^{1/4}\right)
    &=\error_d+\frac{4\error_d^{3/4}}{3\log d}.
    \label{eqn:calibrated-comparison-scales}
\end{align}
This gap is smaller than the large-constant lower bound in the
uniform window inputs. Lemma~\ref{lem:shrinking-gap-profile-inputs},
proved in Subsubsection~\ref{subsec:semilocalized-shrinking-gap-proofs},
verifies the pruning and profile estimates directly at this scale.

For the right comparison, select every centre in
$\mathcal C_{\mathrm R}=\gR_1\cup\gR_2$. Choosing the standing
constant sufficiently large makes
$\gap=\frac34\error_d^{1/4}\leq\min\{1,4/\qr\}$.
The degree displacements on both right branches exceed one, so
the pruning cutoffs retain all of $\gR_1\cup\gR_2$.
For the global left comparison, use the same profile recipe at every
$x\in\gL_2$ from \eqref{eqn:left-outlier-profile-set}. With
the pruning thresholds in \eqref{eqn:calibrated-comparison-scales}, abbreviate
$Q=Q_{\gL_2}^{(\vtau)}$ and
$\hatH_{\mathrm L,\mathrm{edge}}=\hatH_{\gL_2}$.
This centre set has no fixed positive lower profile cutoff and differs
from $\gV_{\mathrm L}^\star$ in the window comparison.
Parity balance gives $Q=Q_1\oplus Q_2$, with
$\operatorname{rank}(\id_m-Q_2)=|\gL_2|$ and
$\dim\operatorname{Ran}Q_2=m-|\gL_2|$. Its rectangular
complementary block is
\[
    X_Q\coloneqq Q_1\rmX^{(\vtau)}Q_2:
    \operatorname{Ran}Q_2\longrightarrow\operatorname{Ran}Q_1.
\]

\begin{corollary}[Location scale comparison]
\label{cor:semilocalized-calibrated-comparison}
Use the pruning scale and error in \eqref{eqn:calibrated-comparison-scales}.
\begin{enumerate}
    \item \emph{Right edge.}
    For every $\nu>0$, choose the standing constant $\gC_0$ in
    \eqref{eqn:subcritical_regime} sufficiently large depending on
    $(B,\qr,\nu)$. With probability at least
    $1-C_{B,\qr,\nu}N^{-\nu}$,
    \begin{align}
        \|\scA-\hatH_{\mathrm R}\!\left(\frac34\error_d^{1/4}\right)\|
        &\leq K_{\mathrm{bl}}\left(\error_d+\frac{4\error_d^{3/4}}{3\log d}\right),
        \label{eqn:right-location-block-error}\\
        \left\|\overline\Pi_{\mathrm R}^{(\vtau)}
          \scA^{(\vtau)}\overline\Pi_{\mathrm R}^{(\vtau)}\right\|
        &\leq\redge+\frac34\error_d^{1/4}+K_{\mathrm{comp}}\left(\error_d+\frac{4\error_d^{3/4}}{3\log d}\right).
        \label{eqn:right-fixed-gap-complement}
    \end{align}
    Here $K_{\mathrm{bl}},K_{\mathrm{comp}}$ depend on $(B,\qr,\nu)$,
    and the profiles are selected at all $x\in\gR_1\cup\gR_2$.
    \item \emph{Global left edge.}
    Under critical sparsity, with $\qr>1$ and $b>\qr^{-2}$, there
    are constants $\nu_0,C,K_{\mathrm{bl}}>0$ depending on
    $(b,\qr)$ such that, with probability at least $1-CN^{-\nu_0}$,
    the profiles at all $x\in\gL_2$ are orthonormal and
    \begin{align}
        \|\scA-\hatH_{\mathrm L,\mathrm{edge}}\|
        &\leq K_{\mathrm{bl}}\left(\error_d+\frac{4\error_d^{3/4}}{3\log d}\right),
        \label{eqn:left-location-block-error}\\
        \|X_Qv\|&\geq(\ledge-\frac34\error_d^{1/4})\|v\|,
        \qquad v\in\operatorname{Ran}Q_2.
        \label{eqn:left-edge-complement-column-bound}
    \end{align}
    The comparison retains all $m$ column singular-value positions.
\end{enumerate}
\end{corollary}
The proof is given in Subsubsection~\ref{subsec:semilocalized-calibrated-proofs}.

\subsection{Proofs of outlier locations}
\label{subsec:semilocalized-eigenvalue-locations}

\begin{proof}[Proof of Theorems~\ref{thm:right_edge_behavior}
and~\ref{thm:left_edge_behavior}]
\emph{Right edge.}
Fix $\nu>0$ and work on the event in
Corollary~\ref{cor:semilocalized-calibrated-comparison}.
Write $r_{\mathrm R}=|\gR_1|+|\gR_2|$ and order the positive
profile values as $\Lambda_{\mathrm R,1}\geq\cdots\geq
\Lambda_{\mathrm R,r_{\mathrm R}}$.
Apply \eqref{eqn:ordered-profile-transfer} to the comparison operator
in \eqref{eqn:right-location-block-error}, obtaining distinct
positions $j_1<\cdots<j_{r_{\mathrm R}}\leq m$, with multiplicity:
\begin{align}
    |\lambda_{j_k}(\scA)-\Lambda_{\mathrm R,k}|
    \leq K_{\mathrm{bl}}\left(\error_d+\frac{4\error_d^{3/4}}{3\log d}\right),
    \qquad 1\leq k\leq r_{\mathrm R}.
    \label{eqn:right-fixed-gap-profile-matching}
\end{align}
By \eqref{eqn:calibrated-comparison-scales}, this matching has the error
scale in \eqref{eqn:right-outlier-locations}.

Only the $r_{\mathrm R}$ positive profile values can exceed the
complementary spectral bound in \eqref{eqn:right-fixed-gap-complement}.
If $r_{\mathrm R}<m$, min--max and Weyl's inequality give
\begin{align}
    \lambda_{r_{\mathrm R}+1}(\scA)
    \leq\redge+\frac34\error_d^{1/4}
      +(K_{\mathrm{comp}}+K_{\mathrm{bl}})\left(\error_d+\frac{4\error_d^{3/4}}{3\log d}\right).
    \label{eqn:right-fixed-gap-final-comparison}
\end{align}
The ratio in \eqref{eqn:location-error-gap-ratio} is uniformly small
when $\gC_0$ is large. Choose it so that the combined block and
complementary error is at most $\error_d^{1/4}/4$.
This proves \eqref{eqn:right-no-additional-outliers}.

\emph{Left edge.}
Assume $b>\qr^{-2}$ and use the global left-edge event in
Corollary~\ref{cor:semilocalized-calibrated-comparison}.
Write $r=|\gL_2|$. The comparison operator
$\hatH_{\mathrm L,\mathrm{edge}}$ has the $r$ positive profile values
and the $m-r$ complementary column singular values of $X_Q$.
Their order may interleave, but \eqref{eqn:left-edge-complement-column-bound}
places every complementary value above $\ledge-\frac34\error_d^{1/4}$.
If $r<m$, min--max and \eqref{eqn:left-location-block-error} yield
\begin{align}
    \lambda_{m-r}(\scA)
    \geq\ledge-\frac34\error_d^{1/4}-K_{\mathrm{bl}}\left(\error_d+\frac{4\error_d^{3/4}}{3\log d}\right)
    \geq\ledge-\error_d^{1/4}
\end{align}
for large $N$, by the critical limit in
\eqref{eqn:location-error-gap-ratio}.
This proves \eqref{eqn:left-no-additional-outliers}.

Apply \eqref{eqn:ordered-profile-transfer} to the $r$ ordered
profile values $\Lambda_{\mathrm L,k}$ and all $m$ column positions
of $\hatH_{\mathrm L,\mathrm{edge}}$. With multiplicity and possible
interleaving, their distinct positions $j_1<\cdots<j_r\leq m$ satisfy
\begin{align}
    |\lambda_{j_k}(\scA)-\Lambda_{\mathrm L,k}|
    \leq K_{\mathrm{bl}}\left(\error_d+\frac{4\error_d^{3/4}}{3\log d}\right),
    \qquad 1\leq k\leq r.
    \label{eqn:left-fixed-gap-profile-matching}
\end{align}
Under critical sparsity,
\[
    \error_d^{1/4}\log d
    =\frac{(\log d)^{5/4}}{b^{1/4}(\log N)^{1/8}}
    \longrightarrow0.
\]
Thus $\error_d=o(\error_d^{3/4}/\log d)$ and
$\error_d+\frac{4\error_d^{3/4}}{3\log d}=O(\error_d^{3/4}/\log d)$, which is the error scale in
\eqref{eqn:left-outlier-locations}.
\end{proof}

\subsection{Proofs of eigenvector semilocalization}
\label{subsec:semilocalized-eigenvector-proofs}

\subsubsection{Resonant profiles and spectral projection}
\label{subsec:semilocalized-eigenvector-semilocalization}

For $\lambda>0$ and $\delta>0$, define the resonant sets
\begin{subequations}
\begin{align}
    \gW_{\lambda,\delta}^{(1)}
    \coloneqq{}&
    \{x\in\gV_1:\lambda>\redge,\ \alpha_x>\qr^{-2}+1,
      \ |\lambda-\Lambda_{\qr^{-1}}(\alpha_x)|\leq\delta\},
    \label{eqn:resonant_vertices_1}\\
    \gW_{\lambda,\delta}^{(2,\mathrm{R})}
    \coloneqq{}&
    \{x\in\gV_2:\lambda>\redge,\ \alpha_x>\qr^2+1,
      \ |\lambda-\Lambda_{\qr}(\alpha_x)|\leq\delta\},
    \label{eqn:resonant_vertices_2_right}\\
    \gW_{\lambda,\delta}^{(2,\mathrm{L})}
    \coloneqq{}&
    \{x\in\gV_2:0<\lambda<\ledge,\ 0<\alpha_x<\qr^2-1,
      \ |\lambda-\Lambda_{\qr}(\alpha_x)|\leq\delta\},
    \label{eqn:resonant_vertices_2_left}
\end{align}
\end{subequations}

Using the gaps in \eqref{eqn:right-gap} and \eqref{eqn:left-gap}, set
\begin{align}
    (\mathrm E,\gap,\mathcal W_{\lambda,\delta})
    \coloneqq
    \begin{cases}
      \bigl(\mathrm R,\gapR,
      \gW_{\lambda,\delta}^{(1)}\cup
      \gW_{\lambda,\delta}^{(2,\mathrm{R})}\bigr),
      &\lambda>\redge,\\
      \bigl(\mathrm L,\gapL,
      \gW_{\lambda,\delta}^{(2,\mathrm{L})}\bigr),
      &0<\lambda<\ledge.
    \end{cases}
\end{align}
Let $(\lambda,\rvw)$ be a normalized eigenpair of $\scA$, with
$\lambda\geq\redge+C\error_d^{1/4}$ on the right or
$\lambda\in\sK$ on the left. On the window-comparison and profile
events, the window restrictions
$0<\delta\leq\gapR/2$ on the right and
$0<\delta\leq\frac12\min\{\lambda,\gapL\}$ on the left, together
with the degree--eigenvalue check in
\eqref{eqn:profile-energy-degree-separation}, give the required
degree cutoffs. On the left, these restrictions also place every
resonant value between $\inf\sK/2$ and
$\ledge-2\error_d^{1/4}$ for sufficiently large $N$. Thus
$\mathcal W_{\lambda,\delta}\subset\gV_{\mathrm E}^\star$.
Only positive profiles are resonant. Define
\begin{align}
    \Pi_{\lambda,\delta}^{(\vtau)}
    &\coloneqq
    \sum_{x\in\mathcal W_{\lambda,\delta}}
      \hatvtau_+(x)\hatvtau_+(x)^*,
    &\overline\Pi_{\lambda,\delta}^{(\vtau)}
    &\coloneqq\id_N-\Pi_{\lambda,\delta}^{(\vtau)}.
    \label{eqn:resonant-profile-projection}
\end{align}
This projection commutes with $\hatH_{\mathrm E}(\gap)$.
If $\delta\leq C(\error_d+\error_{\mathrm{prof}}(\gap))$, the
mass bounds below are trivial after increasing their constants.
Assume otherwise.

Every positive profile eigenvalue of $\hatH_{\mathrm E}(\gap)$
outside the resonant projection is more than $\delta$ from $\lambda$,
and every negative profile eigenvalue is at least $\lambda$ away.
At the right edge, choose the error constants so that
$C(\error_d+\error_{\mathrm{prof}}(\gap))\leq\gap/4$.
Equation~\eqref{eqn:right-complement-spectrum} separates the
remaining block from $\lambda=\redge+\gap$ by at least
$\gap-s/2-C(\error_d+\error_{\mathrm{prof}}(\gap))
\geq\gap/4\geq\delta/2$.
At the left edge, \eqref{eqn:left-complement-spectrum} separates its
spectrum above $3\inf\sK/4$ by at least $c\gap$; the remaining
positive spectrum is at least $\inf\sK/4$ away, and the zero and
negative spectra are at least $\lambda\geq2\delta$ away.
Consequently,
\begin{align}
    \operatorname{dist}\left(
      \lambda,\operatorname{Spec}\!\left(
        \hatH_{\mathrm E}(\gap)|_{\operatorname{Ran}
        \overline\Pi_{\lambda,\delta}^{(\vtau)}}\right)
    \right)
    \geq c\delta.
    \label{eqn:resonant-complement-gap}
\end{align}
Apply the eigenvector-mass clause of Lemma~\ref{lem:spectral-counting}
with $\rmS=\scA$, $\rmT=\hatH_{\mathrm E}(\gap)$,
$P=\Pi_{\lambda,\delta}^{(\vtau)}$, and gap $c\delta$.
The full block approximation and
\eqref{eqn:general-comparison-profile-mass} yield
\begin{align}
    \|\overline\Pi_{\lambda,\delta}^{(\vtau)}\rvw\|
    \leq C\frac{\error_d+\error_{\mathrm{prof}}(\gap)}{\delta}.
    \label{eqn:semilocalized-complement-mass}
\end{align}
By orthonormality, for every $\nu>0$ there is a constant
$\gC_{B,\nu}>0$ such that, with probability at least
$1-\gC_{B,\nu}N^{-\nu}$,
\begin{align}
    \sum_{x\in \gW_{\lambda,\delta}^{(1)}\cup
    \gW_{\lambda,\delta}^{(2,\mathrm{R})}}
    |\<\hatvtau_+(x),\rvw\>|^{2}
    \geq
    1-\gC_{B,\nu}
    \left(
      \frac{\error_d+\error_{\mathrm{prof}}(\gapR)}{\delta}
    \right)^{2}.
    \label{eqn:right-semilocalized-mass}
\end{align}
This holds simultaneously for every normalized eigenpair
$(\lambda,\rvw)$ with
$\lambda\geq\redge+\gC_{B,\nu}\error_d^{1/4}$ and every
$0<\delta\leq\gapR/2$.

On the fixed left interval $\sK$, the gap satisfies
$\gapL\geq\ledge-\sup\sK>0$, so
\eqref{eqn:profile-error-scale} gives
$\error_d+\error_{\mathrm{prof}}(\gapL)\leq C_{\sK}\error_d$.
Under the standing left-edge assumptions, there are constants
$\nu_0=\nu_0(b,\qr,\sK)>0$ and
$\gC=\gC(b,\qr,\sK)>0$ such that, with probability at least
$1-\gC N^{-\nu_0}$,
\begin{align}
    \sum_{x\in\gW_{\lambda,\delta}^{(2,\mathrm{L})}}
    |\<\hatvtau_+(x),\rvw\>|^{2}
    \geq
    1-\gC
    \left(
      \frac{\error_d}{\delta}
    \right)^{2}.
    \label{eqn:left-semilocalized-mass}
\end{align}
This holds simultaneously for every normalized eigenpair
$(\lambda,\rvw)$ with $\lambda\in\sK$ and every
$0<\delta\leq\frac12\min\{\lambda,\gapL\}$.
The input events can be intersected using larger failure exponents
on the right and, if necessary, a smaller $\nu_0$ on the left.

\subsubsection{From resonant profiles to semilocalization}

The degree counts in
Section~\ref{sec:transition-from-delocalization-to-localization} combine
with \eqref{eqn:profile-support-size} to bound resonant support.

\paragraph{Resonance width and half mass.}
Choose $K_{\mathrm R}$ and $K_{\sK}$ sufficiently large for the
respective constants in \eqref{eqn:right-semilocalized-mass} and
\eqref{eqn:left-semilocalized-mass}, and set
\begin{align}
 \delta_\lambda\coloneqq
 \begin{cases}
 K_{\mathrm R}(\error_d+\error_{\mathrm{prof}}(\gapR)),
     &\lambda>\redge,\\
 K_{\sK}\error_d,&\lambda\in\sK.
 \end{cases}
 \label{eqn:lp-resonance-width}
\end{align}
For any prescribed right failure exponent choose $K_{\mathrm R}$
on an input event with a larger exponent; $K_{\sK}$ uses the fixed
left exponent. Once the corresponding width restriction is checked,
those mass estimates give, simultaneously for all eigenpairs under
consideration,
\begin{align}
 \sum_{x\in\mathcal W_{\lambda,\delta_\lambda}}
 |\langle\hatvtau_+(x),\rvw\rangle|^2\geq\frac12.
 \label{eqn:right-uniform-resonant-mass}
\end{align}

\paragraph{Bessel--H\"older conversion.}
Let
$\{\phi_x:x\in\mathcal W\}$ be a finite nonempty family of orthonormal
profiles, let $S=\bigcup_{x\in\mathcal W}\operatorname{supp}\phi_x$,
and suppose that $\sum_{x\in\mathcal W}|\langle\phi_x,\rvw\rangle|^2
\geq M$. Bessel's inequality applied to the coordinate restriction
$\rvw|_S$ gives
\begin{align}
    M
    &\leq\sum_{x\in\mathcal W}|\langle\phi_x,\rvw\rangle|^2
      =\sum_{x\in\mathcal W}|\langle\phi_x,\rvw|_S\rangle|^2
      \leq\|\rvw|_S\|_2^2.
    \label{eqn:right-uniform-bessel-step}
\end{align}
For $2<p<\infty$, H\"older's inequality gives
\begin{align}
    \|\rvw|_S\|_2^2
    &\leq\left(\sum_{y\in S}|\rvw_y|^p\right)^{2/p}
       |S|^{1-2/p}
      \leq\|\rvw\|_p^2|S|^{1-2/p}.
    \label{eqn:right-uniform-holder-step}
\end{align}
The same bound holds at $p=2$ by equality and at $p=\infty$ by
$\sum_{y\in S}|\rvw_y|^2\leq|S|\|\rvw\|_\infty^2$.
Consequently, $\|\rvw\|_p^2\geq M|S|^{2/p-1}$ for every
$2\leq p\leq\infty$.

\begin{proof}[Proof of the noncritical bound
\eqref{eqn:right-lp-semilocalization-sublog}]
Fix $\nu>0$ and use input events with failure exponent larger than
$\nu$. Write $A\coloneqq\error_d$, $K=K_{\mathrm R}=K_{B,\nu}$,
and use the right width in \eqref{eqn:lp-resonance-width}.
We check the restriction $\delta_\lambda\leq\gapR/2$.
By \eqref{eqn:profile-error-scale},
$\error_{\mathrm{prof}}(\gapR)
=A/(\gapR\log d)$.
If $\gapR\geq C A^{1/4}$, then
\begin{align}
    \frac{\delta_\lambda}{\gapR}
    &=K\left(\frac{A}{\gapR}
       +\frac{A}{\gapR^2\log d}\right)
      \leq
      K\left(\frac{A^{3/4}}{C}
       +\frac{A^{1/2}}{C^2\log d}\right).
    \label{eqn:right-uniform-lp-window-check}
\end{align}
The lower bound on $d$ in \eqref{eqn:subcritical_regime} gives
$A\leq C_B/\gC_0$ for large $N$. Taking the standing constant
$\gC_0$ and then the spectral threshold $C$ sufficiently large
makes the right-hand side of
\eqref{eqn:right-uniform-lp-window-check} at most $1/2$.
Equation~\eqref{eqn:right-uniform-resonant-mass} now applies
simultaneously throughout the stated spectral range.

We next count resonant centres without restricting $\lambda$ to a
bounded interval. Put $s_1=\qr^{-1}$ and $s_2=\qr$ and use
Lemma~\ref{lem:right-uniform-degree-tail-count}.

Set $t_\lambda\coloneqq\lambda-\delta_\lambda\geq\redge$.
On the right branches the location maps are increasing. Hence
resonance implies
$\alpha_x\geq\invR{s_i}(t_\lambda)$
for every $x\in\mathcal W_{\lambda,\delta_\lambda}\cap\gV_i$.
Applying \eqref{eqn:right-uniform-degree-tail-count} to both
classes and using \eqref{eqn:right-degree-rate} yields
\begin{align}
    |\mathcal W_{\lambda,\delta_\lambda}|
    &\leq C_{\qr,\nu}
       \bigl(N\mathrm e^{-d\rateR(t_\lambda;\qr)}
             +\log N\bigr)
    \label{eqn:right-uniform-resonant-count-general}
\end{align}
simultaneously for all $\lambda$ in the stated range.

Let $Q_N\coloneqq\sqrt{d\log d}+\log\log N$ and let
$S_\lambda$ be the union of the supports of the resonant profiles.
The estimates \eqref{eqn:profile-support-size} and
\eqref{eqn:right-uniform-resonant-count-general} give
\begin{align}
    |S_\lambda|
    &\leq C_{\qr,\nu}\mathrm e^{C_{B,\qr,\nu}Q_N}
       \bigl(N\mathrm e^{-d\rateR(t_\lambda;\qr)}
             +\log N\bigr)
      \nonumber\\
    &\leq
       \exp\left[
         \bigl(\log N-d\rateR(t_\lambda;\qr)\bigr)_+
         +C_{B,\qr,\nu}Q_N
       \right].
    \label{eqn:right-uniform-resonant-support}
\end{align}
For the second inequality, put
$h=\log N-d\rateR(t_\lambda;\qr)$ and use
$\log(\mathrm e^h+\log N)
\leq h_++\log(1+\log N)$; the last logarithm and the
multiplicative constant are absorbed into $C_{B,\qr,\nu}Q_N$.

Apply \eqref{eqn:right-uniform-bessel-step} and
\eqref{eqn:right-uniform-holder-step} with
$\phi_x=\hatvtau_+(x)$, $S=S_\lambda$, and $M=1/2$.
Together with \eqref{eqn:right-uniform-resonant-support},
this gives
\begin{align}
    \|\rvw\|_p^2
    &\geq\frac12|S_\lambda|^{2/p-1}
      \nonumber\\
    &\geq\frac12 N^{2/p-1}\exp\left[
       \left(1-\frac2p\right)
        \left\{
          \min\bigl\{
            d\rateR(t_\lambda;\qr),\log N
          \bigr\}
          -C_{B,\qr,\nu}Q_N
        \right\}
      \right],
    \label{eqn:right-uniform-lp-conclusion}
\end{align}
Here the second line follows from
$(\log N-d\rateR(t_\lambda;\qr))_+
=\log N-\min\{d\rateR(t_\lambda;\qr),\log N\}$
and $2/p-1=-(1-2/p)$. This is
\eqref{eqn:right-lp-semilocalization-sublog}.
The simplified subcritical consequence is uniform in the whole stated
spectral range. For $s\in\{\qr^{-1},\qr\}$ and $0<h\leq1$,
\eqref{eqn:Lambda-inverse-right-edge} gives
\[
 \invR{s}(\redge+h)-(s^2+1)\geq2\sqrt h.
\]
Since $I_s'(\alpha)\geq\log(1+s^{-2})$ for
$\alpha\geq s^2+1$, both branches, and hence their minimum,
increase by at least $c_\Gamma\sqrt h$.
The half-gap bound gives $t_\lambda-\redge\geq c\error_d^{1/4}$;
monotonicity also handles $t_\lambda-\redge>1$. Thus
\[
 d\bigl[\rateR(t_\lambda;\qr)-\rateR(\redge;\qr)\bigr]
 \geq c d\error_d^{1/8}\gg Q_N,
\]
because
\[
 \frac{d\error_d^{1/8}}{\sqrt{d\log d}}
 =\frac{(\log N)^{1/16}d^{3/8}}{(\log d)^{3/8}}
 \longrightarrow\infty,
 \qquad
 \frac{d\error_d^{1/8}}{\log\log N}\longrightarrow\infty.
\]
If the minimum in \eqref{eqn:right-uniform-lp-conclusion} is
$\log N$, then $d=o(\log N)$ gives
$d\rateR(\redge;\qr)+C Q_N=o(\log N)$.
The two cases prove the stated bound at $p=\infty$.

\end{proof}

\begin{proof}[Proof of the critical case of
Theorem~\ref{thm:right_vector_semilocalization} and
Theorem~\ref{thm:left_vector_semilocalization}]
Assume critical sparsity and fix $\epsilon>0$. Put $\eta=\epsilon/4$.
Use the resonance widths in \eqref{eqn:lp-resonance-width}.

\paragraph{Right edge.}
Fix $\nu>0$ and intersect the right window-comparison event with those in
Lemmas~\ref{lem:bipartite-max-degree}
and~\ref{lem:semilocalized-profile-properties}, using failure
exponents larger than $\nu$.
At the deterministic gap $\gap=1$, the critical maximum-degree
bound makes every ideal profile value bounded. The window comparison
bounds its complementary spectrum, and its full block error is
$O(\error_d)$. Hence $\|\scA\|\leq E_{\max}$ for a fixed
$E_{\max}=E_{\max}(b,\qr,\nu)$, and all eigenvalues under
consideration lie in $[\redge,E_{\max}]$. Write $A=\error_d$, which
tends to zero at critical sparsity. For every fixed $c>0$,
\eqref{eqn:right-uniform-lp-window-check}, with $K=K_{\mathrm R}$
and $C=c$, gives
$\delta_\lambda/\gapR\to0$ uniformly when $\gapR\geq cA^{1/4}$.
The additional width condition follows from
\begin{align}
    \frac{\delta_\lambda}{\sqrt{\gapR}}
    &\leq K_{\mathrm R}
       \left(\frac{A^{7/8}}{\sqrt c}
       +\frac{A^{5/8}}{c^{3/2}\log d}\right)\longrightarrow0.
    \label{eqn:critical-right-lp-window-check}
\end{align}
Thus $\delta_\lambda\leq\gapR/2$ for all sufficiently large $N$,
and the hypotheses \eqref{eqn:right-uniform-resonance-window}
hold for $\gamma_N=c\error_d^{1/4}$.
Apply the right clause of
Corollary~\ref{cor:critical-resonant-upper-counts} with accuracy
$\eta$ and failure exponent larger than $\nu$, choosing $c$ to be
the spectral threshold in the theorem. The resulting intersection
has probability at least
$1-\gC_{B,\epsilon,\nu}N^{-\nu}$, and both the mass estimate and
the count estimate hold simultaneously throughout the stated
right spectral range.

\paragraph{Left edge.}
On $\sK$, the width $\delta_\lambda=K_{\sK}\error_d$ tends to zero,
whereas $\lambda$ and $\gapL$ are bounded below by positive
constants. Hence
$\delta_\lambda\leq\frac12\min\{\lambda,\gapL\}$ for all
sufficiently large $N$, and
\eqref{eqn:left-uniform-resonance-window} holds with $K_0=\sK$.
Apply the left clause of
Corollary~\ref{cor:critical-resonant-upper-counts} with accuracy
$\eta$ and a failure exponent larger than $\nu_0$.
Intersect this event with the left window-comparison and profile
events. After decreasing $\nu_0$ to absorb the profile event,
the intersection has probability at least
$1-\gC_\epsilon N^{-\nu_0}$.
Both the mass estimate and the count estimate hold simultaneously
for every $\lambda\in\sK$.

The two clauses of the general critical count estimate give,
on the respective events,
\begin{align}
    |\mathcal W_{\lambda,\delta_\lambda}|
    \leq N^{\rho_{\mathrm E}(\lambda,b)+\eta},
    \qquad \mathrm E\in\{\mathrm L,\mathrm R\}.
    \label{eqn:semilocalized-critical-upper-count}
\end{align}

\paragraph{From resonant mass to $\ell^p$ mass.}
Combine the half-mass bound
\eqref{eqn:right-uniform-resonant-mass} with the centre count
\eqref{eqn:semilocalized-critical-upper-count}.
Each profile has support of size $N^{o(1)}$ by
\eqref{eqn:profile-support-size}. Hence
\begin{align}
    S_\lambda\coloneqq
      \bigcup_{x\in\mathcal W_{\lambda,\delta_\lambda}}
      \operatorname{supp}\hatvtau_+(x)
    \qquad\text{satisfies}\qquad
    |S_\lambda|\leq N^{\rho_{\mathrm E}(\lambda,b)+2\eta}.
    \label{eqn:lp-profile-support}
\end{align}
Apply \eqref{eqn:right-uniform-bessel-step} and
\eqref{eqn:right-uniform-holder-step} with
$\phi_x=\hatvtau_+(x)$, $M=1/2$, and $S=S_\lambda$.
For every $2\leq p\leq\infty$, they give
\begin{align}
    \|\rvw\|_p^2
    \geq |S_\lambda|^{2/p-1}\|\rvw|_{S_\lambda}\|_2^2
    \geq\frac12N^{(2/p-1)(\rho_{\mathrm E}(\lambda,b)+2\eta)}
    \geq N^{(2/p-1)\rho_{\mathrm E}(\lambda,b)-\epsilon}
    \label{eqn:profile-mass-to-lp}
\end{align}
for large $N$, uniformly in $\lambda$ and $p$. This proves
\eqref{eqn:right-lp-semilocalization-critical} and
\eqref{eqn:left-lp-semilocalization}, including $p=\infty$.
\end{proof}

%% file: transition.tex
\section{Proofs of Transition Theorems}\label{sec:transition-from-delocalization-to-localization}

This section establishes the degree-counting and spectral-comparison
lemmas used in the transition proofs and in
Section~\ref{sec:semilocalized_phase}. We then identify the spectral
boundaries and prove the transition theorems. The degree laws and graph
estimates are collected in
Section~\ref{sec:bipartite-random-graph-inputs}; the auxiliary proofs are
deferred to Section~\ref{sec:transition-deferred-proofs}.

\subsection{General degree estimates}
\label{subsec:transition-general-degree-inputs}

For $i\in\{1,2\}$, let $\rD^{(i)}$ be the degree of a fixed vertex in
$\gV_i$, so that $\rD^{(1)}\sim\Bin(m,p)$ and
$\rD^{(2)}\sim\Bin(n,p)$. Write $d_i\coloneqq\E\rD^{(i)}$, and set
$s_1=\qr^{-1}$ and $s_2=\qr$. For an interval $J\subset[0,\infty)$,
including a singleton or a one-sided interval, define
\begin{align}
    X_i(J)&\coloneqq\#\{x\in\gV_i:\alpha_x\in J\},
    &\mu_i(J)&\coloneqq|\gV_i|\P(\rD^{(i)}/d\in J).
    \label{eqn:general-degree-counts}
\end{align}
\begin{proposition}[Degree counts and binomial tails]
\label{prop:uniform-degree-counts}
Suppose Assumption~\ref{ass:proportional_regime} holds and $0<p\leq1$.
\begin{enumerate}
    \item For every $\nu>0$, with probability at least
    $1-C_\nu N^{-\nu}$, simultaneously for $i\in\{1,2\}$ and every
    interval $J\subset[0,\infty)$,
    \begin{align}
        |X_i(J)-\mu_i(J)|
        &\leq C_\nu\bigl(\sqrt{\mu_i(J)\log N}+\log N\bigr),
        \label{eqn:uniform-degree-count-concentration}\\
        X_i(J)&\leq C_\nu\bigl(\mu_i(J)+\log N\bigr).
        \label{eqn:uniform-degree-count-upper}
    \end{align}
    In particular, $X_i(J)/\mu_i(J)=1+o(1)$ uniformly over any
    deterministic family of intervals whose smallest mean is
    $\gg\log N$.
    \item For each deterministic interval $J$,
    \begin{align}
        \P(X_i(J)\geq1)&\leq\mu_i(J),
        &\P(X_i(J)=0)&\leq\exp[-\mu_i(J)].
        \label{eqn:degree-window-presence-absence}
    \end{align}
    \item For any binomial random variable $Z$ with mean $\mu>0$, put
    \[
        \Psi_\mu(k)\coloneqq k\log(k/\mu)-k+\mu,
        \qquad k\geq0,
    \]
    with $0\log0=0$. For every nonnegative integer $k$,
    \begin{align}
        \P(Z\geq k)&\leq\mathrm e^{-\Psi_\mu(k)},
            &&k\geq\mu,
        \label{eqn:degree-chernoff-upper}\\
        \P(Z\leq k)&\leq\mathrm e^{-\Psi_\mu(k)},
            &&k\leq\mu.
        \label{eqn:degree-chernoff-lower}
    \end{align}
    In particular, these bounds apply to $Z=\rD^{(i)}$ with $\mu=d_i$.
    Write
    \begin{align}
        \Psi_i(k)\coloneqq\Psi_{d_i}(k)
        =k\log(k/d_i)-k+d_i,
        \qquad k\geq0.
        \label{eqn:exact-degree-tail-rate}
    \end{align}
\end{enumerate}
\end{proposition}

For a bounded degree window $J$, set
$J_d\coloneqq J\cap d^{-1}\mathbb Z_{\geq0}$.

\begin{lemma}[Sparse degree-window means]
\label{lem:sparse-degree-window-means}
Under Assumption~\ref{ass:proportional_regime}, suppose
$1\ll d\leq B\log N$ for a fixed $B>0$.
For each fixed $A>0$, uniformly over $i\in\{1,2\}$ and intervals
$J\subset[0,A]$ with $J_d\ne\varnothing$,
\begin{align}
    \log\mu_i(J)
    =\log|\gV_i|-d\min_{a\in J_d}I_{s_i}(a)
      +O_{\qr,A,B}(\log d).
    \label{eqn:sparse-degree-window-mean}
\end{align}
An empty lattice window has $\mu_i(J)=X_i(J)=0$.
Moreover, for every fixed $a>0$ on the indicated side of the mean,
\begin{align}
    \mu_i([a,\infty))&\leq C_\qr N\mathrm e^{-dI_{s_i}(a)},
        &&a>s_i^2,
        \notag\\
    \mu_i([0,a])&\leq C_\qr N\mathrm e^{-dI_{s_i}(a)},
        &&0<a<s_i^2.
    \label{eqn:sparse-degree-tail-mean-bound}
\end{align}
The last bounds hold uniformly when $a$ ranges over a fixed compact
subset of the respective open intervals.
\end{lemma}
\begin{lemma}[Uniform right-edge degree tails]
\label{lem:right-uniform-degree-tail-count}
Under Assumption~\ref{ass:proportional_regime}, suppose
$1\ll d\leq B\log N$. For every $\nu>0$, with probability at least
$1-C_\nu N^{-\nu}$, simultaneously for $i\in\{1,2\}$ and every
$a\geq s_i^2+1$,
\begin{align}
    \#\{x\in\gV_i:\alpha_x\geq a\}
    \leq C_{\qr,\nu}\bigl(N\mathrm e^{-dI_{s_i}(a)}+\log N\bigr).
    \label{eqn:right-uniform-degree-tail-count}
\end{align}
\end{lemma}
The uniform tail bound supplies the count estimate for noncritical
right-edge semilocalization. The proofs of the general degree estimates
are given in Subsection~\ref{subsec:transition-general-degree-proofs}.

\subsection{Critical degree input}

At critical sparsity, inverse profile maps identify spectral resonance
windows with degree windows. Uniform upper counts give semilocalization
bounds; two-sided counts at fixed exterior parameters with positive
count exponent give the transition estimates.

\begin{corollary}[Uniform upper bounds for critical resonant counts]
\label{cor:critical-resonant-upper-counts}
Under Assumptions~\ref{ass:proportional_regime}
and~\ref{ass:critical_sparsity}, the following statements hold.
\begin{enumerate}[label=\textnormal{(\roman*)}]
    \item \emph{Right edge.}
    Fix $L>\redge$. Let $0<\gamma_N<L-\redge$ and let
    $\delta_{\mathrm{R},N}:[\redge+\gamma_N,L]\to(0,\infty)$
    be deterministic and satisfy
    \begin{align}
        \sup_{\lambda\in[\redge+\gamma_N,L]}
        \frac{\delta_{\mathrm{R},N}(\lambda)}{\gapR}
        &\leq\frac12,
        &
        \sup_{\lambda\in[\redge+\gamma_N,L]}
        \frac{\delta_{\mathrm{R},N}(\lambda)}{\sqrt{\gapR}}
        &\longrightarrow0.
        \label{eqn:right-uniform-resonance-window}
    \end{align}
    For every fixed $\varepsilon,\nu>0$, there is a constant
    $C_{b,\qr,L,\varepsilon,\nu}>0$ such that, with probability at
    least $1-C_{b,\qr,L,\varepsilon,\nu}N^{-\nu}$,
    \begin{align}
        |\mathcal W_{\lambda,\delta_{\mathrm{R},N}(\lambda)}|
        \leq N^{\rhoR(\lambda,b)+\varepsilon}
        \label{eqn:right-uniform-resonant-upper-bound}
    \end{align}
    simultaneously for every $\lambda\in[\redge+\gamma_N,L]$.

    \item \emph{Left edge.}
    Suppose that $\qr>1$ and fix a compact interval
    $K_0\subset(0,\ledge)$. Let
    $\delta_{\mathrm{L},N}:K_0\to(0,\infty)$ be deterministic and
    satisfy
    \begin{align}
        \sup_{\lambda\in K_0}\delta_{\mathrm{L},N}(\lambda)
        &\longrightarrow0,
        \notag\\
        \delta_{\mathrm{L},N}(\lambda)
        &\leq\frac12\min\{\lambda,\ledge-\lambda\},
        \qquad\lambda\in K_0,
        \label{eqn:left-uniform-resonance-window}
    \end{align}
    for all sufficiently large $N$.
    For every fixed $\varepsilon,\nu>0$, there is a constant
    $C_{b,\qr,K_0,\varepsilon,\nu}>0$ such that, with probability at
    least $1-C_{b,\qr,K_0,\varepsilon,\nu}N^{-\nu}$,
    \begin{align}
        |\gW_{\lambda,\delta_{\mathrm{L},N}(\lambda)}^{(2,\mathrm{L})}|
        \leq N^{\rhoL(\lambda,b)+\varepsilon}
        \label{eqn:left-uniform-resonant-upper-bound}
    \end{align}
    simultaneously for every $\lambda\in K_0$.
\end{enumerate}
\end{corollary}

\begin{proposition}[Resonant counts at a fixed spectral parameter]
\label{prop:fixed-exterior-resonant-counts}
Under Assumptions~\ref{ass:proportional_regime}
and~\ref{ass:critical_sparsity}, fix an edge
$\mathrm E\in\{\mathrm L,\mathrm R\}$ and a spectral parameter
$\lambda>\redge$ if $\mathrm E=\mathrm R$, or
$0<\lambda<\ledge$ with $\qr>1$ if $\mathrm E=\mathrm L$.
Suppose that $\rho_{\mathrm E}(\lambda,b)>0$.
For $\sigma\in\{-1,+1\}$, let $\delta_{\sigma,N}>0$ be deterministic
and satisfy
\begin{align}
    \delta_{\sigma,N}&\longrightarrow0,
    &\frac{d\delta_{\sigma,N}}{\log\log N}
    &\longrightarrow\infty.
    \label{eqn:fixed-exterior-resonant-widths}
\end{align}
There is a deterministic sequence $\eta_N\downarrow0$ with
$\eta_N\log N\to\infty$ such that, for every $\nu>0$, there is a
constant $C_{b,\qr,\lambda,\nu}>0$ for which, with probability at
least $1-C_{b,\qr,\lambda,\nu}N^{-\nu}$ for all sufficiently large $N$,
\begin{align}
    N^{\rho_{\mathrm E}(\lambda,b)-\eta_N}
    \leq |\mathcal W_{\lambda,\delta_{\sigma,N}}|
    \leq N^{\rho_{\mathrm E}(\lambda,b)+\eta_N},
    \qquad\sigma\in\{-1,+1\}.
    \label{eqn:fixed-exterior-resonant-counts}
\end{align}
\end{proposition}
The counting arguments for both statements are given in
Subsection~\ref{subsec:transition-critical-count-proofs}.

The spectral-window condition \eqref{eqn:transition-spectral-window}
implies the width condition in
Proposition~\ref{prop:fixed-exterior-resonant-counts}, since
\begin{align}
    \frac{d\delta}{\log\log N}
    =\frac{\delta}{\error_d}
      \frac{\sqrt{\log N}\log d}{\log\log N}
    \longrightarrow\infty.
    \label{eqn:transition-degree-window}
\end{align}

\subsection{Degree roots and spectral boundaries}
\label{subsec:transition-degree-roots}

For $s>0$, extend \eqref{eqn:right-explicit-degree-root} by replacing
$\qr$ with $s$. The upper root $a_{\mathrm R}(b;s)$ is characterized by
\begin{align}
    I_s(a_{\mathrm R}(b;s))=b^{-1},
    \qquad a_{\mathrm R}(b;s)>s^2.
    \label{eqn:upper-extreme-degree-root}
\end{align}
Indeed, $I_s'(a)=\log(a/s^2)>0$ for $a>s^2$, with $I_s(s^2)=0$
and $I_s(a)\to\infty$ as $a\to\infty$.
For $\qr>1$, the lower root in \eqref{eqn:left-explicit-degree-root}
has the characterization
\begin{align}
    a_{\mathrm L}(b;\qr)
    =\inf\left\{
      a\in[0,\qr^2]:bI_{\qr}(a)\leq1
    \right\}.
    \label{eqn:lower-extreme-degree-root}
\end{align}
Since $I_{\qr}$ decreases strictly from $I_{\qr}(0)=\qr^2$ to zero on
$[0,\qr^2]$, this infimum is zero for $b\leq\qr^{-2}$ and is the
unique solution of $I_{\qr}(a)=b^{-1}$ in $(0,\qr^2)$ otherwise.

To obtain the explicit formulas, set $y=a/s^2$ and $w=\log y-1$.
For $a>0$, the degree-rate equation becomes
\begin{align}
    y\log y-y+1=\frac1{bs^2}
    \quad\Longleftrightarrow\quad
    w\mathrm e^w=\frac{1/(bs^2)-1}{\mathrm e}.
    \label{eqn:degree-root-lambert-equation}
\end{align}
On the upper branch, $a>s^2$ gives $w>-1$, so the principal branch
$W_0$ from Definition~\ref{def:lambert-w} yields
\begin{align}
    a_{\mathrm R}(b;s)
    =s^2\exp\!\left[
      1+W_0\!\left(\frac{1/(bs^2)-1}{\mathrm e}\right)
    \right].
    \label{eqn:upper-degree-root-explicit-general}
\end{align}
For the positive lower root, $s=\qr$ and $0<a<\qr^2$ give $w<-1$.
When $b>\qr^{-2}$, the Lambert argument lies in $(-\mathrm e^{-1},0)$,
and $W_{-1}$ gives the second case of
\eqref{eqn:left-explicit-degree-root}. The first case follows from
\eqref{eqn:lower-extreme-degree-root}; the positive root tends to zero
as $b\downarrow\qr^{-2}$.

The roots give the limiting extremal normalized degrees.

\begin{corollary}[Extremal normalized degrees]
\label{cor:critical-extreme-degrees}
Suppose that Assumptions~\ref{ass:proportional_regime}
and~\ref{ass:critical_sparsity} hold. Then
\begin{subequations}\label{eqn:extreme-normalized-degrees}
\begin{align}
    \max_{x\in\gV_1}\alpha_x
    &\xrightarrow{\P}a_{\mathrm R}(b;\qr^{-1}),
    \label{eqn:extreme-degree-V1}\\
    \max_{x\in\gV_2}\alpha_x
    &\xrightarrow{\P}a_{\mathrm R}(b;\qr).
    \label{eqn:extreme-degree-V2}
\end{align}
If $\qr>1$, then
\begin{align}
    \min_{x\in\gV_2}\alpha_x
    \xrightarrow{\P}a_{\mathrm L}(b;\qr).
    \label{eqn:extreme-degree-V2-lower}
\end{align}
\end{subequations}
\end{corollary}

Applying the location maps to the extreme degree roots, with the outlier thresholds as cutoffs, gives the spectral boundaries. The two right inverse degrees in \eqref{eqn:Lambda-inverse-right-edge} have the common form $\invR{s}(\lambda)=s^2+u(\lambda)$, where
\begin{align}
    u(\lambda)
    \coloneqq\frac12\left(
      \lambda^2-\qr^2-\qr^{-2}+\Delta_{\Lambda}(\lambda)
    \right)\geq1,
    \qquad\lambda>\redge.
\end{align}
For $r=s^2>0$ and $u>0$,
\begin{align}
    I_s(s^2+u)
    &=(r+u)\log(1+u/r)-u,
    \nonumber\\
    \frac{\partial}{\partial r}
      \bigl[(r+u)\log(1+u/r)-u\bigr]
    &=\log(1+u/r)-u/r<0.
\end{align}
Since $\qr^2\geq\qr^{-2}$, the right degree rate associated with
$\gV_2$ is no larger at every fixed $\lambda>\redge$. Thus
\begin{align}
    \bcurveR(\lambda;\qr)
    =\max_{i\in\{1,2\}}\bcurveRi{i}(\lambda;\qr)
    =\bcurveRi{2}(\lambda;\qr).
\end{align}

The function $\bcurveR(\,\cdot\,;\qr)$ is decreasing on $(\redge,\infty)$, and
\begin{align}
    \lim_{\lambda\downarrow\redge}\bcurveR(\lambda;\qr)
    =\bcritR(\qr)
\end{align}
by
\begin{align}
    \invR{\qr^{-1}}(\redge)
    &=\qr^{-2}+1,
    &
    \invR{\qr}(\redge)
    &=\qr^2+1.
\end{align}
The upper degree root crosses $\qr^2+1$ at $b=\bcritR(\qr)$.
Thus \eqref{eqn:global-right-transition-boundary} agrees with
$\Lambda_{\qr}(a_{\mathrm R}(b;\qr))$ below this threshold and with
$\redge$ at and above it. In particular,
\begin{align}
    \bigl[\bcurveR(\,\cdot\,;\qr)\bigr]^{-1}(b)
    =\Lambda_{\qr}(a_{\mathrm R}(b;\qr)),
    \qquad 0<b<\bcritR(\qr).
\end{align}

For the profile maps $\Lambda_s$ in \eqref{eqn:Lambda-maps}, define the
continuous, nondecreasing upper-branch extensions
\begin{align}
    \Gamma_s^{\mathrm{R}}(a)
    &\coloneqq
    \begin{cases}
        \redge,&a\leq s^2+1,\\
        \Lambda_s(a),&a>s^2+1,
    \end{cases}
    \qquad s\in\{\qr^{-1},\qr\}.
    \label{eqn:extended-upper-profile-map}
\end{align}
Monotonicity of $I_s$ and of the upper inverse branches gives
\begin{align}
    \redge(b;\qr)
    =\max\left\{
      \Gamma_{\qr^{-1}}^{\mathrm{R}}(a_{\mathrm R}(b;\qr^{-1})),
      \Gamma_{\qr}^{\mathrm{R}}(a_{\mathrm R}(b;\qr))
    \right\}.
    \label{eqn:right-boundary-extreme-degree-formula}
\end{align}
For $\lambda>\redge$, the inequality $b<\bcurveR(\lambda;\qr)$
is equivalent to $I_s(\invR{s}(\lambda))<b^{-1}$ for at least one
$s\in\{\qr^{-1},\qr\}$, and hence to
$\invR{s}(\lambda)<a_{\mathrm R}(b;s)$ for that class.

For the left edge, $\invL(\ledge)=\qr^2-1$ and
$\bcurveL(\,\cdot\,;\qr)$ increases from $\qr^{-2}$ to
$\bcritL(\qr)$ on $(0,\ledge)$. Hence
\begin{align}
    \bigl[\bcurveL(\,\cdot\,;\qr)\bigr]^{-1}(b)
    =\Lambda_{\qr}(a_{\mathrm L}(b;\qr)),
    \qquad \bcritIso(\qr)<b<\bcritL(\qr).
    \label{eqn:left-boundary-extreme-degree-formula}
\end{align}
The cutoff in \eqref{eqn:global-left-transition-boundary} gives
$\ledge$ for $b\geq\bcritL(\qr)$, while
\eqref{eqn:lower-extreme-degree-root} gives a zero boundary for
$b\leq\bcritIso(\qr)$. At $b=\bcritIso(\qr)$, the count of isolated vertices in $\gV_2$
is tight by Lemma~\ref{lem:bipartite-isolated-column-threshold}.

For each exterior spectral parameter, the positive count exponent
therefore has the characterization
\begin{align}
    \rhoR(\lambda,b)>0
    &\quad\Longleftrightarrow\quad
    \redge<\lambda<\redge(b;\qr),
    &&\lambda>\redge,\notag\\
    \rhoL(\lambda,b)>0
    &\quad\Longleftrightarrow\quad
    \ledge(b;\qr)<\lambda<\ledge,
    &&0<\lambda<\ledge.
    \label{eqn:transition-positive-count-intervals}
\end{align}
The second equivalence assumes $\qr>1$.

\subsection{Spectral comparison and profile-to-spectrum transfer}

The following comparison is the shared spectral input for the transition
counts and the eigenvector-mass estimates in
Section~\ref{sec:semilocalized_phase}.

\begin{lemma}[Spectral comparison]
\label{lem:spectral-counting}
Let $\rmS$ and $\rmT$ be self-adjoint matrices of the same size, with
$\|\rmS-\rmT\|\leq\varepsilon$.
\begin{enumerate}
    \item For every $\lambda\in\R$ and $\delta>\varepsilon$,
    \begin{align}
        \mathrm{N}_{\lambda,\delta-\varepsilon}(\rmT)
        \leq\mathrm{N}_{\lambda,\delta}(\rmS)
        \leq\mathrm{N}_{\lambda,\delta+\varepsilon}(\rmT).
        \label{eqn:transition-counting-sandwich-general}
    \end{align}
    \item Let $P$ be an orthogonal projection commuting with $\rmT$,
    put $P^\perp=\id-P$, and suppose
    \begin{align}
        \operatorname{dist}\bigl(
          \lambda,\operatorname{Spec}(\rmT|_{\operatorname{Ran}P^\perp})
        \bigr)\geq g>0.
    \end{align}
    Every normalized eigenvector $\rmS u=\lambda u$ satisfies
    \begin{align}
        \|P^\perp u\|\leq\frac{\varepsilon}{g},
        \qquad
        \|Pu\|^2\geq1-\left(\frac{\varepsilon}{g}\right)^2.
        \label{eqn:general-comparison-profile-mass}
    \end{align}
\end{enumerate}
\end{lemma}
Its proof is given in
Subsection~\ref{subsec:transition-spectral-comparison-proof}.

In fixed exterior windows separated from the complementary spectrum,
the comparison eigenvalue counts equal the resonant-vertex counts.

\begin{proposition}[Profile-to-spectrum counting transfer]
\label{prop:transition-counting-transfer}
Suppose that Assumptions~\ref{ass:proportional_regime}
and~\ref{ass:critical_sparsity} hold, and let
$\delta=\delta_N>0$ be a deterministic sequence satisfying
\eqref{eqn:transition-spectral-window}.
\begin{enumerate}
    \item Fix $B>b$ and a fixed $\lambda>\redge$. For every $\nu>0$,
    there exist constants $C_{B,\qr,\nu}>0$ and
    $\gC_{B,\qr,\nu}>0$ such that, for the deterministic sequence
    \begin{align}
        \varepsilon_{\mathrm{R},N}
        \coloneqq
        C_{B,\qr,\nu}\error_d
        =o(\delta)
        \label{eqn:transition-block-error-right}
    \end{align}
    the following holds with probability at least
    $1-\gC_{B,\qr,\nu}N^{-\nu}$,
    \begin{align}
        &|\gW_{\lambda,\delta-\varepsilon_{\mathrm{R},N}}^{(1)}|
        +|\gW_{\lambda,\delta-\varepsilon_{\mathrm{R},N}}^{(2,\mathrm{R})}|
        \leq \mathrm{N}_{\lambda,\delta}(\scA)
        \notag\\
        &\hspace{27mm}\leq
        |\gW_{\lambda,\delta+\varepsilon_{\mathrm{R},N}}^{(1)}|
        +|\gW_{\lambda,\delta+\varepsilon_{\mathrm{R},N}}^{(2,\mathrm{R})}|.
        \label{eqn:transition-counting-sandwich-right}
    \end{align}

    \item Suppose in addition that $\qr>1$. Fix
    $0<\lambda<\ledge$.
    For any compact
    interval $K\subset(0,\ledge)$ whose interior contains $\lambda$,
    there are constants $\nu_0=\nu_0(b,\qr,K)>0$ and
    $\gC_{b,\qr,K}>0$, and a deterministic sequence
    \begin{align}
        \varepsilon_{\mathrm{L},N}
        \coloneqq
        C_{b,\qr,K}\error_d
        =o(\delta)
        \label{eqn:transition-block-error-left}
    \end{align}
    such that, with probability at least
    $1-\gC_{b,\qr,K}N^{-\nu_0}$,
    \begin{align}
        |\gW_{\lambda,\delta-\varepsilon_{\mathrm{L},N}}^{(2,\mathrm{L})}|
        \leq \mathrm{N}_{\lambda,\delta}(\scA)
        \leq
        |\gW_{\lambda,\delta+\varepsilon_{\mathrm{L},N}}^{(2,\mathrm{L})}|.
        \label{eqn:transition-counting-sandwich-left}
    \end{align}
\end{enumerate}
\end{proposition}

The proof is given in
Subsection~\ref{subsec:transition-spectral-transfer-proof}.

\subsection{Proofs of the transition theorems}

\begin{proof}[Proof of Theorem~\ref{thm:right-edge-transition}]
Choose a fixed $B>b$. By Assumption~\ref{ass:critical_sparsity},
\eqref{eqn:subcritical_regime} holds for all sufficiently large $N$,
so all right-edge inputs to
Proposition~\ref{prop:transition-counting-transfer} are available.
We use the boundary identities from
Subsection~\ref{subsec:transition-degree-roots}.

\emph{(i) Spectral boundary.}
Fix $\epsilon>0$ and, for $s\in\{\qr^{-1},\qr\}$, put
\begin{align}
    A_{s,\epsilon}
    \coloneqq
    \invR{s}\left(
      \redge(b;\qr)+\frac{\epsilon}{3}
    \right).
\end{align}
Equation \eqref{eqn:right-boundary-extreme-degree-formula} implies
$A_{s,\epsilon}>a_{\mathrm R}(b;s)$. Define
\begin{align}
    c_{\mathrm{R}}(\epsilon)
    \coloneqq
    \frac12\left[
      b\rateR\left(\redge(b;\qr)+\frac{\epsilon}{3};\qr\right)-1
    \right]
    >0.
    \label{eqn:rightmost-success-exponent}
\end{align}
The shared tail-mean and presence bounds
\eqref{eqn:sparse-degree-tail-mean-bound} and
\eqref{eqn:degree-window-presence-absence} give
\begin{align}
    \P\left(
      \max_{x\in\gV_i}\alpha_x>A_{s_i,\epsilon}
      \text{ for some }i\in\{1,2\}
    \right)
    &\leq
    \sum_{i=1}^2
    N^{1-bI_{s_i}(A_{s_i,\epsilon})+o(1)}
    \notag\\
    &\leq 2N^{-c_{\mathrm{R}}(\epsilon)}
    \label{eqn:rightmost-extreme-degree-probability}
\end{align}
for all sufficiently large $N$.

Choose a fixed $\gap_0>0$ so small that
$(\gap_0\wedge1)/2<\epsilon/4$. Construct
$\hatH_{\mathrm{R}}(\gap_0)$ at the right edge.  By
\eqref{eqn:right-complement-spectrum}, its complementary block has
norm at most $\redge+\epsilon/3$ for all large $N$. On the complement
of the event in \eqref{eqn:rightmost-extreme-degree-probability}, every
positive profile eigenvalue is at most
$\redge(b;\qr)+\epsilon/3$. Invoke the complementary-spectrum
and block-approximation estimates with target failure exponent
$c_{\mathrm{R}}(\epsilon)+1$. Since
\eqref{eqn:full-block-approximation-right} has $o(1)$ error, it follows
that
\begin{align}
    \P\bigl(
      \lambda_{1}(\scA)>
      \redge(b;\qr)+\epsilon
    \bigr)
    \leq \gC_\epsilon N^{-c_{\mathrm{R}}(\epsilon)}.
    \label{eqn:rightmost-eigenvalue-upper-bound}
\end{align}

If $\redge(b;\qr)>\redge$, choose
\begin{align}
    \max\{\redge,\redge(b;\qr)-\epsilon\}
    <\lambda<\redge(b;\qr).
\end{align}
Then $b<\bcurveR(\lambda;\qr)$ and $\rhoR(\lambda,b)>0$.
Take $\delta_N=\error_d^{1/2}$, which satisfies
\eqref{eqn:transition-spectral-window}, and apply
Propositions~\ref{prop:transition-counting-transfer} and
\ref{prop:fixed-exterior-resonant-counts} with target failure exponent
$c_{\mathrm R}(\epsilon)+1$. The lower buffered resonant count is
at least $N^{\rhoR(\lambda,b)-o(1)}>0$, so
\eqref{eqn:transition-counting-sandwich-right} produces an eigenvalue
near $\lambda$ with failure probability at most
$\gC_\epsilon N^{-c_{\mathrm R}(\epsilon)-1}$.
Since $\delta_N\to0$ and
$\lambda>\redge(b;\qr)-\epsilon$, the resulting eigenvalue is
larger than $\redge(b;\qr)-\epsilon$ for all sufficiently
large $N$.

It remains to quantify the lower bound when
$\redge(b;\qr)=\redge$. Choose
\begin{align}
    0<\kappa<\min\left\{\frac{\epsilon}{3},\frac{\qr^{-1}}3\right\},
    \qquad
    E\coloneqq\redge-2\kappa,
    \qquad
    \eta\coloneqq N^{-1+\kappa}.
\end{align}
Then $E$ belongs to the positive $\kappa$-bulk and
$E>\redge-\epsilon$. By
\eqref{eqn:base-stieltjes-lower} and
\eqref{eqn:stieltjes-Lipschitz-alpha}, there are constants
$r_0,c_0>0$, depending only on $(\epsilon,\qr)$, such that
\begin{align}
    |\alpha-\qr^{-2}|\leq r_0
    \quad\Longrightarrow\quad
    \Im\stj^{(1)}_\alpha(E+\ii\eta)\geq c_0.
    \label{eqn:rightmost-typical-transform-lower-bound}
\end{align}
For $x\in\gV_1$, the graph row-square parameter satisfies
\begin{align}
    \zeta_x
    =(1-2p)\alpha_x+\frac{mp^2}{d}.
\end{align}
The one-vertex bounds in Proposition~\ref{prop:uniform-degree-counts}
show that the degree interval corresponding to
$|\zeta_x-\qr^{-2}|\leq r_0$ has expected count $(1-o(1))n$.
The simultaneous concentration estimate in that proposition therefore
gives at least $n/2$ such vertices, with failure probability at most
$C_\epsilon N^{-c_{\mathrm R}(\epsilon)-1}$. Invoke the local law in
Theorem~\ref{thm:bipartite-sparse-local-law} at $z=E+\ii\eta$ with
target failure exponent $c_{\mathrm{R}}(\epsilon)+1$. On this event,
\eqref{eqn:rightmost-typical-transform-lower-bound} and the
nonnegativity of every diagonal resolvent imaginary part yield
\begin{align}
    \frac1N\Im\Tr(\scA-z)^{-1}\geq c_{\epsilon,\qr}>0
    \label{eqn:rightmost-bulk-trace-lower-bound}
\end{align}
with failure probability at most
$\gC_\epsilon N^{-c_{\mathrm{R}}(\epsilon)-1}$. On the other hand, if
$\lambda_{1}(\scA)\leq\redge-\epsilon$, the spectral decomposition
gives
\begin{align}
    \frac1N\Im\Tr(\scA-z)^{-1}
    \leq
    \frac{\eta}{(E-\redge+\epsilon)^2}
    =o(1),
\end{align}
which contradicts \eqref{eqn:rightmost-bulk-trace-lower-bound} for all
large $N$. Thus the lower-deviation probability is smaller than
$\gC_\epsilon N^{-c_{\mathrm{R}}(\epsilon)}$ in both cases. Combining
this with \eqref{eqn:rightmost-eigenvalue-upper-bound} proves
\eqref{eqn:rightmost-eigenvalue-transition} in both regimes, with
$c_\epsilon=c_{\mathrm{R}}(\epsilon)$.

\emph{(ii) Outlier counts.}
Suppose first that $b>\bcritR(\qr)$, so $\redge(b;\qr)=\redge$.
Fix $\lambda>\redge$ and apply
\eqref{eqn:rightmost-eigenvalue-transition} with
$\epsilon=\gapR/3$. For all sufficiently large $N$, one has
$\delta_N\leq\gapR/3$, and on the
event in \eqref{eqn:rightmost-eigenvalue-transition},
\begin{align}
    \lambda_{1}(\scA)
    \leq\redge+\frac{\gapR}{3}
    <\lambda-\delta_N.
\end{align}
This proves \eqref{eqn:right-high-eigenvalue-absence} with
$c_\lambda=c_{\mathrm{R}}(\gapR/3)$.

Assume now that $b<\bcritR(\qr)$ and fix
$\lambda\in(\redge,\redge(b;\qr))$. By
\eqref{eqn:transition-positive-count-intervals}, $\rhoR(\lambda,b)>0$.
Fix $\nu>0$ and invoke Proposition~\ref{prop:transition-counting-transfer} with target
failure exponent $\nu+1$. Both buffered widths
\begin{align}
    \delta_\sigma
    \coloneqq
    \delta+\sigma\varepsilon_{\mathrm{R},N},
    \qquad \sigma\in\{-1,+1\},
\end{align}
are $(1+o(1))\delta$ and satisfy
\eqref{eqn:transition-degree-window}. Applying
Proposition~\ref{prop:fixed-exterior-resonant-counts} to both buffered
widths gives a deterministic sequence $\eta_N\downarrow0$
such that, with failure probability at most
$\gC_{\lambda,\nu}N^{-\nu}$,
\begin{align}
    N^{\rhoR(\lambda,b)-\eta_N}
    \leq |\mathcal W_{\lambda,\delta_\sigma}|
    \leq N^{\rhoR(\lambda,b)+\eta_N},
    \qquad \sigma\in\{-1,+1\}.
    \label{eqn:right-fixed-window-count-concentration}
\end{align}
Combining this event with
\eqref{eqn:transition-counting-sandwich-right} proves
\eqref{eqn:right-transition-eigenvalue-count} with the stated success
probability.

\end{proof}

\begin{proof}[Proof of Theorem~\ref{thm:left-edge-transition}]
We use the left comparison in
Proposition~\ref{prop:transition-counting-transfer}, under
Assumptions~\ref{ass:proportional_regime} and
\ref{ass:critical_sparsity} with $\qr>1$. The boundary identities
in Subsection~\ref{subsec:transition-degree-roots} apply throughout.

\emph{(i) Spectral absence below the boundary.}
When $\ledge(b;\qr)=0$, the compact-absence assertion is vacuous.
Suppose therefore that $\ledge(b;\qr)>0$, and fix a nonempty
compact set $K\subset(0,\ledge(b;\qr))$.  Choose $\eta>0$ such that
\begin{align}
    K^{(2\eta)}
    \coloneqq
    \{t>0:\operatorname{dist}(t,K)\leq2\eta\}
    \subset(0,\ledge(b;\qr)).
\end{align}
Choose a compact interval
$K_0\subset(0,\ledge(b;\qr))$ containing $K^{(2\eta)}$,
and put
\begin{align}
    t_0&\coloneqq\sup K_0,
    &a_0&\coloneqq\invL(t_0),
    &c_{\deg}&\coloneqq bI_{\qr}(a_0)-1>0.
\end{align}
The last inequality follows from
\eqref{eqn:global-left-transition-boundary} and the monotonicity of
$I_{\qr}$ on the lower branch.
The tail-mean bound \eqref{eqn:sparse-degree-tail-mean-bound} and the
presence bound \eqref{eqn:degree-window-presence-absence}, applied
to $[0,a_0]$, give
\begin{align}
    \P\bigl(\exists x\in\gV_2:\alpha_x\leq a_0\bigr)
    \leq C_\qr N^{-c_{\deg}}
    \leq N^{-c_{\deg}/2}
\end{align}
for all sufficiently large $N$. With probability at least
$1-N^{-c_{\deg}/2}$, no positive lower-branch profile value lies in
$K_0$. Set $\gap_0=\ledge-t_0$.
Apply the left constructions with $\sK=K_0$ and $\gap=\gap_0$.
The block-approximation and complementary-spectrum events hold with
probability at least
$1-C_{b,\qr,K_0}N^{-\nu_0}$ for some $\nu_0>0$. On the intersection
with the degree event, \eqref{eqn:left-complement-spectrum} excludes the
complementary positive spectrum from $K_0$ for all sufficiently large
$N$, while the profile spectrum also misses $K_0$. Zero and the negative
spectrum are separated from $K$. Hence
$\operatorname{Spec}(\hatH_{\mathrm{L}}(\gap_0))\cap K^{(2\eta)}=\varnothing$.
Since $\|\scA-\hatH_{\mathrm{L}}(\gap_0)\|=o(1)$,
Lemma~\ref{lem:weyl} yields
\begin{align}
    \P\bigl(\operatorname{Spec}(\scA)\cap K=\varnothing\bigr)
    \geq1-\gC_KN^{-c_{\mathrm{L},K}}
    \label{eqn:left-compact-spectral-absence}
\end{align}
for some $\gC_K>0$, with
$c_{\mathrm{L},K}=\min\{c_{\deg}/2,\nu_0\}>0$.
This proves item~\textnormal{(i)}.

\emph{(ii) Outlier counts.}
The following arguments also cover $\ledge(b;\qr)=0$.
Suppose first that $b>\bcritL(\qr)$. For fixed
$\lambda\in(0,\ledge)$, choose a compact interval
$K_\lambda\subset(0,\ledge)$ whose interior contains $\lambda$. Since
$\ledge(b;\qr)=\ledge$, \eqref{eqn:left-compact-spectral-absence}
applies to $K_\lambda$. For all sufficiently large $N$,
$[\lambda-\delta_N,\lambda+\delta_N]\subset K_\lambda$, which proves
\eqref{eqn:left-high-eigenvalue-absence}.

Suppose now that $0<b<\bcritL(\qr)$ and fix
$\lambda\in(\ledge(b;\qr),\ledge)$. By
\eqref{eqn:transition-positive-count-intervals}, $\rhoL(\lambda,b)>0$.
Choose a compact interval $K_0\subset(0,\ledge)$ whose interior contains $\lambda$, and
invoke Proposition~\ref{prop:transition-counting-transfer}\textnormal{(ii)}.
For $\sigma\in\{-1,+1\}$, set
$\delta_\sigma=\delta+\sigma\varepsilon_{\mathrm{L},N}$.
Both buffered widths are $(1+o(1))\delta$ and satisfy
\eqref{eqn:transition-degree-window}.
Proposition~\ref{prop:fixed-exterior-resonant-counts} gives a deterministic
sequence $\eta_N\downarrow0$ such that, with failure probability smaller
than any fixed negative power of $N$,
\begin{align}
    N^{\rhoL(\lambda,b)-\eta_N}
    \leq |\gW_{\lambda,\delta_\sigma}^{(2,\mathrm L)}|
    \leq N^{\rhoL(\lambda,b)+\eta_N},
    \qquad \sigma\in\{-1,+1\}.
    \label{eqn:left-fixed-window-count-concentration}
\end{align}
Combining this with \eqref{eqn:transition-counting-sandwich-left} proves
\eqref{eqn:left-transition-eigenvalue-count}; the fixed polynomial
success exponent is inherited from the left spectral-transfer event.

\end{proof}

%% file: appendix.tex

\section{Stieltjes Transform of the Symmetric Marchenko-Pastur Distribution}\label{sec:proof-of-lemma-stieltjes-smp-biregular-tree}
\begin{proof}[Proof of Lemma~\ref{lem:stieltjes-smp}]
    Assume first that $\ratio>1$. Recall that $\ratio=\qr^{4}$, $\lambda_{\mathrm L}=\qr-\qr^{-1}$, and $\lambda_{\mathrm R}=\qr+\qr^{-1}$. The Stieltjes transform of \eqref{eqn:symmetric_MP_measure} can be written as
        \begin{align}
            \stj(z) =&\, \int_{\R} \frac{1}{x-z}\,\mu_{\rm SMP}(dx) = \frac{1}{\ratio + 1}\int_{-\lambda_{\mathrm R}}^{-\lambda_{\mathrm L}} \frac{f(-x)}{x - z} dx + \frac{\ratio - 1}{\ratio + 1}\int \delta_{0}(x) \frac{1}{x-z} dx + \frac{1}{\ratio + 1}\int_{\lambda_{\mathrm L}}^{\lambda_{\mathrm R}} \frac{f(x)}{x - z} dx\\
            &= \frac{1-\ratio}{(1+\ratio)z} + \frac{1}{1+\ratio}I(z), \label{eqn:stieltjes-smp-decomposition}
        \end{align}
    where we denote $h(x)\coloneqq \sqrt{(x^{2}-\lambda_{\mathrm L}^{2})(\lambda_{\mathrm R}^{2}-x^{2})}$ and
    \begin{align}
        I(z)
        &\coloneqq
        -\frac{\qr^{2}}{\pi}
        \int_{-\lambda_{\mathrm R}}^{-\lambda_{\mathrm L}}
        \frac{h(x)}{x(x-z)}\,dx
        +\frac{\qr^{2}}{\pi}
        \int_{\lambda_{\mathrm L}}^{\lambda_{\mathrm R}}
        \frac{h(x)}{x(x-z)}\,dx \label{eqn:stieltjes-smp-I}
    \end{align}
    Furthermore, let
    \[
    \widetilde h(\zeta)\coloneqq
        \sqrt{(\zeta^{2}-\lambda_{\mathrm L}^{2})
        (\zeta^{2}-\lambda_{\mathrm R}^{2})}
    \]
    denote the branch holomorphic on
    $
        \C\setminus\big([-\lambda_{\mathrm R},-\lambda_{\mathrm L}]
        \cup[\lambda_{\mathrm L},\lambda_{\mathrm R}]\big)
    $
    which satisfies $\widetilde h(\zeta)\sim\zeta^{2}$ as
    $|\zeta|\to\infty$. In particular,
    $
        \widetilde h(0)=-\lambda_{\mathrm L}\lambda_{\mathrm R}.
    $
    For fixed $z\in\C\setminus\R$, define
    \[
        g(\zeta)\coloneqq
        \frac{\qr^{2}}{\pi}
        \frac{\widetilde h(\zeta)}{\zeta(\zeta-z)}.
    \]
    The function $g$ is meromorphic in this domain, with only the poles
    $\zeta=0$ and $\zeta=z$. Their residues are
    \begin{align}
        \operatorname{Res}(g,0)
        &=\frac{\qr^{2}}{\pi}
        \frac{\lambda_{\mathrm L}\lambda_{\mathrm R}}{z},&
        \operatorname{Res}(g,z)
        &=\frac{\qr^{2}}{\pi}
        \frac{\widetilde h(z)}{z}. \label{eqn:stieltjes-smp-residues}
    \end{align}

    Choose $R>|z|+\lambda_{\mathrm R}+1$ and $0<\eta<\min\{\lambda_{\mathrm L}/2,|\Im z|/2\}$. Let $\Gamma$ be the positively oriented boundary of the disk $|\zeta|<R$ with the closures of the $\eta$-neighborhoods of the two cuts removed. Thus the outer circle is counterclockwise and the boundaries of the two slits are clockwise. See Figure~\ref{fig:contour-smp}.
    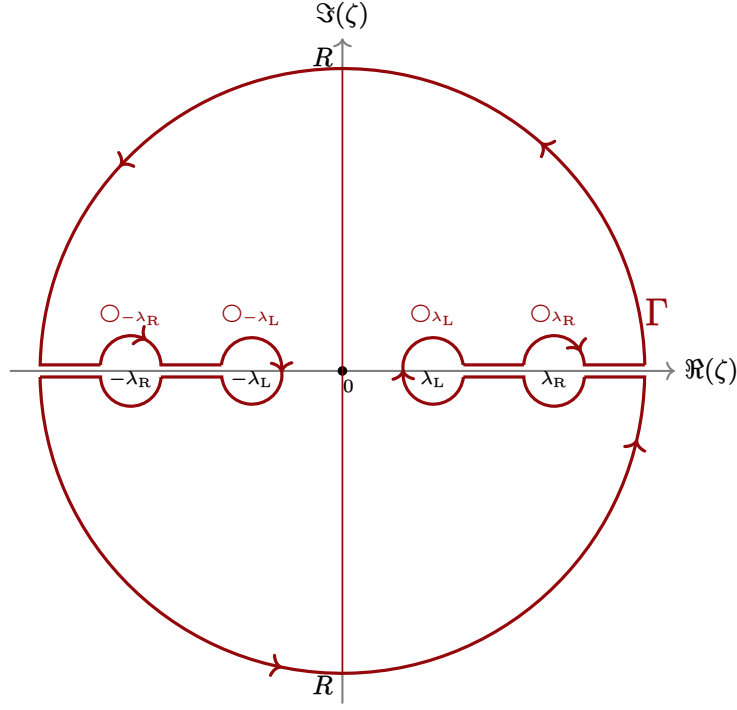
\begin{figure}[h]
        \centering
        \begin{tikzpicture}[
            scale=0.8,
            contour/.style={
                line width=1.2pt,
                ThemeColor,
                decoration={
                    markings,
                    mark=at position 0.045 with {\arrow{>}},
                    mark=at position 0.15 with {\arrow{>}},
                    mark=at position 0.30 with {\arrow{>}},
                    mark=at position 0.410 with {\arrow{>}},
                    mark=at position 0.470 with {\arrow{>}},
                    mark=at position 0.70 with {\arrow{>}},
                    mark=at position 0.85 with {\arrow{>}},
                    mark=at position 0.970 with {\arrow{>}}
                },
                postaction={decorate}
            },
            branchcut/.style={line width=2pt, ThemeColor, dashed, cap=round}
        ]
        \def\largeRadius{5}   
        \def\smallRadius{0.5} 
        \def\offsetDistance{0.1} 
        \def\leftEnd{1.5}
        \def\rightEnd{3.5}
        \draw[->, gray, line width=0.8pt] (-5.5,0) -- (5.5,0) node[right, black]{$\Re(\zeta)$};
        \draw[->, gray, line width=0.8pt] (0,-5.5) -- (0,5.5) node[above, black]{$\Im(\zeta)$};


        \pgfmathsetmacro{\endAngleL}{asin(\offsetDistance/(\leftEnd - \smallRadius))}
        \pgfmathsetmacro{\startAngleL}{180 - \endAngleL}
        \pgfmathsetmacro{\endAngleR}{asin(\offsetDistance/(\rightEnd - \smallRadius))}
        \pgfmathsetmacro{\startAngleR}{180 - \endAngleR}
        \pgfmathsetmacro{\startAngleLarge}{asin(\offsetDistance/(\largeRadius))}
        \pgfmathsetmacro{\endAngleLarge}{180 - \startAngleLarge}

        \draw[contour]
            ({\leftEnd + \smallRadius}, \offsetDistance)
            -- ({\rightEnd - \smallRadius}, \offsetDistance)
            (\rightEnd - \smallRadius, \offsetDistance) arc [start angle= \startAngleR, end angle= \endAngleR, radius=\smallRadius]
            -- (\largeRadius, \offsetDistance)
            (\largeRadius, \offsetDistance) arc [start angle= \startAngleLarge, end angle= \endAngleLarge, radius=\largeRadius]
            -- ({-\rightEnd - \smallRadius}, \offsetDistance)
            (-\rightEnd - \smallRadius, \offsetDistance) arc [start angle= \startAngleR, end angle= \endAngleR, radius=\smallRadius]
            -- ({ -\leftEnd - \smallRadius}, \offsetDistance)
            (-\leftEnd - \smallRadius, \offsetDistance) arc [start angle= \startAngleL, end angle= \endAngleL, radius=\smallRadius]
             -- (-\leftEnd + \smallRadius, {-\offsetDistance})
            (-\leftEnd + \smallRadius, -\offsetDistance) arc [start angle= -\endAngleL, end angle= -\startAngleL, radius=\smallRadius]
            -- ({ -\rightEnd + \smallRadius}, -\offsetDistance)
            (-\rightEnd + \smallRadius, -\offsetDistance) arc [start angle= -\endAngleR, end angle= -\startAngleR, radius=\smallRadius]
            --({-\largeRadius}, - \offsetDistance)
            (-\largeRadius, -\offsetDistance) arc [start angle= -\endAngleLarge, end angle= -\startAngleLarge, radius=\largeRadius]
            --({\rightEnd + \smallRadius}, - \offsetDistance)
            (\rightEnd + \smallRadius, -\offsetDistance) arc [start angle = -\endAngleR, end angle = -\startAngleR, radius=\smallRadius]
            -- ({ \leftEnd + \smallRadius}, -\offsetDistance)
            (\leftEnd + \smallRadius, -\offsetDistance) arc [start angle= -\endAngleL, end angle= -\startAngleL, radius=\smallRadius]
             -- (\leftEnd - \smallRadius, {\offsetDistance})
            (\leftEnd - \smallRadius, \offsetDistance) arc [start angle= \startAngleL, end angle= \endAngleL, radius=\smallRadius]
                ;
        \node[ThemeColor, above] at (\largeRadius + 0.2, 6 *\offsetDistance) {\Large{$\Gamma$}};

        \node[black, below] at (\leftEnd, 0.15){\scriptsize{$\lambda_{\mathrm L}$}};

        \node[black, below] at (\rightEnd, 0.15){\scriptsize{$\lambda_{\mathrm R}$}};

        \node[black, below] at (-\leftEnd, 0.15){\scriptsize{$-\lambda_{\mathrm L}$}};

        \node[black, below] at (-\rightEnd, 0.15){\scriptsize{$-\lambda_{\mathrm R}$}};

        \filldraw[black] (0,0) circle (2pt);
        \node[black, below] at (0.1, 0){\scriptsize{$0$}};

        \draw[ThemeColor]  (0,0) -- (0, \largeRadius);
        \node[black, left] at (0, {\largeRadius + 0.2}){$R$};
        \draw[ThemeColor]  (0,0) -- (0, -\largeRadius);
        \node[black, left] at (0, {-\largeRadius - 0.2}){$R$};

        \node[ThemeColor, above] at (\leftEnd, \smallRadius + \offsetDistance){\scriptsize{$\bigcirc_{\lambda_{\mathrm L}}$}};
        \node[ThemeColor, above] at (\rightEnd, \smallRadius + \offsetDistance){\scriptsize{$\bigcirc_{\lambda_{\mathrm R}}$}};
        \node[ThemeColor, above] at (-\leftEnd, \smallRadius + \offsetDistance){\scriptsize{$\bigcirc_{-\lambda_{\mathrm L}}$}};
        \node[ThemeColor, above] at (-\rightEnd, \smallRadius + \offsetDistance){\scriptsize{$\bigcirc_{-\lambda_{\mathrm R}}$}};
        \end{tikzpicture}
        \caption{Schematic contour $\Gamma$}\label{fig:contour-smp}
        \end{figure}
    According to the residue theorem, we then have
    \begin{align}
         \oint_{\Gamma} g(\zeta) d\zeta = 2\pi \ii \cdot \big( \mathrm{Res}(g(\zeta), 0) +  \mathrm{Res}(g(\zeta), z) \big).\label{eqn:stieltjes-smp-residues-evaluation}
    \end{align}
    To evaluate the two clockwise slit boundaries, let $x$ lie in the
    interior of either cut and use $h(x)$ defined above.
    The branch choice gives
    \begin{align}
        \widetilde h(x+\ii0)&=\sign(x)\ii h(x),&
        \widetilde h(x-\ii0)&=-\sign(x)\ii h(x).
    \end{align}
    Hence their upper-minus-lower jumps are
    \begin{align}
        g(x+\ii0)-g(x-\ii0)
        &=2\ii\frac{\qr^2}{\pi}\frac{h(x)}{x(x-z)}
          =2\ii\frac{f(x)}{x-z},&&x>0,\\
        g(x+\ii0)-g(x-\ii0)
        &=-2\ii\frac{\qr^2}{\pi}\frac{h(x)}{x(x-z)}
          =2\ii\frac{f(-x)}{x-z},&&x<0.
    \end{align}
    At each endpoint, $\widetilde h(\zeta)=O(\sqrt\eta)$ on an arc
    of length $O(\eta)$, while $\zeta(\zeta-z)$ stays bounded away from
    zero. Thus the four endpoint arcs contribute $O(\eta^{3/2})$. For an interval $J$, let $U_\eta(J)$ be its open $\eta$-neighborhood.
    Passing to the limit along the slit sides yields
    \begin{align}
        &\lim_{\eta\downarrow0}
        \sum_{J\in\{[-\lambda_{\mathrm R},-\lambda_{\mathrm L}],
                     [\lambda_{\mathrm L},\lambda_{\mathrm R}]\}}
        \oint_{\partial U_\eta(J)}g(\zeta)\,d\zeta \nonumber\\
        =&\int_{-\lambda_{\mathrm R}}^{-\lambda_{\mathrm L}}
           \bigl(g(x+\ii0)-g(x-\ii0)\bigr)\,dx
          +\int_{\lambda_{\mathrm L}}^{\lambda_{\mathrm R}}
           \bigl(g(x+\ii0)-g(x-\ii0)\bigr)\,dx
          =2\ii I(z).
          \label{eqn:contiour-smp-1st-to-I}
    \end{align}
    Here each $\partial U_\eta(J)$ has the clockwise orientation
    inherited from $\Gamma$. On the outer circle, the Laurent expansion
    $\widetilde h(\zeta)=\zeta^2+O(1)$ gives
    $g(\zeta)=\qr^2\pi^{-1}(1+z/\zeta+O(\zeta^{-2}))$, so
    \begin{align}
        \oint_{|\zeta|=R}g(\zeta)\,d\zeta
        =2\ii\qr^2z.
        \label{eqn:stieltjes-smp-outer-contribution}
    \end{align}

    Combining \eqref{eqn:stieltjes-smp-residues-evaluation}, \eqref{eqn:contiour-smp-1st-to-I}, and
    \eqref{eqn:stieltjes-smp-outer-contribution}, we have
    \[
        2\ii I(z)+2\ii\qr^{2}z
        =
        2\pi\ii\cdot \big( \mathrm{Res}(g(\zeta), 0) +  \mathrm{Res}(g(\zeta), z) \big).
    \]
    Then by substituting the residues in \eqref{eqn:stieltjes-smp-residues}, we have
    \begin{align}
        I(z)
        &=
        \frac{\qr^{2}}{z}
        \left(-z^{2}+\widetilde h(z)
        +\lambda_{\mathrm L}\lambda_{\mathrm R}\right). \label{eqn:stieltjes-smp-I-explicit}
    \end{align}
    Since $\qr^{2}\lambda_{\mathrm L}\lambda_{\mathrm R}=\qr^{4}-1=\ratio-1$, substituting \eqref{eqn:stieltjes-smp-I-explicit} into
    \eqref{eqn:stieltjes-smp-decomposition} yields
    \[
        \stj(z)
        =
        \frac{\sqrt{\ratio}}{(1+\ratio)z}
        \left(-z^{2}
        +\sqrt{(z^{2}-\lambda_{\mathrm L}^{2})
        (z^{2}-\lambda_{\mathrm R}^{2})}\right).
    \]
    Both sides are holomorphic away from the two cuts and the origin, so
    the identity extends from $\C\setminus\R$ to their common domain by
    analytic continuation.
    The case $\ratio=1$ follows by taking the limit $\ratio\downarrow1$.

    For completeness, we check that the resulting formula is consistent
    with Stieltjes inversion. For $x$ in the interior of either cut with
    $x\neq0$, the Stieltjes inversion formula gives
    \begin{align}
        \frac{1}{\pi}\lim_{\eta\downarrow0}\Im\stj(x+\ii\eta)
        &=
        \frac{\qr^{2}}{\pi(1+\qr^{4})}
        \frac{\sqrt{(x^{2}-\lambda_{\mathrm L}^{2})
        (\lambda_{\mathrm R}^{2}-x^{2})}}{|x|},
    \end{align}
    which is the absolutely continuous density in
    \eqref{eqn:symmetric_MP_measure}. When $\ratio=1$, the right-hand side
    extends continuously to $x=0$. To recover the atom at zero, note that
    $\widetilde h(\ii\eta)
    =-\sqrt{(\eta^2+\lambda_{\mathrm L}^2)
             (\eta^2+\lambda_{\mathrm R}^2)}$, direct substitution gives
    \begin{align}
        \lim_{\eta\downarrow0}\eta\Im\stj(\ii\eta)
        &=\frac{\qr^2}{1+\ratio}
          \lim_{\eta\downarrow0}
          \left(\sqrt{(\eta^2+\lambda_{\mathrm L}^2)
                       (\eta^2+\lambda_{\mathrm R}^2)}-\eta^2\right)
          =\frac{\qr^2\lambda_{\mathrm L}\lambda_{\mathrm R}}{1+\ratio}
          =\frac{\ratio-1}{\ratio+1}.
    \end{align}
    \end{proof}


\section{Analysis of the Infinite Rooted \texorpdfstring{$(q, d_1, d_2)$}{(q, d1, d2)}-Biregular Tree}\label{sec:analysis-of-biregular-tree}
We focus on the infinite rooted $(q, d_1, d_2)$-biregular tree defined below.
\begin{definition}\label{def:biregular_tree}
Let $x$ be the root vertex of the tree $\tree$. The $(q, d_1, d_2)$-biregular tree $\tree$ is defined as follows.
\begin{enumerate}
    \item The root $x$ has $q$ children, so $|\layer_1^{\tree}(x)|=q$.
    \item $x\in \gV_1$: vertices in $\layer_{2j - 1}^{\tree}(x)$ have $d_2$ children, while vertices in $\layer_{2j}^{\tree}(x)$ have $d_1$ children.
    \item $x\in \gV_2$: vertices in $\layer_{2j - 1}^{\tree}(x)$ have $d_1$ children, while vertices in $\layer_{2j}^{\tree}(x)$ have $d_2$ children.
\end{enumerate}
\end{definition}
The main tool we use is the tridiagonalization technique, which was originally developed in numerical linear algebra for matrix computations and has found numerous applications in random matrix theory, e.g., the construction of general beta ensembles in \cite{dumitriua2002matrix}.

The section is organized as follows. We first present basic properties in
Section~\ref{sec:basic-properties-biregular-tree} and then study the
Stieltjes transforms of the spectral measures in
Section~\ref{sec:stieltjes-transform-spectral-measures-biregular-tree}.
The scalar Stieltjes transform is computed in
Appendix~\ref{sec:proof-of-lemma-stieltjes-smp-biregular-tree}.

\subsection{Basic properties}\label{sec:basic-properties-biregular-tree}
Below, Lemma~\ref{lem:tridiagrep} is collected from our previous paper \cite{dumitriu2025singular} regarding the tridiagonal representation of the normalized adjacency matrix of $\tree$. Here, we denote $d = \sqrt{d_{1} d_{2}}$.

\begin{lemma}[{\cite[Lemma C.1]{dumitriu2025singular}}]\label{lem:tridiagrep}
Let $\rvs_0, \rvs_1, \rvs_2, \ldots, \rvs_{r}$ denote the vectors obtained from $\ones_{x}$, $\rmA\ones_x$, $\rmA^2\ones_x$, $\ldots$, $\rmA^{r}\ones_x$ after Gram-Schmidt orthonormalization. Then the following holds.
\begin{enumerate}
    \item For all $j = 0, 1, \ldots, r$, we have
        \begin{align}
            \rvs_j = |\layer_j(x)|^{-1/2} \ones_{\layer_j(x)}.
        \end{align}
    \item Let $\rvs_{r + 1}, \ldots, \rvs_{N-1}$ be a completion of $\rvs_0, \ldots, \rvs_{r}$ such that $\{\rvs_j\}_{j=0}^{N-1}$ form an orthonormal basis of $\R^{N}$. Define the basis $\rmS \coloneqq [\rvs_0, \rvs_1, \ldots, \rvs_{N-1}]$ and let
    \begin{align}
        \rmZ \coloneqq \frac{1}{\sqrt{d}} \rmS^{\sT} \rmA \rmS
    \end{align}
    denote the representation of $\rmA/\sqrt{d}$ in the basis $\rmS$. Then the upper-left $(r +1)\times (r +1)$ block $\rmZ_{\llbracket r \rrbracket}$ of $\rmZ$ has the following tridiagonal form
    \begin{align}
        \rmZ_{\llbracket r \rrbracket}(\alpha) = \indi{x\in \gV_1} \rmZ^{(1)}(\alpha_x ) + \indi{x\in \gV_2} \rmZ^{(2)}(\alpha_x),
    \end{align}
    where for $r$ being even, $\rmZ^{(1)}(\alpha)$ and $\rmZ^{(2)}(\alpha)$ are $(r+1)\times (r+1)$ matrices defined through
\begin{align}\label{eqn:Z1Z2}
    \rmZ^{(1)}(\alpha) \coloneqq
    \begin{bmatrix}
   0  & \sqrt{\alpha} & \\
   \sqrt{\alpha} & 0 & \qr\\
       &  \qr &  0 &  \qr^{-1}\\
        & &\qr^{-1} & 0 & \ddots\\
        & &  &  \ddots &  \ddots & \qr \\
        & &  & &   \qr & 0
    \end{bmatrix}, \,\,
    \rmZ^{(2)}(\alpha)=\begin{bmatrix}
   0  & \sqrt{\alpha} & \\
   \sqrt{\alpha} & 0 & \qr^{-1}\\
       &  \qr^{-1} &  0 &  \qr\\
        & &\qr & 0 & \ddots\\
        & &  &  \ddots &  \ddots & \qr^{-1} \\
        & &  & &   \qr^{-1} & 0
    \end{bmatrix}.
\end{align}
The bottom right corner of $\rmZ^{(1)}(\alpha)$ (resp. $\rmZ^{(2)}(\alpha)$) changes to $\qr^{-1}$ (resp. $\qr$) when $r$ is odd.
\end{enumerate}
\end{lemma}

The spectral properties of $\rmZ^{(1)}(\alpha)$ and $\rmZ^{(2)}(\alpha)$ have been studied in \cite{dumitriu2025singular}. We summarize the relevant results in the following lemma.

\begin{lemma}[{\cite[Lemma C.2]{dumitriu2025singular}}]\label{lem:eigenZ1Z2}
Let $\lambda^{(1)}_{1} \geq \cdots \geq \lambda^{(1)}_{r + 1}$ and $\lambda^{(2)}_{1} \geq \cdots \geq \lambda^{(2)}_{r + 1}$ denote the eigenvalues of $\rmZ^{(1)}(\alpha)$ and $\rmZ^{(2)}(\alpha)$, respectively.  The following holds.
\begin{enumerate}
    \item Bulk eigenvalues:
    \begin{align}
        \big\{ \{\lambda^{(1)}_{l}\}_{l=2}^{r} \cup \{\lambda^{(2)}_{l}\}_{l=2}^{r} \big\} \subset [-\qr -\qr^{-1},-\qr +\qr^{-1}]\cup [\qr-\qr^{-1}, \qr +\qr^{-1}]
    \end{align}
    \item Extreme eigenvalues:
    \begin{enumerate}
        \item If $\alpha>\qr^{-2}+1$, $\lambda^{(1)}_{1}$ and $\lambda^{(1)}_{r + 1}$ converge to $\Lambda_{\qr^{-1}}(\alpha)$ and $ -\Lambda_{\qr^{-1}}(\alpha)$, respectively, as $r\to\infty$.
        \item If $\alpha>\qr^2+1$, $\lambda^{(2)}_{1}$ and $\lambda^{(2)}_{r + 1}$ converge to $ \Lambda_{\qr}(\alpha)$ and $ -\Lambda_{\qr}(\alpha)$, respectively, as $r\to\infty$.
        \item If $0<\alpha<\qr^2-1$, the smallest positive eigenvalue and largest negative eigenvalue of $\rmZ^{(2)}(\alpha)$ converge to $\Lambda_{\qr}(\alpha)$ and $-\Lambda_{\qr}(\alpha)$, respectively, as $r\to\infty$.
    \end{enumerate}

    \item Approximate eigenvectors:
        \begin{enumerate}
            \item For $x\in \gV_{1}$ and $\alpha > \qr^{-2} + 1$, let $\rvu_{+} = \{\ervu_{j}\}_{j=0}^{r}$, $\rvu_{-} = \{(-1)^{j}\ervu_{j}\}_{j=0}^{r}$ have the components
            \begin{align}\label{eqn:Z1Eigenvec}
                     \ervu_0 \in \R \setminus \{0\}, \quad \ervu_1 &\,= \alpha^{-1/2} \Lambda_{\qr^{-1}} (\alpha)\ervu_0, \quad \ervu_2 = \alpha^{1/2}\qr(\alpha - \qr^{-2})^{-1} \ervu_0\,,\\
                      \ervu_{2j+1} &\, = (\alpha - \qr^{-2})^{-j} \ervu_1, \quad \ervu_{2j+2}= (\alpha -\qr^{-2})^{-j} \ervu_2, \quad 1 \leq j \leq \lfloor \frac{r-3}{2} \rfloor.
            \end{align}
            Here, $\ervu_0$ is chosen to ensure $\|\rvu_{+}\|_2 = 1$. Then as $r\to \infty$, $\rvu_{+}$ and $\rvu_{-}$ are the approximate eigenvectors of $\rmZ^{(1)}(\alpha)$ corresponding to $\Lambda_{\qr^{-1}}(\alpha)$ and $ -\Lambda_{\qr^{-1}}(\alpha)$, respectively.
            \item For $x\in \gV_{2}$, and $\alpha > \qr^{2} + 1$ or $0<\alpha < \qr^{2} - 1$, define $\rvu_{+} = \{\ervu_{j}\}_{j=0}^{r}$ and $\rvu_{-} = \{(-1)^{j}\ervu_{j}\}_{j=0}^{r}$ by
                \begin{align}\label{eqn:Z2Eigenvec}
                    \ervu_0 \in \R \setminus \{0\}, \quad \ervu_1 &\,= \alpha^{-1/2} \Lambda_{\qr} (\alpha)\ervu_0, \quad \ervu_2 = \alpha^{1/2}\qr^{-1}(\alpha - \qr^2)^{-1} \ervu_0\,,\\
                     \ervu_{2j+1} &\,= (\alpha - \qr^2)^{-j} \ervu_1, \quad \ervu_{2j+2} = (\alpha -\qr^2)^{-j} \ervu_2, \quad 1 \leq j \leq \lfloor (r-3)/2 \rfloor.
                \end{align}
            Then as $r\to \infty$, $\rvu_{+}$ and $\rvu_{-}$ are the approximate eigenvectors of $\rmZ^{(2)}(\alpha)$ corresponding to $\Lambda_{\qr}(\alpha)$ and $ -\Lambda_{\qr}(\alpha)$, respectively.
        \end{enumerate}
\end{enumerate}
\end{lemma}

\begin{lemma}\label{lem:tree-norm}
    Let $\tree$ be a bipartite graph whose vertices have degree at most $q+1$ for some $q\geq 1$. Then $\|\mA^{\tree}\|\leq q+1$. If, in addition, $\tree$ is a tree, then $\|\mA^{\tree}\| \leq 2\sqrt{q}$.
\end{lemma}
\begin{proof}[Proof of Lemma~\ref{lem:tree-norm}]
Every row and every column of $\mA^{\tree}$ contains at most $q+1$
nonzero entries, all of which are equal to one. Hence its maximum row
sum and maximum column sum are both at most $q+1$. Lemma
\ref{lem:schur_test} therefore gives
\begin{align}
    \|\mA^{\tree}\|
    \leq\sqrt{(q+1)(q+1)}
    =q+1.
\end{align}

Suppose now that $\tree$ is a tree. Choose a root $o$ and orient every
edge away from $o$, so that $x\to y$ means that $x$ is the parent of
$y$. Let $\rvf=(\ervf_x)_{x\in \mathrm V(\tree)}$ be finitely
supported. For every oriented edge $x\to y$, we have
\begin{align}
    0
    &\leq
    \left(q^{-1/4}|\ervf_x|-q^{1/4}|\ervf_y|\right)^2
    =
    q^{-1/2}|\ervf_x|^2+q^{1/2}|\ervf_y|^2
    -2|\ervf_x\ervf_y|.
\end{align}
Rearranging this inequality yields
\begin{align}
    2|\ervf_x\ervf_y|
    \leq
    q^{-1/2}|\ervf_x|^2+q^{1/2}|\ervf_y|^2
    \label{eqn:tree-norm-oriented-edge-bound}
\end{align}
Since every unoriented edge appears exactly once in the oriented edge
set, summing \eqref{eqn:tree-norm-oriented-edge-bound} over all edges
gives
\begin{align}
    \left|\langle \rvf,\mA^{\tree}\rvf\rangle\right|
    &\leq
    2\sum_{x\to y}|\ervf_x\ervf_y|
    \leq
    \sum_{x\to y}
      \left(q^{-1/2}|\ervf_x|^2+q^{1/2}|\ervf_y|^2\right).
    \label{eqn:tree-norm-oriented-sum}
\end{align}
We now collect the terms in the last sum according to their vertices.
The root appears only as a parent, whereas every other vertex appears
once as a child and once for each of its own children. Thus the
right-hand side of \eqref{eqn:tree-norm-oriented-sum} equals
\begin{align}
    &\frac{1}{\sqrt q}
      \#\{y:o\to y\}\,|\ervf_o|^2
    +
    \sum_{x\in \mathrm V(\tree)\setminus\{o\}}
      \left(
        \sqrt q+\frac{1}{\sqrt q}
        \#\{y:x\to y\}
      \right)|\ervf_x|^2.
    \label{eqn:tree-norm-vertex-collection}
\end{align}
The root has at most $q+1$ children. Since $q\geq1$, its coefficient in
\eqref{eqn:tree-norm-vertex-collection} satisfies
\begin{align}
    \frac{1}{\sqrt q}\#\{y:o\to y\}
    \leq
    \frac{q+1}{\sqrt q}
    \leq2\sqrt q.
\end{align}
Every vertex $x\neq o$ has exactly one parent. It therefore has at most
$q$ children, because its total degree is at most $q+1$. Consequently,
its coefficient in \eqref{eqn:tree-norm-vertex-collection} satisfies
\begin{align}
    \sqrt q+\frac{1}{\sqrt q}\#\{y:x\to y\}
    \leq
    \sqrt q+\frac{q}{\sqrt q}
    =2\sqrt q.
\end{align}
Combining these coefficient bounds with
\eqref{eqn:tree-norm-oriented-sum} gives
\begin{align}
    \left|\langle \rvf,\mA^{\tree}\rvf\rangle\right|
    \leq2\sqrt q\,\|\rvf\|^2.
\end{align}
Since $\mA^{\tree}$ is self-adjoint, taking the supremum over unit
vectors proves $\|\mA^{\tree}\|\leq2\sqrt q$. For an infinite graph,
both bounds follow by applying the corresponding estimates to finitely
supported vectors and then using density.
\end{proof}

\begin{lemma}
    Let $\tree$ be the infinite $(q, d_1, d_2)$-biregular tree rooted at $x$ as in Definition~\ref{def:biregular_tree}. Let $\mA$ denote the adjacency matrix of $\tree$ and $d = \sqrt{d_{1} d_{2}}$. Denote $\alpha = q /d$ where $d^{2} = d_{1} d_{2}$. Then the following holds.
    \begin{align}
        \|\mA /\sqrt{d}\| \leq \Lambda_{\qr^{-1}}(\alpha \vee \qr^{2}) \indi{x\in \gV_{2}} + \Lambda_{\qr}(\alpha \vee \qr^{-2}) \indi{x\in \gV_{1}}.
    \end{align}
\end{lemma}
\begin{proof}
    It follows directly from Lemma~\ref{lem:tridiagrep} and Lemma~\ref{lem:eigenZ1Z2}.
\end{proof}

\subsection{Stieltjes transforms of spectral measures}\label{sec:stieltjes-transform-spectral-measures-biregular-tree}

The comparison profiles in \eqref{eqn:stieltjes-alpha-inverse} admit a
root-resolvent interpretation.

\begin{lemma}[Root spectral interpretation]
\label{lem:spectra-biregular-tree}
Let $\ve_0\in\ell^2(\N)$ be the first coordinate vector.  For
$j\in\{1,2\}$ and $\alpha\geq0$,
\begin{align}
    \stj_\alpha^{(j)}(z)
    =
    \left\langle
      \ve_0,(\mZ^{(j)}(\alpha)-z)^{-1}\ve_0
    \right\rangle,
    \qquad z\in\C_+,
    \label{eqn:comparison-transform-root-resolvent}
\end{align}
where $\mZ^{(j)}(\alpha)$ is the infinite Jacobi matrix determined by
\eqref{eqn:Z1Z2}.  Hence $\stj_\alpha^{(j)}$ is the Stieltjes transform
of its spectral measure at the root.

Moreover, let $\tree$ be the infinite $(q,d_1,d_2)$-biregular tree
rooted at $x\in\gV_j$, let $\mA$ be its adjacency matrix, and set
$d=\sqrt{d_1d_2}$ and $\alpha=q/d$.  Then
\begin{align}
    \stj_\alpha^{(j)}(z)
    =
    \left\langle
      \ones_x,(\mA/\sqrt d-z)^{-1}\ones_x
    \right\rangle.
    \label{eqn:comparison-transform-tree-root}
\end{align}
\end{lemma}

\begin{proof}[Proof of Lemma~\ref{lem:spectra-biregular-tree}]
Set
\begin{align}
    g_\alpha^{(j)}(z)
    \coloneqq
    \left\langle
      \ve_0,(\mZ^{(j)}(\alpha)-z)^{-1}\ve_0
    \right\rangle.
\end{align}
Schur complementation at the first coordinate gives
\begin{align}
    g_\alpha^{(1)}
    &=\bigl(-z-\alpha g_{\qr^2}^{(2)}\bigr)^{-1},
    &
    g_\alpha^{(2)}
    &=\bigl(-z-\alpha g_{\qr^{-2}}^{(1)}\bigr)^{-1}.
    \label{eqn:tree-root-schur-system}
\end{align}
At $\alpha=\qr^{-2}$ in the first identity and
$\alpha=\qr^2$ in the second, this coupled system, together with the
Stieltjes asymptotic $g_\alpha^{(j)}(z)\sim-z^{-1}$ at infinity, selects
the branches in \eqref{eqn:stieltjes-V1} and
\eqref{eqn:stieltjes-V2}.  Thus
$g_{\qr^{-2}}^{(1)}=\stj^{(1)}$ and
$g_{\qr^2}^{(2)}=\stj^{(2)}$.  The general identities in
\eqref{eqn:tree-root-schur-system} and the definition
\eqref{eqn:stieltjes-alpha-inverse} then give
$g_\alpha^{(j)}=\stj_\alpha^{(j)}$ for all $\alpha\geq0$.

The tree identity follows from the tridiagonalization in
Lemma~\ref{lem:tridiagrep}.
\end{proof}

The following bounds on the comparison functions will be used both in the
local-law proof and in the deduction of eigenvector delocalization.  They
also make explicit that the estimates are uniform in the possibly
atypical row-square parameter \(\alpha\).

\begin{lemma}[Bounds for the comparison Stieltjes transforms]
\label{lem:stieltjes-bound}
Recall the Stieltjes transforms
$\{\stj^{(j)}_{\alpha}(z)\}_{j=1,2}$ in
\eqref{eqn:stieltjes-alpha-inverse}. Fix
$0<\kappa<\qr^{-1}$, $L\geq1$, and $K>0$. The following statements
hold.
\begin{enumerate}
    \item There are constants
    $C_{\rm St}=C_{\rm St}(\kappa,\qr,L)>0$ and
    $c_{\rm St,K}=c_{\rm St,K}(\kappa,\qr,K,L)>0$ such that
    \begin{align}
        \sup_{\alpha\geq0}\max_{j=1,2}
        |\stj^{(j)}_{\alpha}(z)|
        &\leq C_{\rm St},
        &&z\in\specdomain,\notag\\
        \inf_{0\leq\alpha\leq K}\min_{j=1,2}
        |\stj^{(j)}_{\alpha}(z)|
        &\geq c_{\rm St,K},
        &&z\in\specdomain.
        \label{eqn:stieltjes-bound-asymp}
    \end{align}
    \item Uniformly for $z\in\specdomain$ and $0\leq\alpha\leq K$,
    \begin{align}
        |\stj^{(1)}_{\alpha}(z)-\stj^{(1)}_{\qr^{-2}}(z)|
        &\leq C|\alpha-\qr^{-2}|,\notag\\
        |\stj^{(2)}_{\alpha}(z)-\stj^{(2)}_{\qr^{2}}(z)|
        &\leq C|\alpha-\qr^{2}|,
        \label{eqn:stieltjes-Lipschitz-alpha}
    \end{align}
    where $C=C(\qr,\kappa,K,L)>0$.
    \item For every $\alpha\geq0$ and $z\in\C_{+}$,
    \begin{align}
        \max_{j=1,2}|\stj_{\alpha}^{(j)}(z)|
        \leq \frac{1}{\Im z}.
        \label{eqn:stieltjes-uniform-resolvent-bound}
    \end{align}
    Moreover, uniformly for all $\alpha\geq0$ and
    $z,w\in\specdomain$,
    \begin{align}
        \max_{j=1,2}
        |\stj_{\alpha}^{(j)}(z)-\stj_{\alpha}^{(j)}(w)|
        \leq N^{2}|z-w|.
        \label{eqn:stieltjes-Lipschitz-z}
    \end{align}
\end{enumerate}
\end{lemma}

\begin{proof}[Proof of Lemma~\ref{lem:stieltjes-bound}]
    \noindent\textbf{Part (1): uniform bounds on the Stieltjes transforms.}
Due to
\[
    \stj^{(1)}_{\alpha}(z) =\, (-z - \alpha \stj^{(2)}(z))^{-1}, \quad \stj^{(2)}_{\alpha}(z) =\, (-z - \alpha \stj^{(1)}(z))^{-1},
\]
it suffices to prove the upper and lower bounds on $\stj^{(1)}(z)$ and $\stj^{(2)}(z)$. We proceed with the lower bound first. We write $z=E+\ii\eta \in \specdomain$ in \eqref{eqn:spectral_region}, where $E \in\bulk$ and $\eta\in[0,L]$, then
\begin{align}
    |E|\geq\kappa,\qquad |E|-\ledge\geq\kappa,\qquad \redge-|E|\geq\kappa,
\end{align}
where we denote $\lambda_{\pm} = \qr \pm \qr^{-1}$. Hence the function $h(E)$ satisfies
\begin{align}
    h(E) \coloneqq \sqrt{(E^{2}-\ledge^{2})(\redge^{2}-E^{2})} \geq \kappa^{2}.
\end{align}
The square root term in $\stj^{(1)}(z)$ and $\stj^{(2)}(z)$ is the branch satisfying $\sqrt{(z^{2}-\redge^{2})(z^{2}-\ledge^{2})}\asymp z^{2}$ as $|z|\to\infty$ in $\C_{+}$. Equivalently, it can be written as the product
    \[
        \sqrt{z^{2}-\redge^{2}}\sqrt{z^{2}-\ledge^{2}},
    \]
    where each factor is the branch asymptotic to $z$ at infinity. Hence, for $E\in\bulk$,
    \[
        \sqrt{((E+\ii0)^{2}-\redge^{2})((E+\ii0)^{2}-\ledge^{2})}
        =
        \begin{cases}
            \ii h(E),& E>0,\\
            -\ii h(E),& E<0.
        \end{cases}
    \]
    By the explicit formulas in \eqref{eqn:stieltjes-V1} and
    \eqref{eqn:stieltjes-V2}, the boundary imaginary parts are
    \begin{align}
        \lim_{\eta\downarrow0}\Im\stj^{(1)}(E+\ii\eta)
        &= \frac{h(E)}{2\qr^{2}|E|}
        \geq \frac{\kappa^{2}}{2\qr^{2}\redge},
        \label{eqn:m1-boundary-lower}\\
        \lim_{\eta\downarrow0}\Im\stj^{(2)}(E+\ii\eta)
        &= \frac{\qr^{2}h(E)}{2|E|}
        \geq \frac{\qr^{2}\,\kappa^{2}}{2\redge}.
        \label{eqn:m2-boundary-lower}
    \end{align}
For $\eta>0$, both Stieltjes transforms have strictly positive
imaginary part. Together with
\eqref{eqn:m1-boundary-lower}--\eqref{eqn:m2-boundary-lower} and
continuity on the compact set
\[
    \{E+\ii\eta:E\in\bulk,\ 0\leq\eta\leq L\},
\]
this yields a constant $c_0=c_0(\kappa,\qr,L)>0$ such that
    \begin{align}
        \inf_{z\in\specdomain}
        \bigl(\Im\stj^{(1)}(z)\wedge\Im\stj^{(2)}(z)\bigr)
        \geq c_0.
        \label{eqn:base-stieltjes-lower}
    \end{align}
For the upper bounds, since $|E|\geq\kappa$, $|E|\leq\redge$, and $h(E)\leq\redge^{2}$ on $\bulk$, the same explicit formulas give
    \begin{align}
        |\stj^{(1)}(E+\ii0)|
        &\leq
        \frac{\redge^{2}+\qr^{2}+\qr^{-2}+\redge^{2}}
        {2\qr^{2}\kappa},\notag\\
        |\stj^{(2)}(E+\ii0)|
        &\leq
        \frac{\qr^{2}(\redge^{2}+\qr^{2}+\qr^{-2}
        +\redge^{2})}{2\kappa}.
    \end{align}
    By continuity of the explicit formulas on the same compact set
    $\{E+\ii\eta:E\in\bulk,\ 0\leq\eta\leq L\}$, after extending by the
    boundary values above, there is a constant $C_{0}>0$, depending
    only on $\kappa$, $\qr$, and $L$, such that for all
    $z\in\specdomain$,
    \begin{align}
        |\stj^{(1)}(z)|+|\stj^{(2)}(z)|\leq C_{0}.
        \label{eqn:base-stieltjes-upper}
    \end{align}
We now control the denominators in the inverse relationship explicitly. For
$j,k\in\{1,2\}$ with $j\neq k$, write
\[
    D_{\alpha}^{(j)}(z)\coloneqq -z-\alpha \stj^{(k)}(z),
    \qquad
    \stj_{\alpha}^{(j)}(z)=D_{\alpha}^{(j)}(z)^{-1}.
\]
If $0\leq \alpha\leq \kappa/(2C_{0})$, then
\[
    |D_{\alpha}^{(j)}(z)|
    \geq |z|-\alpha |\stj^{(k)}(z)|
    \geq |E|-\alpha C_{0}
    \geq \frac{\kappa}{2}.
\]
If $\alpha>\kappa/(2C_{0})$, then \eqref{eqn:base-stieltjes-lower}
gives
\[
    |D_{\alpha}^{(j)}(z)|
    \geq |\Im D_{\alpha}^{(j)}(z)|
    =
    \eta+\alpha\Im \stj^{(k)}(z)
    \geq
    \frac{\kappa c_{0}}{2C_{0}}.
\]
Therefore,
\[
    |\stj_{\alpha}^{(j)}(z)|
    \leq
    \max\left\{\frac{2}{\kappa},\frac{2C_{0}}{\kappa c_{0}}\right\}.
\]
On the other hand, since $\alpha\leq K$ and $z\in\specdomain$,
\[
    |D_{\alpha}^{(j)}(z)|
    \leq |z|+\alpha |\stj^{(k)}(z)|
    \leq \sqrt{\redge^{2}+L^{2}}+KC_{0},
\]
and hence
\[
    |\stj_{\alpha}^{(j)}(z)|
    \geq
    \left(\sqrt{\redge^{2}+L^{2}}+KC_{0}\right)^{-1}.
\]
Then \eqref{eqn:stieltjes-bound-asymp} is established by renaming the constants, respectively.

\noindent\textbf{Part (2): Lipschitz continuity with respect to $\alpha$.}
    By \eqref{eqn:stieltjes-alpha-inverse}, we obtain
    \begin{align}
        \stj_{\alpha}^{(1)}(z)-\stj^{(1)}(z)
        &=(\alpha-\qr^{-2})\stj^{(2)}(z) \stj_{\alpha}^{(1)}(z)\stj^{(1)}(z),\\
        \stj_{\alpha}^{(2)}(z)-\stj^{(2)}(z)
        &=(\alpha-\qr^{2})\stj^{(1)}(z)\stj_{\alpha}^{(2)}(z)\stj^{(2)}(z).
    \end{align}
    Then \eqref{eqn:stieltjes-Lipschitz-alpha} can be established by using \eqref{eqn:stieltjes-bound-asymp} and \eqref{eqn:base-stieltjes-upper}.

    \noindent\textbf{Part (3): uniform resolvent and Lipschitz bounds.}
    Fix any $\alpha\geq0$. For $\mZ^{(j)}(\alpha)$ in \eqref{eqn:Z1Z2}, define
    \[
        \mR_{\alpha}^{(j)}(z)
        \coloneqq(\mZ^{(j)}(\alpha)-z)^{-1}, \qquad j\in \{1, 2\}.
    \]
    By Lemma~\ref{lem:spectra-biregular-tree},
    $\stj_\alpha^{(j)}(z)=
    \langle\ve_0,\mR_\alpha^{(j)}(z)\ve_0\rangle$.
    The resolvent identity gives
    \[
        \mR_{\alpha}^{(j)}(z)-\mR_{\alpha}^{(j)}(w)
        =(z-w)\mR_{\alpha}^{(j)}(z)\mR_{\alpha}^{(j)}(w).
    \]
    Since $\mZ^{(j)}(\alpha)$ is self-adjoint, the spectral theorem gives
    \[
        \|\mR_{\alpha}^{(j)}(z)\|\leq \frac{1}{\Im z},
        \qquad
        \|\mR_{\alpha}^{(j)}(w)\|\leq \frac{1}{\Im w}.
    \]
    Taking the $\ve_{0}$-$\ve_{0}$ matrix element in the first of these bounds proves \eqref{eqn:stieltjes-uniform-resolvent-bound}. Taking the same matrix element in the resolvent identity, we obtain
    \begin{align}
        |\stj_{\alpha}^{(j)}(z)-\stj_{\alpha}^{(j)}(w)|
        &\, =
        |\<\ve_{0},(\mR_{\alpha}^{(j)}(z)-\mR_{\alpha}^{(j)}(w))\ve_{0}\>|\notag\\
        &\, \leq |z-w|\, \|\mR_{\alpha}^{(j)}(z)\|\, \|\mR_{\alpha}^{(j)}(w)\|
        \leq \frac{|z-w|}{\Im z\,\Im w} \leq N^{2}|z-w|,
    \end{align}
    where the last inequality uses $\Im z,\Im w\geq N^{-1}$. Since the right-hand side is independent of $\alpha$, taking the supremum over all $\alpha\geq0$ proves \eqref{eqn:stieltjes-Lipschitz-z}.
\end{proof}



\section{Properties of Bipartite Erd\H{o}s--R\'enyi Graphs}
\label{sec:bipartite-random-graph-inputs}

Under Assumption~\ref{ass:proportional_regime}, we give degree and
component estimates, determine the kernel and singular vectors on small
components, and construct exact localized singular vectors from
star-tuning forks in the giant. We then compare the centered matrix
with its giant block. Degree-window concentration is stated in
Proposition~\ref{prop:uniform-degree-counts}.
All centered restrictions use the original parameters $p,d$.

We begin with the factorial and binomial-coefficient estimates used below.

\begin{lemma}[Stirling and binomial-coefficient estimates]
\label{lem:stirling-formula}
The following estimates hold.
\begin{enumerate}
    \item[(i)] For every integer $n\geq1$,
    \begin{align}
        n!
        =
        \sqrt{2\pi n}
        \left(\frac n{\mathrm e}\right)^n
        \exp\bigl(O(n^{-1})\bigr).
    \end{align}

    \item[(ii)] Define the falling factorial $(n)_k\coloneqq n!/(n-k)!$. Then uniformly for $1\leq k\leq n$,
    \begin{align}
        (n)_k
        =
        n^k\exp\left(O\left(\frac{k^2}{n}\right)\right).
    \end{align}

    \item[(iii)] Uniformly for $1\leq k\leq n$,
    \begin{align}
        \left(\frac nk\right)^k
        \leq
        \binom nk
        \leq
        \left(\frac{\mathrm e n}{k}\right)^k.
    \end{align}
\end{enumerate}
\end{lemma}
\begin{proof}[Proof of Lemma~\ref{lem:stirling-formula}]
Item~(i) is Stirling's formula with its standard remainder estimate.
For item~(ii), if $k\leq n/2$, expanding
$\sum_{j=0}^{k-1}\log(1-j/n)$ shows that the logarithm of
$(n)_k/n^k$ is $O(k^2/n)$.  If $k>n/2$, then
$n!/n^n\leq(n)_k/n^k\leq1$; item~(i) shows that the logarithm of the
lower bound is $-O(n)=-O(k^2/n)$.

For the lower bound in item~(iii), write
$\binom nk=\prod_{j=0}^{k-1}(n-j)/(k-j)$ and bound each factor from
below by $n/k$.  For the upper bound, use
$\binom nk\leq n^k/k!$ together with
$k!\geq(k/\mathrm e)^k$.
\end{proof}

\subsection{Degree probabilities and maximum degree}
\label{subsec:bipartite-degree-probabilities}

We now specialize to the two degree laws.  Whenever $d=O(\log N)$ is
assumed below, constants may depend on its fixed implicit upper bound. Define the binomial random variables $\rD^{(1)}$ and $\rD^{(2)}$ by
\begin{align}
    \rD^{(1)}\sim\Bin(m,p),
    \qquad
    \rD^{(2)}\sim\Bin(n,p),
    \label{eqn:bipartite-degree-laws}
\end{align}
where $p = d/\sqrt{mn}$ and $d_{1} \coloneqq \E\rD^{(1)} = d\qr^{-2}+O_{\qr}(d/N)$,
$d_{2} \coloneqq \E\rD^{(2)} = d\qr^{2}+O_{\qr}(d/N)$.

Recall $I_s$ from \eqref{eqn:transition-rate-function}. For
$s\in\{\qr^{-1},\qr\}$, $d>0$, and $\alpha>0$, set
\begin{align}
    f_{s,d}(\alpha)
    &\coloneqq
    dI_s(\alpha)+\frac12\log(2\pi\alpha d),
    \label{eqn:bipartite-degree-rate}
\end{align}
For every integer $r\geq0$, write
\[
\alpha_r\coloneqq r/d.
\]
Whenever a real endpoint occurs in a degree event, we round $\alpha d$
to the indicated adjacent integer. Errors from this convention are
absorbed in the displayed bounds.

\begin{lemma}[Degree probabilities: Poisson and Stirling forms]
\label{lem:bipartite-poisson-approx}
\label{cor:bipartite-degree-stirling}
Suppose that $d=O(\log N)$.  Uniformly for $0\leq k\leq m$ and
$0\leq \ell\leq n$,
\begin{subequations}
\begin{align}
    \P(\rD^{(1)}=k)
    &=
    \frac{d_1^k}{k!}\mathrm e^{-d_1}
    \exp\left(
      O_{\qr}\left(\frac{k^2+d^2}{N}\right)
    \right),
    \label{eqn:poisson-approx-V1}\\
    \P(\rD^{(2)}=\ell)
    &=
    \frac{d_2^\ell}{\ell!}\mathrm e^{-d_2}
    \exp\left(
      O_{\qr}\left(\frac{\ell^2+d^2}{N}\right)
    \right).
    \label{eqn:poisson-approx-V2}
\end{align}
\end{subequations}

For positive degrees $1\leq k\leq m$ and $1\leq\ell\leq n$,
the Stirling forms are
\begin{subequations}
\begin{align}
    \P(\rD^{(1)}=k)
    &=
    \exp\left(
      -f_{\qr^{-1},d}(\alpha_k)
      +O_{\qr}\left(\frac1k+\frac{k^2+d^2}{N}\right)
    \right),
    \label{eqn:degree-stirling-V1}\\
    \P(\rD^{(2)}=\ell)
    &=
    \exp\left(
      -f_{\qr,d}(\alpha_\ell)
      +O_{\qr}\left(\frac1\ell+\frac{\ell^2+d^2}{N}\right)
    \right).
    \label{eqn:degree-stirling-V2}
\end{align}
\end{subequations}
\end{lemma}
\begin{proof}
For $i\in\{1,2\}$, put $L_i=|\gV_{3-i}|$, $s_1=\qr^{-1}$ and
$s_2=\qr$, so $d_i=L_ip=ds_i^2+O_{\qr}(d/N)$. The exact mass is
\[
 \P(\rD^{(i)}=k)=\frac{(L_i)_k}{k!}p^k(1-p)^{L_i-k},
 \qquad 0\leq k\leq L_i,
\]
with $(L_i)_0=1$. Lemma~\ref{lem:stirling-formula} and the expansion
of $\log(1-p)$ give the Poisson forms
\eqref{eqn:poisson-approx-V1}--\eqref{eqn:poisson-approx-V2}:
the logarithmic error is
$O(k^2/L_i+L_ip^2+kp)=O_{\qr}((k^2+d^2)/N)$, using
$kp\leq k^2/L_i+L_ip^2$. This also covers $k=0$.
For $k\geq1$, Stirling's formula gives
\[
 -\log\P(\rD^{(i)}=k)
 =k\log(k/d_i)-k+d_i+\frac12\log(2\pi k)
 +O_{\qr}\left(\frac1k+\frac{k^2+d^2}{N}\right).
\]
Replacing $d_i$ by $ds_i^2$ costs $O_{\qr}((k+d)/N)$, which is absorbed
in the error. The resulting main term is $f_{s_i,d}(k/d)$, proving
\eqref{eqn:degree-stirling-V1}--\eqref{eqn:degree-stirling-V2}.
\end{proof}

We close the subsection with a uniform upper bound on all degrees.  Its
scale interpolates between the dense range $d\gtrsim\log N$ and the
sparser range in which the largest degree is substantially larger than
its mean.

\begin{lemma}[Maximum degree]\label{lem:bipartite-max-degree}
For every $\nu>0$ there is a constant $C_{\qr,\nu}>0$ such that, with probability at least $1-C_{\qr,\nu}N^{-\nu}$, the following holds:
\begin{align}
    \max_{x\in\gV}\rD_x
    \leq
    C_{\qr,\nu}\left(
      d+\frac{\log N}
      {\log(1+\log N/d)}
    \right).
    \label{eqn:bipartite-max-degree}
\end{align}
\end{lemma}
\begin{proof}[Proof of Lemma~\ref{lem:bipartite-max-degree}]
For $i\in\{1,2\}$, a degree in $\gV_i$ has law
$\Bin(|\gV_{3-i}|,p)$ and mean $d_i\leq C_{\qr}d$.
The Chernoff upper bound \eqref{eqn:degree-chernoff-upper} gives, for
$k\geq\mathrm e d_i$,
\begin{align}
    \P(\rD^{(i)}\geq k)
    \leq
    \exp\left[-d_i\benrate\left(\frac{k}{d_i}-1\right)\right]
    \leq
    \exp\left[-k\log\left(\frac{k}{\mathrm e d_i}\right)\right],
    \label{eqn:maximum-degree-bennett}
\end{align}
where the second inequality follows from
$d_i\benrate(k/d_i-1)=k\log(k/d_i)-k+d_i$.

Put $t=\log N/d$ and
$k_*=C_{\qr,\nu}\bigl(d+\log N/\log(1+t)\bigr)$.
For $C_{\qr,\nu}$ sufficiently large, $k_*\geq\mathrm e d_i$ and,
uniformly in $i\in\{1,2\}$,
\begin{align}
    k_*\log\left(\frac{k_*}{\mathrm e d_i}\right)
    \geq(\nu+2)\log N.
    \label{eqn:maximum-degree-exponent}
\end{align}
Indeed, if $t\leq2$, then $d\geq\frac12\log N$ and
$\log(k_*/(\mathrm e d_i))\geq1$ after increasing $C_{\qr,\nu}$.
If $t>2$, then $\log(1+t)\leq\sqrt t$, so
$k_*/d\geq C_{\qr,\nu}\sqrt t$. Thus
$\log(k_*/(\mathrm e d_i))\geq c\log(1+t)$ for an absolute $c>0$
once $C_{\qr,\nu}$ is large enough. Combining this with
$k_*\geq C_{\qr,\nu}\log N/\log(1+t)$ proves
\eqref{eqn:maximum-degree-exponent} in the second case as well.
The Chernoff bound and a union bound now give
\begin{align}
    \P\left(\max_{x\in\gV}\rD_x\geq k_*\right)
    \leq N\max_{i\in\{1,2\}}
    \exp\left[-k_*\log\left(\frac{k_*}{\mathrm e d_i}\right)\right]
    \leq N^{-\nu-1},
\end{align}
which proves \eqref{eqn:bipartite-max-degree}.
\end{proof}

\subsection{Isolated vertices and component counts}
\label{subsec:bipartite-vertex-counts}

We count isolated vertices and components, obtaining a common event
for the giant and its complement and sharp counts for each fixed tree.
All critical-sparsity statements below use fixed $\qr$.

\paragraph{Isolated vertices.}
\begin{lemma}[Isolated-vertex counts]
\label{lem:bipartite-isolated-vertices}
Suppose that Assumptions~\ref{ass:proportional_regime} and
\ref{ass:critical_sparsity} hold with fixed $\qr$. Let $\rZ_i$ count
isolated vertices in $\gV_i$. Then
\begin{align}
 \E\rZ_1
 &=\left(\frac{\qr^4}{1+\qr^4}+o(1)\right)N^{1-b\qr^{-2}},
 &\E\rZ_2
 &=\left(\frac1{1+\qr^4}+o(1)\right)N^{1-b\qr^2}.
 \label{eqn:bipartite-isolated-expectations}
\end{align}
For every $\varepsilon,\nu>0$ and all sufficiently large $N$, if
$b<\qr^2$ then
\begin{align}
 \P\left(
 \left|N^{-(1-b\qr^{-2})}\rZ_1-\frac{\qr^4}{1+\qr^4}\right|
 \leq\varepsilon\right)
 \geq1-C_{b,\qr,\varepsilon,\nu}N^{-\nu}.
 \label{eqn:bipartite-isolated-V1-count}
\end{align}
If $b<\qr^{-2}$ then
\begin{align}
 \P\left(
 \left|N^{-(1-b\qr^2)}\rZ_2-\frac1{1+\qr^4}\right|
 \leq\varepsilon\right)
 \geq1-C_{b,\qr,\varepsilon,\nu}N^{-\nu}.
 \label{eqn:bipartite-isolated-V2-count}
\end{align}
At the respective thresholds,
\begin{align}
 \rZ_1&\xrightarrow{\mathrm d}
 \Pois\left(\frac{\qr^4}{1+\qr^4}\right)
 &&\textnormal{if }b=\qr^2,
 \nonumber\\
 \rZ_2&\xrightarrow{\mathrm d}
 \Pois\left(\frac1{1+\qr^4}\right)
 &&\textnormal{if }b=\qr^{-2}.
 \label{eqn:bipartite-isolated-threshold-poisson}
\end{align}
\end{lemma}
\begin{proof}[Proof of Lemma~\ref{lem:bipartite-isolated-vertices}]
Isolation events within a vertex class involve disjoint edge sets, so
\begin{align}
 \rZ_1&\sim\Bin\bigl(n,(1-p)^m\bigr),
 &\rZ_2&\sim\Bin\bigl(m,(1-p)^n\bigr).
 \label{eqn:bipartite-isolated-binomial-laws}
\end{align}
Lemma~\ref{lem:bipartite-poisson-approx} at degree zero gives the
isolation probabilities
$N^{-b\qr^{-2}}(1+o(1))$ and $N^{-b\qr^2}(1+o(1))$.
This proves \eqref{eqn:bipartite-isolated-expectations}.
Below either threshold, write $\mu_i=\E\rZ_i=(c_i+o(1))N^{a_i}$,
where $a_i>0$ is fixed and $c_i$ is the corresponding constant in
\eqref{eqn:bipartite-isolated-expectations}. Bennett's inequality
(Lemma~\ref{lem:Bennett}) for the binomial law gives, for fixed
$\delta>0$,
\[
 \P(|\rZ_i-\mu_i|>\delta\mu_i)
 \leq2\exp\{-c_\delta\mu_i\}.
\]
Taking $\delta$ sufficiently small in terms of $\varepsilon$ proves
\eqref{eqn:bipartite-isolated-V1-count} and
\eqref{eqn:bipartite-isolated-V2-count}.

At the threshold for class $i$, write $\rZ_i\sim\Bin(L_i,q_{i,N})$
using \eqref{eqn:bipartite-isolated-binomial-laws}. The means above give
$L_iq_{i,N}\to c_i$ and $L_iq_{i,N}^2\to0$, where
$c_1=\qr^4/(1+\qr^4)$ and $c_2=(1+\qr^4)^{-1}$. For each fixed
$k\in\N_0$, the same falling-factorial and logarithmic expansion gives
\[
 \P(\rZ_i=k)
 =\frac{(L_iq_{i,N})^k}{k!}\exp\{-L_iq_{i,N}+o(1)\}
 \longrightarrow\mathrm e^{-c_i}\frac{c_i^k}{k!},
\]
proving \eqref{eqn:bipartite-isolated-threshold-poisson}.
\end{proof}

\begin{proof}[Proof of Lemma~\ref{lem:bipartite-isolated-column-threshold}]
Since $|\sI_2|=\rZ_2$, items~(ii) and~(iii) are the second line of
\eqref{eqn:bipartite-isolated-threshold-poisson} and
\eqref{eqn:bipartite-isolated-V2-count} in
Lemma~\ref{lem:bipartite-isolated-vertices}. For item~(i), Markov's
inequality and \eqref{eqn:bipartite-isolated-expectations} give
$\P(\rZ_2>0)\leq\E\rZ_2=\bigl((1+\qr^4)^{-1}+o(1)\bigr)N^{1-b\qr^2}$.
\end{proof}

\paragraph{Small-component counts.}
A component of size at least two contains a spanning tree and has no
boundary edges.

\begin{lemma}[Component counting]\label{lem:bipartite-component-counting}
Suppose that Assumption~\ref{ass:proportional_regime} holds and
$d\geq1$. Let $\rW_k$ count connected components of size $k$, and
let $\widehat{\rW}_k$ count those that are not trees. There are
$c_{\qr},C_{\qr}>0$ such that, for $1\leq k\leq N/2$,
\begin{align}
 \E\rW_k
 &\leq C_{\qr}N\exp\{-k(c_{\qr}d-\log d-C_{\qr})\},
 \label{eqn:bipartite-component-count}\\
 \E\widehat{\rW}_k
 &\leq C_{\qr}\exp\{-k(c_{\qr}d-\log d-C_{\qr})\}.
 \label{eqn:bipartite-nontree-component-count}
\end{align}
Writing $\tau=c_{\qr}d-\log d-C_{\qr}$, whenever $\tau\geq1$,
\begin{align}
 \sum_{k=K}^{\lfloor N/2\rfloor}\E\rW_k
 &\leq C_{\qr}N\mathrm e^{-K\tau}\quad(K\geq1),
&
 \sum_{k=1}^{\lfloor N/2\rfloor}k\E\rW_k
 &\leq C_{\qr}N\mathrm e^{-\tau},
&
 \sum_{k=1}^{\lfloor N/2\rfloor}\E\widehat{\rW}_k
 &\leq C_{\qr}\mathrm e^{-\tau}.
 \label{eqn:bipartite-component-tails}
\end{align}
\end{lemma}
\begin{proof}[Proof of Lemma~\ref{lem:bipartite-component-counting}]
For $k=1$, the exact isolation laws give
$\E\rW_1\leq N\mathrm e^{-c_{\qr}d}$ and
$\widehat{\rW}_1=0$. For $k=a+c\geq2$, choose $a\geq1$ rows and
$c\geq1$ columns and a spanning tree. The boundary has
$a(m-c)+c(n-a)$ possible edges. The bipartite Cayley formula
\cite{AbuSbeih1990OnTN} gives $a^{c-1}c^{a-1}$ labelled spanning trees
of $K_{a,c}$, hence
\begin{align}
 \E\rW_k
 \leq\sum_{\substack{a+c=k\\a,c\geq1}}
 \binom na\binom mc a^{c-1}c^{a-1}p^{k-1}
 (1-p)^{a(m-c)+c(n-a)}.
 \label{eqn:component-count-first}
\end{align}
Put $u=a/n$ and $v=c/m$. Since $k\leq N/2$, at least one of
$u,v$ is at most $1/2$; in either case
$u(1-v)+v(1-u)\geq(u+v)/4$. As $n\asymp_{\qr}m$, this yields
\begin{align}
 a(m-c)+c(n-a)\geq c_{\qr}k\sqrt{mn},
 \qquad
 (1-p)^{a(m-c)+c(n-a)}\leq\mathrm e^{-c_{\qr}dk}.
 \label{eqn:bipartite-boundary-suppression}
\end{align}
Using $a^{c-1}c^{a-1}\leq k^{k-2}$, Vandermonde's identity,
$\binom Nk\leq(\mathrm eN/k)^k$, and $Np\leq C_{\qr}d$, we obtain
\[
 \E\rW_k
 \leq\binom Nk k^{k-2}p^{k-1}\mathrm e^{-c_{\qr}dk}
 \leq\frac N{k^2}\mathrm e^k(Np)^{k-1}\mathrm e^{-c_{\qr}dk},
\]
which proves \eqref{eqn:bipartite-component-count} after adjusting the
constants. A non-tree component has, in addition, an edge outside the
chosen tree. A union bound supplies the factor $acp\leq k^2p$, hence
\[
 \E\widehat{\rW}_k
 \leq\binom Nk k^kp^k\mathrm e^{-c_{\qr}dk}
 \leq(\mathrm eNp)^k\mathrm e^{-c_{\qr}dk},
\]
proving \eqref{eqn:bipartite-nontree-component-count}.
Summing the geometric series, and its derivative for the weighted
sum, proves \eqref{eqn:bipartite-component-tails}.
\end{proof}

Let $\gG_{\rm{giant}}$ denote a largest component and call all other
components \emph{small}. Let $\chi_{\rm{giant}}$ denote the coordinate
projection onto $\gG_{\rm{giant}}$.
\begin{corollary}[Giant and small components]
\label{cor:components-size}
Suppose that Assumption~\ref{ass:proportional_regime} holds and
$d\gg1$. For every $\nu>0$, there are $c_{\qr}>0$ and $C_{\qr,\nu}>0$ such that, with probability at least
$
 1-C_{\qr,\nu}\left(\mathrm e^{-c_{\qr}d}+N^{-\nu}\right)
$,
the following hold simultaneously:
\begin{enumerate}
 \item $|\gG_{\rm{giant}}|\geq N(1-\mathrm e^{-c_{\qr}d})$;
 \item every small component has fewer than
 $\lceil C_{\qr,\nu}\log N/d\rceil$ vertices;
 \item every small component is a tree;
 \item the raw normalized adjacency operator on their union satisfies
 \begin{align}
  \left\|(\id-\chi_{\rm{giant}})\frac{\rmA}{\sqrt d}
   (\id-\chi_{\rm{giant}})\right\|
  \leq C_{\qr,\nu}\frac{\sqrt{\log N}}d.
  \label{eqn:bipartite-small-component-norm}
 \end{align}
 At critical sparsity $d=b\log N$, this bound is
 $O_{b,\qr,\nu}((\log N)^{-1/2})=o(1)$.
\end{enumerate}
\end{corollary}
\begin{proof}[Proof of Corollary~\ref{cor:components-size}]
For $d\to\infty$, the parameter $\tau$ in
\eqref{eqn:bipartite-component-tails} satisfies
$\tau\geq2c_{\qr}d$ after decreasing $c_{\qr}$. If $\rS$ counts
vertices in components of size at most $N/2$, the weighted tail bound
and Markov's inequality give
\[
 \E\rS\leq C_{\qr}N\mathrm e^{-2c_{\qr}d},
 \qquad
 \P(\rS>N\mathrm e^{-c_{\qr}d})
 \leq C_{\qr}\mathrm e^{-c_{\qr}d}.
\]
On the complementary event $\rS<N$, a component has more than $N/2$
vertices and its complement has size $\rS$, proving the first claim.
Every other component has size at most $N/2$. Choose
$k_0=\lceil C_{\qr,\nu}\log N/d\rceil$ with
$2c_{\qr}C_{\qr,\nu}>\nu+2$. The remaining two tail bounds give
\begin{align}
 \sum_{k=k_0}^{\lfloor N/2\rfloor}\E\rW_k
 &\leq C_{\qr}N\mathrm e^{-2c_{\qr}dk_0}
 \leq C_{\qr,\nu}N^{-\nu},
 &
 \sum_{k=1}^{\lfloor N/2\rfloor}\E\widehat{\rW}_k
 &\leq C_{\qr}\mathrm e^{-2c_{\qr}d}.
\end{align}
Markov's inequality and a union bound prove the first three assertions,
including $k_0=1$, when there are no small components.
For the fourth assertion, an isolated vertex has zero adjacency
operator, while Lemma~\ref{lem:tree-norm} bounds the adjacency norm
of a tree with $k\geq2$ vertices by $2\sqrt{k-1}$. Since $k<k_0$
implies $k<C_{\qr,\nu}\log N/d$, normalization by $\sqrt d$
gives the stated bound after increasing $C_{\qr,\nu}$. The restriction
to the union is a direct sum, whose norm is the maximum of the
component norms.
\end{proof}

\begin{lemma}[Fixed tree components]
\label{lem:fixed-tree-components}
Suppose that Assumptions~\ref{ass:proportional_regime} and
\ref{ass:critical_sparsity} hold with fixed $\qr$. Let $T$ be a fixed
bipartite tree with $a\geq1$ rows, $c\geq1$ columns, and $e=a+c-1$
edges. Let $Z_T$ count components isomorphic to $T$, preserving the
vertex classes, and let $\operatorname{aut}_{\rm bip}(T)$ count its
class-preserving automorphisms. Then
\begin{align}
 \E Z_T
 &=\frac{(n)_a(m)_c}{\operatorname{aut}_{\rm bip}(T)}
 p^e(1-p)^{am+cn-ac-e}\notag\\
 &=\bigl(C_{T,\qr}+o(1)\bigr)
 N d^e\exp\{-d(a\qr^{-2}+c\qr^2)\},
 \label{eqn:fixed-tree-component-mean}
\end{align}
where
\[
 C_{T,\qr}=\frac{\qr^{2(a-c+1)}}
 {(1+\qr^4)\operatorname{aut}_{\rm bip}(T)}>0.
\]
Writing $\mu_T=\E Z_T$, for $0<\varepsilon<1$ and sufficiently large
$N$,
\begin{align}
 \P(|Z_T-\mu_T|>\varepsilon\mu_T)
 \leq C_T\varepsilon^{-2}(\mu_T^{-1}+p).
 \label{eqn:fixed-tree-component-concentration}
\end{align}
Consequently, $Z_T/\E Z_T\to1$ in probability if
$b(a\qr^{-2}+c\qr^2)\leq1$, whereas $Z_T=0$ with probability tending
to one if this inequality is reversed.
\end{lemma}
\begin{proof}[Proof of Lemma~\ref{lem:fixed-tree-components}]
Each embedding requires the $e$ tree edges and forbids the other
$am+cn-ac-e$ edges incident to its vertices, proving the mean formula.
Distinct overlapping component copies cannot both occur. For disjoint
copies, the joint probability divided by the product of their
probabilities is $(1-p)^{-2ac}=1+O_T(p)$, so
\begin{align}
 \operatorname{Var}Z_T\leq\mu_T+C_Tp\mu_T^2.
 \label{eqn:fixed-tree-component-variance}
\end{align}
Chebyshev's inequality gives
\eqref{eqn:fixed-tree-component-concentration}. The mean diverges
below the threshold and also at equality, where
$\mu_T=(C_{T,\qr}+o(1))d^e\to\infty$. Above the threshold it tends to
zero, and Markov's inequality applies.
\end{proof}

\subsection{Centering and small-component spectra}
\label{subsec:bipartite-near-zero}
\label{subsec:zero-and-near-zero-modes}

Centering acts through coordinate sums. We first determine its effect
on kernels and component singular subspaces, then bound the least
positive singular value contributed by centering the zero columns of
$\widetilde{\rmA}$.

Since $n\geq m$, the least of the $m$ singular values of $\mX$ is
$\lambda_m(\scA)$; it vanishes exactly when $\ker\mX\neq\{0\}$.
Write $s_{\min}^{+}(\mX)$ for the least positive singular value and
$\sI_2$ for the isolated vertices in $\gV_2$, so $|\sI_2|=\rZ_2$.

A column with one neighbor is called pendant.
The centered action is
\begin{align}
 \sqrt d\,\mX\vv
 =\widetilde{\rmA}\vv-p(\ones_{\gV_2}^{\sT}\vv)\ones_{\gV_1},
 \qquad p^2n/d=d/m.
 \label{eqn:near-zero-centered-action}
\end{align}

\subsubsection{Column rank deficiency from isolated vertices and components}
\label{subsec:bipartite-exact-zeros}

\begin{lemma}[Raw and centered kernels]
\label{lem:centered-raw-kernel}
Let $p>0$.
\begin{enumerate}
 \item[(i)] $\ker\widetilde{\rmA}\cap\ones_{\gV_2}^{\perp}\subset\ker\mX$.
 Hence $\dim\ker\mX\geq\dim\ker\widetilde{\rmA}-1$.
 \item[(ii)] If $\widetilde{\rmA}$ has an isolated row, then
 $\ker\mX=\ker\widetilde{\rmA}\cap\ones_{\gV_2}^{\perp}$. Consequently,
 $\dim\ker\mX=\dim\ker\widetilde{\rmA}-1$ if some vector of
 $\ker\widetilde{\rmA}$ has nonzero coordinate sum, and
 $\dim\ker\mX=\dim\ker\widetilde{\rmA}$ otherwise. If there is also a
 zero column of $\widetilde{\rmA}$, then
 $\operatorname{rank}\mX=\operatorname{rank}\widetilde{\rmA}+1$.
\end{enumerate}
\end{lemma}
\begin{proof}
Equation~\eqref{eqn:near-zero-centered-action} gives (i), since imposing
one linear functional on $\ker\widetilde{\rmA}$ reduces its dimension
by at most one. If a row is isolated, its coordinate in $\mX\vv=0$
forces $\ones^{\sT}\vv=0$, and the remaining coordinates then give
$\widetilde{\rmA}\vv=0$. This proves (ii) and its dimension formulas.
A zero column of $\widetilde{\rmA}$ supplies a kernel vector of nonzero sum, so
rank--nullity gives the final rank identity.
\end{proof}

\begin{proof}[Proof of Lemma~\ref{lem:bipartite-isolated-column-zero-modes}]
For $\vv$ supported on $\sI_2$, \eqref{eqn:near-zero-centered-action}
reduces to $\sqrt d\,\mX\vv=-p(\sum_{y\in\sI_2}\evv_y)\ones$.
Thus its kernel is exactly the zero-sum subspace, of dimension
$(\rZ_2-1)_+$, including the empty case.
\end{proof}

A matching is a set of edges with no shared endpoints. For a finite
graph $F$, let $\mat(F)$ denote its maximum matching size. In particular,
\begin{align}
 \mat(\gG)\coloneqq
 \max\bigl\{|M|:M\subseteq\gE\textnormal{ is a matching}\bigr\}.
 \label{eqn:matching-number}
\end{align}
The column matching deficiency $m-\mat(\gG)$ counts the columns left
unmatched by a maximum matching.

For distinct components $C_1,\ldots,C_z$ containing columns, let $W_i$
be their column sets and $A_i$ their bi-adjacency matrices. For a
component consisting of a single isolated vertex in $\gV_2$, $A_i$
has zero rows and $\ker A_i=\R$. With
$W=\bigcup_iW_i$, absence of edges between components gives
\[
 \sK_W\coloneqq\{\vv:\operatorname{supp}\vv\subset W,\
                   \widetilde{\rmA}\vv=0\}
 =\bigoplus_i\ker A_i,
\]
where each summand is extended by zero. By
Lemma~\ref{lem:centered-raw-kernel}(i), its zero-sum part lies in
$\ker\mX$ and has dimension at least $\sum_i\dim\ker A_i-1$.
For trees these dimensions are determined by the matching number
in \eqref{eqn:matching-number}.

\begin{lemma}[Tree kernels and matching]
\label{lem:component-kernels}
If $C$ is a bipartite tree with $c$ columns and bi-adjacency matrix
$A_C$, then
\[
 \operatorname{rank}A_C=\mat(C),\qquad
 \dim\ker A_C=c-\mat(C).
\]
\end{lemma}
\begin{proof}
Remove a leaf $y$ and its neighbor $x$ from a
bipartite forest $F$. After transposing if necessary, $y$ is a column.
Its single nonzero entry permits column elimination at row $x$, giving
$\operatorname{rank}A_F=1+\operatorname{rank}A_{F-\{x,y\}}$.
Some maximum matching contains $xy$: if $x$ is matched elsewhere,
replace that edge by $xy$, and if $x$ is unmatched, add $xy$.
Thus the matching number satisfies the same recursion. Induction gives
the rank identity and rank--nullity gives the kernel dimension.
\end{proof}

For example, a component formed by one row and $c\geq2$ pendant columns
has matching number one, hence $c-1$ kernel directions, all with zero
coordinate sum. By Lemma~\ref{lem:fixed-tree-components} with $a=1$,
such components have diverging mean count exactly when
$b(\qr^{-2}+c\qr^2)\leq1$; the first is $c=2$, at
$b=1/(\qr^{-2}+2\qr^2)$.

\subsubsection{Singular vectors from small components}
\label{subsec:bipartite-compact-modes}

The raw matrix is block diagonal across components, preserving their
singular-value multiplicities. Centering couples the blocks through
the row and column coordinate sums.

\begin{lemma}[Component singular subspaces under centering]
\label{lem:small-component-singular-subspaces}
Let $d=p\sqrt{nm}>0$, and let $A_C$ be the unnormalized bi-adjacency
matrix of a connected component $C$ with $a_C$ row vertices and
$c_C$ column vertices.
\begin{enumerate}
\item A unit singular pair $(u,v)$ of $A_C$ at $\sigma>0$ extends
by zero to a singular pair of $\widetilde{\rmA}/\sqrt d$ at
$s=\sigma/\sqrt d$. For $\psi=2^{-1/2}(u,v)$, extended by zero,
\begin{align}
 \left\|(\scA-s)\psi\right\|^2
 &=\frac{p^2}{2d}
   \left[n(\ones^{\sT}v)^2+m(\ones^{\sT}u)^2\right]
 \leq\frac d2\left(\frac{c_C}{m}+\frac{a_C}{n}\right).
 \label{eqn:small-component-centered-quasimode}
\end{align}
The pair is preserved exactly by centering if and only if both
coordinate sums vanish. If the residual $\|(\scA-s)\psi\|$, the square
root of \eqref{eqn:small-component-centered-quasimode}, is less than $s$,
there is a positive singular value of $\mX$ within this distance of $s$.
Under the proportional and critical-sparsity assumptions, a fixed
component has residual $O_{C,\qr}(\sqrt{d/N})=o(d^{-1/2})$.

\item Fix $\sigma>0$ and any collection of distinct components
$C_1,\ldots,C_z$.
Write $A_i=A_{C_i}$ and
\[
 E_i(\sigma)=\ker(A_i^{\sT}A_i-\sigma^2\id),\qquad
 r_i(\sigma)=\dim E_i(\sigma),\qquad
 M_\sigma=\sum_{i=1}^z r_i(\sigma).
\]
On $\bigoplus_i E_i(\sigma)$ define the two-sum map
\begin{align}
 S_\sigma((v_i)_{i=1}^z)
 =\left(\sum_i\ones^{\sT}A_i v_i/\sigma,\,
         \sum_i\ones^{\sT}v_i\right).
 \label{eqn:component-singular-sum-map}
\end{align}
Its kernel gives an exact centered singular subspace at
$\sigma/\sqrt d$ of dimension
\begin{align}
 M_\sigma-\operatorname{rank}S_\sigma
 \geq (M_\sigma-2)_+.
 \label{eqn:component-centered-exact-multiplicity}
\end{align}
Thus at most two dimensions of this raw singular subspace fail the
exact cancellation condition, irrespective of the number of components.

\item If the collection consists of $z\geq1$ class-preserving copies
of one component $T$, let $r_T(\sigma)$ be the full multiplicity of
$\sigma$ as a singular value of $A_T$, and put
\[
 S_{T,\sigma}:v\in E_T(\sigma)\longmapsto
 (\ones^{\sT}A_Tv/\sigma,\ones^{\sT}v),\qquad
 \kappa_T(\sigma)=\operatorname{rank}S_{T,\sigma}.
\]
Then $0\leq\kappa_T(\sigma)\leq\min\{2,r_T(\sigma)\}$, and the
dimension in \eqref{eqn:component-centered-exact-multiplicity} is
\begin{align}
 z r_T(\sigma)-\kappa_T(\sigma)
 \geq (z-1)r_T(\sigma).
 \label{eqn:repeated-component-centered-multiplicity}
\end{align}
The constructed subspace has generators supported on at most two
copies of $T$; an orthonormal basis can involve more copies.
\end{enumerate}
\end{lemma}
\begin{proof}[Proof of Lemma~\ref{lem:small-component-singular-subspaces}]
Zero extension preserves the raw pair equations. Applying
\eqref{eqn:near-zero-centered-action} and its transpose gives
\eqref{eqn:small-component-centered-quasimode}; Cauchy--Schwarz bounds
the two coordinate sums. The spectral theorem
places an eigenvalue within the residual; a residual smaller than $s$
keeps it positive.

On $\bigoplus_iE_i(\sigma)$, the map $v_i\mapsto A_iv_i/\sigma$ is an
isometry. Hence $\ker S_\sigma$ gives orthonormal centered singular
pairs, and rank--nullity gives (ii). For identical copies the sum map
factors through the surjection $(v_i)_i\mapsto\sum_i v_i$, so its rank
is $\kappa_T(\sigma)$, proving (iii). Generators consist of vectors
from $\ker S_{T,\sigma}$ on individual copies and differences of each
remaining profile on two copies.
\end{proof}

\begin{corollary}[Exact singular vectors from fixed tree components]
\label{cor:fixed-tree-centered-modes}
Under the assumptions of Lemma~\ref{lem:fixed-tree-components}, fix a
positive singular value $\sigma$ of $A_T$ with multiplicity $r_T(\sigma)$.
Lemma~\ref{lem:small-component-singular-subspaces}(iii) gives the exact
centered multiplicity lower bound at $\sigma/\sqrt d$:
\begin{align}
 \bigl(Z_T r_T(\sigma)-\kappa_T(\sigma)\bigr)_+
 \geq (Z_T-1)_+r_T(\sigma).
 \label{eqn:fixed-tree-centered-multiplicity}
\end{align}
Let $\mu_T=\E Z_T$, as in \eqref{eqn:fixed-tree-component-mean}. If $b(a\qr^{-2}+c\qr^2)\leq1$, then for every fixed $0<\varepsilon<1$, with probability at least $1-C_T\varepsilon^{-2}(\mu_T^{-1}+p)$, the first lower bound is at least
\begin{align}
 (1-\varepsilon)r_T(\sigma)\mu_T-\kappa_T(\sigma).
 \label{eqn:fixed-tree-centered-probability-bound}
\end{align}
In particular the lower bound is $(1-o_{\mathbb P}(1))r_T(\sigma)\mu_T$, including equality at the
threshold, where $\mu_T\asymp_{T,\qr}d^e\to\infty$.
With probability tending to one, all these copies are small
components.
\end{corollary}
\begin{proof}[Proof of Corollary~\ref{cor:fixed-tree-centered-modes}]
Apply Lemma~\ref{lem:small-component-singular-subspaces} (iii) to the $Z_T$ copies and Lemma~\ref{lem:fixed-tree-components} to their count. Since $\mu_T\to\infty$, the finite loss $\kappa_T(\sigma)$ is negligible. Corollary~\ref{cor:components-size} puts every fixed-size copy outside the giant component with probability tending to one.
\end{proof}

\subsubsection{The least positive singular value created by centering}
\label{subsec:centering-least-singular-value}

Columns of $\mX$ corresponding to isolated vertices in $\gV_2$ are identical and nonzero, each equal to
$g\coloneqq-p\ones_{\gV_1}/\sqrt d$, with $\|g\|^2=d/m$.
Their contribution has scale $\rZ_2\|g\|^2=d\rZ_2/m$; whether it adds
a positive singular value depends on the range of the remaining columns.

\begin{lemma}[Centering-induced singular value]
\label{lem:centering-induced-singular-value}
Suppose $\rZ_2\geq1$. Order the columns so that
$\mX=[\,g\ones_{\sI_2}^{\sT}\ \ R\,]$, where $R$ consists of the columns
of $\mX$ outside $\sI_2$. Let $P_R$ be the orthogonal projection onto
$\ker R^{\sT}$, put $q=\|P_Rg\|$ and $\kappa=q/\|g\|\in[0,1]$, and let
$s_R$ be the least positive singular value of $R$, with $s_R=+\infty$ if
$R=0$.
\begin{enumerate}
 \item[(i)] $\operatorname{rank}\mX=\operatorname{rank}R+1$ if and only
 if $q>0$.
 \item[(ii)] If $q>0$, then
 \begin{align}
  \kappa\sqrt{\frac{d\rZ_2}{m}}
  \left(1+\frac{d\rZ_2}{ms_R^2}\right)^{-1/2}
  \leq s_{\min}^{+}(\mX)
  \leq\kappa\sqrt{\frac{d\rZ_2}{m}}.
  \label{eqn:centering-induced-scale}
 \end{align}
 \item[(iii)] Suppose that $\widetilde{\rmA}$ has an isolated row, and
 let $\widetilde{\rmA}_{\rm rest}$ consist of the columns of
 $\widetilde{\rmA}$ outside $\sI_2$. Then $q>0$ if and only if every
 vector of $\ker\widetilde{\rmA}_{\rm rest}$ has zero coordinate sum.
 In that case,
 $\operatorname{rank}\mX=\operatorname{rank}\widetilde{\rmA}+1$ and, if
 $\widetilde{\rmA}\neq0$,
 \begin{align}
  s_{\min}^{+}(\mX)
  \leq\min\left\{\sqrt{\frac{d\rZ_2}{m}},\
  \frac{\sigma_{\min}^{+}(\widetilde{\rmA})}{\sqrt d}\right\},
  \label{eqn:centering-induced-upper}
 \end{align}
 where $\sigma_{\min}^{+}$ denotes the least positive singular value.
\end{enumerate}
\end{lemma}
\begin{proof}[Proof of Lemma~\ref{lem:centering-induced-singular-value}]
The range identity
$\operatorname{Ran}\mX=\operatorname{Ran}R+\operatorname{span}\{g\}$
proves (i). For (ii), suppose first that $R\neq0$, put
$e=P_Rg/q$ and $h=g-P_Rg$, and restrict
$M=RR^{\sT}+\rZ_2gg^{\sT}$ to its range
$\operatorname{Ran}R\oplus\operatorname{span}\{e\}$.
The unit vector $e$ has Rayleigh quotient $\rZ_2q^2$.
For a unit vector $x=a+te$ in this range,
\[
 x^{\sT}Mx\geq s_R^2\|a\|^2+
 \rZ_2\bigl(\langle h,a\rangle+tq\bigr)^2.
\]
If $h\neq0$, split $a$ along $h$ and $h^\perp$. On the span of $h,e$
the lower quadratic form has matrix
\[
 \begin{bmatrix}
 s_R^2+\rZ_2\|h\|^2&\rZ_2\|h\|q\\
 \rZ_2\|h\|q&\rZ_2q^2
 \end{bmatrix}.
\]
Its least eigenvalue is at least its determinant divided by its trace,
namely $s_R^2\rZ_2q^2/(s_R^2+\rZ_2\|g\|^2)$; this bound is also at
most $s_R^2$, the bound on the remaining directions. The case $h=0$
follows from the orthogonal sum of $RR^{\sT}$ and $\rZ_2gg^{\sT}$.
If $R=0$, the only positive eigenvalue is $\rZ_2\|g\|^2$.
We have proved
\begin{align}
 \frac{s_R^2\,\rZ_2q^2}{s_R^2+\rZ_2\|g\|^2}
 \leq s_{\min}^{+}(\mX)^2\leq\rZ_2q^2,
 \label{eqn:centering-induced-two-sided}
\end{align}
with the limiting convention for $s_R=+\infty$.
Substituting $q=\kappa\|g\|$ and $\|g\|^2=d/m$ gives
\eqref{eqn:centering-induced-scale}.

For (iii), an isolated-row coordinate in $Ru=g$ forces
$\ones^{\sT}u=1$, and the other coordinates give
$\widetilde{\rmA}_{\rm rest}u=0$. Conversely, any raw remainder-kernel
vector of nonzero sum can be rescaled to solve $Ru=g$. Thus $q>0$
is equivalent to the stated zero-sum condition. The first upper bound
follows from (ii). Put $r=\operatorname{rank}\widetilde{\rmA}\geq1$.
Lemma~\ref{lem:centered-raw-kernel}(ii) gives
$s_{\min}^{+}(\mX)=\sigma_{r+1}(\mX)$. Choose an $(m-r)$-dimensional
space orthogonal to $\ones$ and the top $r-1$ raw right singular vectors.
On this space centering vanishes and the raw operator norm is at most
$\sigma_r(\widetilde{\rmA})/\sqrt d$. Min--max proves the second bound.
\end{proof}

\begin{corollary}[Least singular values near the isolation threshold in $\gV_2$]
\label{cor:least-singular-value-regimes}
Suppose that Assumptions~\ref{ass:proportional_regime} and
\ref{ass:critical_sparsity} hold with fixed $\qr$. Put $c_2=(1+\qr^4)^{-1}$, and let
$\mathcal E_\perp$ be the event that $\widetilde{\rmA}$ has an isolated
row and every vector of $\ker\widetilde{\rmA}_{\rm rest}$ has zero
coordinate sum.
\begin{enumerate}
 \item[(i)] If $b=\qr^{-2}$, then
 \[
  \liminf_{N\to\infty}\P\bigl(\lambda_m(\scA)=0\bigr)
  \geq\lim_{N\to\infty}\P(\rZ_2\geq2)=1-(1+c_2)\mathrm e^{-c_2},
 \]
 $\P(\rZ_2=1)\to c_2\mathrm e^{-c_2}$, and
 $s_{\min}^{+}(\mX)\leq\sqrt{d/m}$ on $\{\rZ_2=1\}\cap\mathcal E_\perp$.
 \item[(ii)] If $b<\qr^{-2}$, then for every $0<\varepsilon\leq c_2$
 and $\nu>0$, with probability at least
 $1-C_{b,\qr,\varepsilon,\nu}N^{-\nu}$,
 \[
  \lambda_m(\scA)=0,\qquad
  \dim\ker\mX\geq(c_2-\varepsilon)N^{1-b\qr^2}-1,
 \]
 and, on $\mathcal E_\perp$,
 $s_{\min}^{+}(\mX)\leq\sqrt{d\rZ_2/m}\leq2\sqrt d\,N^{-b\qr^2/2}$.
\end{enumerate}
\end{corollary}
\begin{proof}
Combine the Poisson limit and concentration in
Lemma~\ref{lem:bipartite-isolated-vertices} with the exact kernel
formula in Lemma~\ref{lem:bipartite-isolated-column-zero-modes}.
Below the threshold, use \eqref{eqn:bipartite-isolated-V2-count}
with $\varepsilon\wedge(c_2/2)$ to obtain
$(c_2-\varepsilon)N^{1-b\qr^2}\leq\rZ_2\leq2c_2N^{1-b\qr^2}$
and $\rZ_2\geq2$ for large $N$.
On $\mathcal E_\perp$, \eqref{eqn:centering-induced-upper} gives the
positive-value bounds; $m\geq c_2N/2$ yields
$d\rZ_2/m\leq4dN^{-b\qr^2}$ below the threshold.
\end{proof}

\begin{remark}
\label{rem:centering-induced-sharpness}
The event $\mathcal E_\perp$ holds when $\widetilde{\rmA}$ has an isolated
row and $\ker\widetilde{\rmA}_{\rm rest}$ is spanned by differences of
pendant columns of a common row.

For $q>0$, \eqref{eqn:centering-induced-scale} gives
$s_{\min}^{+}(\mX)\asymp\kappa\sqrt{d\rZ_2/m}$ if
$d\rZ_2/(ms_R^2)=O(1)$, and
$s_{\min}^{+}(\mX)=(1+o(1))\kappa\sqrt{d\rZ_2/m}$ if that ratio is
$o(1)$. Below $b=\qr^{-2}$, the concentration event in
Corollary~\ref{cor:least-singular-value-regimes} gives
$d\rZ_2/m=O(dN^{-b\qr^2})$, so $s_R\geq cd^{-1/2}$ implies
$d\rZ_2/(ms_R^2)=O(d^2N^{-b\qr^2})=o(1)$.
If $\kappa$ is also bounded away from zero, the value has order
$\sqrt{d\rZ_2/m}$.

In simulations with $\ratio=4$, $m=1500$, and
$b\in[0.22,0.45]$, every sample with $\rZ_2\geq1$ lay in
$\mathcal E_\perp$, and $s_{\min}^{+}(\mX)/\sqrt{d\rZ_2/m}$ lay in
$[0.80,0.84]$. The centering-induced singular value is then far below
the fixed-tree and fork singular values of order $d^{-1/2}$ in
Subsection~\ref{subsec:bipartite-compact-modes} and
Appendix~\ref{subsec:bipartite-component-structure}.
\end{remark}

\subsection{Star-tuning forks on the giant component}
\label{subsec:bipartite-component-structure}

We prove the fork results summarized in
Section~\ref{subsec:intro-zero-modes}. All threshold statements use
fixed $\qr$.

\subsubsection{Exact singular vectors and kernel vectors}
\label{subsec:bipartite-pendant-configurations}

\begin{proof}[Proof of Lemma~\ref{lem:bipartite-star-tuning-fork-modes}]
For $F=0$, the hubs have the sole neighbor $o$, so their contributions
in $\rvw_0=(\ones_{x_1}-\ones_{x_2})/\sqrt2$ cancel and
$\rmA\rvw_0=0$. Its coordinate sum vanishes, so centering preserves
this identity. The vector lies in the right kernel when $o\in\gV_1$
and in the left kernel otherwise.

For $F\geq1$, extend the profiles in
\eqref{eqn:star-tuning-fork-singular-profiles} by zero. The degree
conditions give, at each leaf and hub respectively,
\[
 \rmA\vh_F=\sqrt F\,\vll_F,\qquad
 \rmA\vll_F=\sqrt F\,\vh_F;
\]
the hub contributions at the base cancel. Both profiles have zero
coordinate sum in their vertex class, so centering and normalization
give the stated singular pair at $s_F=\sqrt{F/d}$ and the orthonormal
signed eigenvectors
\begin{align}
 \rvw_\sigma&=
 \frac12(\ones_{x_1}-\ones_{x_2})+
 \frac{\sigma}{2\sqrt F}(\ones_{R_1}-\ones_{R_2}),\notag\\
 \scA\rvw_\sigma&=\sigma\sqrt{F/d}\,\rvw_\sigma,
 \qquad \sigma\in\{-1,1\}.
 \label{eqn:star-tuning-fork-eigenvalue}
\end{align}
A zero-sum vector on $R_i$ is annihilated by both the raw adjacency
and the mean matrix. This single constraint gives dimension $F-1$.
The two leaf-kernel spaces have disjoint supports and are orthogonal
to $\vll_F$. Their vertex classes give the kernel orientations in
item~(iii), and the displayed profiles give the stated support sizes
and supremum norms.

If the base lies in the giant, the entire fork lies there. Every
constructed eigenvector therefore satisfies
$\chi_{\rm giant}\rvw=\rvw$, and its eigenvalue is preserved by
$\chi_{\rm giant}\scA\chi_{\rm giant}$.
\end{proof}

\subsubsection{Counts, thresholds, and occurrence in the giant}

\begin{lemma}[Fork counts and moments]
\label{lem:bipartite-star-tuning-fork-counts}
Suppose that Assumptions~\ref{ass:proportional_regime} and
\ref{ass:critical_sparsity} hold, and fix $F\in\N_0$.
Let $N_{i,F}$ count $F$-forks with base in $\gV_i$, and put
\begin{align}
 a_1&=\qr^2+F\qr^{-2},& a_2&=\qr^{-2}+F\qr^2,\\
 t_i&=d^{F+1}\mathrm e^{-da_i},& \mu_i&=\E N_{i,F}.
\end{align}
Then
\begin{subequations}
\begin{align}
 \mu_1&=
 \left(\frac{\qr^{4F}b^{2F+2}}
 {2(1+\qr^4)(F!)^2}+o(1)\right)
 (\log N)^{2F+2}N^{1-2ba_1},
 \label{eqn:star-tuning-fork-expectations-V1}\\
 \mu_2&=
 \left(\frac{\qr^{4-4F}b^{2F+2}}
 {2(1+\qr^4)(F!)^2}+o(1)\right)
 (\log N)^{2F+2}N^{1-2ba_2}.
 \label{eqn:star-tuning-fork-expectations-V2}
\end{align}
\end{subequations}
For $k\in\{2,3,4\}$, let $\rT_{i,k}$ count configurations with a
common base in $\gV_i$ and $k$ unordered arms satisfying the fork
degree conditions, with pairwise disjoint leaf sets. Thus
$\rT_{i,2}=N_{i,F}$, and $\E\rT_{i,k}\asymp_{\qr,F}Nt_i^k$.
Moreover,
\begin{align}
 \Var(N_{i,F})\leq(1+o(1))\mu_i,
 \qquad i\in\{1,2\}.
 \label{eqn:star-tuning-fork-variance}
\end{align}
If $b>\bcritF{F,i}(\qr)$, then $\P(N_{i,F}>0)\sim\mu_i$, with the
leading constants in the displayed means.
\end{lemma}
\begin{proof}[Proof of Lemma~\ref{lem:bipartite-star-tuning-fork-counts}]
For orientation $1$, choose the base, the unordered hubs, and their
leaf sets in $n\binom mk(n-1)_{kF}/(F!)^k$ ways, with $(n-1)_0=1$.
Among the $kn+kFm-k^2F$ pairs incident to hubs or leaves,
$k(F+1)$ must be present and the rest absent. Hence
\begin{align}
 \E\rT_{1,k}
 &=n\binom mk\frac{(n-1)_{kF}}{(F!)^k}
 p^{k(F+1)}(1-p)^{kn+kFm-k(k+1)F-k}\notag\\
 &=\left(\frac{\qr^{4+2k(F-1)}}
 {(1+\qr^4)k!(F!)^k}+o(1)\right)
 Nd^{k(F+1)}\mathrm e^{-kda_1}
 \asymp_{\qr,F}Nt_1^k.
 \label{eqn:star-tuning-fork-multiple-arm-count}
\end{align}
The expansion follows from the fixed-$F$ factorial estimates in
Lemma~\ref{lem:stirling-formula} and $\log(1-p)=-p+O(p^2)$.
Taking $k=2$ and $d=b\log N$ gives
\eqref{eqn:star-tuning-fork-expectations-V1}. Exchanging the vertex
classes gives the second orientation, with numerator
$\qr^{2k(1-F)}$ in the leading constant, proving
\eqref{eqn:star-tuning-fork-expectations-V2} and its arm-count estimates.

Two compatible intersecting forks of one orientation have the same
base: a base has degree at least two and cannot be a leaf; a shared
leaf forces a shared hub; and a shared hub determines its base and
leaf set. For $F=0$ its neighbor is unique, while for $F\geq1$ the
base is its unique neighbor of degree at least two. Thus such a pair
coincides, shares one complete arm, or shares only its base.
Three- and four-arm configurations each give six ordered pairs of
the latter two kinds. For vertex-disjoint forks, the forbidden-edge
sets overlap in $O_F(1)$ cross edges, multiplying the product of
their probabilities by at most $1+O_F(p)$. Consequently,
\begin{align}
 \E N_{i,F}^2
 \leq(1+O_F(p))\mu_i^2+\mu_i
       +6\E\rT_{i,3}+6\E\rT_{i,4}.
 \label{eqn:star-tuning-fork-second-moment-upper-bound}
\end{align}
Since $b>0$ is fixed, $t_i\to0$, and the arm counts give
\begin{align}
 \frac{\E\rT_{i,3}}{\mu_i}&=O_{\qr,F}(t_i)=o(1),
 &\frac{\E\rT_{i,4}}{\mu_i}&=O_{\qr,F}(t_i^2)=o(1),\notag\\
 p\mu_i&=O_{\qr,F}\bigl(d^{2F+3}\mathrm e^{-2da_i}\bigr)=o(1).
 \label{eqn:star-tuning-fork-overlap-ratios}
\end{align}
Subtracting $\mu_i^2$ proves the variance bound. This is the
bipartite counterpart of the second-moment argument in
\cite[Lemma~A.19]{alt2021delocalization}.
For $b>\bcritF{F,i}(\qr)$, the mean tends to zero. Markov gives
$\P(N_{i,F}>0)\leq\mu_i$, while Cauchy--Schwarz gives
\[
 \P(N_{i,F}>0)\geq\frac{\mu_i^2}{\E N_{i,F}^2}
 \geq\frac{\mu_i}{\mu_i+1+o(1)}=(1+o(1))\mu_i.
\]
Thus $\P(N_{i,F}>0)\sim\mu_i$, with the leading constants in
\eqref{eqn:star-tuning-fork-expectations-V1}--
\eqref{eqn:star-tuning-fork-expectations-V2}.
\end{proof}

\begin{proof}[Proof of Lemma~\ref{lem:bipartite-star-tuning-fork-threshold}]
The notation in Lemma~\ref{lem:bipartite-star-tuning-fork-counts}
satisfies $a_i=[2\bcritF{F,i}(\qr)]^{-1}$, so
\begin{align}
 \mu_i\asymp_{b,\qr,F}Nt_i^2
 \asymp_{b,\qr,F}
 (\log N)^{2F+2}N^{1-b/\bcritF{F,i}(\qr)}.
 \label{eqn:star-tuning-fork-mean-threshold-scale}
\end{align}
If $b\leq\bcritF{F,i}(\qr)$, then $\mu_i\to\infty$, including at
equality, where $\mu_i\asymp_{b,\qr,F}d^{2F+2}$. The variance bound
and Chebyshev's inequality give
\begin{align}
 N_{i,F}/\mu_i&\longrightarrow1
 \quad\text{in probability},
 \label{eqn:star-tuning-fork-count-concentration}\\
 \P(N_{i,F}=0)&\leq\frac{1+o(1)}{\mu_i}
 \lesssim_{b,\qr,F}
 N^{b/\bcritF{F,i}(\qr)-1}(\log N)^{-2F-2}.
 \label{eqn:star-tuning-fork-emergence-probability-rate}
\end{align}
The failure bound is polynomial below the threshold and
$O_{b,\qr,F}((\log N)^{-2F-2})$ at equality.
If $b>\bcritF{F,i}(\qr)$, then $\mu_i\to0$ and Markov's inequality
gives
\begin{align}
 \P(N_{i,F}>0)\leq\mu_i
 \lesssim_{b,\qr,F}
 N^{1-b/\bcritF{F,i}(\qr)}(\log N)^{2F+2}.
 \label{eqn:star-tuning-fork-absence-probability-rate}
\end{align}

For fixed $F$, these bounds and a union bound give the threshold
$\max_i\bcritF{F,i}(\qr)$ for either orientation. For $F\geq1$ and
$\qr\geq1$, $a_2-a_1=(F-1)(\qr^2-\qr^{-2})\geq0$.
Hence the larger threshold is $\bcritF{F,1}(\qr)$; it decreases with
$F$ and is largest at $F=1$, where the orientations agree, giving
\eqref{eqn:nonzero-star-tuning-fork-threshold}. For absence of all
positive forks above this threshold, the bounded-$F$ conclusion of
Corollary~\ref{cor:bipartite-star-tuning-fork-zero-scale} below reduces
the union over $F$ to a finite one. The thresholds for fork-induced
column rank deficiency follow from the kernel orientations in
Lemma~\ref{lem:bipartite-star-tuning-fork-modes}.
\end{proof}

\begin{corollary}[Occurrence of fixed forks in the giant]
\label{cor:bipartite-star-tuning-fork-giant-emergence}
Suppose that Assumptions~\ref{ass:proportional_regime} and
\ref{ass:critical_sparsity} hold with fixed $b>0$, $\qr$, and
$F\in\N_0$. Let $N_{i,F}^{\rm giant}$ count $F$-forks with base in
$\gV_i\cap\gG_{\rm giant}$. If $b\leq\bcritF{F,i}(\qr)$, then
\[
 N_{i,F}^{\rm giant}=(1-o_{\mathbb P}(1))\E N_{i,F}.
\]
In particular, a fork of this orientation occurs in the giant with
probability $1-o(1)$, including equality at the threshold.
\end{corollary}
\begin{proof}[Proof of
Corollary~\ref{cor:bipartite-star-tuning-fork-giant-emergence}]
On the common event of Corollary~\ref{cor:components-size}, every
small component has at most a fixed $K=K(b,\qr)$ vertices.
Conditional on a prescribed labelled fork, the base's additional
neighbors have law $\Bin(m-2,p)$ in orientation $1$ and $\Bin(n-2,p)$
in orientation $2$. These edges are unconstrained by the fork, and
their means are asymptotic to $d\qr^{-2}$ and $d\qr^2$. The probability
$\eta_{i,N}$ of fewer than $K$ such neighbors therefore tends to zero.

Let $B_{i,F}$ count forks failing this degree condition. Then
$\E B_{i,F}=\eta_{i,N}\mu_i=o(\mu_i)$, so Markov's inequality gives
$B_{i,F}/\mu_i\to0$ in probability. On the common component event,
a base with at least $K$ additional neighbors lies in the giant,
giving $0\leq N_{i,F}-N_{i,F}^{\rm giant}\leq B_{i,F}$.
Combine this with \eqref{eqn:star-tuning-fork-count-concentration}.
The quantitative rates in Lemma~\ref{lem:bipartite-star-tuning-fork-threshold}
concern the global count.
\end{proof}

\subsubsection{Uniform confinement}

\begin{lemma}[Uniform confinement of star-tuning-fork singular values]
\label{lem:bipartite-tuning-forks}
Suppose that $d^2/\log N\to\infty$, and let $L_x$ count the pendant
neighbors of $x$. For some $c_{\qr}>0$ and all sufficiently large $N$,
\begin{align}
 \P\left(\max_{x\in\gV}L_x\geq\sqrt{\log N}\right)
 \leq\exp(-c_{\qr}d\sqrt{\log N}).
 \label{eqn:pendant-neighbour-tail}
\end{align}
On the complementary event, simultaneously for every $F\in\N_0$
and every $F$-fork rooted in $\gG$,
$|R_1|=|R_2|=L_{x_1}=L_{x_2}=F<\sqrt{\log N}$.
For $F\geq1$, its constructed nonzero eigenvalues satisfy
\begin{align}
 |\lambda|=\sqrt{F/d}
 <\left(\frac{\log N}{d^2}\right)^{1/4}=o(1).
 \label{eqn:star-tuning-fork-eigenvalue-tail}
\end{align}
The constructed kernel vectors have eigenvalue zero and satisfy the
same bound.
\end{lemma}
\begin{proof}[Proof of Lemma~\ref{lem:bipartite-tuning-forks}]
For fixed $x\in\gV_1$, the pendant-neighbor events depend on disjoint
edge sets and each has probability $p(1-p)^{n-1}$. Thus
\[
 L_x\sim
 \begin{cases}
 \Bin(m,p(1-p)^{n-1}),&x\in\gV_1,\\
 \Bin(n,p(1-p)^{m-1}),&x\in\gV_2,
 \end{cases}
\]
and
\begin{align}
 \mu_*\coloneqq\max_{x\in\gV}\E L_x
 &=\max\{mp(1-p)^{n-1},np(1-p)^{m-1}\}\notag\\
 &\leq C_{\qr}d\mathrm e^{-c_{\qr}d}.
 \label{eqn:pendant-neighbour-mean}
\end{align}
Chernoff's inequality and a union bound give
\begin{align}
 \P(\max_{x\in\gV}L_x\geq r)
 \leq N\left(\frac{\mathrm e\mu_*}{r}\right)^r,
 \qquad r>0,\quad r\geq\mu_*.
 \label{eqn:pendant-neighbour-uniform-tail}
\end{align}
If $\mu_*=0$, this probability is zero. Otherwise set
$r_N=\sqrt{\log N}=o(d)$. Since
$\log(C_{\qr}d/r_N)=o(d)$ and $\log N=o(dr_N)$, the logarithm of
the right-hand side is at most $-c_{\qr}dr_N$, after decreasing
$c_{\qr}$. This proves \eqref{eqn:pendant-neighbour-tail}.
The fork degree conditions give $L_{x_1}=L_{x_2}=F$; apply
\eqref{eqn:star-tuning-fork-eigenvalue} for the spectral conclusion.
\end{proof}

\begin{corollary}[Star-tuning-fork eigenspaces on the giant component]
\label{cor:bipartite-star-tuning-fork-zero-scale}
Suppose that Assumptions~\ref{ass:proportional_regime} and
\ref{ass:critical_sparsity} hold with fixed $\qr$.
For some $c_{b,\qr}>0$ and all sufficiently large $N$, with probability
at least $1-\exp[-c_{b,\qr}(\log N)^{3/2}]$, every $F$-fork rooted in
$\gG$ has $F<\sqrt{\log N}$ and, for $F\geq1$,
\begin{align}
 \lambda_\pm=\pm s_F,\qquad
 s_F=\sqrt{\frac{F}{b\log N}}
 <b^{-1/2}(\log N)^{-1/4}=o(1).
 \label{eqn:critical-star-tuning-fork-zero-scale}
\end{align}
For forks based in the giant, all constructed eigenvectors retain
their eigenvalues under $\chi_{\rm giant}\scA\chi_{\rm giant}$.
When $\qr>1$, none of these eigenvalues lies in a fixed compact
interval $\sK\subset(0,\ledge)$ for sufficiently large $N$ on this event.

For every fixed $\nu>0$, choose an integer $r=r_{b,\qr,\nu}$ with
$rb\qr^{-2}>\nu+1$. For some $C_{b,\qr,\nu}>0$ and all sufficiently
large $N$, with probability at least $1-C_{b,\qr,\nu}N^{-\nu}$,
$\max_xL_x<r$. On this event every $F$-fork rooted in $\gG$ has $F<r$
and, for $F\geq1$,
\begin{align}
 |\lambda_\pm|=s_F\leq\sqrt{\frac{r-1}{d}}
 =O_{b,\qr,\nu}(d^{-1/2}).
 \label{eqn:critical-star-tuning-fork-bounded-size}
\end{align}
\end{corollary}
\begin{proof}[Proof of
Corollary~\ref{cor:bipartite-star-tuning-fork-zero-scale}]
The first probability bound is Lemma~\ref{lem:bipartite-tuning-forks}
with $d=b\log N$. Compression preserves the constructed vectors by
the singular-vector proof above, and $\inf\sK>0$ gives the exclusion.
For the stronger bound, the exact mean in
\eqref{eqn:pendant-neighbour-mean} gives
$\mu_*\leq C_{\qr}d\mathrm e^{-d\qr^{-2}}
=C_{\qr}dN^{-b\qr^{-2}}$. Apply
\eqref{eqn:pendant-neighbour-uniform-tail} at the fixed integer $r$:
\[
 \P(\max_xL_x\geq r)
 \leq C_{b,\qr,\nu}(\log N)^rN^{1-rb\qr^{-2}}
 \leq C_{b,\qr,\nu}N^{-\nu}.
\]
The strict inequality $rb\qr^{-2}>\nu+1$ absorbs the logarithmic
factor, and $F\leq r-1$ gives
\eqref{eqn:critical-star-tuning-fork-bounded-size}.
\end{proof}

\subsection{Reduction to the giant component}
\label{sec:giant-component-reduction-and-spectrum-near-zero}

The raw adjacency decomposes across the giant and its complement.
Keeping the small-component block isolates the effect of centering;
discarding it gives the comparison on fixed positive spectral windows.
We first record the deterministic identities and eigenpair estimates.

\begin{lemma}[Component decomposition and eigenpair comparison]
\label{lem:giant-component-deterministic-comparison}
Let $d=p\sqrt{nm}>0$, $P=\chi_{\rm{giant}}$, $Q=\id-P$, and
$\mathcal S_i=\gV_i\setminus(\gG_{\rm{giant}}\cap\gV_i)$. Put
\[
 \rmH_{\rm g}=P\scA P,\qquad
 \rmA_{\rm s}=Q\rmA Q/\sqrt d,\qquad
 \rmE_{\rm g}=-\bigl(\E\rmA-P(\E\rmA)P\bigr)/\sqrt d.
\]
Write $\beta_N=\|\rmA_{\rm s}\|$ and $\delta_N=\|\rmE_{\rm g}\|$.
Then the following hold for every realization.
\begin{enumerate}[label=\textnormal{(\roman*)}]
\item The block decomposition is
\begin{align}
 \scA&=\rmH_{\rm g}+\rmA_{\rm s}+\rmE_{\rm g},
 &P\rmE_{\rm g}P&=0,
 &\operatorname{rank}\rmE_{\rm g}&\leq4.
 \label{eqn:giant-component-block-decomposition}
\end{align}
Consequently,
\[
 \|\scA-\rmH_{\rm g}\|\leq\beta_N+\delta_N,
 \qquad
 \delta_N\leq\frac p{\sqrt d}
 \left(\sqrt{|\mathcal S_1|m}+\sqrt{n|\mathcal S_2|}\right).
\]
\item For every unit eigenpair $\scA\rvw=\lambda\rvw$,
\[
 \|\rmH_{\rm g}(P\rvw)-\lambda P\rvw\|
 \leq\delta_N\|Q\rvw\|.
\]
If $\lambda\neq0$, the coarse bounds are
\begin{align}
 \|Q\rvw\|&\leq\frac{\beta_N+\delta_N}{|\lambda|},
 \label{eqn:eigenvector-mass-off-giant}\\
 \|\rmH_{\rm g}(P\rvw)-\lambda P\rvw\|
 &\leq\delta_N.
 \label{eqn:giant-component-eigenvector-defect}
\end{align}
For $Q\neq0$, define
\[
 \gamma_\lambda=
 \operatorname{dist}\bigl(\lambda,
 \operatorname{Spec}(\rmA_{\rm s}|_{\operatorname{Ran}Q})\bigr).
\]
If $\gamma_\lambda>0$, then
\begin{align}
 \|Q\rvw\|&\leq\frac{\delta_N}{\gamma_\lambda},
 \notag\\
 \|\rmH_{\rm g}(P\rvw)-\lambda P\rvw\|
 &\leq\frac{\delta_N^2}{\gamma_\lambda}.
 \label{eqn:giant-component-sharp-eigenpair-comparison}
\end{align}
When $Q=0$, take $\gamma_\lambda=+\infty$; both quantities vanish.
If $|\lambda|>\beta_N$, the denominator $\gamma_\lambda$ may be
replaced by $|\lambda|-\beta_N$.
\end{enumerate}
\end{lemma}
\begin{proof}
There are no raw edges between distinct components, giving the
decomposition. The rectangular block of
$\E\rmA-P(\E\rmA)P$ is the difference of two rank-one matrices,
so its Hermitian linearization has rank at most four. The same block is
\[
 p\left(\ones_{\mathcal S_1}\ones_{\gV_2}^{\sT}
   +\ones_{\gG_{\rm{giant}}\cap\gV_1}
      \ones_{\mathcal S_2}^{\sT}\right),
\]
whose norm gives the bound on $\delta_N$. The triangle inequality gives
the comparison norm.

Projecting the eigenvector equation gives
\[
 (\lambda-\rmA_{\rm s})Q\rvw=Q\rmE_{\rm g}\rvw,\qquad
 (\rmH_{\rm g}-\lambda)P\rvw=-P\rmE_{\rm g}Q\rvw.
\]
The second identity bounds the defect. Since $Q\rmH_{\rm g}=0$,
$\lambda Q\rvw=Q(\scA-\rmH_{\rm g})\rvw$ gives the coarse mass
bound. When $\gamma_\lambda>0$, inversion of the first identity on
$\operatorname{Ran}Q$ gives the sharp mass bound, and substitution
into the defect estimate gives its sharp counterpart. Finally,
$\operatorname{Spec}(\rmA_{\rm s}|_{\operatorname{Ran}Q})
\subset[-\beta_N,\beta_N]$, so
$\gamma_\lambda\geq|\lambda|-\beta_N$ whenever
$|\lambda|>\beta_N$.
\end{proof}

\begin{corollary}[Reduction to the giant component]
\label{cor:bipartite-giant-component-reduction}
Suppose that $d=b\log N+c_N$, where $b>0$ is fixed and
$0\leq c_N\leq\log N$. Use the notation of
Lemma~\ref{lem:giant-component-deterministic-comparison}.
There are constants $c,\nu,C>0$, depending only on $b$ and $\qr$, and
an event $\mathcal E_N$ of probability at least $1-CN^{-\nu}$ such that
\begin{enumerate}[label=\textnormal{(\roman*)}]
\item On $\mathcal E_N$, the two error terms satisfy
\begin{align}
 \beta_N&\leq C\frac{\sqrt{\log N}}d,
 &\delta_N&\leq C\sqrt d\,\mathrm e^{-cd}.
 \label{eqn:centering-off-giant-bound}
\end{align}
The centered giant block, extended by zero, therefore satisfies
\begin{align}
 \|\scA-P\scA P\|
 &\leq\beta_N+\delta_N
 \leq C\left(\frac{\sqrt{\log N}}d+\sqrt d\,\mathrm e^{-cd}\right)
 \leq C\error_d,
 \label{eqn:giant-component-centered-comparison}
\end{align}
where $\error_d$ is defined in \eqref{eqn:error}. In particular,
every unit eigenpair with $\lambda\neq0$ satisfies
$\|Q\rvw\|\leq C\error_d/|\lambda|$ and
$\|\rmH_{\rm g}(P\rvw)-\lambda P\rvw\|\leq C\error_d$.
\item On $\mathcal E_N$, for every fixed $a_0>0$ and all sufficiently
large $N$, the condition $\min_{y\in\gV_2}\alpha_y\geq a_0$
implies that every column belongs to the giant and its complement
consists only of isolated rows. Thus compression to the giant retains
all $m$ column coordinates, and
\begin{align}
 \|\scA-\rmH_{\rm g}\|
 =p\sqrt{m\rZ_1/d}
 =\sqrt{\frac{d\rZ_1}{n}}
 \leq C\sqrt d\,\mathrm e^{-cd}.
 \label{eqn:column-preserving-giant-reduction}
\end{align}
\item If $\qr$ is fixed and $b>\qr^{-2}$, there are
$a_0\in(0,\qr^2)$ and $\nu_0>0$, depending only on $b$ and $\qr$, such
that, with probability at least $1-CN^{-\nu_0}$, $\mathcal E_N$ and
\begin{align}
 \min_{y\in\gV_2}\alpha_y\geq a_0
 \label{eqn:giant-column-degree-floor}
\end{align}
hold simultaneously. The constant $a_0$ may be chosen arbitrarily small.
\end{enumerate}
\end{corollary}
\begin{proof}[Proof of
Corollary~\ref{cor:bipartite-giant-component-reduction}]
Take $\mathcal E_N$ to be the event in
Corollary~\ref{cor:components-size}. Since $d\geq b\log N$, its failure
probability is $O_{b,\qr}(N^{-\nu})$ for some $\nu>0$. On this event,
$|\mathcal S_1|+|\mathcal S_2|\leq N\mathrm e^{-c_{\qr}d}$ and
\eqref{eqn:bipartite-small-component-norm} bounds $\beta_N$.
The deterministic estimate for $\delta_N$ gives
\[
 \delta_N\leq C_{\qr}\sqrt d
 \left(\sqrt{\frac{|\mathcal S_1|}{n}}
       +\sqrt{\frac{|\mathcal S_2|}{m}}\right),
\]
which yields \eqref{eqn:centering-off-giant-bound}. The comparison and
eigenpair bounds in (i) now follow from
Lemma~\ref{lem:giant-component-deterministic-comparison}, using
$d\asymp_b\log N$ for the comparison with $\error_d$.

For (ii), a small component containing a column would have at least
$a_0d+1$ vertices, contradicting its bounded size on $\mathcal E_N$ for
large $N$. Hence every column belongs to the giant and its complement
consists of the $\rZ_1$ isolated rows. Its raw adjacency is zero, and
the rectangular block of $\scA-\rmH_{\rm g}$ is
$-p\ones_{\mathcal S_1}\ones_{\gV_2}^{\sT}/\sqrt d$.
Its norm gives \eqref{eqn:column-preserving-giant-reduction}.

For (iii), $I_{\qr}(a)\to\qr^2$ as $a\downarrow0$, so choose $a_0>0$
arbitrarily small with $bI_{\qr}(a_0)>1$.
The sparse tail mean bound \eqref{eqn:sparse-degree-tail-mean-bound}
and the presence bound \eqref{eqn:degree-window-presence-absence} give
\[
 \P\bigl(\min_{y\in\gV_2}\alpha_y<a_0\bigr)
 \leq\mu_2([0,a_0])
 \leq C_{b,\qr}N\mathrm e^{-dI_{\qr}(a_0)}
 \leq C_{b,\qr}N^{1-bI_{\qr}(a_0)}.
\]
Intersecting with $\mathcal E_N$ proves (iii), with some $\nu_0>0$.
\end{proof}

At critical sparsity, the common component event restricts the small
trees to finitely many types. Their positive raw normalized singular
values are positive constants times $d^{-1/2}$, whereas the centering
coupling has the smaller scale $\delta_N=O_{b,\qr}(\sqrt d\,\mathrm e^{-cd})$.
The comparison norm is
$O_{b,\qr}((\log N)^{-1/2})=o(1)$, while
$\error_d=O_b(\log\log N/\sqrt{\log N})=o(1)$.
The norm bound permits reduction on fixed positive windows.
On shrinking windows, the sharp eigenpair comparison applies when the
distance from the raw small-component spectrum is positive; a distance
large relative to $\delta_N$ makes the mass outside the giant small.
Near zero, retain the raw small-component block and the centering term
separately, as in Subsections~\ref{subsec:bipartite-compact-modes}
and~\ref{subsec:centering-least-singular-value}.

%% file: app-delocalized.tex
\section{Deferred Proofs in Section~\ref{sec:delocalized_phase}}
\label{sec:delocalized-deferred-proofs}

We prove the estimates deferred from Section~\ref{sec:delocalized_phase}.
Throughout, the hypotheses of
Theorem~\ref{thm:bipartite-sparse-local-law} and the definitions in
Subsection~\ref{sec:local-law-typical-vertices} are in force, and
$z\in\specdomain$ is deterministic. Polynomial union bounds are absorbed
by increasing the failure exponent in the corresponding input.

\subsection{Resolvent identities and concentration}
\label{sec:delocalized-local-law-inputs}

\begin{lemma}[Resolvent identities, {\cite{benaych2016lectures}}]
\label{lem:resolvent-identities}
For $\setT\subset[N]$, the following identities hold:
\begin{enumerate}[label=\textup{(\arabic*)}]
    \item For $x\in[N]$,
    \begin{align}
        \frac1{\eres_{xx}}=\emM_{xx}-z
        -\sum_{a,b}^{(x)}\emM_{xa}\eres_{ab}^{(x)}\emM_{bx}.
        \label{eqn:schur-complement-diagonal}
    \end{align}
    \item For $x\notin\setT$, Ward's identity is
    \begin{align}
        \sum_y^{(\setT)}|\eres_{xy}^{(\setT)}|^2
        =\frac{\Im\eres_{xx}^{(\setT)}}{\Im z}.
        \label{eqn:ward-identity}
    \end{align}
    \item For distinct $x,y\notin\setT$,
    \begin{align}
        \eres_{xy}^{(\setT)}
        =-\eres_{xx}^{(\setT)}
          \sum_b^{(\setT x)}\emM_{xb}\eres_{by}^{(\setT x)}
        =-\eres_{yy}^{(\setT)}
          \sum_a^{(\setT y)}\eres_{xa}^{(\setT y)}\emM_{ay}.
        \label{eqn:GxyT-expansion}
    \end{align}
    \item For $u\notin\setT$ and $x,y\notin\setT u$,
    \begin{align}
        \eres_{xy}^{(\setT u)}=\eres_{xy}^{(\setT)}
        -\frac{\eres_{xu}^{(\setT)}\eres_{uy}^{(\setT)}}
        {\eres_{uu}^{(\setT)}}.
        \label{eqn:GxyTu-equality}
    \end{align}
    \item If $\Im z,\Im w\geq N^{-1}$, then
    \begin{align}
        |\eres_{xy}(z)-\eres_{xy}(w)|\leq N^2|z-w|,
        \qquad x,y\in[N].
        \label{eqn:resolvent-Lipschitz}
    \end{align}
\end{enumerate}
\end{lemma}

\begin{proof}
Apply the block inverse formula in Lemma~\ref{lem:shur-complement},
isolating one index, to obtain \textup{(1)}, \textup{(3)}, and
\textup{(4)}. Since $\mM^{(\setT)}$ is Hermitian,
$\Im\res^{(\setT)}=(\Im z)\res^{(\setT)}(\res^{(\setT)})^*$;
its diagonal entries give \textup{(2)}. Finally,
$\res(z)-\res(w)=(z-w)\res(z)\res(w)$ and
$\|\res(z)\|\leq(\Im z)^{-1}$ give \textup{(5)}.
\end{proof}

\begin{proposition}[Sparse large deviations;
{\cite[Propositions 3.1, 3.2 and 3.5]{he2019local}}]
\label{prop:sparse-large-deviation}
Let $q\in\N$ be even and $1\leq d\leq N$. Suppose that
$(\rX_j)_{j=1}^N$ are independent real random variables satisfying
\begin{align}
    \E\rX_j=0,\qquad
    \E|\rX_j|^k\leq\frac1{Nd^{(k-2)/2}},\qquad 2\leq k\leq q.
\end{align}
Let $(\rA_j)$ and $(\rB_{ij})$ be deterministic complex coefficients.
Suppose that $\Gamma,\psi>0$ satisfy
\begin{align}
    \left(\frac1N\sum_j|\rA_j|^2\right)^{1/2}&\leq\Gamma,
    &\frac1{\sqrt d}\max_j|\rA_j|&\leq\psi,\notag\\
    \left(\max_j\frac1N\sum_i|\rB_{ij}|^2\right)^{1/2}
    \vee\left(\max_i\frac1N\sum_j|\rB_{ij}|^2\right)^{1/2}
    &\leq\Gamma,
    &\frac1d\max_{i,j}|\rB_{ij}|&\leq\psi.
\end{align}
Writing $\|\rX\|_q=(\E|\rX|^q)^{1/q}$, one has
\begin{align}
    \left\|\sum_i\rA_i\rX_i\right\|_q
    &\leq\left(\frac{2q}{1+2(\log(\psi/\Gamma))_+}\vee2\right)
    (\Gamma\vee\psi),
    \label{eqn:sparse-large-deviation-mean}\\
    \left\|\sum_i\rA_i(|\rX_i|^2-\E|\rX_i|^2)\right\|_q
    &\leq2\left(1+\frac{2d}{N}\right)\max_i|\rA_i|
    \left(\frac qd\vee\sqrt{\frac qd}\right),
    \label{eqn:sparse-large-deviation-variance}\\
    \left\|\sum_{i\neq j}\rB_{ij}\rX_i\rX_j\right\|_q
    &\leq\left(\frac{4q}{1+(\log(\psi/\Gamma))_+}\vee4\right)^2
    (\Gamma\vee\psi).
    \label{eqn:sparse-large-deviation-covariance}
\end{align}
\end{proposition}

\begin{corollary}\label{cor:sparse-Hanson-Wright-Estimates}
Under the assumptions of Proposition~\ref{prop:sparse-large-deviation},
suppose that $\psi/\Gamma\geq N^{\kappa/4}$. For each $\nu>0$, if
$q\geq C_{\kappa,\nu}\log N$, then with probability at least
$1-C_{\kappa,\nu}N^{-\nu}$,
\begin{align}
    \left|\sum_j\rA_j\rX_j\right|\leq C_{\kappa,\nu}\psi,
    \qquad
    \left|\sum_{i\neq j}\rB_{ij}\rX_i\rX_j\right|
    \leq C_{\kappa,\nu}\psi.
    \label{eqn:sparse-Hanson-Wright-Estimates}
\end{align}
\end{corollary}

\begin{proof}
Apply Markov's inequality to
\eqref{eqn:sparse-large-deviation-mean} and
\eqref{eqn:sparse-large-deviation-covariance} at an even moment of order
$\log N$. The logarithmic denominators keep both moment bounds of order
$\psi$; increasing the constant in the tail threshold gives the desired
exponent. See also \cite[Corollary A.21]{alt2021delocalization}.
\end{proof}

We apply these estimates conditionally on a minor, so the remaining row
entries are independent of the coefficients. The moment bounds for
$\emH_{xy}$ follow from Definition~\ref{def:sparse-random-rectangular-matrix}
after rescaling by a fixed constant.

\subsection{Minor-resolvent comparison}
\label{sec:delocalized-minor-resolvent-proof}

\begin{proof}[Proof of Lemma~\ref{lem:res-max-difference-bound}]
Enumerate $\setT=\{u_1,\ldots,u_j\}$ and put
$\setT_\ell=\{u_1,\ldots,u_\ell\}$, with $\setT_0=\emptyset$.
Set
\begin{align}
    R_\ell\coloneqq1\vee
    \max_{x,y\notin\setT_\ell}|\eres_{xy}^{(\setT_\ell)}|.
\end{align}
On $\{\theta_1^{(\emptyset)}=1\}$, we have $R_0\leq\resmax$.
We show that deleting one index costs a relative error $O(\resmax^2/d)$.

\paragraph{One deletion.}
For this step, write $\setS=\setT_\ell$, $u=u_{\ell+1}$, and
$\eta=\Im z$. For $x\notin\setS u$, define the noise and signal sums
\begin{align}
    A_x\coloneqq\sum_a^{(\setS u)}\eres_{xa}^{(\setS u)}\emH_{au},
    \qquad
    D_x\coloneqq\frac{2\signal}{N}
    \sum_a^{(\setS u)}\indi{\evsigma_a\neq\evsigma_u}
    \eres_{xa}^{(\setS u)}.
\end{align}
By \eqref{eqn:GxyT-expansion} and \eqref{eqn:GxyTu-equality},
\begin{align}
    \eres_{xy}^{(\setS u)}-\eres_{xy}^{(\setS)}
    &=\eres_{uy}^{(\setS)}(A_x+D_x)
    =-\eres_{uu}^{(\setS)}(A_x+D_x)(A_y+D_y).
    \label{eqn:resolvent-deletion-expansion}
\end{align}
Conditionally on $\mH^{(\setS u)}$, Ward's identity gives
\begin{align}
    \frac1N\sum_a^{(\setS u)}|\eres_{xa}^{(\setS u)}|^2
    \leq N^{-\kappa}R_{\ell+1}.
\end{align}
Use Corollary~\ref{cor:sparse-Hanson-Wright-Estimates} with
$\Gamma=N^{-\kappa/2}R_{\ell+1}^{1/2}$ and
$\psi=R_{\ell+1}/\sqrt d$. Its ratio hypothesis holds since
$R_{\ell+1}\geq1$ and $d=O(\log N)$. Cauchy--Schwarz and Ward's
identity give the signal bound deterministically. Thus
\begin{align}
    |A_x|\leq C\frac{R_{\ell+1}}{\sqrt d},
    \qquad
    |D_x|\leq\frac{2\signal}{\sqrt{N\eta}}R_{\ell+1}^{1/2}
    \leq2N^{-\kappa/3}R_{\ell+1}^{1/2}.
    \label{eqn:minor-row-sum-bounds}
\end{align}
All these bounds hold simultaneously over the entries and deletion
steps with failure probability at most $CN^{-\nu}$.

\paragraph{Iteration.}
The first equality in \eqref{eqn:resolvent-deletion-expansion} implies
\begin{align}
    R_{\ell+1}\leq R_\ell
    +CR_\ell R_{\ell+1}(d^{-1/2}+N^{-\kappa/3}).
\end{align}
As long as $R_\ell\leq2\resmax$, absorption gives
$R_{\ell+1}\leq2R_\ell$. The second equality in
\eqref{eqn:resolvent-deletion-expansion} then improves this to
\begin{align}
    R_{\ell+1}\leq R_\ell
      \left(1+\frac{CR_{\ell+1}^2}{d}\right)
    \leq R_\ell\left(1+\frac{C_1\resmax^2}{d}\right),
\end{align}
since $N^{-\kappa/3}=o(d^{-1/2})$. Choosing the cutoff $c_*$ in
\eqref{eqn:setT-size} so that $C_1c_*\leq1/2$, induction gives
\begin{align}
    R_\ell\leq\resmax
    \left(1+\frac{C_1\resmax^2}{d}\right)^\ell
    \leq\sqrt e\,\resmax<2\resmax,
    \qquad 0\leq\ell\leq j.
    \label{eqn:minor-induction-bound}
\end{align}
This proves the first estimate in
\eqref{eqn:resolvent-entry-upper-bound}. Reserve one additional deletion
by decreasing $c_*$ if needed. Applying
\eqref{eqn:minor-row-sum-bounds} and the second equality in
\eqref{eqn:resolvent-deletion-expansion} to that deletion proves
\eqref{eqn:resolvent-schur-difference-upper-bound}.
Finally, \eqref{eqn:GxyT-expansion} with $\setT=\emptyset$ and the
same row-sum bounds give $|\eres_{xy}|\lesssim d^{-1/2}$ for $x\neq y$.
A union bound over these entries proves the remaining assertion.
\end{proof}

\begin{lemma}[Nested-minor comparison]
\label{lem:nested-minor-resolvent-bound}
Fix $c\geq1$ and $\nu>0$. There are constants
$c_{\rm nest},C>0$, depending only on $c,\kappa,\qr,L,\nu$, such that
for deterministic disjoint $\setS,\setT\subset[N]$ with
$|\setS\cup\setT|+1\leq c_{\rm nest}d$, and deterministic
$v\notin\setS\setT$, with probability at least $1-CN^{-\nu}$,
\begin{align}
    \theta_c^{(\setS)}
    \max_{a,b\notin\setS\setT}|\eres_{ab}^{(\setS\setT)}|
    &\leq2c\resmax,
    \label{eqn:nested-minor-entry-bound}\\
    \theta_c^{(\setS)}
    \max_{a,b\notin\setS\setT v}
    |\eres_{ab}^{(\setS\setT)}-\eres_{ab}^{(\setS\setT v)}|
    &\leq Cd^{-1}.
    \label{eqn:nested-minor-difference-bound}
\end{align}
In particular, outside the same exceptional event,
\begin{align}
    \theta_c^{(\setS)}\leq\theta_{2c}^{(\setS\setT)}.
    \label{eqn:nested-minor-indicator-comparison}
\end{align}
\end{lemma}

\begin{proof}
Repeat the preceding proof from $\res^{(\setS)}$, with $c\resmax$ in
place of $\resmax$. Deleting deterministic indices preserves the
independence and entry bounds used in the conditional estimates.
The constants are uniform in the sets, and choosing $c_{\rm nest}$
sufficiently small gives \eqref{eqn:minor-induction-bound} with the new
initial bound.
\end{proof}

For subsequent applications at a given failure exponent, decrease the
cutoff $c_*$ in \eqref{eqn:setT-size} so that
$c_*\resmax^{-2}\leq c_{\rm nest}/2$ for $1\leq c\leq4$.
Then \eqref{eqn:setT-size} permits up to two additional deletions in
Lemma~\ref{lem:nested-minor-resolvent-bound} for large $N$.

\subsection{Concentration and common-minor decoupling}
\label{sec:delocalized-typical-set-proof}

We first establish conditional concentration and one-deletion bounds
for the typicality parameters.

\begin{lemma}[Concentration of the typicality parameters]
\label{lem:Phi-Psi-concentration}
Let $\sU\subset[N]$ be deterministic and satisfy \eqref{eqn:setT-size}.
There are constants $0<\fkb<1$ and $C_{\resmax}>0$ such that, for
$1\leq c\leq4$, $x\in\sU$, and $0<\delta\leq1$ with
$\delta^2d\geq C_{\resmax}$, almost surely,
\begin{align}
    \theta_c^{(\sU)}\P\bigl(|\Phi_x^{(\sU)}|>\delta
        \mid\mH^{(\sU)}\bigr)
    \vee
    \theta_c^{(\sU)}\P\bigl(|\Psi_x^{(\sU)}|>\delta
        \mid\mH^{(\sU)}\bigr)
    \leq e^{-128\fkb\delta^2d}.
    \label{eqn:Phi-Psi-smallness}
\end{align}
For every $x\in[N]$ and $v\notin\sU$, with overwhelming probability,
\begin{align}
    |\Phi_x^{(\sU v)}-\Phi_x^{(\sU)}|\lesssim d^{-1},
    \qquad
    \theta_c^{(\sU)}|\Psi_x^{(\sU v)}-\Psi_x^{(\sU)}|
    \lesssim d^{-1}(1+\zeta_x).
    \label{eqn:diff-Phi-Psi-smallness}
\end{align}
The constants are uniform for $1\leq c\leq4$ and depend only on the
fixed entry bounds and $\resmax$.
\end{lemma}

\begin{proof}
For \eqref{eqn:Phi-Psi-smallness}, condition on $\mH^{(\sU)}$.
Since $x\in\sU$, the coefficients in $\Psi_x^{(\sU)}$ are
$\eres_{yy}^{(\sU)}$. On $\{\theta_c^{(\sU)}=1\}$ they are bounded
by $4\resmax$. Equation~\eqref{eqn:sparse-large-deviation-variance}
therefore gives, for even $q\leq d$,
\begin{align}
    \theta_c^{(\sU)}
    \left(\E[|\Psi_x^{(\sU)}|^q\mid\mH^{(\sU)}]\right)^{1/q}
    \leq C\resmax\sqrt{q/d}.
\end{align}
Take $q=2\lceil64\fkb\delta^2d\rceil$. Choosing $\fkb$ sufficiently
small and $C_{\resmax}$ sufficiently large ensures
$128\fkb\delta^2d\leq q\leq256\fkb\delta^2d\leq d$ and
$C\resmax\sqrt{q/d}\leq\delta/e$. Markov's inequality gives
\eqref{eqn:Phi-Psi-smallness}; for $\Phi_x^{(\sU)}$, use coefficients
identically one.

For \eqref{eqn:diff-Phi-Psi-smallness}, the case $v=x$ is immediate.
Otherwise, subtracting the definitions gives
\begin{align}
    \Psi_x^{(\sU v)}-\Psi_x^{(\sU)}
    ={}&\sum_y^{(\sU v x)}(|\emH_{xy}|^2-\E|\emH_{xy}|^2)
       (\eres_{yy}^{(\sU v x)}-\eres_{yy}^{(\sU x)})\notag\\
    &-(|\emH_{xv}|^2-\E|\emH_{xv}|^2)\eres_{vv}^{(\sU x)}.
\end{align}
Lemma~\ref{lem:nested-minor-resolvent-bound} bounds the resolvent
difference by $Cd^{-1}$ and the remaining entry by $2c\resmax$ on
$\{\theta_c^{(\sU)}=1\}$. Since
$\sum_y\E|\emH_{xy}|^2=O(1)$ and $|\emH_{xv}|^2\lesssim d^{-1}$,
the claimed bound follows. The difference of the two $\Phi$ parameters
is just $-(|\emH_{xv}|^2-\E|\emH_{xv}|^2)$.
\end{proof}

\begin{lemma}[Uniform row-square bound]
\label{lem:uniform-row-square-bound}
For each $\nu>0$, with probability at least $1-C_{\qr,\nu}N^{-\nu}$,
\begin{align}
    \max_x\zeta_x\leq C_{\qr,\nu}\left(1+\frac{\log N}{d}\right).
    \label{eqn:uniform-row-square-bound}
\end{align}
Consequently, for every fixed $K,\fka>0$, with overwhelming probability,
\begin{align}
    K\max_x\frac{1+\zeta_x}{d}\leq\varphi_{\fka/2}.
    \label{eqn:row-square-absorbed}
\end{align}
\end{lemma}

\begin{proof}
For each row, the variables $Y_y=|\emH_{xy}|^2$ are independent,
$0\leq Y_y\leq C/d$, and
$\sum_y\Var(Y_y)\leq(C/d)\E\zeta_x\lesssim d^{-1}$ by
\eqref{eqn:zeta-class-means}. Bennett's inequality
(Lemma~\ref{lem:Bennett}) gives
\begin{align}
    \P(\zeta_x-\E\zeta_x\geq t)
    \leq\exp\left(-\frac{cdt^2}{1+t}\right),\qquad t>0.
\end{align}
Take $t=A_{\qr,\nu}(1+\log N/d)$ and sum over $x$ to obtain
\eqref{eqn:uniform-row-square-bound}. On this event,
\begin{align}
    \frac{K(1+\max_x\zeta_x)}{d\varphi_{\fka/2}}
    \lesssim(d\log N)^{-1/3}
       +\left(\frac{\log N}{d^2}\right)^{2/3}=o(1),
\end{align}
by Assumption~\ref{ass:subcritical_regime}, proving
\eqref{eqn:row-square-absorbed}.
\end{proof}

\begin{corollary}\label{cor:PhiPsi-difference}
For fixed $\fka>0$, with overwhelming probability, simultaneously for
all $x,y\in[N]$,
\begin{align}
    \theta_1^{(\emptyset)}|\Phi_y-\Phi_y^{(x)}|
    \vee\theta_1^{(\emptyset)}|\Psi_y-\Psi_y^{(x)}|
    \leq\varphi_{\fka/2}.
    \label{eqn:single-vertex-typicality-transfer}
\end{align}
\end{corollary}

\begin{proof}
Use \eqref{eqn:diff-Phi-Psi-smallness} with $\sU=\emptyset$, $c=1$,
and a union bound over $x,y$. Equation~\eqref{eqn:row-square-absorbed}
absorbs both errors into $\varphi_{\fka/2}$.
\end{proof}

For the remaining proofs, fix $\fkb$ from
Lemma~\ref{lem:Phi-Psi-concentration} and write
\begin{align}
    \omega\coloneqq\fkb\varphi_{\fka}^2d.
    \label{eqn:atypical-decay-scale}
\end{align}
For every fixed $\fka>0$, the sparsity regime gives
$\omega\gtrsim(\log N)^{1/3}$ and $\omega=o(\log N)$.

\begin{lemma}[Joint atypicality]
\label{lem:typical-vertices-probability}
For each $\nu>0$ and fixed $\fka>0$, there is $C>0$, depending only
on $\nu$ and the fixed model parameters, such that every deterministic
$\sU\subset[N]$ with $|\sU|\leq\varphi_{\fka}d/C$ satisfies
\begin{align}
    \P(\sU\subset\Vtyp_{\fka/2}^{\complement},\,
        \theta_1^{(\emptyset)}=1)
    &\leq e^{-4\omega|\sU|}+CN^{-\nu},
    \label{eqn:typical-vertices-probability-1}\\
    \P(\sU\subset(\Vtyp_{\fka/2}^{(v)})^{\complement},\,
        \theta_2^{(v)}=1)
    &\leq e^{-4\omega|\sU|}+CN^{-\nu},\qquad v\notin\sU.
    \label{eqn:typical-vertices-probability-2}
\end{align}
\end{lemma}

\begin{proof}
We prove both estimates with
$(\setS,c)=(\emptyset,1)$ or $(\{v\},2)$, respectively.
The empty-set case $\sU=\emptyset$ is immediate. Otherwise, enlarge
$C$ so that all minors below satisfy \eqref{eqn:setT-size} and the
nested-minor size condition.

\paragraph{Transfer to a common minor.}
Delete the vertices of $\sU$ one at a time. Equations
\eqref{eqn:diff-Phi-Psi-smallness} and
\eqref{eqn:nested-minor-indicator-comparison}, with a union bound over
the intermediate minors, give on $\{\theta_c^{(\setS)}=1\}$,
\begin{align}
    |\Phi_x^{(\setS)}-\Phi_x^{(\setS\sU)}|
    &\lesssim |\sU|d^{-1},
    &
    |\Psi_x^{(\setS)}-\Psi_x^{(\setS\sU)}|
    &\lesssim |\sU|(1+\zeta_x)d^{-1},
    \label{eqn:joint-typicality-transfer}
\end{align}
simultaneously for $x\in\sU$, outside an event of probability
$CN^{-\nu}$. If $|\Phi_x^{(\setS)}|\leq\varphi_{\fka/2}$, then
\begin{align}
    \zeta_x\leq|\Phi_x^{(\setS)}|+\E\zeta_x
       +\sum_{y\in\setS}|\emH_{xy}|^2\lesssim1,
\end{align}
since $|\setS|\leq1$. Thus, for sufficiently large $C$, the first
transfer error is at most $\varphi_{\fka/4}$, and the second is also
at most $\varphi_{\fka/4}$ whenever $\Phi_x^{(\setS)}$ is small.
It follows that atypicality of $x$ implies the event
\begin{align}
    \sF_x\coloneqq
    \{|\Phi_x^{(\setS\sU)}|\vee|\Psi_x^{(\setS\sU)}|
       >\varphi_{\fka/4}\}.
\end{align}

\paragraph{Conditional product bound.}
Given $\mH^{(\setS\sU)}$, the events $(\sF_x)_{x\in\sU}$ are independent:
their random entries lie in distinct rows outside the common minor,
and edges within $\sU$ have been removed. The nested-minor comparison
gives $\theta_c^{(\setS)}\leq\theta_{2c}^{(\setS\sU)}$ with failure
probability $O(N^{-\nu})$. Since $2c\leq4$ and
$\varphi_{\fka/4}^2d\to\infty$, \eqref{eqn:Phi-Psi-smallness} yields
\begin{align}
    \theta_{2c}^{(\setS\sU)}
    \P(\sF_x\mid\mH^{(\setS\sU)})
    \leq2e^{-8\omega}\leq e^{-4\omega}.
\end{align}
Combining the transfer and indicator comparisons, then conditioning,
we obtain
\begin{align}
    \P\bigl(\sU\subset(\Vtyp_{\fka/2}^{(\setS)})^{\complement},\,
        \theta_c^{(\setS)}=1\bigr)
    &\leq\E\left[\theta_{2c}^{(\setS\sU)}
       \prod_{x\in\sU}\P(\sF_x\mid\mH^{(\setS\sU)})\right]
       +CN^{-\nu}\notag\\
    &\leq e^{-4\omega|\sU|}+CN^{-\nu}.
\end{align}
These are the two claimed bounds.
\end{proof}

\subsection{Size and edge weight of the atypical set}
\label{sec:delocalized-atypical-set-proof}

We now convert the joint atypicality bound into counting estimates.

\begin{lemma}[Counting atypical vertices]
\label{lem:atypical-set-size}
For each $\nu>0$, one can choose a fixed $\fka>0$ such that, for every
deterministic $\sU\subset[N]$, the following estimates fail on
$\{\theta_1^{(\emptyset)}=1\}$ with probability at most $CN^{-\nu}$:
\begin{enumerate}[label=\textup{(\arabic*)}]
    \item $|\sU\cap\Vtyp_{\fka/2}^{\complement}|
        \leq e^\omega+|\sU|e^{-2\omega}$.
    \item If $|\sU|\leq e^{2\omega}$, then
        $|\sU\cap\Vtyp_{\fka/2}^{\complement}|\leq\varphi_{\fka}d$.
\end{enumerate}
For fixed $v\in[N]$, the same estimates hold with
$(\Vtyp_{\fka/2}^{(v)})^{\complement}$ in place of
$\Vtyp_{\fka/2}^{\complement}$ for any possibly random
$\sU\subset[N]\setminus\{v\}$ independent of $\mH^{(v)}$.
They hold on either $\{\theta_2^{(v)}=1\}$ or
$\{\theta_1^{(\emptyset)}=1\}$, with the same failure bound.
\end{lemma}

\begin{proof}[Proof of Lemma~\ref{lem:atypical-set-size}]
Set $\nu_0=20(\nu+1)$ and let $C_0\geq1$ be the constant in
Lemma~\ref{lem:typical-vertices-probability} for exponent $\nu_0$.
Choose
\begin{align}
    \fka^3=\frac{4C_0(\nu+1)}{\fkb},
    \qquad \ell=\left\lfloor\frac{\varphi_{\fka}d}{C_0}\right\rfloor,
    \qquad Z=|\sU\cap\Vtyp_{\fka/2}^{\complement}|.
    \label{eqn:atypical-counting-parameters}
\end{align}
The identity $\varphi_{\fka}^3d^2=\fka^3\log N$ and the sparsity
regime imply
\begin{align}
    \omega\ell=(4(\nu+1)+o(1))\log N,
    \qquad \log\ell=o(\omega).
    \label{eqn:omega-ell}
\end{align}

\paragraph{Small sets.}
If $|\sU|\leq e^{2\omega}$ and $Z\geq\ell$, some $\ell$-element
subset of $\sU$ is entirely atypical. Equation
\eqref{eqn:typical-vertices-probability-1} and a union bound give
\begin{align}
    \P(Z\geq\ell,\,\theta_1^{(\emptyset)}=1)
    &\leq\binom{|\sU|}{\ell}
       (e^{-4\omega\ell}+C_0N^{-\nu_0})
    \leq e^{-2\omega\ell}+C_0N^{-\nu_0}e^{2\omega\ell}
    \leq CN^{-\nu}.
\end{align}
Since $\ell\leq\varphi_{\fka}d$, this proves \textup{(2)}.

\paragraph{Arbitrary sets.}
Expand $Z^\ell$ and group the tuples according to their $m$ distinct
indices. There are at most
$\binom{\ell}{m}\ell^{\ell-m}|\sU|^m$ such tuples. Using
\eqref{eqn:typical-vertices-probability-1} for each set of distinct
indices yields
\begin{align}
    \E[\theta_1^{(\emptyset)}Z^\ell]
    &\leq\sum_{m=1}^{\ell}\binom{\ell}{m}\ell^{\ell-m}|\sU|^m
       (e^{-4\omega m}+C_0N^{-\nu_0})\notag\\
    &\leq(\ell+|\sU|e^{-4\omega})^\ell
       +C_0N^{-\nu_0}(\ell+|\sU|)^\ell.
    \label{eqn:atypical-count-moment}
\end{align}
For $t=e^\omega+|\sU|e^{-2\omega}$,
\begin{align}
    \frac{\ell+|\sU|e^{-4\omega}}{t}
    \leq e^{-(1-o(1))\omega},
    \qquad
    \frac{\ell+|\sU|}{t}\leq e^{(2+o(1))\omega}.
\end{align}
Markov's inequality and \eqref{eqn:omega-ell} therefore give
\begin{align}
    \P(Z\geq t,\,\theta_1^{(\emptyset)}=1)
    \leq e^{-(1-o(1))\omega\ell}
       +C_0N^{-\nu_0}e^{(2+o(1))\omega\ell}
    \leq CN^{-\nu},
\end{align}
which proves \textup{(1)}.

\paragraph{Independent random sets in a minor.}
For deterministic $\sU$, the same proof uses
\eqref{eqn:typical-vertices-probability-2} to give both estimates on
$\{\theta_2^{(v)}=1\}$. If $\sU$ is random and independent of
$\mH^{(v)}$, condition on $\sU$; the bound is uniform over its
realizations. Taking expectations proves the random-set version on
$\{\theta_2^{(v)}=1\}$. By Lemma~\ref{lem:res-max-difference-bound},
$\theta_1^{(\emptyset)}\leq\theta_2^{(v)}$ outside an event of
probability $O(N^{-\nu})$, which transfers both estimates to
$\{\theta_1^{(\emptyset)}=1\}$.
\end{proof}

\begin{proof}[Proof of Proposition~\ref{prop:size-typical-vertices}]
Use Lemma~\ref{lem:atypical-set-size} with exponent $\nu+3$.
Since $\Vtyp_{\fka}^{\complement}\subset\Vtyp_{\fka/2}^{\complement}$,
its first estimate with $\sU=[N]$ gives the sharper bound
\begin{align}
    \frac{|\Vtyp_{\fka}^{\complement}|}{N}
    \leq N^{-1}e^\omega+e^{-2\omega}.
    \label{eqn:atypical-set-cardinality-sharp}
\end{align}

For the edge weight, \eqref{eqn:single-vertex-typicality-transfer}
implies, simultaneously for every $x$,
\begin{align}
    \Vtyp_{\fka}^{\complement}\setminus\{x\}
    \subset(\Vtyp_{\fka/2}^{(x)})^{\complement}
    \quad\text{on }\{\theta_1^{(\emptyset)}=1\},
    \label{eqn:Vtyp-inclusion}
\end{align}
outside an event of probability $O(N^{-\nu})$.
Fix $C_H$ with $|\emH_{xy}|^2\leq C_H/d$. For a fixed row $x$, put
\begin{align}
    \sU_k\coloneqq\{y\neq x:|\emH_{xy}|^2\geq C_Hd^{-k-2}\},
    \qquad 0\leq k\leq\lfloor\log N\rfloor.
\end{align}
These sets are independent of $\mH^{(x)}$. Decomposing the weights
into levels with ratio $d$, and using $\emH_{xx}=0$, gives
\begin{align}
    \sum_{y\in(\Vtyp_{\fka/2}^{(x)})^{\complement}}|\emH_{xy}|^2
    \lesssim\sum_{k=0}^{\lfloor\log N\rfloor}d^{-k-1}
       |\sU_k\cap(\Vtyp_{\fka/2}^{(x)})^{\complement}|+N^{-1}.
    \label{eqn:weight-level-decomposition}
\end{align}
By \eqref{eqn:uniform-row-square-bound} and $\log N\leq d^2$,
$\zeta_x\leq C_{\qr,\nu}d$ with overwhelming probability. Hence
\begin{align}
    |\sU_k|\leq C_{\qr,\nu}d^{k+3}.
    \label{eqn:weight-level-size}
\end{align}
The random-set version of Lemma~\ref{lem:atypical-set-size} applies
simultaneously to all $\sU_k$ and all rows after a union bound over
$O(N\log N)$ pairs.

Let $K$ be the largest integer at most $\lfloor\log N\rfloor$ with
$C_{\qr,\nu}d^{K+3}\leq e^{2\omega}$. Such an integer exists since
$\omega\gg\log d$. For $k\leq K$, the small-set bound gives
\begin{align}
    \sum_{k=0}^{K}d^{-k-1}
    |\sU_k\cap(\Vtyp_{\fka/2}^{(x)})^{\complement}|
    \leq\varphi_{\fka}\sum_{k=0}^{K}d^{-k}\leq2\varphi_{\fka}.
\end{align}
For $k>K$, the arbitrary-set bound and
\eqref{eqn:weight-level-size} give
\begin{align}
    \sum_{k=K+1}^{\lfloor\log N\rfloor}d^{-k-1}
    |\sU_k\cap(\Vtyp_{\fka/2}^{(x)})^{\complement}|
    &\lesssim e^\omega d^{-K-2}+d^2e^{-2\omega}\log N\notag\\
    &\lesssim d^2e^{-\omega}+d^2e^{-2\omega}\log N
    \lesssim d^4e^{-\omega}.
\end{align}
Here maximality gives $d^{-K-2}\lesssim d^2e^{-2\omega}$; if
$K=\lfloor\log N\rfloor$, the tail is empty. Combining these estimates
with \eqref{eqn:Vtyp-inclusion} and
\eqref{eqn:weight-level-decomposition} yields
\begin{align}
    \max_x\sum_{y\in\Vtyp_{\fka}^{\complement}}|\emH_{xy}|^2
    \lesssim\varphi_{\fka}+d^4e^{-\omega}.
    \label{eqn:atypical-set-weight-sharp}
\end{align}
Finally, \eqref{eqn:atypical-decay-scale} and the sparsity regime give
\begin{align}
    N^{-1}e^\omega+e^{-2\omega}+d^4e^{-\omega}=o(d^{-1}).
    \label{eqn:atypical-errors-absorbed}
\end{align}
Since $d^{-1}=o(\varphi_{\fka})$, the two sharper bounds imply
\eqref{eqn:atypical-set-cardinality} and
\eqref{eqn:atypical-set-weight}. Proportional class sizes and
$\varphi_{\fka}=o(1)$ also give
$|\Vtyp_{j,\fka}|\geq|\gV_j|/2$ for large $N$.
\end{proof}

\subsection{Schur expansions and error estimates}
\label{sec:delocalized-schur-proof}
\label{sec:delocalized-atypical-schur-proof}

We expand the residuals defined by
\eqref{eqn:res-inverse-expansion-raw} and
\eqref{eqn:res-inverse-expansion-error-epsilon}. Define the variance
defects by
\begin{align}
    r_{xy}^{\rm var}\coloneqq\E|\emH_{xy}|^2
       -\frac{\indi{\evsigma_x\neq\evsigma_y}}{\sqrt{mn}},
    \qquad
    \rho_j\coloneqq\frac{|\gV_{3-j}|}{\sqrt{mn}}-c_j.
    \label{eqn:variance-defect-parameters}
\end{align}
Definition~\ref{def:sparse-random-rectangular-matrix} and the proportional
regime imply
\begin{align}
    \max_{x,y}|r_{xy}^{\rm var}|\leq\frac{Cd}{N\sqrt{mn}},
    \qquad \max_{j=1,2}|\rho_j|\leq\frac{C_{\qr}}N.
    \label{eqn:variance-defect-bounds}
\end{align}

\begin{proof}[Schur expansions for Lemma~\ref{lem:approximate-self-consistent}]
Insert \eqref{eqn:two-spike-model} into
\eqref{eqn:schur-complement-diagonal}. Since $\emM_{xx}=\emH_{xx}=0$,
separating the diagonal noise terms and combining the mixed signal
terms by symmetry gives
\begin{align}
    \rY_x={} &
    -\sum_{a\neq b}^{(x)}\emH_{xa}\eres_{ab}^{(x)}\emH_{bx}
    -\frac{4\signal}{N}\sum_{a,b}^{(x)}
      \indi{\evsigma_a\neq\evsigma_x}\eres_{ab}^{(x)}\emH_{bx}
    -\frac{4\signal^2}{N^2}\sum_{a,b}^{(x)}
      \indi{\evsigma_a\neq\evsigma_x}
      \indi{\evsigma_b\neq\evsigma_x}\eres_{ab}^{(x)},
    \label{eqn:error-Y-x}
\end{align}

For $x\in\gV_j$, centering the row weights gives
\begin{align}
    \sum_y^{(x)}|\emH_{xy}|^2\eres_{yy}^{(x)}
    =\Psi_x+\sum_y^{(x)}r_{xy}^{\rm var}\eres_{yy}^{(x)}
       +\frac{c_j+\rho_j}{|\gV_{3-j}|}
          \sum_{y\in\gV_{3-j}}\eres_{yy}^{(x)}.
\end{align}
Substituting into \eqref{eqn:res-inverse-expansion-raw} and adding and
subtracting $c_j\overline G_{3-j}$ yields
\begin{align}
    \widehat\epsilon_x=\rY_x-\Psi_x
       -\sum_y^{(x)}r_{xy}^{\rm var}\eres_{yy}^{(x)}
       +c_j\overline G_{3-j}
       -\frac{c_j+\rho_j}{|\gV_{3-j}|}
          \sum_{y\in\gV_{3-j}}\eres_{yy}^{(x)}.
    \label{eqn:error-epsilon-x}
\end{align}
\end{proof}

\begin{proof}[Proof of Lemma~\ref{lem:schur-complement-error}]
We work on $\{\theta_1^{(\emptyset)}=1\}$. By
Lemma~\ref{lem:res-max-difference-bound}, the bounds
$\theta_2^{(x)}=1$ hold simultaneously for all singleton minors outside
an event of overwhelmingly small probability.

\paragraph{The remainder $\rY_x$.}
Condition on $\mH^{(x)}$. For the off-diagonal quadratic form in
\eqref{eqn:error-Y-x}, Ward's identity permits
$\Gamma=O(N^{-\kappa/2})$ and $\psi=O(d^{-1})$ on
$\{\theta_2^{(x)}=1\}$. Corollary~\ref{cor:sparse-Hanson-Wright-Estimates}
therefore bounds it by $Cd^{-1}$.
For the mixed signal term, put
\begin{align}
    a_b\coloneqq\frac{4\signal}{N}
       \sum_a^{(x)}\indi{\evsigma_a\neq\evsigma_x}\eres_{ab}^{(x)}.
\end{align}
Cauchy--Schwarz and Ward's identity give
\begin{align}
    \max_b|a_b|\lesssim\frac{\signal}{\sqrt{N\Im z}}
    \lesssim N^{-\kappa/3}.
\end{align}
Apply \eqref{eqn:sparse-large-deviation-mean} conditionally, with an
even moment of order $\log N$. Since the coefficient norm $\Gamma$
can be bounded by $\max_b|a_b|$, this gives
\begin{align}
    \left|\sum_b^{(x)}a_b\emH_{bx}\right|
    \lesssim N^{-\kappa/3}\log N=o(d^{-1})
\end{align}
with overwhelming probability. The quadratic signal term
satisfies
\begin{align}
    \frac{4\signal^2}{N^2}
    \left|\sum_{a,b}^{(x)}
       \indi{\evsigma_a\neq\evsigma_x}
       \indi{\evsigma_b\neq\evsigma_x}\eres_{ab}^{(x)}\right|
    &\leq\frac{4\signal^2}{N}
       \left(\sum_{a,b}^{(x)}|\eres_{ab}^{(x)}|^2\right)^{1/2}\notag\\
    &\lesssim\frac{\signal^2}{\sqrt{N\Im z}}
    \leq N^{-\kappa/6}=o(d^{-1}).
\end{align}
Here we used $\signal\leq N^{\kappa/6}$, $N\Im z\geq N^\kappa$, and
$d=O(\log N)$. A union bound over $x$ proves the first estimate in
\eqref{eqn:schur-complement-error-bounds}.

\paragraph{The residual in the typical-class equation.}
Proposition~\ref{prop:size-typical-vertices} ensures
$|\Vtyp_{j,\fka}|\geq|\gV_j|/2$ for large $N$.
For $x\in\Vtyp_{j,\fka}$, \eqref{eqn:error-epsilon-x} gives
\begin{align}
    |\widehat\epsilon_x|
    \leq{}&|\rY_x|+|\Psi_x|
       +\sum_y^{(x)}|r_{xy}^{\rm var}||\eres_{yy}^{(x)}|
       +\frac{|\rho_j|}{|\gV_{3-j}|}
          \sum_{y\in\gV_{3-j}}|\eres_{yy}^{(x)}|\notag\\
    &+\frac{c_j}{|\gV_{3-j}|}\sum_{y\in\gV_{3-j}}
       |\eres_{yy}^{(x)}-\eres_{yy}|
       +c_j\left|\overline G_{3-j}
          -\frac1{|\gV_{3-j}|}\sum_{y\in\gV_{3-j}}\eres_{yy}\right|.
    \label{eqn:epsilon-decomposition-bound}
\end{align}
The first two terms are $O(d^{-1})$ and at most $\varphi_{\fka}$,
respectively. The variance defects contribute $O(d/N+N^{-1})$ by
\eqref{eqn:variance-defect-bounds}; the minor difference contributes
$O(d^{-1})$ by \eqref{eqn:resolvent-schur-difference-upper-bound}.
Finally, boundedness of the diagonal entries gives
\begin{align}
    \left|\overline G_{3-j}
       -\frac1{|\gV_{3-j}|}\sum_{y\in\gV_{3-j}}\eres_{yy}\right|
    \leq\frac{2\resmax|\gV_{3-j}\setminus\Vtyp_{3-j,\fka}|}
       {|\gV_{3-j}|}
    \lesssim\varphi_{\fka}.
\end{align}
All terms are $O(\varphi_{\fka})$ by
\eqref{eqn:atypical-set-cardinality} and the sparsity regime. This proves
the second estimate in \eqref{eqn:schur-complement-error-bounds}.
\end{proof}

\subsection{Instability of the neighborhood equation}
\label{sec:delocalized-neighborhood-instability-proof}

Let $\mA_*$ be the adjacency matrix of a finite bipartite graph on
$N$ vertices with
degree $D_j=dc_j\in\N$, $D_j\geq2$, at every vertex in class $j$.
Assume that the ball of radius $r\in\N$, $r\geq2$, around some vertex is a tree,
and set $\mS=\mA_*/d$ and $m_j(z)=\stj^{(j)}(z)$.
The neighborhood equation and its linearization at the class profiles
are
\begin{align}
    \frac1{g_x}
    &=-z-\sum_y S_{xy}g_y+\rho_x,
    \label{eqn:bipartite-neighborhood-equation}\\
    \mathcal L_{\rm nb}(z)
    &\coloneqq I-\operatorname{diag}\bigl(m_{\evsigma_x}(z)^2\bigr)\mS.
    \label{eqn:bipartite-neighborhood-linearization}
\end{align}
At $\eta=0$, the profiles are interpreted as boundary values from
$\C_+$. There are constants $\eta_0,c>0$, depending only on $\kappa,\qr$,
such that
\begin{align}
    \|\mathcal L_{\rm nb}(E+\ii\eta)^{-1}\|_{\infty\to\infty}
    \geq\frac{c}{r^{-1}+d^{-1}+\eta},
    \qquad E\in\bulk,\quad 0\leq\eta\leq\eta_0.
    \label{eqn:bipartite-neighborhood-instability}
\end{align}
\begin{proof}[Proof of \eqref{eqn:bipartite-neighborhood-instability}]
Write $\beta(z)=m_1(z)m_2(z)$ and define
\begin{align}
    P_{xy}\coloneqq\frac{(\mA_*)_{xy}}{D_{\evsigma_x}},
    \qquad
    \mT(z)\coloneqq
       \operatorname{diag}\bigl(\sqrt{c_{\evsigma_x}}\,m_{\evsigma_x}(z)\bigr).
\end{align}
The matrix $\mP$ is stochastic and is similar to the symmetric matrix
$\mS$, since $c_1c_2=1$. A direct calculation gives
\begin{align}
    \mT(z)^{-1}\mathcal L_{\rm nb}(z)\mT(z)=I-\beta(z)\mP.
    \label{eqn:neighborhood-stability-similarity}
\end{align}
The explicit profiles \eqref{eqn:stieltjes-V1}--\eqref{eqn:stieltjes-V2}
are bounded and bounded away from zero on a sufficiently small strip
above the bulk. Hence $\|\mT(z)\|_{\infty\to\infty}
\|\mT(z)^{-1}\|_{\infty\to\infty}\leq C_{\kappa,\qr}$ there.

\paragraph{The bulk phase.}
The coupled profile equations imply
\begin{align}
    z^2\beta=1+(c_1+c_2)\beta+\beta^2.
\end{align}
For $E\in\bulk$, put $\beta_0=\beta(E+\ii0)$. Then
\begin{align}
    \beta_0+\beta_0^{-1}=E^2-(c_1+c_2)\in(-2,2).
\end{align}
Thus $|\beta_0|=1$ and $\beta_0$ is separated from $\pm1$.
Since $\mP$ has real spectrum, $I-\beta_0\mP$ is invertible.
For $\eta>0$, taking imaginary parts of the profile equations gives
\begin{align}
    \Im m_j=|m_j|^2\bigl(\eta+c_j\Im m_{3-j}\bigr),
    \qquad j=1,2.
\end{align}
Multiplying the resulting strict inequalities yields $|\beta(z)|<1$,
so $I-\beta(z)\mP$ is invertible as well.

\paragraph{A radial approximate null vector.}
Let $o$ be the center of the tree ball and write
$k=\operatorname{dist}(o,x)$ for $x$ in this ball. Define
\begin{align}
    f_k=(1-k/r)_+,\qquad u_x=\beta_0^{-k}f_k.
    \label{eqn:neighborhood-instability-test-vector}
\end{align}
Set $u_x=0$ outside the ball. Then $\|\vu\|_\infty=1$, and $\vu$
vanishes at distance $r$ and beyond.
For $1\leq k\leq r-1$, a vertex $x$ of class $j$ has one inward
neighbor and $D_j-1$ outward neighbors. Therefore
\begin{align}
    \bigl((I-\beta_0\mP)\vu\bigr)_x
    =\beta_0^{-k}
      \left(f_k-f_{k+1}
       +\frac{f_{k+1}-\beta_0^2f_{k-1}}{D_j}\right).
\end{align}
Its modulus is at most $r^{-1}+2D_j^{-1}$.
At the root the residual is $r^{-1}$, and at distance $r$ its
modulus is $(rD_j)^{-1}$. Outside the ball it vanishes, because
an outside vertex has no neighbor at distance less than $r$.
Consequently,
\begin{align}
    \|(I-\beta_0\mP)\vu\|_\infty
    \leq r^{-1}+2D_{\min}^{-1}
    \lesssim_{\qr} r^{-1}+d^{-1}.
    \label{eqn:neighborhood-instability-radial-residual}
\end{align}
Here $D_{\min}=\min\{D_1,D_2\}$.

\paragraph{Positive imaginary part.}
The explicit profile formulas are uniformly Lipschitz in $\eta$ on
a sufficiently small bulk strip. In particular,
$|\beta(E+\ii\eta)-\beta_0|\leq C_{\kappa,\qr}\eta$.
Since $\|\mP\|_{\infty\to\infty}=1$, the same vector satisfies
\begin{align}
    \|(I-\beta(E+\ii\eta)\mP)\vu\|_\infty
    \lesssim_{\kappa,\qr}r^{-1}+d^{-1}+\eta.
\end{align}
Set $\vv=\mT(z)\vu/\|\mT(z)\vu\|_\infty$.
The uniform bounds on $\mT$ and its inverse, together with
\eqref{eqn:neighborhood-stability-similarity}, give
\begin{align}
    \|\vv\|_\infty=1,\qquad
    \|\mathcal L_{\rm nb}(z)\vv\|_\infty
    \lesssim_{\kappa,\qr}r^{-1}+d^{-1}+\eta.
    \label{eqn:neighborhood-instability-normalized-residual}
\end{align}
Applying $\mathcal L_{\rm nb}(z)^{-1}$ proves
\eqref{eqn:bipartite-neighborhood-instability}.
\end{proof}

\paragraph{A logarithmic benchmark.}
The scales $d\asymp\log N$ and $r\asymp\log N/\log d$ can be
realized by finite biregular graphs. For even $r$, take $D_2$ copies
of the alternating tree of depth $r$ rooted in class $1$, and $D_1$
copies rooted in class $2$. Their leaves belong to opposite classes,
and both total degree deficits equal
$D_1D_2[(D_1-1)(D_2-1)]^{r/2}$.
Join the leaf classes by a simple
$(D_1-1,D_2-1)$-biregular graph to complete the degrees.
This preserves the radius-$r$ tree ball around each root.
For bounded $\qr$ and $r\asymp d/\log d$, the resulting graph has
$\log N\asymp r\log d\asymp d$, as required.

At local spectral scales, \eqref{eqn:bipartite-neighborhood-instability}
then gives an inverse norm at least of order $\log N/\log\log N$,
while $\varphi_{\fka}\asymp(\log N)^{-1/3}$ for fixed $\fka$.

\paragraph{Failure of a uniform nonlinear stability estimate.}
For this benchmark, take $z=E+\ii N^{-1+\kappa}$ and the normalized
vector $\vv$ from \eqref{eqn:neighborhood-instability-normalized-residual}.
Put $t=\sqrt{\varphi_{\fka}}$, $h_x=tv_x$, and
$g_x=m_{\evsigma_x}+h_x$. Since the bulk profiles have imaginary
parts bounded below, $g_x\in\C_+$ for large $N$. Moreover,
$\max_x|g_x-m_{\evsigma_x}|=t=o(1)$, so these vectors satisfy any
fixed small a priori proximity threshold.
Define $\rho_x$ by \eqref{eqn:bipartite-neighborhood-equation}.
Using the exact profile equation and the biregular row sums gives
\begin{align}
    \rho_x
    &=-\frac{(\mathcal L_{\rm nb}(z)\vh)_x}{m_j^2}
      +\frac{h_x^2}{m_j^2(m_j+h_x)},\qquad x\text{ in class }j,\\
    \|\rho\|_\infty
    &\lesssim_{\kappa,\qr}
      t(r^{-1}+d^{-1}+N^{-1+\kappa})+t^2.
\end{align}
Here $r^{-1}+d^{-1}+N^{-1+\kappa}
=O(\log\log N/\log N)=o(\sqrt{\varphi_{\fka}})$.
It follows that
\begin{align}
    \|\rho\|_\infty\lesssim_{\kappa,\qr}\varphi_{\fka},
    \qquad
    \max_x|g_x-m_{\evsigma_x}|=\sqrt{\varphi_{\fka}}
       \gg\varphi_{\fka}.
    \label{eqn:neighborhood-nonlinear-instability}
\end{align}
Thus weak proximity and a residual of size $O(\varphi_{\fka})$ do not
imply an error of the same size through a graph-uniform deterministic
estimate for general approximate solutions of the neighborhood
equation. This example makes no assertion about the stochastic
structure of genuine Schur residuals.

\subsection{Bulk stability}
\label{sec:delocalized-stability-proof}

\begin{proof}[Proof of Lemma~\ref{lem:stability-self-consistent}]
For $x\in\sX_j$, put $\evh_x=\evg_x-\stj^{(j)}$ and let
$\overline h_j=|\sX_j|^{-1}\sum_{x\in\sX_j}\evh_x$.
Subtracting \eqref{eqn:stieltjes-m1-m2-definition} from
\eqref{eqn:self-consistent-equation-stability} gives
\begin{align}
    (\mB\vh)_x&\coloneqq\evh_x-c_j(\stj^{(j)})^2\overline h_{3-j}
    =c_j\stj^{(j)}\evh_x\overline h_{3-j}
       -\evg_x\stj^{(j)}\epsilon_x.
    \label{eqn:stability-operator-B}
\end{align}

\paragraph{Invertibility of the linear part.}
For $\vu=\mB\vh$, class averaging gives
\begin{align}
    \begin{bmatrix}
        1&-c_1(\stj^{(1)})^2\\
        -c_2(\stj^{(2)})^2&1
    \end{bmatrix}
    \begin{bmatrix}\overline h_1\\\overline h_2\end{bmatrix}
    =\begin{bmatrix}\overline u_1\\\overline u_2\end{bmatrix}.
\end{align}
Its determinant is $1-(\stj^{(1)}\stj^{(2)})^2$, since $c_1c_2=1$.
Writing $p=\stj^{(1)}\stj^{(2)}$, the Dyson system gives
\begin{align}
    z^2p=1+(\qr^2+\qr^{-2})p+p^2.
    \label{eqn:stability-product-identity}
\end{align}
Thus a zero of the determinant would have $p=\pm1$ and
$z^2=(\qr\pm\qr^{-1})^2$, which are the spectral edges.
These points are separated from $\Re z\in\bulk$. The profiles in
\eqref{eqn:stieltjes-V1}--\eqref{eqn:stieltjes-V2} are continuous up to
the real axis inside the bulk and tend to zero
as $\Im z\to\infty$. Compactness therefore gives
$|1-p^2|\geq c_{\kappa,\qr}>0$ uniformly in this region.
Solving the averaged system and returning to
\eqref{eqn:stability-operator-B} yields
\begin{align}
    \|\mB^{-1}\|_{\infty\to\infty}\leq C_{\kappa,\qr}.
    \label{eqn:stability-inverse-bound}
\end{align}

\paragraph{Absorbing the nonlinear error.}
The profiles are uniformly bounded, and the a priori assumption gives
$\max_x|\evg_x|\lesssim1$. Equations
\eqref{eqn:stability-operator-B} and \eqref{eqn:stability-inverse-bound}
then imply
\begin{align}
    \|\vh\|_\infty\leq C_{\kappa,\qr}
    \left(\|\vh\|_\infty^2+\max_x|\epsilon_x|\right).
\end{align}
Choose $\epsilon_*$ so that the quadratic term is at most
$\|\vh\|_\infty/2$. Absorption proves
\eqref{eqn:stability-self-consistent}.
\end{proof}

%% file: app-semilocalized.tex
\section{Deferred Proofs in Section~\ref{sec:semilocalized_phase}}
\label{sec:semilocalized-deferred-proofs}

\subsection{Pruning and layer concentration}
\label{subsec:semilocalized-pruning-proofs}

\subsubsection{Cutoffs and admissible scales}
Write $L_N=\log N$, $\xi=\error_d$, and $R_0=8\radius+2$.
Fix $a>0$ and restrict to $a\xi^{1/4}\leq s=s(\gap)\leq1$.
Constants may depend on $(B,\qr,\nu,a)$ and on any specified root
degree floor. Throughout this subsection $N$ is sufficiently large.
The sparsity assumptions and \eqref{eqn:error} give
\begin{align}
 \xi&\leq C_B/\gC_0,& d&\gg\log L_N,\nonumber\\
 \frac{ds^2}{R_0\log L_N}&\longrightarrow\infty,
 &\frac{L_N}{d^2s^3}&\longrightarrow0,
 &ds^2&\gg\log d.
 \label{eqn:pruning-scale-comparison}
\end{align}
Indeed, $ds^2\geq a^2\sqrt d\,L_N^{1/4}\sqrt{\log d}$,
whereas $R_0\log L_N=O_\qr(\sqrt{d/\log d}\log L_N)$.
The second ratio is at most $a^{-3}\xi^{5/4}/(\log d)^2$.

Each exceptional set depends on three integer degree cutoffs.
The affine thresholds in \eqref{eqn:energy-dependent-thresholds}
have $O_\qr(d)$ deterministic breakpoints and hence, including
endpoint conventions, $O_\qr(d)$ possible triples. Choose a
deterministic representative in each cutoff interval and include
the endpoint triples. The intervals have width $O_\qr(d^{-1})$,
so $s'/s=1+o(1)$ uniformly because $da\xi^{1/4}\to\infty$.
The construction depends on the triple; its $s$ budgets extend over
the interval after adjusting constants. Prove the estimates at all
representatives with an input failure exponent larger than $\nu$,
then intersect their events. This deterministic-cutoff intersection
precedes any eigenvalue-dependent choice of gap and requires no
independence between events. For $\gap\geq1$, the pruning is fixed
at $s=1$. Let $\eta_\gap=\xi+\error_{\mathrm{prof}}(\gap)$. The truncation
$s=\gap\wedge1$ gives
\begin{align}
 \frac{\sqrt{L_N}}{ds}
 \leq\xi+\error_{\mathrm{prof}}(\gap)=\eta_\gap.
 \label{eqn:truncated-gap-error}
\end{align}
Since $d^{-1}\leq\xi$, $\xi^{-1}\leq d$, and
$s\radius/\log d\to\infty$, for every fixed $c>0$ and
$K>0$ the comparison budgets are
\begin{align}
 d^{-1}+\frac{L_N}{d^2s^2}&\leq C_a\eta_\gap,\nonumber\\
 (1+K L_N/d)s^{-4}e^{-cs\radius}&=o(\eta_\gap),\nonumber\\
 \frac{\eta_\gap}{s}
 &\leq a^{-1}\xi^{3/4}
      +a^{-2}\frac{\xi^{1/2}}{\log d},
       &&\gap\leq1,\nonumber\\
 \frac{\eta_\gap}{s}&=O(\xi),&&\gap\geq1.
 \label{eqn:pruning-generic-comparison-budgets}
\end{align}
These ratios are uniformly small after choosing $\gC_0$, even when
$\xi$ does not tend to zero at the lowest sparsity. At critical sparsity,
$d=bL_N$ and $\xi\asymp_b\log L_N/\sqrt{L_N}$, so the ratios tend to zero.

At the location gap $s=3\xi^{1/4}/4$, choose $\gC_0$ so that
$s\leq\min\{1,4/\qr\}$ and $\xi\leq1$, and set
$\eta_{\rm loc}=\xi+4\xi^{3/4}/(3\log d)$. Then
\begin{align}
 \frac{ds^2}{\log d}&\longrightarrow\infty,
 &\frac{s\radius}{\log d}&\longrightarrow\infty,\nonumber\\
 \frac{L_N}{d^2s^2}
 &=\frac{16\xi^{3/2}}{9(\log d)^2}=o(s),
 &\frac{\sqrt{L_N}}{ds}
 &=\frac{4\xi^{3/4}}{3\log d}.
 \label{eqn:location-pruning-scale-identities}
\end{align}
The error-to-gap ratios are
\begin{align}
 \frac{\eta_{\rm loc}}{\xi^{1/4}}
 &=\xi^{3/4}+\frac{4\xi^{1/2}}{3\log d},
 &\frac{\eta_{\rm loc}}s
 &=\frac43\xi^{3/4}+\frac{16\xi^{1/2}}{9\log d}.
 \label{eqn:location-error-gap-ratio}
\end{align}
They are uniformly as small as needed on the right by increasing
$\gC_0$; at critical sparsity they tend to zero. The weighted
square of the root-loss scale is
\begin{align}
 \xi^{-1}\frac{L_N^2}{d^4s^4}
 =\frac{256\xi^2}{81(\log d)^4}.
 \label{eqn:location-weighted-root-loss-square}
\end{align}
For every fixed $c,K>0$, the exponential and column-loss budgets are
\begin{align}
 \frac{256}{81\xi}e^{-\frac{3c}{4}\xi^{1/4}\radius}
 &=o(\eta_{\rm loc}),
 &
 (1+K L_N/d)s^{-4}e^{-cs\radius}&=o(\eta_{\rm loc}),
 &
 d^{-1}+\frac{16L_N}{9d^2\xi^{1/2}}
 &\leq C\eta_{\rm loc}.
 \label{eqn:location-exponential-column-budgets}
\end{align}
The first two bounds follow from $s\radius/\log d\to\infty$;
the last follows from $d^{-1}\leq\xi$ and
\eqref{eqn:location-pruning-scale-identities}. These calculations
verify the location gap directly, independently of the larger
constant required for the window statements.

\subsubsection{Exploration, exceptional vertices and collisions}
For every root $x$, fix a breadth-first spanning tree of the induced
unpruned ball $\ball_{R_0}^{\gG}(x)$, using a deterministic ordering
to select the parent whenever more than one is possible. For
$A\subset\layer_1^{\gG}(x)$, write $W_j^x(A)$ for the vertices in
layer $j$ whose first tree ancestor after $x$ lies in $A$.
Let $\mathcal F_j(x)$ reveal the entire first $j$ breadth-first
steps, starting with all edges incident to $x$. If $\mathcal U_j(x)$
is the unvisited opposite class, then, conditionally on
$\mathcal F_{j-1}(x)$,
\begin{align}
 Z_j^x(y)&=\sum_{z\in\mathcal U_j(x)}\indi{\{y,z\}\in\gE}
       \sim\Bin(|\mathcal U_j(x)|,p),\nonumber\\
 \{Z_j^x(y):y\in\layer_{j-1}^{\gG}(x)\}
 &\text{ are independent},\nonumber\\
 |\layer_j^{\gG}(x)|
 &\sim\Bin\bigl(|\mathcal U_j(x)|,
           1-(1-p)^{|\layer_{j-1}^{\gG}(x)|}\bigr).
 \label{eqn:bfs-conditional-laws}
\end{align}
Set $\overline d_j(x)$ equal to $d_1$ when the frontier
$\layer_{j-1}^{\gG}(x)$ lies in $\gV_1$, and to $d_2$ otherwise.
Thus, for $j\geq2$,
\begin{align}
 \overline d_j(x)
 ={}&d_1\bigl(\indi{x\in\gV_1,\ j\ {\rm odd}}
              +\indi{x\in\gV_2,\ j\ {\rm even}}\bigr)\nonumber\\
 &+d_2\bigl(\indi{x\in\gV_1,\ j\ {\rm even}}
              +\indi{x\in\gV_2,\ j\ {\rm odd}}\bigr).
 \label{eqn:bfs-branching-means-unpruned}
\end{align}

\begin{lemma}[Unpruned breadth-first estimates]
\label{lem:semilocalized-unpruned-bfs}
Fix $c_0>0$ and $\nu>0$. With probability at least $1-N^{-\nu}$,
simultaneously for all roots and $2\leq j\leq R_0$, the following
hold. Every subset $A\subset\layer_1^{\gG}(x)$ satisfies
\begin{align}
 |W_j^x(A)|\leq C_{B,\qr,\nu}
       (d|A|+L_N)d^{j-2}.
 \label{eqn:bfs-aggregate-branch-growth}
\end{align}
If $\rD_x\geq c_0d$, then
\begin{align}
 \left|\frac{|\layer_j^{\gG}(x)|}
                    {|\layer_{j-1}^{\gG}(x)|}
                -\overline d_j(x)\right|
 &\leq C_{B,\qr,c_0,\nu}
          \sqrt{\frac{dL_N}{|\layer_{j-1}^{\gG}(x)|}},
 \label{eqn:unpruned-layer-ratio-concentration}\\
 |\layer_j^{\gG}(x)|
 &=\rD_xd^{j-1}\vartheta_j(x)
       \left(1+O_{B,\qr,c_0,\nu}
          \left(\sqrt{\frac{L_N}{d\rD_x}}\right)\right),
 \label{eqn:unpruned-layer-product-concentration}
\end{align}
where $\vartheta_j=1$ for odd $j$ and
$\vartheta_j=(d_2/d)\indi{x\in\gV_1}
 +(d_1/d)\indi{x\in\gV_2}$ for even $j$.
For any fixed $K>0$, with
$\zeta_j=p|\mathcal U_j(x)|$, the same event satisfies
\begin{align}
 \sum_{y\in\layer_{j-1}^{\gG}(x)}
   |Z_j^x(y)-\zeta_j|^2\indi{Z_j^x(y)\leq Kd}
 \leq C_{B,\qr,\nu,K}d
       \bigl(|\layer_{j-1}^{\gG}(x)|+L_N\bigr).
 \label{eqn:bfs-truncated-square-sum}
\end{align}
These events do not depend on the pruning cutoff.
\end{lemma}

\begin{proof}
By Lemma~\ref{lem:bipartite-max-degree}, using $d\leq BL_N$ and
failure exponent $\nu+11$, choose $D_{\rm cap}=C_{B,\qr,\nu}L_N$
so that $\max_x\rD_x\leq D_{\rm cap}$ with failure at most
$N^{-\nu-10}$ for large $N$. Every $R_0$-ball then has at most
\begin{align}
 V_N=(R_0+1)(1+D_{\rm cap})^{R_0}
     =\exp\bigl(O_{B,\qr,\nu}(R_0\log L_N)\bigr)=N^{o(1)}
 \label{eqn:bfs-ball-size-bootstrap}
\end{align}
vertices. Stop each exploration if this bound is exceeded; this
stopping does not occur on the maximum-degree event.

For a fixed first-neighbour subset $A$, the set
$W_{j-1}^x(A)$ is $\mathcal F_{j-1}(x)$-measurable. Its number of
new tree descendants is at most the number of edges from it to
$\mathcal U_j(x)$, which is conditionally binomial with mean at most
$\overline d_j|W_{j-1}^x(A)|$. Conditional Bernstein, with exponent
$t=C_{B,\qr,\nu}L_N$, gives
\begin{align}
 |W_j^x(A)|\leq
 \overline d_j|W_{j-1}^x(A)|
 +C\sqrt{d|W_{j-1}^x(A)|L_N}+CL_N.
 \label{eqn:bfs-branch-recursion}
\end{align}
The constant in $t$ is chosen after the one in $D_{\rm cap}$. Conditional on
$\mathcal F_1(x)$ there are at most $2^{D_{\rm cap}}$ subsets $A$. Therefore
the union over roots, layers, and all these subsets still has failure
at most $N^{-\nu-8}$.
For $j=2$, \eqref{eqn:bfs-branch-recursion} gives
$|W_2^x(A)|\leq C(d|A|+L_N)$. Starting the deterministic upper
recursion at this value proves \eqref{eqn:bfs-aggregate-branch-growth}:
the successive relative increments are bounded by
$C\sqrt{L_N/(dE_j)}+CL_N/(dE_j)$, where
$E_j\geq c(d|A|+L_N)d^{j-2}$ is the deterministic upper envelope.
Their sum from $j=2$ is $O(d^{-1/2})$. The factors
$\overline d_j/d$ cancel in pairs because $d_1d_2=d^2$; their
partial products are bounded above and below by constants depending
on $\qr$. This proves a constant, rather than an exponential-in-$j$,
loss in \eqref{eqn:bfs-aggregate-branch-growth}.

For the full frontier, the last law in
\eqref{eqn:bfs-conditional-laws} has conditional mean
\begin{align}
 \overline d_j|\layer_{j-1}^{\gG}(x)|
       \bigl(1+O_{\qr}(pV_N+V_N/N)\bigr).
 \label{eqn:bfs-layer-conditional-mean}
\end{align}
For $\rD_x\geq c_0d$, the first relative Bernstein error is
$O(\sqrt{L_N/(d\rD_x)})=o(1)$. Induction gives geometrically
decreasing errors in subsequent generations. Their sum is bounded
by a constant times the first error. The depletion errors are
$N^{-1+o(1)}$ and remain negligible after $R_0$ steps. This proves
\eqref{eqn:unpruned-layer-ratio-concentration} and
\eqref{eqn:unpruned-layer-product-concentration}. A union bound over
the $NR_0=N^{1+o(1)}$ rooted layers makes the statements simultaneous.

Finally, conditional Bernstein for a single forward count implies
\[
 \P(|Z-\zeta_j|\geq t\mid\mathcal F_{j-1})
 \leq2\exp\left(-c_{\qr}\min\{t^2/d,t\}\right).
\]
On the truncation $Z\leq Kd$, one has $t\leq C_{\qr,K}d$.
Consequently
$Y=(Z-\zeta_j)^2\indi{Z\leq Kd}/d$ satisfies
$\P(Y\geq u\mid\mathcal F_{j-1})\leq2e^{-c_{\qr,K}u}$,
and $\E[e^{c'Y}\mid\mathcal F_{j-1}]\leq C_{\qr,K}$ for a
fixed $c'>0$. Multiplying these conditional exponential moments over
the independent frontier counts gives
\[
 \P\left(\sum Y\geq C(n+t)
          \,\middle|\,\mathcal F_{j-1}\right)\leq e^{-t},
 \qquad n=|\layer_{j-1}^{\gG}(x)|.
\]
Take $t=C_{B,\qr,\nu,K}L_N$ and union over rooted layers.
This proves \eqref{eqn:bfs-truncated-square-sum}.
\end{proof}

\paragraph{Exceptional vertices and collision witnesses.}
\begin{lemma}[Local exceptional counts and excess]
\label{lem:semilocalized-pruning-witnesses}
On an event of probability at least $1-N^{-\nu}$, simultaneously
over the deterministic admissible cutoffs,
\begin{align}
 \max_x|\gV^{(\vtau)}\cap
          (\ball_{R_0}^{\gG}(x)\setminus\{x\})|
 &\leq C_{B,\qr,\nu,a}\frac{L_N}{ds^2},
 \label{eqn:pruning-local-exceptional-count}\\
 \max_x\operatorname{exc}
      (\gG[\ball_{R_0}^{\gG}(x)])
 &\leq C_{B,\qr,\nu}.
 \label{eqn:pruning-local-excess}
\end{align}
For the graph $H_s$ consisting of all edges incident to an
exceptional vertex, for every fixed integer $\ell\geq1$,
\begin{align}
 \P\left(\max_{C\in\operatorname{Comp}(H_s)}
                     \operatorname{exc}(C)\geq\ell\right)
 \leq N^{1-\ell+o(1)}.
 \label{eqn:pruning-surplus-witness}
\end{align}
\end{lemma}

\begin{proof}[Proof of Lemma~\ref{lem:semilocalized-pruning-witnesses}]
A degree has binomial mean $d\qr^{-2}+O(d/N)$ or
$d\qr^2+O(d/N)$. The three cutoff displacements have size
at least $c_{\qr}ds$, so the Chernoff bounds
\eqref{eqn:degree-chernoff-upper}--\eqref{eqn:degree-chernoff-lower} prove
\eqref{eqn:pruning-single-vertex-tail}.
Here is the conditional version needed for witnesses. Force the
edges of a deterministic witness graph. Suppose its maximum degree
at a required exceptional vertex is $o(ds)$. If $f_v$ incident
edges have been forced, its conditional degree is exactly
$f_v+\Bin(m-f_v,p)$ or $f_v+\Bin(n-f_v,p)$.
The residual mean changes by $f_vp=o(ds)$ and the forced contribution
is $f_v=o(ds)$, regardless of the total witness size. Applying the
same Chernoff bounds with the residual binomial mean gives
\begin{align}
 \P(v\in\gV^{(\vtau)}\mid\hbox{witness edges present})
 \leq2e^{-c_{\qr}ds^2}.
 \label{eqn:pruning-conditional-degree-tail}
\end{align}
For required vertices in one bipartition class the residual edge
families are disjoint. Apply the product bound to the larger
bipartition subset of any $k$ required exceptional vertices to obtain
the joint bound $2^ke^{-c_{\qr}kds^2/2}$.

If $\ball_{R_0}(x)\setminus\{x\}$ contains $k$ exceptional
vertices, a minimal subtree of the breadth-first tree connecting
them to $x$ has at most $kR_0+1$ vertices and maximum degree at most
$k+1$. Take
$k=\lceil C_{B,\qr,\nu,a}L_N/(ds^2)\rceil$.
By \eqref{eqn:pruning-scale-comparison}, $k=o(ds)$, so
\eqref{eqn:pruning-conditional-degree-tail} applies to every such
witness. For a fixed root, a plane rooted tree with $v$ vertices
has at most $4^v$ shapes, at most $2^v$ markings, and at most
$N^{v-1}$ embeddings. Its forced edges have probability $p^{v-1}$.
Summing over $v\leq kR_0+1$ and using the joint tail gives
\[
 \P\bigl(|\gV^{(\vtau)}\cap(B_{R_0}(x)\setminus\{x\})|
                       \geq k\bigr)
 \leq\exp\bigl(C_{\qr}kR_0\log d-c_{\qr}kds^2\bigr).
\]
The first term in the exponent is $o(kds^2)$. Increasing $C$ and
unioning over roots and deterministic cutoffs proves
\eqref{eqn:pruning-local-exceptional-count}.

For \eqref{eqn:pruning-local-excess}, perform the breadth-first
exploration stopped at the degree and volume bounds $D_{\rm cap},V_N$ in the
preceding proof, and expose all remaining edges of the induced ball.
Collision tests include edges to already discovered vertices during
exploration, including multiple-parent edges, as well as the
remaining edges exposed after the ball is determined. There are at
most $V_N^2$ such tests. At each test an
unexposed edge is Bernoulli with conditional probability $p$; selecting
any $h$ successful tests and taking factorial moments gives
\[
 \P(\operatorname{exc}(\gG[B_{R_0}(x)])\geq h,
                    \max_v\rD_v\leq D_{\rm cap})
 \leq\binom{V_N^2}{h}p^h
 \leq N^{-h+o(1)}
\]
for fixed $h$. Tests which introduce a new vertex are the tree edges;
every additional present edge counts as a collision. Take
$h>\nu+3$ and union over $x$, adding the maximum-degree failure.

For \eqref{eqn:pruning-surplus-witness}, a connected component with
excess at least $\ell$ contains a connected subgraph of excess exactly
$\ell$. Delete its leaves and suppress its degree-two paths. For
$\ell\geq2$ the resulting kernel has at most $2\ell-2$ vertices
and $3\ell-3$ edges; for $\ell=1$ it is a cycle. There are finitely many kernels for fixed $\ell$; only subdivisions
that form simple bipartite graphs are counted, so every forced edge
is distinct.
Every edge of a subdivision witness must have at least one
exceptional endpoint. Along each subdivided kernel edge, group the
internal degree-two vertices into disjoint blocks of four, leaving
at most three ungrouped vertices. The middle edge of each block
requires one of its two endpoints to be exceptional. For $q$ blocks,
there are at most $2^q$ choices of these required endpoints. Their
larger bipartition group has at least $q/2$ members, with independent
residual degree sums after the witness edges are forced. The forced
degree at such vertices is two, and at kernel vertices is $O(\ell)$.
Hence the joint tail bounds the block contribution by
$2^qe^{-c_{\qr}qds^2}$. A block adds four vertices and four
edges, whose embedding and edge costs are at most $(C_{\qr}d)^4$.
The subdivision sum is therefore a geometric series with ratio
\[
 2(C_{\qr}d)^4e^{-c_{\qr}ds^2}=o(1).
\]
The kernel and the ungrouped vertices contribute only $O(\ell)$
vertices. For any connected excess-$\ell$ witness with $v$ vertices,
its embedding and forced-edge cost is
$N^vp^{v-1+\ell}=N^{1-\ell}(C_{\qr}d)^{v-1+\ell}$.
The ungrouped costs are $d^{O(\ell)}=N^{o(1)}$; multiplying the geometric series over the $O(\ell)$ kernel edges
and summing the finitely many kernels proves
\eqref{eqn:pruning-surplus-witness}. For a cycle the same argument
uses the cyclic blocks of four, leaving at most a fixed number of
ungrouped vertices. The conditional tail applies to every subdivision
length because the incident forced degrees remain two at internal
vertices and $O(\ell)$ at kernel vertices.
\end{proof}

\subsubsection{Pruned construction and retained-layer concentration}
\begin{proof}[Proof of Lemma~\ref{lem:existence_pruned_graph}]
Work on the simultaneous events in the two preceding lemmas and the
maximum-degree event. At each exceptional vertex $x$, use its fixed
breadth-first tree of $\gG[\ball_{R_0}(x)]$. Call a first neighbour
$y$ bad if its tree branch contains an exceptional vertex other than
$x$, or contains an endpoint of a nontree edge of this induced ball.
Delete $xy$ for all bad first neighbours at all exceptional roots.
Denote the resulting graph by $\gG^{(\vtau)}$.

A good branch contains no exceptional vertex through depth $R_0$,
and no nontree edge of the induced ball has an endpoint in it. All
its tree edges through depth $R_0$ survive: an edge deleted at another
exceptional vertex in this ball would require such a vertex in the
branch. A path of length at most $R_0$ from $x$ stays in its original
$R_0$-ball. It cannot leave a good branch for a bad one, since that
would use a nontree edge incident to the good branch. Consequently
the retained spheres through depth $R_0$ are exactly the unions of
the good canonical tree branches.

This identity proves the sphere inclusions and the retained
neighbour identity \eqref{eqn:retained-branch-neighbours}. In the
latter, any original neighbour of a retained nonroot vertex in an
original sphere of depth at most $2\radius$ belongs to the same
good branch (or is $x$), so both the edge and the endpoint survive.
A retained path of length at most $4\radius$ joining $x$ to another
exceptional root would lie in a branch declared bad at $x$.
Thus distinct exceptional roots have distance at least $4\radius+1$.
No nontree edge survives in the retained $4\radius$-ball, which is
therefore a tree. This proves items (1), (2), and (4).

At an exceptional root, each nearby exceptional vertex marks at most
one first branch and each nontree edge marks at most two. Thus
\[
 m_x\leq|\gV^{(\vtau)}\cap(B_{R_0}(x)\setminus\{x\})|
            +2\operatorname{exc}(\gG[B_{R_0}(x)]).
\]
An edge to another exceptional first neighbour is already marked at
$x$, so deletion at other roots does not add a further root loss.
At a nonexceptional vertex every deleted neighbour is exceptional,
so its deleted degree is at most the same local exceptional-count
bound. By \eqref{eqn:pruning-local-exceptional-count},
\eqref{eqn:pruning-local-excess}, and $d\leq BL_N$, this gives
$\max_xm_x\leq C L_N/(ds^2)$.

All deleted edges lie in $H_s$. Take
$\ell>\nu+4$ in \eqref{eqn:pruning-surplus-witness} and union
over the deterministic cutoffs. On the resulting event every
component of $\gD^{(\vtau)}$ has excess at most $\ell-1$.
Select a spanning tree of each component and let their union be
$\gF^{(\vtau)}$. The complement $\gU^{(\vtau)}$ has at most
$\ell-1$ edges in each component, so its maximum degree is bounded
by $\ell-1$. The forest is a subgraph of the deleted graph and has
maximum degree at most $C L_N/(ds^2)$. This proves item (3),
including \eqref{eqn:removed-graph-degrees}.

Finally, let $A_x$ be the bad first-neighbour set at an exceptional
root. We proved that $|A_x|=m_x$ and that every missing sphere
vertex through depth $R_0$ is precisely a descendant of $A_x$ in
the original breadth-first tree. Applying
\eqref{eqn:bfs-aggregate-branch-growth} to this graph-dependent
subset gives
$\ell_j(x)\leq C(dm_x+L_N)d^{j-2}$ for $j\geq2$.
The first-layer identity is $\ell_1=m_x$. This proves item (5).
The maximum-degree event proves
\eqref{eqn:pruning-crude-ball-size}; the deterministic-cutoff
intersection proves simultaneity in the gap.
\end{proof}

\paragraph{Retained layers and child-count second moments.}
\begin{lemma}[Breadth-first concentration]
\label{lem:semilocalized-layer-concentration}
Fix $c_0>0$, $a>0$, and $\nu>0$. With probability at least
$1-N^{-\nu}$, simultaneously over all cutoffs with
$s(\gap)\geq a\error_d^{1/4}$ and all exceptional roots with
$\rD_x\geq c_0d$, the retained layers through $\radius+1$ are
nonempty and
\begin{align}
 c_{\qr,c_0}\rD_xd^{j-1}
 \leq|\layer_j^{(\vtau)}(x)|
 \leq C_{B,\qr,c_0,\nu,a}\rD_xd^{j-1},
 \qquad1\leq j\leq\radius+1.
 \label{eqn:pruned-layer-two-sided-growth}
\end{align}
With $\overline d_j$ as in
\eqref{eqn:bfs-branching-means-unpruned}, define
\begin{align}
 \sigma_x^{(\vtau)}=
       \frac{m_x}{\rD_x}+\frac{L_N}{d\rD_x}.
 \label{eqn:bfs-branching-means}
\end{align}
For $2\leq j\leq\radius+1$,
\begin{align}
 \left|\frac{|\layer_j^{(\vtau)}(x)|}
                    {|\layer_{j-1}^{(\vtau)}(x)|}
                 -\overline d_j(x)\right|
 \leq C\left(
       \sqrt{\frac{dL_N}{|\layer_{j-1}^{(\vtau)}(x)|}}
       +d\sigma_x^{(\vtau)}\right),
 \label{eqn:pruned-layer-ratio-concentration}
\end{align}
and
\begin{align}
 &\frac1{|\layer_{j-1}^{(\vtau)}(x)|}
 \sum_{y\in\layer_{j-1}^{(\vtau)}(x)}
 \left(|\layer_1^{(\vtau)}(y)\cap\layer_j^{(\vtau)}(x)|
            -\frac{|\layer_j^{(\vtau)}(x)|}
                        {|\layer_{j-1}^{(\vtau)}(x)|}\right)^2
\leq C d\left(1+\frac{L_N}{\rD_x}\right)
              +C d^2\bigl(\sigma_x^{(\vtau)}\bigr)^2.
 \label{eqn:pruned-child-square-concentration}
\end{align}
The weaker right-hand side obtained by multiplying the first term by
$\log d$ may therefore be used in the later profile calculation.
Every profile and its residual have support of size $N^{o(1)}$ by
\eqref{eqn:pruning-crude-ball-size}.
\end{lemma}

\begin{proof}[Proof of Lemma~\ref{lem:semilocalized-layer-concentration}]
Extract the unpruned event of
Lemma~\ref{lem:semilocalized-unpruned-bfs} before selecting any
pruning. Equations \eqref{eqn:pruned-layer-loss} and
\eqref{eqn:unpruned-layer-product-concentration} imply
\begin{align}
 \frac{\ell_j(x)}{|\layer_j^{\gG}(x)|}
 \leq C\sigma_x^{(\vtau)},\qquad1\leq j\leq\radius+1.
 \label{eqn:bfs-relative-pruned-loss}
\end{align}
The degree bound gives
\begin{align}
 \sigma_x^{(\vtau)}
 \leq C\frac{L_N}{d^2s^2}
       +C\frac{L_N}{d^2}
 \leq C_a\frac{\xi^{3/2}}{(\log d)^2}=o(1).
 \label{eqn:bfs-pruned-loss-small}
\end{align}
Here $\xi$ need not tend to zero at the lowest sparsity; its uniform
boundedness and $\log d\to\infty$ suffice. This proves
\eqref{eqn:pruned-layer-two-sided-growth}, including nonemptiness.
Comparing the full and retained sphere ratios adds at most
$Cd\sigma_x^{(\vtau)}$, proving
\eqref{eqn:pruned-layer-ratio-concentration}.

For a retained parent $y$ and $j\geq2$, separation implies that
$y$ is nonexceptional. Its degree is at most $K_{\qr}d$.
The exact neighbour identity gives
\begin{align}
 |\layer_1^{(\vtau)}(y)\cap\layer_j^{(\vtau)}(x)|
       =Z_j^x(y).
 \label{eqn:bfs-retained-child-identity}
\end{align}
Here $Z_j^x(y)$ counts all original forward edges to unvisited
vertices, rather than only the chosen canonical tree children.
For a retained parent, the absence of a nontree endpoint makes these
counts equal. Thus the sum of squared deviations from $\zeta_j$ over the retained
parents is bounded by the truncated sum in
\eqref{eqn:bfs-truncated-square-sum}. Divide by the retained parent
count, which is comparable to the full count, and use
$|\layer_{j-1}^{\gG}(x)|\geq c\rD_x$.
The difference between $\zeta_j$ and the retained ratio is at most
the right-hand side of
\eqref{eqn:pruned-layer-ratio-concentration}, plus the negligible
depletion term. The inequality $(u+v)^2\leq2u^2+2v^2$ now proves
\eqref{eqn:pruned-child-square-concentration}. This transfer uses
the unconditional unpruned square sum, so no conditional independence
is assumed after the graph-dependent selection of retained parents.
\end{proof}

\subsubsection{Norm of the pruning error}
\begin{proof}[Proof of Lemma~\ref{lem:semilocalized-pruning-input}]
The forest norm bound and the maximum-row-sum bound for the bounded
degree remainder give
\begin{align}
 \frac1{\sqrt d}\|\rmA-\rmA^{(\vtau)}\|
 &\leq\frac2{\sqrt d}\sqrt{\frac{CL_N}{ds^2}}
                  +\frac C{\sqrt d}
 \leq C\left(\frac{\sqrt{L_N}}{ds}+\xi\right)
 \leq C(\xi+\error_{\mathrm{prof}}(\gap)).
 \label{eqn:adjacency-pruning-norm}
\end{align}
Put $\mathcal U=\bigcup_{x\in\gV^{(\vtau)}}
\ball_{2\radius}^{(\vtau)}(x)$.
The one-vertex estimate \eqref{eqn:pruning-single-vertex-tail} bounds
the mean exceptional count on side $i$ by
$2|\gV_i|e^{-c_{\qr}ds^2}$. Choose
$0<c_{B,\qr}<\min\{c_{\qr}/2,1/(4B)\}$.
Proposition~\ref{prop:uniform-degree-counts}, applied to the defining
degree tails with a sufficiently large failure exponent, gives
$|\gV^{(\vtau)}\cap\gV_i|
\leq|\gV_i|e^{-c_{B,\qr}ds^2}$.
The smaller exponential rate absorbs the constants multiplying the
mean, and the threshold $c_{\qr}N^{3/4}$ absorbs the additive $\log N$.
The crude ball-size estimate therefore gives
\begin{align}
 \frac{|\mathcal U|}{N}
 \leq \exp\bigl(-c_{B,\qr}ds^2+\log V_N\bigr)
 \leq e^{-c_{B,\qr}ds^2/2}.
 \label{eqn:pruned-union-size}
\end{align}
Writing $P_{\mathcal U}=\id_N-\chi^{(\vtau)}$, the constant off-diagonal
entries of $\E\rmA$ imply
\begin{align}
 \|\E\rmA-\chi^{(\vtau)}(\E\rmA)\chi^{(\vtau)}\|
 \leq2\|P_{\mathcal U}\E\rmA\|
 \leq C_{\qr}p\sqrt{N|\mathcal U|}
 \leq C_{\qr}d e^{-c_{B,\qr}ds^2/4}.
 \label{eqn:centering-pruning-norm}
\end{align}
After division by $\sqrt d$ this is $o(\xi)$ because
$ds^2\gg\log d$. The intersection over deterministic cutoff triples
preserves the arbitrary prescribed right failure exponent by taking
larger input exponents. The constants and $\gC_0$ may depend on
$(B,\qr,\nu,a)$.

Under the standing left-edge assumptions, the giant comparison
\eqref{eqn:giant-component-centered-comparison} adds
$O_{b,\qr}(\xi)$. Intersect its event with the simultaneous pruning
events and retain the prescribed fixed left failure exponent.
This gives the stated probability and constants depending on
$(b,\qr,a)$.
\end{proof}

\subsection{Radial profiles and block approximation}
\label{subsec:semilocalized-profile-proofs}

\subsubsection{Selected centres and signed profiles}
Select a set $\mathcal C\subset\gV^{(\vtau)}$ of positive-branch
centres: $x\in\gR_1\cup\gR_2$, or $x\in\gV_2$ with
$0<\alpha_x<\qr^2-1$. Define the local residual relative to this
selection by
\begin{align}
 \varepsilon_{x,\mathcal C}
 \coloneqq\max_{\sigma=\pm}
   \|\bigl(\scA^{(\vtau)}-\sigma\Lambda_x\bigr)
          \hatvtau_\sigma(x)\|,
 \qquad x\in\mathcal C.
 \label{eqn:local-profile-error}
\end{align}
For the window sets, write
$\varepsilon_{x,\mathrm E}=\varepsilon_{x,\gV_{\mathrm E}^\star}$
for $x\in\gV_{\mathrm E}^\star$, $\mathrm E\in\{\mathrm L,\mathrm R\}$.
Maxima over an empty selected set are zero.

\begin{lemma}[Selected-profile estimates]
\label{lem:semilocalized-selected-profile-estimates}
Let $s=s(\gap)\geq a\xi^{1/4}$ and work on the common pruning and
retained-layer event, with $\xi=\error_d\leq1$. Suppose, for fixed
$c_0,c_1>0$, that every $x\in\mathcal C$ satisfies
$\rD_x\geq c_0d$ and $1-|q_x|\geq c_1\xi^{1/2}$.
The selected signed profiles are orthonormal. Their residuals are
supported in $\ball_{\radius+1}^{(\vtau)}(x)$, with disjoint
residual supports between distinct centres, and satisfy
\begin{align}
 \varepsilon_{x,\mathcal C}^2
 \leq{}&C\frac{\log d}{d}\left(1+\frac{L_N}{\rD_x}\right)
       +C\frac{m_x^2}{d\rD_x}
       +C\frac{L_N^2}{d^2\rD_x^2}
       +Ce^{-c_{\qr,c_1}\xi^{1/2}\radius},
 \label{eqn:weighted-local-profile-error}
\end{align}
where $L_N=\log N$ and $m_x=\rD_x^{\gD^{(\vtau)}}$.
In particular,
\begin{align}
 \xi^{-1}\max_{x\in\mathcal C}
       (1+\alpha_x)\varepsilon_{x,\mathcal C}^2
 &\leq C_{B,\qr,\nu,a,c_0,c_1}\xi
 \leq C_{B,\qr,\nu,a,c_0,c_1}
             \bigl(\xi+\error_{\mathrm{prof}}(\gap)\bigr),
 \label{eqn:pruning-weighted-profile-acceptance}
\end{align}
and $\max_{x\in\mathcal C}\varepsilon_{x,\mathcal C}\leq C\xi$.
The tail bound \eqref{eqn:profile-tail-bound} holds for
$0\leq r\leq\radius$ and for either sign. Every layer up to depth
$\radius$ is nonempty, and the quantitative support bound
\eqref{eqn:profile-support-size} holds with the same parameter
dependence as the residual estimate.
\end{lemma}

\begin{proof}[Proof of Lemma~\ref{lem:semilocalized-selected-profile-estimates}]
\emph{Parity balance and orthogonality.}
For $x\in\mathcal C$, set $\rho_x=q_x^2\in(0,1)$.
Substituting \eqref{eqn:signed-profile-heuristic} gives
\begin{align}
 \ervu_1^2-\ervu_2^2=(1-\rho_x)\ervu_0^2.
 \label{eqn:profile-first-pair-identity}
\end{align}
The geometric tail therefore gives, for every $k\geq0$,
\begin{align}
 \sum_{j=0}^{2k}(-1)^j\ervu_j^2
 =\ervu_0^2\rho_x^k>0.
 \label{eqn:profile-alternating-partial-sum}
\end{align}
For $J=2k+1$, this proves positivity of the terminal square root in
\eqref{eqn:explicit-terminal-coefficient} and the exact parity
balance \eqref{eqn:profile-parity-balance}. Its size is
$O_\qr(\rho_x^{k/2})$, the geometric order of the omitted tail.
The two signs satisfy
\[
 \langle\hatvtau_+(x),\hatvtau_-(x)\rangle
 =\sum_j(-1)^ja_j(x)^2=0.
\]
Disjoint balls separate different centres, so the signed family is
orthonormal.

\paragraph{Residual decomposition.}
Fix $x\in\mathcal C$ and abbreviate
\begin{align}
 \layer_j&\coloneqq\layer_j^{(\vtau)}(x),
 &\rvs_j&\coloneqq|\layer_j|^{-1/2}\ones_{\layer_j},
 &\rvv_{x,r}&\coloneqq\sum_{j=0}^r\ervu_j\rvs_j,
       \qquad0\leq r\leq\radius.
 \label{eqn:unbalanced-radial-profile}
\end{align}
Write
\begin{align}
 \rN_j(y)&\coloneqq|\layer_j\cap\layer_1^{(\vtau)}(y)|,
 &t_j&\coloneqq\sqrt{\frac{|\layer_j|}{d|\layer_{j-1}|}},
       \qquad j\geq1.
 \label{eqn:radial-edge-count}
\end{align}
The centering term vanishes on the profile ball, and each nonroot
vertex has one parent. Grouping by layer and using the ideal
recurrence \eqref{eqn:ideal-radial-recurrence} gives
\begin{align}
 (\scA^{(\vtau)}-\Lambda_x)\rvv_{x,r}
 =\rvw_{\rm ch}+\rvw_{\rm gr}+\rvw_{\rm bd},
 \label{eqn:radial-error-decomposition}
\end{align}
where
\begin{subequations}\label{eqn:radial-error-terms}
\begin{align}
 \rvw_{\rm ch}
 &\coloneqq\frac1{\sqrt d}\sum_{j=1}^r
    \frac{\ervu_j}{\sqrt{|\layer_j|}}
    \sum_{y\in\layer_{j-1}}
    \left(\rN_j(y)-\frac{|\layer_j|}{|\layer_{j-1}|}\right)\ones_y,
 \label{eqn:radial-error-children}\\
 \rvw_{\rm gr}
 &\coloneqq\sum_{j=1}^r(t_j-\omega_j)
      (\ervu_j\rvs_{j-1}+\ervu_{j-1}\rvs_j),
 \label{eqn:radial-error-growth}\\
 \rvw_{\rm bd}
 &\coloneqq-\omega_{r+1}\ervu_{r+1}\rvs_r
                 +t_{r+1}\ervu_r\rvs_{r+1}.
 \label{eqn:radial-error-boundary}
\end{align}
\end{subequations}
The three terms are child-count fluctuations, layer-growth errors,
and the outer boundary. The growth term at $j=1$ records the loss of
root branches. This is the tree specialization of
\cite[Lemma~4.2]{dumitriu2025singular}; the identity applies to
both bipartition classes, and the negative-profile identity follows
by the bipartite sign involution.

Apply the decomposition at $r=J$ to the normalized truncation of
$\sum_{j=0}^J\ervu_j\rvs_j$. The child term at $j=1$ vanishes.
The retained child-count estimate and layer comparability give
\begin{align}
 \|\rvw_{\rm ch}\|^2
 \leq C\frac{\log d}{d}
       \left(1+\frac{L_N}{\rD_x}\right)
       +C(\sigma_x^{(\vtau)})^2.
 \label{eqn:profile-child-error-from-bfs}
\end{align}
For the growth term, $t_1=\sqrt{(\rD_x-m_x)/d}$ and
\[
 |t_1-\sqrt{\alpha_x}|^2
 =\frac{m_x^2}{d(\sqrt{\rD_x-m_x}+\sqrt{\rD_x})^2}
 \leq C\frac{m_x^2}{d\rD_x}.
\]
For $j\geq2$, \eqref{eqn:pruned-layer-ratio-concentration} and the
square-root inequality give
$|t_j-\omega_j|\leq C\sqrt{L_N/(d|\layer_{j-1}|)}
+C\sigma_x^{(\vtau)}+O(N^{-1})$.
The normalized ideal coefficients have squared sum at most one,
and each contributes to at most two consecutive layers. Hence
\begin{align}
 \|\rvw_{\rm gr}\|^2
 &\leq C\frac{L_N}{d\rD_x}+C(\sigma_x^{(\vtau)})^2
                +C\frac{m_x^2}{d\rD_x}\nonumber\\
 &\leq C\frac{L_N}{d\rD_x}+C\frac{m_x^2}{d\rD_x}
                +C\frac{L_N^2}{d^2\rD_x^2},
 \label{eqn:profile-growth-error-from-bfs}
\end{align}
where the last bound uses $\rD_x\geq c_0d$. Thus
\begin{align}
 \|\rvw_{\rm ch}\|^2+\|\rvw_{\rm gr}\|^2
 \leq C\frac{\log d}{d}\left(1+\frac{L_N}{\rD_x}\right)
       +C\frac{m_x^2}{d\rD_x}
       +C\frac{L_N^2}{d^2\rD_x^2}.
 \label{eqn:radial-profile-error}
\end{align}

\paragraph{Terminal control.}
The boundary formula and geometric decay give
\begin{align}
 \|\rvw_{\rm bd}\|\leq C_\qr\rho_x^{(J-1)/4}.
 \label{eqn:radial-boundary-error}
\end{align}
The terminal adjustment in \eqref{eqn:explicit-terminal-coefficient}
changes the normalized profile by at most the same order.
Lemma~\ref{lem:tree-norm} and the maximum-degree bound control its
image under $\scA^{(\vtau)}-\sigma\Lambda_x$ by a polynomial in
$L_N$. Consequently
\begin{align}
 \varepsilon_{x,\mathcal C}
 \leq C\bigl(\|\rvw_{\rm ch}\|+\|\rvw_{\rm gr}\|\bigr)
                +CL_N^C\rho_x^{(J-1)/4}.
 \label{eqn:profile-residual-before-tail}
\end{align}
The decay assumption and
\begin{align}
 \xi^{1/2}\radius\geq c_\qr L_N^{1/4}
 \label{eqn:pruning-profile-tail-scale}
\end{align}
absorb the polynomial factor into the exponential:
\begin{align}
 L_N^C\rho_x^{(J-1)/4}
 \leq e^{-c_{\qr,c_1}\xi^{1/2}\radius}\leq\xi,
 \label{eqn:profile-terminal-negligible}
\end{align}
uniformly over the sparsity range. Squaring
\eqref{eqn:profile-residual-before-tail} and using
\eqref{eqn:radial-profile-error} proves
\eqref{eqn:weighted-local-profile-error}. The decomposition also
places the residual in $\ball_{\radius+1}^{(\vtau)}(x)$.
Distinct retained roots have distance at least $4\radius+1$, so
these residual supports are disjoint for large $N$.

\paragraph{Weighted acceptance.}
Use $\rD_x\geq c_0d$, $\max_x\rD_x\leq CL_N$, and
$m_x\leq CL_N/(ds^2)$. The algebraic contributions in
\eqref{eqn:weighted-local-profile-error} satisfy
\begin{align}
 \xi^{-1}\left(1+\frac{\rD_x}{d}\right)
   \frac{\log d}{d}\left(1+\frac{L_N}{\rD_x}\right)
 &\leq C\frac{\sqrt{L_N}}d=C\frac{\xi}{\log d},\nonumber\\
 \xi^{-1}\left(1+\frac{\rD_x}{d}\right)
   \frac{m_x^2}{d\rD_x}
 &\leq C\frac{\xi^3}{(\log d)^4s^4}
       \leq C_a\frac{\xi^2}{(\log d)^4},\nonumber\\
 \xi^{-1}\left(1+\frac{\rD_x}{d}\right)
   \frac{L_N^2}{d^2\rD_x^2}
 &\leq C\frac{\xi^3}{(\log d)^4}.
 \label{eqn:pruning-weighted-residual-budget}
\end{align}
The weighted exponential is $o(\xi)$ uniformly in $d$, by
\eqref{eqn:pruning-profile-tail-scale}. These estimates prove
\eqref{eqn:pruning-weighted-profile-acceptance}; since
$1+\alpha_x\geq1$, they also give
$\varepsilon_{x,\mathcal C}\leq C\xi$.

\paragraph{Spatial tail and support.}
Put $g_x=|q_x|^{-1}>1$. For even $r=2h$, summing the geometric
series in \eqref{eqn:Z1Eigenvec} and \eqref{eqn:Z2Eigenvec},
including the adjusted terminal term, bounds the normalized tail
by $g_x^{-2h}$. For odd $r=2h+1$, the first omitted even
coefficient satisfies
\[
 \ervu_{2h+2}^2\leq C_\qr g_x^{-2h-1}
                          \sum_j\widetilde a_j(x)^2.
\]
The remaining coefficients form a geometric series of ratio
$g_x^{-2}$, so the tail is at most $C_\qr g_x^{-r}$ in both
cases. This proves \eqref{eqn:profile-tail-bound}; the signs have
the same coordinate magnitudes. The two-sided retained-layer
bound gives nonemptiness and
$|\ball_{\radius}^{(\vtau)}(x)|\leq C\rD_xd^{\radius}$.
Its logarithm is $O(\sqrt{d\log d}+\log\log N)$ because
$\max_x\rD_x\leq C\log N$, proving
\eqref{eqn:profile-support-size}.
\end{proof}

\begin{proof}[Proof of Lemma~\ref{lem:semilocalized-profile-properties}]
Choose the simultaneous pruning and retained-layer events before
selecting the gap. For the right window centres,
$\rD_x\geq c_\qr d$; the right profile sets and the
degree--eigenvalue separation \eqref{eqn:profile-energy-degree-separation}
(with its first-side analogue) give
$1-|q_x|\geq c_\qr\xi^{1/8}$. For the left window centres,
the positive lower cutoff in \eqref{eqn:V2-left-profile-centres}
and monotonicity of $\invL$ give
$\rD_x\geq c_{\qr,\sK}d$, while the upper cutoff
$\Lambda_x\leq\ledge-2\xi^{1/4}$ gives the same decay bound.
After choosing $\gC_0$ so that $\xi\leq1$, both sets satisfy the
weaker hypothesis $1-|q_x|\geq c_\qr\xi^{1/2}$.
Apply Lemma~\ref{lem:semilocalized-selected-profile-estimates} to
$\mathcal C=\gV_{\mathrm R}^\star$ and
$\mathcal C=\gV_{\mathrm L}^\star$.
Take the right input failure exponent larger than the prescribed
$\nu$ before intersection. At the left edge, intersect with the
fixed-exponent giant event, retaining probability
$1-CN^{-\nu_0}$ and constants depending on $(b,\qr,\sK)$.
\end{proof}

\begin{samepage}
\paragraph{Block approximation.}
\label{subsec:semilocalized-block-proof}

\begin{proof}[Proof of
Proposition~\ref{prop:semilocalized-block-approximation}]
Let $U$ have the selected signed profiles as columns, and let
$D_{\mathcal C}$ be diagonal with entries $\sigma\Lambda_x$.
The residual operator is
$R=\scA^{(\vtau)}U-UD_{\mathcal C}$.
Residual supports are disjoint between centres; at each centre there
are only two signs. Thus, for every coefficient vector $a$,
\[
    \|Ra\|^2\leq2\varepsilon_{\mathrm{res}}^2\|a\|^2.
\]
Write $\Pi=UU^*$ and $Q=\id_N-\Pi$. Since $U^*U=\id$,
\begin{align}
    \scA^{(\vtau)}-\hatH_{\mathcal C}
    &=RU^*+UR^*-U(U^*R)U^*,\\
    \Pi\scA^{(\vtau)}Q&=UR^*Q.
\end{align}
Their norms are at most $3\|R\|$ and $\|R\|$, respectively.
This proves \eqref{eqn:pruned-profile-block-error}; adding
$\|\scA-\scA^{(\vtau)}\|$ proves
\eqref{eqn:full-profile-block-error}.
\end{proof}
\end{samepage}\subsection{Nonbacktracking specialization and regularized forms}
\label{subsec:semilocalized-loewner-proofs}

\begin{proof}[Proof of Proposition~\ref{prop:semilocalized-nb-radius-input}]
Apply the bipartite clause of the nonbacktracking-radius proposition
in \cite{wang2026regularized}. Center the Bernoulli matrix by its mean
before applying either projection.
The source sparsity parameter is $p\sqrt{nm}=d$, its uniform
sparsity constant is one, and its variance sums satisfy
\[
 \overline v_{\rm r}=mp(1-p)/d,\qquad
 \overline v_{\rm c}=np(1-p)/d,\qquad
 (\overline v_{\rm r}\overline v_{\rm c})^{1/4}
 =(1-p)^{1/2}\leq1.
\]
Apply its degree-cutoff algorithm to the uncentered Bernoulli matrix
with $L_{\rm r}=L_{\rm c}=0$,
$U_{\rm r}=\lceil\tau_1d\rceil-1$, and $U_{\rm c}=n$.
The column cutoff $U_{\rm c}=n$ leaves all columns in place, so the
padded matrix is $\mP^{(\tau_1)}\mX$, with centering preserved on
each retained row.

Write $C_{\rm prof,NB}$ for the fixed degree-profile constant in
the source cutoff condition, distinguishing it from $\gC_0$.
The required margins are
\[
 U_{\rm r}-mp=\frac{ds}{8\qr}+O_\Gamma(1+d/N),
 \qquad U_{\rm c}-np=n(1-p).
\]
Uniform admissibility gives $ds/\sqrt{d\log d}\to\infty$;
both margins therefore exceed
$5\sqrt{C_{\rm prof,NB}d\log d}$ for sufficiently large $N$.
The remaining conditions $d\geq C_{\rm NB}$ and
$d^6(\log N)^3\leq N^{1-\delta_{\rm NB}}$, for a fixed
$0<\delta_{\rm NB}<1$, hold throughout $d\leq B\log N$.
The source radius estimate is uniform over these admissible cutoffs.
There are $O_\Gamma(d)$ possible integer row cutoffs; apply it with
failure exponent $\nu+2$ and take their union. This produces one
event of probability at least $1-N^{-\nu}$ on which the bound holds
for every gap in the deterministic family. At critical sparsity take
$B=2b$. Zero padding preserves the nonzero nonbacktracking spectrum,
so its host-edge operator agrees with the one defined here.
\end{proof}

\begin{proof}[Proof of Proposition~\ref{prop:semilocalized-regularized-loewner-bounds}]
Let $\mQ_{\rm r}$ and $\mQ_{\rm c}$ be the
diagonal matrices of row and column squared norms of
$\mX^{(\tau_1)}$, and put
$\mQ=\operatorname{diag}(\mQ_{\rm r},\mQ_{\rm c})$ and
$L_{\rm ent}=\max_{x,y}|\emX^{(\tau_1)}_{xy}|$.
Separating present and absent Bernoulli edges gives
\begin{align}
    (\mQ_{\rm r})_{xx}
    &=\frac{\indi{x\in I_1}}{d}
      \bigl((1-2p)\rD_x+mp^2\bigr),
    \nonumber\\
    (\mQ_{\rm c})_{yy}
    &=\frac{(1-2p)\emD^{(2,\tau_1)}_{yy}+|I_1|p^2}{d}.
    \label{eqn:semilocalized-original-square-sums}
\end{align}
Since $p=O_\qr(d/N)$, these identities give
$\mQ\preceq d^{-1}\mD^{(\tau_1)}+C_\qr p\id_N$ and
$L_{\rm ent}\leq d^{-1/2}$. By \cite[Proposition~2.4]{wang2026regularized},
for $t>\max\{\rho(\mB^{(\tau_1)}),L_{\rm ent}\}$,
\begin{align}
    \mH^{(\tau_1)}
    \preceq t\id_N+t^{-1}\mQ+\frac{L_{\rm ent}}{t(t-L_{\rm ent})}\mQ.
    \label{eqn:semilocalized-imported-upper-form}
\end{align}
On the event of Proposition~\ref{prop:semilocalized-nb-radius-input}, take
$t=1+(C_{B,\qr,\nu}+1)d^{-1/2}$. Since $p=o(d^{-1/2})$, this
proves \eqref{eqn:upper_bound_H} and
\eqref{eqn:upper-bound-H-loewner-error}, simultaneously for every gap
in the admissible family, with probability at least $1-N^{-\nu}$.

For the left clause, take $B=2b$ in the nonbacktracking input.
At critical sparsity, Lemma~\ref{lem:bipartite-max-degree} gives
$\|\mD^{(\tau_1)}\|/d=O_{b,\qr,\nu}(1)$ with any prescribed
polynomial failure exponent. The square-sum identities give
\[
    \|\mQ_{\rm r}-d^{-1}\mD^{(1,\tau_1)}\|
    +\|\mQ_{\rm c}-d^{-1}\mD^{(2,\tau_1)}\|
    =O_{b,\qr,\nu}(p).
\]
Together with \eqref{eqn:good-row-margin} and
\eqref{eqn:semilocalized-nb-input}, these bounds verify the hypotheses
of \cite[Proposition~2.5]{wang2026regularized} with $t=1$,
$u=C_{b,\qr,\nu}d^{-1/2}$,
$\mW_{\rm r}=d^{-1}\mD^{(1,\tau_1)}$, and
$\mW_{\rm c}=d^{-1}\mD^{(2,\tau_1)}$.
Indeed, the row inverse has a fixed positive margin, and its
square-sum error is $O(p)=o(1)$.
Moreover,
$(\id_n-\mW_{\rm r})^{-1}=\mM^{(\tau_1)}+\id_n-\mP^{(\tau_1)}$,
and $\id_n-\mP^{(\tau_1)}$ annihilates $\mX^{(\tau_1)}$.
The cited lower comparison therefore yields
\eqref{eqn:lower_bound_H}. Increasing the failure exponents before
intersecting the radius and maximum-degree events gives probability
at least $1-N^{-\nu}$, simultaneously over the admissible family
of left-edge gaps.
\end{proof}

\begin{proof}[Proof of Lemma~\ref{lem:semilocalized-column-degree-loss}]
Fix a gap in the family and put $a=ds^2$.
Conditioning on the single Bernoulli edge $\emA_{xy}$ gives
\begin{align}
    \P(\emA_{xy}=1,\ \rD_x\geq U_1)
    &=p\,\P\bigl(\Bin(m-1,p)\geq U_1-1\bigr)
    \leq C_\qr p\exp(-c_\qr a).
\end{align}
Indeed, the row-cutoff margin verified in the nonbacktracking proof
gives a deviation of order $ds$ above the conditional mean. The
Chernoff bound \eqref{eqn:degree-chernoff-upper} has exponent of order
$ds^2$, since $s\leq1$.
For fixed $y$, the summands defining $\ell_y(s)$ are independent over
$x$, so $\ell_y(s)$ is binomial with mean
$\lambda_s\leq C_\qr d\exp(-c_\qr a)$.
Admissibility gives $a/\log d\to\infty$, uniformly over the family,
and hence $\lambda_s\leq\exp(-c_\qr a/2)$ for all large $N$.
Choose
$h=\lceil C_{B,\qr,\nu}(1+\log N/a)\rceil$.
If $\lambda_s=0$, the loss count vanishes almost surely. Otherwise,
the Chernoff bound \eqref{eqn:degree-chernoff-upper} yields
\begin{align}
    \P(\ell_y(s)\geq h)
    \leq\left(\frac{\mathrm e\lambda_s}{h}\right)^h
    \leq\exp(-c_\qr ah/3)
    \leq N^{-\nu-3}
\end{align}
after increasing the constant in $h$.

As the gap varies, the discarded-row event changes only when
$U_1=\lceil\tau_1d\rceil$ changes, giving $O(d)$ possibilities.
For each attained cutoff $k$, use the supremum $s_k\in(0,1]$ producing it, with a one-sided
limit when that endpoint is not attained. Rounding changes the
conditional deviation by $O(1)=o(ds_k)$, the same tail estimate applies to $s_k$, which is at least as large as
any $s$ with that cutoff, and
$1+\log N/(ds_k^2)\leq1+\log N/(ds^2)$.
A union bound over columns and cutoffs proves
\eqref{eqn:regularized-incidence}. Finally,
$\emD^{(2)}_{yy}-\emD^{(2,\tau_1)}_{yy}=\ell_y(s)$, which gives
\eqref{eqn:regularized-degree-loss}.
\end{proof}

\subsection{Local spectra and identification of root eigenvalues}
\label{subsec:semilocalized-local-spectra}

Root identification uses an adjacency envelope on the right and a
squared-operator Schur comparison on the left, with separate degree
hypotheses on positive windows and through zero.

\subsubsection{Finite-forest forms and local operators}

\begin{lemma}[Finite-forest quadratic forms]
\label{lem:semilocalized-finite-forest-forms}
Let $\rmT$ be the adjacency matrix of a finite bipartite forest with
edge weights $d^{-1/2}$, where $d>1$, and let $\rmD_{\rmT}$ be its degree
matrix. Write $\rmT_{12}$ for its $\gV_1\times\gV_2$ block.
Then
\begin{align}
    \rmT\preceq
    (1-d^{-1})\id+d^{-1}\rmD_{\rmT}.
    \label{eqn:finite-tree-upper-form}
\end{align}
If $\rmD_{\rmT}^{(1)}\prec d\id$, then
\begin{align}
    \rmT_{12}^*
      \left(\id-\frac{\rmD_{\rmT}^{(1)}}d\right)^{-1}
      \rmT_{12}
    \succeq
      \frac{\rmD_{\rmT}^{(2)}}d-\id
      -\frac{\rmD_{\rmT}^{(2)}}{d(d+1)}.
    \label{eqn:finite-forest-lower-form}
\end{align}
\end{lemma}

\begin{proof}[Proof of Lemma~\ref{lem:semilocalized-finite-forest-forms}]
The nonbacktracking matrix of a finite forest is nilpotent: a closed
nonbacktracking walk would contain a cycle.

\paragraph{Upper form.}
For $t>d^{-1/2}$, the real symmetric Ihara--Bass matrix has nonzero
determinant and converges to the identity as $t\to\infty$.
Its eigenvalues remain positive throughout this interval by
continuity. Its diagonal and adjacency terms are
$\id+\rmD_{\rmT}/(dt^2-1)$ and
$t\rmT/(t^2-d^{-1})$, respectively. Multiplying their comparison by
$(t^2-d^{-1})/t$ gives
\begin{align}
    \rmT\preceq
    \left(t-\frac1{dt}\right)\id
    +\frac1{dt}\rmD_{\rmT}.
\end{align}
Taking $t=1$ proves \eqref{eqn:finite-tree-upper-form}.

\paragraph{Lower form.}
Define the edge deformation at $\lambda=\ii$ by
\begin{align}
    (\rmT_{12}(\lambda))_{xy}
    \coloneqq
    \frac{\lambda(\rmT_{12})_{xy}}{\lambda^2-d^{-1}}.
\end{align}
The exact Ihara blocks for edge weight $d^{-1/2}$ are
\begin{align}
    \rmM_i=\id-\frac{\rmD_{\rmT}^{(i)}}{d+1},
    \qquad
    \lambda\rmT_{12}(\lambda)=\frac d{d+1}\rmT_{12}.
\end{align}
The first block is positive because
$\rmD_{\rmT}^{(1)}\prec d\id$, and remains positive along the
imaginary axis from $\ii\infty$ to $\ii$.
Nilpotence rules out a zero of the Ihara--Bass determinant.
Continuity and the Schur complement give
\begin{align}
    \left(\frac d{d+1}\right)^2
    \rmT_{12}^*
      \left(\id-\frac{\rmD_{\rmT}^{(1)}}{d+1}\right)^{-1}
      \rmT_{12}
    \succeq\frac{\rmD_{\rmT}^{(2)}}{d+1}-\id.
    \label{eqn:exact-finite-forest-lower-form}
\end{align}
For every diagonal value $0\leq D<d$,
\begin{align}
    \frac d{d-D}
    -\frac{d^2}{(d+1)(d+1-D)}
    =\frac{d(2d+1-D)}{(d-D)(d+1)(d+1-D)}>0.
\end{align}
Therefore
\begin{align}
    \left(\frac d{d+1}\right)^2
      \left(\id-\frac{\rmD_{\rmT}^{(1)}}{d+1}\right)^{-1}
    \preceq
      \left(\id-\frac{\rmD_{\rmT}^{(1)}}d\right)^{-1}.
\end{align}
Insert this comparison into
\eqref{eqn:exact-finite-forest-lower-form} and use
\begin{align}
    \frac{\rmD_{\rmT}^{(2)}}{d+1}-\id
    =\frac{\rmD_{\rmT}^{(2)}}d-\id
      -\frac{\rmD_{\rmT}^{(2)}}{d(d+1)}
\end{align}
to obtain \eqref{eqn:finite-forest-lower-form}.
\end{proof}

For the local analysis, fix a cutoff triple with
$s(\gap)\geq a\error_d^{1/4}$ for a fixed $a>0$, and work on the
common pruning and retained-layer event in
Lemmas~\ref{lem:existence_pruned_graph} and
\ref{lem:semilocalized-layer-concentration}, intersected with the
applicable maximum-degree event. The right analysis uses
the subcritical assumptions; the left analysis uses critical sparsity,
$\qr>1$, $0<\gap<\ledge$ and the giant-component reduction.
The retained-layer estimates require $\rD_x\geq c d$. We obtain
this floor from a positive target eigenvalue or from the assumed global
retained degree floor. These inputs apply at the window gaps and at the
separate location gap $s=3\error_d^{1/4}/4$.
Write $s=s(\gap)$ and fix $x\in\gV^{(\vtau)}$. Let $\rmT_x$ be
the Dirichlet restriction of
$\scA^{(\vtau)}$ to
$\ball_{2\radius-1}^{(\vtau)}(x)$.  The centering in
\eqref{eqn:def_H_tau} vanishes on this ball, so $\rmT_x$ is the
normalized adjacency matrix of a finite tree.
Every nonroot vertex in this ball lies outside
$\gV^{(\vtau)}$. Recall the residuals
$\varepsilon_{x,\mathrm{R}}$ and $\varepsilon_{x,\mathrm{L}}$ from
\eqref{eqn:local-profile-error}.

The retained root degree is
\[
    \bar\alpha_x
    \coloneqq
    \frac{|\layer_1^{(\vtau)}(x)|}{d}.
\]
Equations \eqref{eqn:retained-branch-neighbours} and
\eqref{eqn:removed-graph-degrees} give
\begin{align}
    |\bar\alpha_x-\alpha_x|
    \leq\frac{\rD_x^{\gD^{(\vtau)}}}{d}
    \leq C\frac{\log N}{d^2s^2}
    =o(s).
    \label{eqn:local-root-pruning-error}
\end{align}
The last relation follows from the admissible gap range and
\eqref{eqn:truncated-gap-error}.

\subsubsection{Right local spectrum and root eigenvalues}

\begin{lemma}[Uniform right envelope]
\label{lem:uniform-right-envelope}
Put $a_1=\qr^{-2}$, $a_2=\qr^2$ and, for $0<s\leq1$,
$\beta_1=a_1+s/(8\qr)$, $\beta_2=a_2+\qr s/8$.
Let $\rmT$ be a finite rooted bipartite tree with edge weights
$d^{-1/2}$, retained normalized root degree $A$, and at most
$d\beta_i$ children at each nonroot vertex on side $i$.
For a root on side $i$, set $i'=3-i$ and let $J_i(A;\beta)$
have zero diagonal, first coupling $\sqrt A$, and subsequent
couplings $\sqrt{\beta_{i'}},\sqrt{\beta_i},\sqrt{\beta_{i'}},\ldots$.
Then
\begin{align}
 \sup\operatorname{Spec}\rmT
 &\leq\sup\operatorname{Spec}J_i(A;\beta)
 \leq F_i(A)+s/8,
 \label{eqn:finite-tree-root-envelope}\\
 F_i(A)&=\begin{cases}
 \redge,&0\leq A\leq a_i+1,\\
 \Lambda_{\mathrm{br},i}(A),&A>a_i+1,
 \end{cases}
 \qquad
 \Lambda_{\mathrm{br},1}=\Lambda_{\qr^{-1}},\quad
 \Lambda_{\mathrm{br},2}=\Lambda_\qr.
 \label{eqn:uniform-right-ideal-envelope}
\end{align}
In particular, if $A\leq\alpha_x$, a local eigenvalue lies within
$c_0s$ of $\mu\geq\redge+\kappa s$, and $\kappa>1/4+c_0$, then
$\alpha_x>a_i+1$ and
\begin{align}
 \Lambda_{\mathrm{br},i}(\alpha_x)
 &\geq\redge+(\kappa-1/8-c_0)s,\\
 \Lambda_{\mathrm{br},i}(\alpha_x)
   -(\redge+s/8+C\error_d)
 &\geq(\kappa-1/4-c_0)s-C\error_d.
 \label{eqn:right-envelope-separation}
\end{align}
\end{lemma}

\begin{proof}
Write $h=\sqrt{\beta_i}$ and $k=\sqrt{\beta_{i'}}$.
Rescale Lemma~\ref{lem:eigenZ1Z2} by
$\rho=(\beta_1\beta_2)^{1/4}$, with
$q_\beta=(\beta_2/\beta_1)^{1/4}$ and $\alpha=A/\rho^2$.
For $A>h(h+k)$, its formulas
\eqref{eqn:Z1Eigenvec}--\eqref{eqn:Z2Eigenvec} extend to a
positive square-summable eigenvector $y$, normalized by $y_0=1$, at
\[
 E^2=A+\beta_{i'}+\frac{\beta_i\beta_{i'}}{A-\beta_i},
 \qquad t=\frac{hk}{A-\beta_i}<1,
\]
where $t$ is the ratio between successive coefficients of each parity.
For $0<A\leq h(h+k)$,
take $E=h+k$, the upper spectral edge, with $y_j=1$ for $j\geq1$
and $y_0=h/\sqrt A$.
The first nonroot equation is an equality and the root inequality is
$Ey_0\geq\sqrt A y_1$; all subsequent equations are equalities.
The case $A=0$ consists of an isolated root and regular branches.

For a vertex $v$ at depth $j\geq1$, put
\[
 f(v)=d^{-j/2}\frac{y_j}{\sqrt A\prod_{l=1}^{j-1}
                  \sqrt{\beta_{\operatorname{side}(l)}}},
 \qquad f(x)=y_0.
\]
At an interior vertex the adjacency sum is bounded by the parent
term plus $d\beta_{\operatorname{side}(j)}$ copies of the next-layer
value, so $\rmT f\leq Ef$. At a terminal vertex the omitted
next-layer value is positive; explicitly,
$Ef(v)-d^{-1/2}f(\operatorname{parent}(v))\geq0$.
Thus the Dirichlet boundary preserves the supersolution inequality.
The positive-vector comparison on each connected component gives
$\sup\operatorname{Spec}\rmT\leq E$.
The same supersolution bounds finite compressions of $J_i$;
long tail test vectors approach $h+k$, and in the supercritical case
the displayed $y$ is an eigenvector. These facts identify $E$ with
$\sup\operatorname{Spec}J_i(A;\beta)$, including the subcritical case.

The first coupling does not change when $\beta$ is replaced by $a$.
The maximum absolute row sum of the difference is at most
\begin{align}
 \|J_i(A;\beta)-J_i(A;a)\|
 &\leq(\sqrt{\beta_1}-\sqrt{a_1})
       +(\sqrt{\beta_2}-\sqrt{a_2})\leq s/8.
 \label{eqn:right-jacobi-uniform-perturbation}
\end{align}
Concavity gives $s/16$ for each summand. At $a_1a_2=1$ the
supercritical threshold is $a_i+1$ and the supercritical eigenvalue is
precisely $\Lambda_{\mathrm{br},i}(A)$. This proves the first two
claims uniformly for every $0<s\leq1$.
If a local eigenvalue is at least $\redge+(\kappa-c_0)s$, the
envelope excludes $A\leq a_i+1$. Monotonicity of the ideal upper
branch and $A\leq\alpha_x$ give the final two inequalities.
\end{proof}

\paragraph{Root-zero spectrum and active eigenvalues.}
Set the root coordinate equal to zero and let
$(\mu_{\rm loc},\rvv)$ be a normalized positive eigenpair of the
resulting forest. Bipartite symmetry gives
$\|\rvv^{(1)}\|^2=\|\rvv^{(2)}\|^2=1/2$.
Test \eqref{eqn:finite-tree-upper-form} at
$(\qr\rvv^{(1)},\rvv^{(2)})$. Every remaining degree on side
$i$ is below $\tau_i d$, where $\tau_2=\tau_2^{\mathrm R}$.
Lemma~\ref{lem:semilocalized-scalar-form-algebra}, with
$c=1-d^{-1}$ and $\theta=\epsilon=0$, gives
\begin{align}
    \mu_{\rm loc}
    &\leq\frac12\left((1-d^{-1})\redge
           +\qr\tau_1+\qr^{-1}\tau_2^{\mathrm R}\right)
    \nonumber\\
    &\leq\redge+\frac{s}{8}+O_\qr(d^{-1}).
    \label{eqn:right-root-zero-local-edge}
\end{align}
Restoring the root changes one row or one column of the rectangular
block.  Singular-value interlacing therefore permits at most one
positive eigenvalue above this interval, together with its negative
partner.  Apart from this possible root-created pair, the entire
spectrum is contained in
\begin{align}
    [-\redge-s/8-C_{\qr}\error_d,
      \redge+s/8+C_{\qr}\error_d].
    \label{eqn:right-local-radial-spectrum}
\end{align}
For a root on side $i$, each nonroot vertex on side $l$ has
at most $d(\tau_l-d^{-1})\leq d\beta_l$ children, with
$\tau_2=\tau_2^{\mathrm R}$. Apply
Lemma~\ref{lem:uniform-right-envelope} with $A=\bar\alpha_x$.
A root exception within $c_0s$ of a right target
$\mu\geq\redge+s/2$, where $c_0=1/16$, satisfies
\[
 \Lambda_{\mathrm{br},i}(\alpha_x)
 \geq\redge+5s/16.
\]
This implies activity whenever $s\geq\frac34\error_d^{1/4}$.
Indeed, on side $1$, $\alpha_x-a_1>1>s/(8\qr)$.
On side $2$, put $u=\alpha_x-a_2>1$. The identity
\eqref{eqn:profile-energy-degree-separation} gives
\[
 u\geq\frac{(u-1)^2}{u}
 =\Lambda_\qr(\alpha_x)^2-\redge^2
 \geq\frac58\redge s>\frac\qr8s.
\]
Thus both degree cutoffs hold, including $\qr>8$ and $s=1$;
the profile cutoff $\Lambda_x\geq\redge+\error_d^{1/4}/8$
also holds. Low-degree exceptional roots on side $2$ are subcritical
and cannot yield this exceptional eigenvalue. In addition,
\begin{align}
 \operatorname{dist}\bigl(\Lambda_x,
 [-\redge-s/8-C\error_d,\redge+s/8+C\error_d]\bigr)
 \geq3s/16-C\error_d\geq s/8,
 \label{eqn:right-root-exception-is-active}
\end{align}
provided $C\error_d\leq s/16$. This error budget holds in the
window application by its scale checks and at the location scale
by the standing choice of $\gC_0$.

\subsubsection{Left squared spectrum and tail comparison}
Fix $x\in\Vtwolow{\tau_2^{\mathrm{L}}(\gap)}\subset\gV_2$. To retain a positive spectral
gap at the left edge, work with the square on the root bipartition. Write
$\rmT_{12}$ for the off-diagonal block of the local tree operator and
$\rmD_{\rmT}^{(i)}$ for its degree matrix on $\gV_i$, and put
\begin{align}
    \rmS_x\coloneqq\rmT_{12}^*\rmT_{12}
    \label{eqn:left-local-squared-root-operator}
\end{align}
on the root bipartition.
Every first-side vertex lies outside $\Vonehigh{\tau_1(\gap)}$,
so its degree is below $\tau_1d<d$. Since the Dirichlet radius is odd,
the outer layer lies on side $1$. Every nonroot side-$2$ vertex is
therefore interior and retains all its incident edges by
\eqref{eqn:retained-branch-neighbours}; its degree is at least
$\tau_2^{\mathrm L}d$. Truncating the outer first-side degrees
improves the inverse margin.

For a normalized vector $v$ in the second-side root-zero subspace,
\eqref{eqn:finite-forest-lower-form} and the maximum-degree bound give
\begin{align}
 \tau_2^{\mathrm L}-1-\frac{C_{b,\qr}}d
 &\leq\left\langle\rmT_{12}v,
       \left(\id-\frac{\rmD^{(1)}_{\rmT}}d\right)^{-1}
       \rmT_{12}v\right\rangle
 \leq\frac{\|\rmT_{12}v\|^2}{1-\tau_1}.
 \label{eqn:finite-tree-lower-form}
\end{align}
The left error is the forest correction
$\rmD^{(2)}_{\rmT}/(d(d+1))=O_{b,\qr}(d^{-1})$.
The lower conclusion of Lemma~\ref{lem:semilocalized-scalar-form-algebra}
yields, on this subspace,
\begin{align}
    \rmT_{12}^*\rmT_{12}
    \succeq
      (1-\tau_1)(\tau_2^{\mathrm{L}}-1)\id-O_{b,\qr}(d^{-1})\id.
    \label{eqn:left-root-zero-local-edge}
\end{align}
Restoring the root is a one-dimensional principal extension, so
interlacing permits at most one eigenvalue below this floor. It may be
zero. The upper forest form and the critical maximum-degree bound give
$\|\rmT_x\|\leq C_{b,\qr}$ and hence
$\|\rmS_x\|\leq C_{b,\qr}$. Let $\rmS_{x,\rm reg}$ be the restriction
after removing the eigenspace of the possible low root eigenvalue. The
threshold identity
\begin{align}
    (1-\tau_1(\gap))(\tau_2^{\mathrm{L}}(\gap)-1)
    =(\ledge-s/8)^2
    \label{eqn:left-threshold-identity}
\end{align}
gives the weaker but more convenient bound
\begin{align}
    \operatorname{Spec}(\rmS_{x,\rm reg})
    \subset[(\ledge-s/8-C_{b,\qr}\error_d)^2,C_{b,\qr}].
    \label{eqn:left-local-squared-spectrum}
\end{align}
At the left edge, a scalar Schur equation identifies the degree
associated with an exceptional eigenvalue.
Let $\kappa_{x,\mathrm{L}}^{\rm loc}$ denote the possible nonnegative
eigenvalue of $\rmS_x$ below the squared left regular interval.  With respect to the
root and its orthogonal complement, write
\begin{align}
    \rmS_x
    =
    \begin{bmatrix}
       \bar\alpha_x&c_x^*\\
       c_x&\rmS_x^{(0)}
    \end{bmatrix}.
    \label{eqn:left-root-rest-block}
\end{align}
An exception below the root-zero floor
\eqref{eqn:left-root-zero-local-edge} is the smallest eigenvalue
of $\rmS_x$. If it lies in a fixed positive target neighborhood,
its lower bound and the root Rayleigh quotient give
\begin{align}
 \bar\alpha_x\geq\kappa_{x,\mathrm L}^{\rm loc}
       \geq c_{\qr,\sK}>0.
 \label{eqn:left-positive-root-degree-floor}
\end{align}
This floor permits the retained-layer estimates before a profile
is invoked. The Schur complement gives the exact equation
\begin{align}
    \bar\alpha_x-\kappa_{x,\mathrm{L}}^{\rm loc}
    =\left\langle
       c_x,
       (\rmS_x^{(0)}-\kappa_{x,\mathrm{L}}^{\rm loc})^{-1}c_x
     \right\rangle.
    \label{eqn:left-root-schur-equation}
\end{align}

\emph{Comparison with the regular tail.}
Put $\eta_s=\error_d+\error_{\mathrm{prof}}(\gap)$ and define
\begin{align}
    \mathcal K_{\mathrm L}(s)
    \coloneqq\left\{t^2:
      \frac{3\inf\sK}{4}\leq t\leq
      \ledge-\frac{s}{4}-C\eta_s\right\}.
    \label{eqn:left-squared-target-interval}
\end{align}
For $0\leq z<\ledge^2$, let
\begin{align}
    m_{\mathrm L}(z)
    \coloneqq\frac{\qr^2+\qr^{-2}-z
      -\sqrt{(\qr^2+\qr^{-2}-z)^2-4}}2.
    \label{eqn:left-regular-tail-resolvent}
\end{align}
\begin{lemma}[Two-factor left tail comparison]
\label{lem:left-two-factor-tail-comparison}
Let $x\in\gV_2$ and let $\rmT_x$ be the normalized adjacency matrix
of its retained tree ball of odd radius $2\radius-1$. Suppose that the
layers through depth $\radius+1$ are nonempty, that
$a_*/2\leq\bar\alpha_x\leq M$, and that $\|\rmT_x\|\leq M$.
For $n_j=|\layer_j^{(\vtau)}(x)|$ and
$t_j=(n_j/(dn_{j-1}))^{1/2}$, assume, uniformly for
$2\leq j\leq\radius+1$,
\begin{align}
    |t_j-\omega_j|&\leq C_*\eta_s,
    & \omega_j&=\begin{cases}\qr^{-1},&j\textnormal{ even},\\
                              \qr,&j\textnormal{ odd},\end{cases}
    \label{eqn:left-tail-coupling-input}\\
    \frac{1}{dn_j}
      \sum_{y\in\layer_{j-1}^{(\vtau)}(x)}
      \left(N_j(y)-\frac{n_j}{n_{j-1}}\right)^2
      &\leq C_*^2\eta_s^2,
    &N_j(y)&=|\layer_j^{(\vtau)}(x)\cap
                       \layer_1^{(\vtau)}(y)|.
    \label{eqn:left-tail-fluctuation-input}
\end{align}
Assume $\eta_s$ is sufficiently small depending on the displayed
constants. Let $\rmS_x^{(0)}$ and $c_x$ be the root-zero block and
root coupling in \eqref{eqn:left-root-rest-block}. If $0<s\leq1$,
$0\leq z\leq\ledge^2-\sigma s$, and
\begin{align}
    \operatorname{dist}(z,\operatorname{Spec}(\rmS_x^{(0)}))
        \geq\sigma s,
    \label{eqn:left-tail-resolvent-input}
\end{align}
then a truncated radial vector $\widetilde g_{x,z}$ satisfies
\begin{subequations}\label{eqn:left-truncated-tail-comparison}
\begin{align}
    \| (\rmS_x^{(0)}-z)\widetilde g_{x,z}-c_x\|
    &\leq C\left(s^{-1/4}\eta_s
                     +s^{-1}e^{-cs\radius}\right),
    \label{eqn:left-truncated-tail-residual}\\
    \left|\langle c_x,\widetilde g_{x,z}\rangle
               -\bar\alpha_x\qr^{-2}m_{\mathrm L}(z)\right|
    &\leq C\eta_s.
    \label{eqn:left-truncated-tail-root-value}
\end{align}
\end{subequations}
Consequently,
\begin{align}
    \left|\langle c_x,(\rmS_x^{(0)}-z)^{-1}c_x\rangle
               -\bar\alpha_x\qr^{-2}m_{\mathrm L}(z)\right|
    &\leq\delta_{\rm root}(\gap),
    \nonumber\\
    \delta_{\rm root}(\gap)
    &\leq C\left(s^{-5/4}\eta_s
                       +s^{-2}e^{-cs\radius}\right).
    \label{eqn:left-root-schur-comparison}
\end{align}
Here $C,c>0$ depend only on $(\qr,a_*,M,C_*,\sigma)$.
\end{lemma}

\begin{proof}
Write $\rvs_j=n_j^{-1/2}\ones_{\layer_j^{(\vtau)}(x)}$ and
$R=2K+1$ for the largest odd integer at most $\radius$.
Let $A$ be $\rmT_{12}$ restricted to the root-zero second-side
space, so that $A^*A=\rmS_x^{(0)}$. The even and odd layer
isometries are
\begin{align}
    U_x\rve_j&=\rvs_{2j+2},&&0\leq j<K,
    &V_x\rve_k&=\rvs_{2k+1},&&0\leq k\leq K.
    \label{eqn:left-tail-even-odd-maps}
\end{align}
Define $H_K:\mathbb R^K\to\mathbb R^{K+1}$ by
\begin{align}
    H_K\rve_j=\qr^{-1}\rve_j+\qr\rve_{j+1}.
    \label{eqn:left-tail-bidiagonal-map}
\end{align}
Then $J_K=H_K^*H_K$ has diagonal $\qr^2+\qr^{-2}$
and both off-diagonals equal to one, including at its first and last
coordinates. This factorization keeps the physical parent edge at
layer two and the physical child edge at layer $R-1$.

For $j\geq2$, put
\begin{align}
    F_j=\frac1{\sqrt{dn_j}}
      \sum_{y\in\layer_{j-1}^{(\vtau)}(x)}
      \left(N_j(y)-\frac{n_j}{n_{j-1}}\right)\rve_y,
    \qquad F_1=0.
    \label{eqn:left-tail-child-fluctuation-vector}
\end{align}
Each $F_j$ is orthogonal to $\rvs_{j-1}$, and its squared norm
is the left-hand side of \eqref{eqn:left-tail-fluctuation-input}.
The tree identity, with $\zeta_j=t_j-\omega_j$ for $j\geq2$, gives
\begin{align}
    AU_x&=V_xH_K+E_{\rm e},
    \nonumber\\
    E_{\rm e}a
      &=\sum_{j=0}^{K-1}a_j
        \left(\zeta_{2j+2}\rvs_{2j+1}
             +\zeta_{2j+3}\rvs_{2j+3}+F_{2j+2}\right),
    \label{eqn:left-tail-first-adjacency-factor}\\
    A^*V_x&=U_xH_K^*+E_{\rm o}+B_R,
    \nonumber\\
    E_{\rm o}v
      &=\sum_{k=1}^{K}v_k
           \left(\zeta_{2k+1}\rvs_{2k}+F_{2k+1}\right)
        +\sum_{k=0}^{K-1}v_k\zeta_{2k+2}\rvs_{2k+2},
    \nonumber\\
    B_Rv&=t_{R+1}v_K\rvs_{R+1}.
    \label{eqn:left-tail-second-adjacency-factor}
\end{align}
At $k=0$, root removal deletes the term $t_1v_0\rve_x$;
there is no root child fluctuation because $F_1=0$.
The last displayed term is the outward Dirichlet truncation error
inside the larger local ball. No layer concentration beyond depth
$\radius+1$ has been used. Within each parity, the fluctuation
vectors have disjoint supports and are orthogonal to every radial
layer indicator. The radial coupling errors form two weighted shifts.
Thus
\begin{align}
    \|E_{\rm e}\|+\|E_{\rm o}\|\leq C\eta_s,
    \qquad t_{R+1}\leq C,
    \qquad\|H_K\|\leq\qr+\qr^{-1}.
    \label{eqn:left-tail-factor-error-norms}
\end{align}

Let $J_{\mathrm L}$ be the infinite constant Jacobi matrix on
$\ell^2(\mathbb N_0)$ and set
$g(z)=(J_{\mathrm L}-z)^{-1}\rve_0$.
Writing $m=m_{\mathrm L}(z)$, its recurrence gives
\begin{align}
    g_j=(-1)^jm^{j+1},
    \quad \|g(z)\|^2=\frac{m^2}{1-m^2},
    \quad 1-m^2\geq c\sqrt{s},
    \quad m\leq e^{-c\sqrt{s}}.
    \label{eqn:left-tail-coefficient-decay}
\end{align}
Indeed, $m+m^{-1}=\qr^2+\qr^{-2}-z$ and
$\ledge^2-z=(1-m)^2/m\geq\sigma s$; on $z\geq0$ one
also has $\qr^{-2}\leq m<1$. Let
$g^{[K]}=(g_0,\ldots,g_{K-1})$ and define
$\widetilde g_{x,z}=\|c_x\|U_xg^{[K]}$.
The root coupling is exactly
\begin{align}
    c_x=d^{-1}\ones_{\layer_2^{(\vtau)}(x)},
    \qquad U_x\rve_0=c_x/\|c_x\|,
    \qquad \|c_x\|^2=\bar\alpha_xt_2^2.
    \label{eqn:left-tail-first-coefficient}
\end{align}
The assumptions ensure $\|c_x\|\leq C$. Since
$(J_K-z)g^{[K]}=\rve_0-g_K\rve_{K-1}$, applying
\emph{both} adjacency factors gives the exact residual
\begin{align}
    (\rmS_x^{(0)}-z)\widetilde g_{x,z}-c_x
      =\|c_x\|\bigl(
        A^*E_{\rm e}g^{[K]}
        +E_{\rm o}H_Kg^{[K]}
        -g_K\rvs_{R-1}
        +\qr t_{R+1}g_{K-1}\rvs_{R+1}\bigr).
    \label{eqn:left-tail-two-factor-residual-identity}
\end{align}
The first two terms are bounded by $C\eta_s\|g\|$
using $\|A\|\leq M$. The two terminal terms are bounded by
$Ce^{-c\sqrt{s}\radius}$; $K\geq\radius/4$ for all sufficiently
large $\radius$. With $0<s\leq1$ these bounds imply
\eqref{eqn:left-truncated-tail-residual}.
Only the first even layer pairs with $c_x$, so
\begin{align}
    \langle c_x,\widetilde g_{x,z}\rangle
       =\|c_x\|^2m
       =\bar\alpha_xt_2^2m
       =\bar\alpha_x\qr^{-2}m+O(\eta_s).
\end{align}
This proves \eqref{eqn:left-truncated-tail-root-value}, including
the first coefficient. Finally,
\eqref{eqn:left-tail-resolvent-input} implies
$\|(\rmS_x^{(0)}-z)^{-1}\|\leq(\sigma s)^{-1}$.
Apply this inverse to the residual and pair with $c_x$ to obtain
\eqref{eqn:left-root-schur-comparison}.
\end{proof}

Lemma~\ref{lem:semilocalized-layer-concentration} supplies
\eqref{eqn:left-tail-coupling-input}--\eqref{eqn:left-tail-fluctuation-input}
for roots $\rD_x\geq c d$.
More explicitly, let
$\sigma_x=\rD_x^{\gD^{(\vtau)}}/\rD_x+
\log N/(\rD_xd)$ denote the aggregate branch-loss parameter.
The retained estimates give
\begin{align}
    |t_j-\omega_j|&\leq C\left(
      \sqrt{\frac{\log N}{d\rD_x}}+\sigma_x\right),
    \nonumber\\
    \|F_j\|&\leq C\left(
       \sqrt{\frac{(1+\log N/\rD_x)\log d}{d}}
          +\sigma_x\right).
    \label{eqn:left-tail-input-from-retained-layers}
\end{align}
For $d\leq B\log N$ and $\rD_x\geq c d$, the square-root
terms are at most $C\error_d$.
The removed-degree bound and $\rD_x\geq cd$ give
\begin{align}
    \sigma_x\leq C\frac{\log N}{d^2s^2}
       =C\frac{\error_d^2}{s^2(\log d)^2}.
    \label{eqn:left-tail-aggregate-loss-budget}
\end{align}
For $s\geq a\error_d^{1/4}$, \eqref{eqn:left-tail-aggregate-loss-budget}
gives $\sigma_x\leq C\error_d$. This includes the separate
location scale $s=3\error_d^{1/4}/4$, whose substitution is
\eqref{eqn:location-pruning-scale-identities}. Both factor-error
inputs therefore hold with
$\eta_s=\error_d+\error_{\mathrm{prof}}(\gap)$.

\subsubsection{Left root recovery}

The ideal Schur equation is
\begin{align}
    \alpha-z-\alpha\qr^{-2}m_{\mathrm{L}}(z)=0.
    \label{eqn:left-ideal-root-schur-equation}
\end{align}
For $0<z<\ledge^2$, its degree parameter is
\begin{align}
    \alpha_*(z)
    \coloneqq
    \frac{z}{1-\qr^{-2}m_{\mathrm{L}}(z)}
    =\qr^2-m_{\mathrm{L}}(z)^{-1}.
    \label{eqn:left-ideal-degree-from-energy}
\end{align}
Thus $0<\alpha_*(z)<\qr^2-1$ and
$z=\Lambda_\qr(\alpha_*(z))^2$. Since
$0<m_{\mathrm{L}}(z)\leq1$, the denominator in
\eqref{eqn:left-ideal-degree-from-energy} is at least $1-\qr^{-2}>0$.
If a local root eigenvalue $z$ satisfies the exact Schur equation
and its tail comparison has error $\delta$, division by this
denominator gives
\begin{align}
 |\bar\alpha_x-\alpha_*(z)|\leq C_\qr\delta.
 \label{eqn:left-schur-degree-recovery}
 \label{eqn:left-global-root-parameter-recovery}
\end{align}
For an original degree $0<\alpha_x<\qr^2-1$, the derivative
\[
    \frac{d}{d\alpha}\Lambda_\qr(\alpha)^2
    =1-(\alpha-\qr^2)^{-2}\in[0,1],
    \qquad 0<\alpha<\qr^2-1,
\]
and $z=\Lambda_\qr(\alpha_*(z))^2$ then yield
\begin{align}
 |z-\Lambda_\qr(\alpha_x)^2|
 \leq C_\qr\bigl(\delta+|\bar\alpha_x-\alpha_x|\bigr).
 \label{eqn:left-root-energy-comparison}
\end{align}
\paragraph{Positive windows.}
The root degree floor is \eqref{eqn:left-positive-root-degree-floor},
and the forest upper form gives $M=C_{b,\qr}$.
For $c_0=c_0(\qr,\sK)>0$ sufficiently small, points within
$c_0s$ of \eqref{eqn:left-squared-target-interval} have
$0<z\leq\ledge^2-c_{\qr,\sK}s$ and distance at least
$c_{\qr,\sK}s$ from the root-zero spectrum. To see the latter,
the regular-to-target margin before form errors is
\begin{align}
    (\ledge-s/8)^2-(\ledge-s/4)^2
       =\frac{s}{8}\left(2\ledge-\frac{3s}{8}\right)
       \geq c_\qr s.
    \label{eqn:left-window-squared-buffer-margin}
\end{align}
Nonemptiness of the target interval entails $s<4\ledge$.
The $O(\error_d)$ floor correction and the $c_0s$
neighborhood are absorbed once $\error_d/s$ is sufficiently
small. Lemma~\ref{lem:left-two-factor-tail-comparison} therefore
applies with constants $C_{b,\qr,\sK}$.
In particular, at critical sparsity, for every admissible
$s\geq C\error_d^{1/4}$,
\begin{align}
    \frac{\delta_{\rm root}(\gap)}s
      \leq C\left(
        \frac{\error_d}{s^{9/4}}
        +\frac{\error_d}{s^{13/4}\log d}
        +s^{-3}e^{-cs\radius}\right)=o(1).
    \label{eqn:left-window-root-error-budget}
\end{align}
For the fixed eigenvector window, $s$ is bounded below by a
constant depending on $\sK$, and this bound is $O_{b,\qr,\sK}(\error_d)$.
Suppose now that
$\operatorname{dist}(\kappa_{x,\mathrm{L}}^{\rm loc},\mathcal K_{\mathrm L}(s))<c_0s$.
Apply \eqref{eqn:left-schur-degree-recovery} using
\eqref{eqn:left-root-schur-comparison}. Throughout the target neighborhood,
$\alpha_*(z)\geq c_{\qr,\sK}$ and
$\qr^2-1-\alpha_*(z)\geq c_{\qr,\sK}\sqrt{s}$.  Hence
\eqref{eqn:local-root-pruning-error} and the scale bound in
\eqref{eqn:left-window-root-error-budget} imply
$0<\alpha_x<\qr^2-1$. Apply \eqref{eqn:left-root-energy-comparison}
with $z=\kappa_{x,\mathrm L}^{\rm loc}$ and
$\delta=\delta_{\rm root}(\gap)$, whose constant depends on
$(b,\qr,\sK)$. The $s/4$ target and $s/8$ regular buffers give
\begin{align}
    x\in\gV_{2,\mathrm{L}}^\star,
    \qquad
    \operatorname{dist}\left(
      \Lambda_\qr(\alpha_x)^2,
      \operatorname{Spec}(\rmS_{x,\rm reg})
    \right)
    \geq c_{b,\qr,\sK}s.
    \label{eqn:left-root-exception-is-active}
\end{align}
A zero root exception, possible when $\alpha_x=0$, lies outside this
target neighborhood.

\paragraph{Global comparison with a degree floor.}
\begin{lemma}[Global left tail and root recovery]
\label{lem:global-left-root-recovery}
Assume critical sparsity and $\qr>1$. On the common pruning and
retained-layer event at $\vtau(s)$, take a root
$x\in\Vtwolow{\tau_2^{\mathrm L}(s)}$ with
$\bar\alpha_x\geq a_0/2>0$.
At $s=3\error_d^{1/4}/4$, the Schur comparison through zero,
positive pivot, and root recovery in
\eqref{eqn:left-global-root-schur-comparison}--
\eqref{eqn:appendix-left-root-recovery} hold for sufficiently large
$N$, with constants depending on $(b,\qr,a_0)$.
\end{lemma}
\begin{proof}
The root-zero floor supplies the second spectral-distance input in
Lemma~\ref{lem:left-two-factor-tail-comparison} uniformly on
\begin{align}
    0\leq z\leq(\ledge-s)^2+c_0s,
    \qquad s=\frac34\error_d^{1/4},
    \label{eqn:left-global-tail-domain}
\end{align}
where $c_0=c_0(b,\qr,a_0)>0$ is sufficiently small.
Indeed, for $s\leq\ledge/2$, the ideal squared margin equals
\begin{align}
    (\ledge-s/8)^2-(\ledge-s)^2
       =\frac{7s}{8}\left(2\ledge-\frac{9s}{8}\right)
       \geq c_\qr s.
    \label{eqn:left-global-squared-buffer-margin}
\end{align}
The retained degree floor
$\bar\alpha_x\geq a_0/2$ supplies
$a_*=a_0$, and the retained-layer inputs supply both factor errors.
The lemma consequently gives, throughout
\eqref{eqn:left-global-tail-domain},
\begin{align}
    \left|\langle c_x,(\rmS_x^{(0)}-z)^{-1}c_x\rangle
              -\bar\alpha_x\qr^{-2}m_{\mathrm L}(z)\right|
      &\leq\delta_{\rm root}^{\rm loc},
    \nonumber\\
    \delta_{\rm root}^{\rm loc}
      &\leq C_{b,\qr,a_0}\left[
        \error_d^{-5/16}\left(
           \error_d+\frac{4\error_d^{3/4}}{3\log d}\right)
        +\error_d^{-1/2}e^{-c_\qr\error_d^{1/4}\radius}
        \right],
    \label{eqn:left-global-root-schur-comparison}\\
    \frac{\delta_{\rm root}^{\rm loc}}s
      &\leq C_{b,\qr,a_0}\left(
         \error_d^{7/16}
         +\frac{\error_d^{3/16}}{\log d}\right)+o(1)
        \longrightarrow0.
    \label{eqn:left-global-root-error-budget}
\end{align}
The exponential budget \eqref{eqn:location-exponential-column-budgets}
absorbs the terminal term, including after division by $s$.
The degree floor makes this comparison uniform through $z=0$, where
$m_{\mathrm L}(0)=\qr^{-2}$.

Set $z_0=\min\{a_0(1-\qr^{-2})/8,\ledge^2/4\}>0$.
For every $0\leq z\leq z_0$ in
\eqref{eqn:left-global-tail-domain}, the exact Schur pivot satisfies
\begin{align}
    \bar\alpha_x-z-
       \langle c_x,(\rmS_x^{(0)}-z)^{-1}c_x\rangle
      &\geq \frac{a_0}{2}(1-\qr^{-2})
                -z_0-\delta_{\rm root}^{\rm loc}
       \geq\frac{a_0}{4}(1-\qr^{-2})>0
    \label{eqn:left-global-positive-schur-pivot}
\end{align}
for all sufficiently large $N$. Since the root-zero block is
positive on this interval, the Schur determinant identity excludes
local eigenvalues there, including at zero.

For a root-created eigenvalue $z\geq z_0$ in
\eqref{eqn:left-global-tail-domain}, apply
\eqref{eqn:left-schur-degree-recovery} with
$\delta=\delta_{\rm root}^{\rm loc}$ from
\eqref{eqn:left-global-root-schur-comparison}.
Also $\alpha_*(z)\geq z_0$ and
\begin{align}
    \qr^2-1-\alpha_*(z)
        =m_{\mathrm L}(z)^{-1}-1\geq c_\qr\sqrt{s}.
    \label{eqn:left-global-parameter-branch-margin}
\end{align}
The latter follows from
$\ledge^2-z=(m_{\mathrm L}(z)^{-1}-1)^2m_{\mathrm L}(z)
\geq c_\qr s$ and $m_{\mathrm L}(z)\leq1$.
The location-scale substitution
\eqref{eqn:location-pruning-scale-identities} gives
\begin{align}
    |\bar\alpha_x-\alpha_x|
       \leq C_{b,\qr}\frac{\log N}{d^2s^2}=o(s).
    \label{eqn:left-global-root-pruning-budget}
\end{align}
Together these estimates imply $0<\alpha_x<\qr^2-1$.
Apply \eqref{eqn:left-root-energy-comparison} to obtain
\begin{align}
    |z-\Lambda_\qr(\alpha_x)^2|
       &\leq C_{b,\qr,a_0}\left(
          \delta_{\rm root}^{\rm loc}
          +\frac{16\log N}{9d^2\error_d^{1/2}}\right)
        =o(s).
    \label{eqn:appendix-left-root-recovery}
\end{align}
Since $z\geq z_0$, passing to positive singular values costs only
$C_{b,\qr,a_0}$.
\end{proof}

\subsection{Exceptional-coordinate mass on the profile complement}
\label{subsec:semilocalized-radial-mass-proof}

The local operators $\rmT_x$ and retained root degrees $\bar\alpha_x$
are as in Subsection~\ref{subsec:semilocalized-local-spectra}.
Profile orthogonality controls the local eigenvectors near the target
spectral parameter; a reduced resolvent controls the outer boundary.
Together they bound the exceptional-coordinate mass.

\begin{lemma}[Radial mass estimate]
\label{lem:semilocalized-radial-mass}
The right- and left-edge estimates are as follows.
\begin{enumerate}
    \item Suppose $s(\gap)\geq\gC_{B,\nu}\error_d^{1/4}$, and let
    $(\mu,\rvw)$ be a normalized eigenpair of
    $\overline\Pi_{\mathrm{R}}^{(\vtau)}\scA^{(\vtau)}
    \overline\Pi_{\mathrm{R}}^{(\vtau)}$ with $\mu>0$. Uniformly over these
    gaps, with probability at least $1-C_{B,\qr,\nu}N^{-\nu}$,
    \begin{align}
        \mu\geq\redge+\frac{s(\gap)}2
        +C_{B,\qr,\nu}(\error_d+\error_{\mathrm{prof}}(\gap))
        \quad\Longrightarrow\quad
        \sum_{x\in\gV^{(\vtau)}}(1+\alpha_x)|\ervw_x|^2
        \leq C_{B,\qr,\nu}(\error_d+\error_{\mathrm{prof}}(\gap)).
        \label{eqn:right-radial-mass}
    \end{align}
    \item Under Assumptions~\ref{ass:proportional_regime} and
    \ref{ass:critical_sparsity}, with $\qr>1$, suppose
    $0<\gap<\ledge$ and
    $s(\gap)\geq\gC\error_d^{1/4}$, and let $(\mu,\rvw)$ be a
    normalized eigenpair of
    $\overline\Pi_{\mathrm{L}}^{(\vtau)}\scA^{(\vtau)}
    \overline\Pi_{\mathrm{L}}^{(\vtau)}$ with $\mu>0$. Uniformly over these
    gaps, with probability at least $1-CN^{-\nu_0}$,
    \begin{align}
        &\frac{3\inf\sK}{4}\leq\mu\leq\ledge-\frac{s(\gap)}4
        -C_{b,\qr,\sK}(\error_d+\error_{\mathrm{prof}}(\gap))
        \nonumber\\
        &\quad\Longrightarrow\quad
        \sum_{x\in\Vtwolow{\tau_2^{\mathrm{L}}(\gap)}}(1+\alpha_x)|\ervw_x|^2
        \leq C_{b,\qr,\sK}(\error_d+\error_{\mathrm{prof}}(\gap)).
        \label{eqn:left-radial-mass}
    \end{align}
\end{enumerate}
\end{lemma}

\subsubsection{Reduced resolvents and local root estimates}

The following resolvent bound applies to the local tree operator and
its square. Its constants are uniform in the exceptional spectral
values, which may grow with the root degree at the right edge.

\begin{lemma}[Reduced resolvent for a local root operator]
\label{lem:reduced-local-resolvent}
Let $\rmT$ be self-adjoint on a finite orthogonal layer decomposition
$\bigoplus_{j=0}^{J_{\max}}\mathcal H_j$.
In the right case it has layer bandwidth one, with
$\|P_{j+1}\rmT P_j\|\leq M$ for $j\geq1$; the coupling between
layers zero and one may be arbitrary. Its spectrum is contained in
$I\cup\{-\ell,\ell\}$, where $I=[-M_1,M_1]$, and the two
exceptional values, if present, satisfy $\ell>M_1$.
In the left case it has bandwidth at most two and
$\|P_j\rmT P_k\|\leq M$ for $j\ne k$;
its spectrum is contained in $I\cup\{\ell\}$, where
$I=[a,b]\subset[-M_1,M_1]$ and $\ell<a$ if present.

Let $0<g\leq1$, $z\in\R$ and let $\mathcal P$ be zero or the
spectral projection of the positive exceptional value in the right
case, or the unique exceptional value in the left case.
Suppose
\[
 \operatorname{dist}(z,I)\geq g,\qquad
 \operatorname{dist}(z,\operatorname{Spec}(\rmT|_{\mathcal P^\perp}))
 \geq g.
\]
When $\mathcal P\ne0$, assume in addition
$\operatorname{dist}(\ell,I)\geq g$ and, in the right case,
$\operatorname{dist}(\ell,-\ell)\geq g$.
Define
\begin{align}
 \mathcal R_{\mathcal P}(z)=\mathcal P^\perp
 (\rmT|_{\mathcal P^\perp}-z)^{-1}\mathcal P^\perp.
 \label{eqn:def-reduced-local-resolvent}
\end{align}
For sums of layer projections $P_A,P_C$ whose layer supports are
at distance at least $L$, the reduced resolvent satisfies
\begin{align}
 \|P_A\mathcal R_{\mathcal P}(z)P_C\|
 \leq C_{M,M_1}g^{-2}e^{-c_{M,M_1}gL}.
 \label{eqn:reduced-local-resolvent-decay}
\end{align}
\end{lemma}
\begin{proof}
Collapse the edge between layer indices zero and one, and let
$\phi(j)$ be the distance from index $j$ to the layer support of
$P_A$ in this collapsed path. Then $\phi=0$ on that support,
$\phi(0)=\phi(1)$, $|\phi(j)-\phi(k)|\leq|j-k|$, and
$\phi\geq L-1$ on the support of $P_C$.
Put $D_\theta=\sum_j e^{\theta\phi(j)}P_j$.
The right root coupling is unchanged by conjugation. All remaining
off-diagonal blocks are bounded, and the bandwidth is at most two,
so the block Schur test gives
\[
 \|D_\theta^{-1}\rmT D_\theta-\rmT\|
 \leq C_M(e^{2\theta}-1).
\]
Choose $\theta=c_M g$ sufficiently small that this bound is at
most $g/8$. Whenever
$\operatorname{dist}(\zeta,\operatorname{Spec}\rmT)\geq g/4$,
the resolvent identity and a Neumann series give
\[
 \|D_\theta^{-1}(\rmT-\zeta)^{-1}D_\theta\|\leq8/g,
 \qquad
 \|P_A(\rmT-\zeta)^{-1}P_C\|
 \leq C_Mg^{-1}e^{-c_MgL}.
\]
If $\mathcal P=0$, apply this directly at $\zeta=z$.

Suppose $\mathcal P\ne0$. Take the positively oriented boundary
of the rectangle whose real sides extend the endpoints of $I$ by
$g/4$ and whose imaginary sides are at $\pm g/4$.
Every point of this contour is at distance at least $g/4$ from the
spectrum and at least $3g/4$ from $z$; its length is bounded in
terms of $M_1$. It encloses all the regular spectrum and excludes
the removed eigenvalue. In the right case also take the circle of
radius $g/4$ around $-\ell$. Symmetry gives
$\operatorname{dist}(-\ell,I)\geq g$; the circle excludes the
positive exceptional value and is at distance at least $3g/4$
from $z$. Its length is $\pi g/2$.
Functional calculus expresses the reduced resolvent as
\[
 \frac1{2\pi\ii}\int_\mathcal C
 \frac1{\zeta-z}(\zeta-\rmT)^{-1}\,d\zeta,
\]
where $\mathcal C$ is the rectangle together with the additional
circle in the right case. The preceding decay bound, the factor
$|\zeta-z|^{-1}\leq4/(3g)$, and the contour lengths give the
claimed $Cg^{-2}e^{-cgL}$ estimate. The constants are independent
of the root coupling and of an exceptional value escaping to infinity.
\end{proof}

\begin{lemma}[Local root estimate]
\label{lem:semilocalized-local-root-estimate}
Use the assumptions and notation of
Lemma~\ref{lem:reduced-local-resolvent}. Let $P_0$ project onto the root
layer, and suppose
\begin{align}
    (\rmT-z)u=f+h,
    \qquad
    h=\sum_{j\geq R}P_jh,
    \qquad R\geq1.
    \label{eqn:local-root-source-equation}
\end{align}
If, for $\varepsilon,W\geq0$ and $K>0$,
\begin{align}
    \|\mathcal Pu\|&\leq K g^{-1}\varepsilon W,
    &\|f\|&\leq K\varepsilon W,
    &\|h\|&\leq KW,
    \label{eqn:local-root-input-bounds}
\end{align}
then
\begin{align}
    \|P_0u\|^2
    \leq C_{M,M_1,K}\left(
      g^{-2}\varepsilon^2+g^{-4}e^{-c_{M,M_1}gR}
    \right)W^2.
    \label{eqn:local-root-estimate}
\end{align}
\end{lemma}
\begin{proof}
Spectral calculus gives
\[
    u=\mathcal Pu+\mathcal R_{\mathcal P}(z)(f+h).
\]
The resolvent norm is at most $g^{-1}$, and
\eqref{eqn:reduced-local-resolvent-decay} controls its action from the
support of $h$ to the root. Thus
\[
    \|P_0u\|
    \leq 2K g^{-1}\varepsilon W
       +C_{M,M_1}K g^{-2}e^{-c_{M,M_1}gR}W.
\]
Squaring and adjusting the constants proves the estimate.
\end{proof}

\subsubsection{Profile orthogonality and cutoff source equations}

For a selected centre set $\mathcal C$, use the residuals
$\varepsilon_{x,\mathcal C}$ from
Lemma~\ref{lem:semilocalized-selected-profile-estimates}, defined only
for $x\in\mathcal C$. On an unselected exceptional ball the profile
compression defect vanishes. Every residual term below carries the
indicator of selection and is read as zero outside $\mathcal C$.
For the windows in Lemma~\ref{lem:semilocalized-radial-mass}, take
$\mathcal C=\gV_{\mathrm E}^\star$ and abbreviate the residual by
$\varepsilon_{x,\mathrm E}$. A normalized complementary eigenpair
$(\mu,\rvw)$ satisfies
\begin{align}
    \scA^{(\vtau)}\rvw&=\mu\rvw+\rvr_{\mathrm E},
    &\rvr_{\mathrm E}&=\Pi_{\mathrm E}^{(\vtau)}\scA^{(\vtau)}\rvw,
    \qquad \mathrm E\in\{\mathrm L,\mathrm R\}.
    \label{eqn:radial-proof-defect}
\end{align}
Disjoint residual supports imply the local defect bound
\begin{align}
    \|\rvr_{\mathrm E}|_{\ball_{\radius+1}^{(\vtau)}(x)}\|
    \leq C\indi{x\in\mathcal C}\varepsilon_{x,\mathrm E}
      \|\rvw|_{\ball_{2\radius}^{(\vtau)}(x)}\|.
    \label{eqn:local-defect-bound}
\end{align}
We use the scalar consequence of
Lemma~\ref{lem:low-column-Schur-spectrum}(ii): if a normalized real
approximate eigenvector of a real symmetric operator has residual at most
$\varepsilon$,
its spectral parameter lies within $\varepsilon$ of a simple eigenvalue,
and that eigenvalue is at distance $g>4\varepsilon$ from the remaining
spectrum, then a choice of sign makes its distance from the normalized
eigenvector at most $C\varepsilon/g$.

\paragraph{Right-edge spectral projection.}
Let $\lambda_{x,\mathrm{R}}^{\rm loc}$ be the possible positive eigenvalue
of $\rmT_x$ above the interval in
\eqref{eqn:right-local-radial-spectrum}. Choose a sufficiently small
constant $c_0=c_0(\qr)>0$, and let
$\mathcal P_{x,\mathrm{R}}(\mu)$ be its spectral projection if
\begin{align}
    |\lambda_{x,\mathrm{R}}^{\rm loc}-\mu|<c_0s;
    \label{eqn:right-dangerous-local-direction}
\end{align}
otherwise set $\mathcal P_{x,\mathrm{R}}(\mu)=0$.
The spectral description \eqref{eqn:right-local-radial-spectrum},
the lower bound on $\mu$, the choice of $c_0$, and
$s\geq\gC_{B,\nu}\error_d^{1/4}$ give
\begin{align}
    \operatorname{dist}\left(
      \mu,\operatorname{Spec}\bigl(
        \rmT_x|_{\mathcal P_{x,\mathrm{R}}(\mu)^\perp}
      \bigr)
    \right)
    \geq c_{\qr}s.
    \label{eqn:right-reduced-local-gap}
\end{align}
When $\mathcal P_{x,\mathrm{R}}(\mu)\neq0$, the same comparison with
\eqref{eqn:right-local-radial-spectrum} and the negative partner gives
\begin{align}
    \operatorname{dist}\left(
      \lambda_{x,\mathrm{R}}^{\rm loc},
      \operatorname{Spec}\bigl(
        \rmT_x|_{\mathcal P_{x,\mathrm{R}}(\mu)^\perp}
      \bigr)
    \right)
    \geq c_{\qr}s.
    \label{eqn:right-removed-local-separation}
\end{align}
If $\mathcal P_{x,\mathrm{R}}(\mu)\neq0$,
\eqref{eqn:right-root-exception-is-active} establishes activity, so the
active profile residual applies. The positive exception is simple
by interlacing, and its uniqueness in this component gives
$|\lambda_{x,\mathrm{R}}^{\rm loc}-\Lambda_x|\leq
C_\qr\varepsilon_{x,\mathrm{R}}$.
The scalar projection estimate and
\eqref{eqn:right-removed-local-separation} therefore give a normalized
local eigenvector $\phi_{x,\mathrm{R}}$, with a suitable sign, such that
\begin{align}
    \|\phi_{x,\mathrm{R}}-\hatvtau_+(x)\|
    \leq
    C_{\qr}s^{-1}\varepsilon_{x,\mathrm{R}}.
    \label{eqn:right-exact-local-profile-comparison}
\end{align}

\paragraph{Left-edge spectral projection.}
Take $x\in\Vtwolow{\tau_2^{\mathrm{L}}(\gap)}$ and $\mathrm E=\mathrm L$.
For the squared operator \eqref{eqn:left-local-squared-root-operator},
let $\mathcal P_{x,\mathrm{L}}(\mu)$ be the spectral projection
associated with the possible low root eigenvalue
$\kappa_{x,\mathrm L}^{\rm loc}$ below the regular interval
\eqref{eqn:left-local-squared-spectrum} when
$|\kappa_{x,\mathrm{L}}^{\rm loc}-\mu^2|$ is less than
$c_0s$, for a sufficiently small
$c_0=c_0(b,\qr,\sK)>0$, and set it equal to zero otherwise. The
regular spectral floor and the target interval in
\eqref{eqn:left-radial-mass} give the reduced gap
\begin{align}
    \operatorname{dist}\left(
      \mu^2,\operatorname{Spec}\bigl(
        \rmS_x|_{\mathcal P_{x,\mathrm{L}}(\mu)^\perp}
      \bigr)
    \right)
    \geq c_{b,\qr,\sK}s.
    \label{eqn:left-reduced-local-gap}
\end{align}
If $\mathcal P_{x,\mathrm{L}}(\mu)\neq0$, its eigenvalue is likewise at distance
at least $c_{b,\qr,\sK}s$ from the rest of
$\operatorname{Spec}(\rmS_x)$.
If $\mathcal P_{x,\mathrm{L}}(\mu)\neq0$,
\eqref{eqn:left-root-exception-is-active} establishes activity, so the
signed profile residuals apply. Interlacing makes the low root eigenvalue
simple. Parity balance normalizes the even component
\begin{align}
    \widehat u_{x,\mathrm{L}}
    \coloneqq
    \sqrt2\,\chi_{\gV_2}\hatvtau_+(x)
    =\frac{\hatvtau_+(x)+\hatvtau_-(x)}{\sqrt2}.
    \label{eqn:left-even-local-profile}
\end{align}
The signed adjacency residuals and the bounded critical local norm
$\|\rmT_x\|\leq C_{b,\qr}$ give a residual at most
$C_{b,\qr,\sK}\varepsilon_{x,\mathrm{L}}$ for the Gram operator
$\rmS_x$ at the spectral parameter $\Lambda_x^2$.
The separation in \eqref{eqn:left-root-exception-is-active}
identifies the low root eigenvalue within that error. The scalar
projection estimate then gives a unit vector
$\phi_{x,\mathrm{L}}$ spanning $\mathcal P_{x,\mathrm{L}}(\mu)$, with a suitable
sign, for which
\begin{align}
    \|\phi_{x,\mathrm{L}}-\widehat u_{x,\mathrm{L}}\|
    \leq
    C_{b,\qr,\sK}s^{-1}\varepsilon_{x,\mathrm{L}}.
    \label{eqn:left-exact-local-profile-comparison}
\end{align}
Since $\rvw$ is orthogonal to both signed profiles, its component on
the root bipartition, namely $\rvw^{(2)}$, is orthogonal to
$\widehat u_{x,\mathrm{L}}$.

\paragraph{Cutoff equations and source bounds.}
Let $\chi_x$ equal one on $\ball_{\radius+1}^{(\vtau)}(x)$, vanish on
layers $j\geq2\radius-2$, and interpolate linearly in the intervening
layers. Set
\[
    u_{x,\mathrm R}=\chi_x\rvw,
    \qquad
    u_{x,\mathrm L}=\chi_x\rvw^{(2)},
    \qquad
    W_x=\|\rvw|_{\ball_{2\radius}^{(\vtau)}(x)}\|.
\]
Because the cutoff vanishes before the Dirichlet boundary,
\eqref{eqn:radial-proof-defect} gives
\begin{align}
    (\rmT_x-\mu)u_{x,\mathrm R}
      &=f_{x,\mathrm R}+h_{x,\mathrm R},
      \label{eqn:right-cutoff-local-equation}\\
    (\rmS_x-\mu^2)u_{x,\mathrm L}
      &=f_{x,\mathrm L}+h_{x,\mathrm L},
      \label{eqn:left-cutoff-local-equation}
\end{align}
where
\begin{align}
    f_{x,\mathrm R}&=\chi_x\rvr_{\mathrm R},
    &h_{x,\mathrm R}&=[\rmT_x,\chi_x]\rvw,\\
    f_{x,\mathrm L}&=\chi_x
      (\mu\rvr_{\mathrm L}^{(2)}+\rmT_{12}^*\rvr_{\mathrm L}^{(1)}),
    &h_{x,\mathrm L}&=[\rmS_x,\chi_x]\rvw^{(2)}.
\end{align}
Both commutators are supported in the outer annulus, at distance at
least $\radius$ from the root. The right commutator involves only edges
with regular endpoints, so the finite-tree bound controls its norm
independently of the root degree. On the left, the upper forest form and
the critical maximum-degree estimate bound the local operator norms.
Together with
\eqref{eqn:local-defect-bound}, this gives
\begin{align}
    \|f_{x,\mathrm E}\|&\leq C\indi{x\in\mathcal C}\varepsilon_{x,\mathrm E}W_x,
    &\|h_{x,\mathrm E}\|&\leq CW_x.
    \label{eqn:local-cutoff-source-bounds}
\end{align}
The cutoff is one on each profile support. Orthogonality and the
profile comparisons
\eqref{eqn:right-exact-local-profile-comparison} and
\eqref{eqn:left-exact-local-profile-comparison} therefore yield
\begin{align}
    \|\mathcal P_{x,\mathrm E}(\mu)u_{x,\mathrm E}\|
    \leq Cs^{-1}\indi{x\in\mathcal C}\varepsilon_{x,\mathrm E}W_x.
    \label{eqn:local-dangerous-direction-overlap}
\end{align}
The spectral descriptions and bounded layer blocks verify
Lemma~\ref{lem:reduced-local-resolvent}: $\rmT_x$ has propagation
one and at most two exceptional values; $\rmS_x$ has propagation two
and at most one. Each nonroot right layer block is a child-incidence
matrix with norm $\sqrt{\max_y N_{j+1}(y)/d}\leq\sqrt{\beta_{\operatorname{side}(j)}}$.
Only its root block can be unbounded. At the critical left edge all
root and nonroot degrees are $O_{b,\qr}(d)$, and
$\|\rmT_x\|+\|\rmS_x\|\leq C_{b,\qr}$, so all squared layer
blocks are bounded. Apply
Lemma~\ref{lem:semilocalized-local-root-estimate} with $g=cs\leq1$
for a sufficiently small edge-dependent constant $c$, and $R=\radius$.
This gives the sharp estimate
\begin{align}
    |\ervw_x|^2
    \leq C\left[
      s^{-2}\indi{x\in\mathcal C}\varepsilon_{x,\mathrm E}^2
      +s^{-4}e^{-cs\radius}
    \right]W_x^2.
    \label{eqn:local-complement-root-mass}
\end{align}
Here $x\in\gV^{(\vtau)}$ on the right and
$x\in\Vtwolow{\tau_2^{\mathrm L}(\gap)}$ on the left, with the
corresponding edge dependence of the constants. For
$s\geq C\error_d^{1/4}$, its residual coefficient can also be
bounded by $C\error_d^{-1}$:
\begin{align}
    |\ervw_x|^2
    \leq C\left[
      \error_d^{-1}\indi{x\in\mathcal C}\varepsilon_{x,\mathrm E}^2
      +s^{-4}e^{-cs\radius}
    \right]W_x^2.
    \label{eqn:local-complement-root-mass-weak}
\end{align}

\subsubsection{Conditional weighted mass and radial-mass probabilities}

\begin{lemma}[Conditional weighted exceptional mass]
\label{lem:semilocalized-weighted-mass}
Let $0<s\leq1$, let $\mathcal D$ be exceptional centres with disjoint
$2\radius$-balls, and select $\mathcal C\subset\mathcal D$.
At each $x\in\mathcal D$, suppose the right local operator or the
left squared local operator satisfies
Lemma~\ref{lem:reduced-local-resolvent} with $g=c_3s\leq1$.
Suppose its cutoff vector has root coordinate $\ervw_x$ and satisfies
\eqref{eqn:local-root-source-equation}--\eqref{eqn:local-root-input-bounds}
with $R=\radius$, uniformly bounded $K$, and
$W_x=\|\rvw|_{\ball_{2\radius}^{(\vtau)}(x)}\|$.
For $x\in\mathcal C$, use $\varepsilon=\varepsilon_{x,\mathcal C}$
in the local source bound. For $x\notin\mathcal C$, require the
dangerous projection and interior source to vanish, and use
$\varepsilon=0$ in that bound. Then
\begin{align}
 \sum_{x\in\mathcal D}(1+\alpha_x)|\ervw_x|^2
 \leq C\left[
   s^{-2}\max_{x\in\mathcal C}(1+\alpha_x)\varepsilon_{x,\mathcal C}^2
   +(1+\max_{x\in\mathcal D}\alpha_x)s^{-4}e^{-cs\radius}
 \right]\|\rvw\|^2.
 \label{eqn:conditional-weighted-exceptional-mass}
\end{align}
Empty maxima and sums are zero. The constants depend only on the
local input constants, including $c_3$, and are independent of an
unbounded right root coupling.
\end{lemma}
\begin{proof}
Lemma~\ref{lem:semilocalized-local-root-estimate} gives
$s^{-2}\varepsilon_{x,\mathcal C}^2$ on selected centres and zero
residual contribution elsewhere. Multiply the sharp estimate by
$1+\alpha_x$ and sum. Disjointness gives
$\sum_{x\in\mathcal D}W_x^2\leq\|\rvw\|^2$, proving
\eqref{eqn:conditional-weighted-exceptional-mass}.
\end{proof}

\begin{proof}[Proof of Lemma~\ref{lem:semilocalized-radial-mass}]
The spectral projections, profile comparisons and cutoff bounds above
verify Lemma~\ref{lem:semilocalized-weighted-mass} with
$\mathcal D=\gV^{(\vtau)}$ on the right and
$\mathcal D=\Vtwolow{\tau_2^{\mathrm L}(\gap)}$ on the left.
Apply \eqref{eqn:pruning-weighted-profile-acceptance} to the window
centre sets verified in the proof of
Lemma~\ref{lem:semilocalized-profile-properties}. This gives
\begin{subequations}
\begin{align}
 \error_d^{-1}\max_{x\in\gV_{\mathrm R}^\star}
   (1+\alpha_x)\varepsilon_{x,\mathrm R}^2
 &\leq C_{B,\qr,\nu}(\error_d+\error_{\mathrm{prof}}(\gap)),
 \label{eqn:right-weighted-profile-errors-negligible}\\
 \error_d^{-1}\max_{x\in\gV_{\mathrm L}^\star}
   (1+\alpha_x)\varepsilon_{x,\mathrm L}^2
 &\leq C_{b,\qr,\sK}(\error_d+\error_{\mathrm{prof}}(\gap)).
 \label{eqn:left-weighted-profile-errors-negligible}
\end{align}
\end{subequations}
Since $s\geq C\error_d^{1/4}$, $s^{-2}\leq C\error_d^{-1}$.
The shared scale budget \eqref{eqn:pruning-generic-comparison-budgets}
controls the second term:
\begin{align}
 (1+\max_x\alpha_x)s^{-4}e^{-c_\qr s\radius}
 \leq C(\error_d+\error_{\mathrm{prof}}(\gap)).
 \label{eqn:weighted-boundary-errors-negligible}
\end{align}
The constants are $C_{B,\qr,\nu}$ on the right and
$C_{b,\qr,\sK}$ on the left. Equation~\eqref{eqn:conditional-weighted-exceptional-mass}
now proves both radial-mass bounds. The deterministic-cutoff input
event makes the argument simultaneous in $\gap$, with arbitrary
prescribed right failure exponent and the stated fixed left exponent.
\end{proof}

\subsection{Complementary-spectrum comparison and scale verification}
\label{subsec:semilocalized-comparison-proofs}

\subsubsection{Conditional comparison and its three domains}
\label{subsec:semilocalized-shared-comparison}

Write $s=\gap\wedge1$ and
$\eta_\gap=\error_d+\error_{\mathrm{prof}}(\gap)$.
Use the profiles, projections and block operator in
\eqref{eqn:profile-projections}--\eqref{eqn:right-hatH}, with the
centre set specified in each clause below.

\begin{proposition}[Conditional comparison]
\label{prop:semilocalized-shared-comparison}
Suppose the pruning and retained-layer conclusions hold, the local
balls are disjoint trees, the profiles have disjoint residual supports,
and the applicable global regularized forms and column-loss bound hold.
Suppose, for fixed input constants,
\begin{align}
 \|\scA-\scA^{(\vtau)}\|&\leq C\eta_\gap,
 &\max_{x\in\mathcal C}\varepsilon_x&\leq C\error_d,
 \label{eqn:shared-comparison-block-budget}\\
 s^{-2}\max_{x\in\mathcal C}(1+\alpha_x)\varepsilon_x^2
 &\leq C\eta_\gap,
 &(1+\max_x\alpha_x)s^{-4}e^{-cs\radius}
 &\leq C\eta_\gap.
 \label{eqn:shared-comparison-mass-budget}
\end{align}
Here $\varepsilon_x=\varepsilon_{x,\mathcal C}$ on $\mathcal C$;
unselected balls contribute no interior compression defect.
Assume that the root of every exceptional local eigenvalue within
$c_0s$ of the target belongs to $\mathcal C$, that the corresponding
ideal spectral value lies at distance at least $c_1s$ from the regular
local spectrum, and that the local profile comparison and reduced-resolvent estimates proved above hold.
For each cutoff equation require a reduced gap at least $c_3s$ and
$\|f_x\|\leq C\varepsilon_x W_x$ for $x\in\mathcal C$,
$f_x=0$ otherwise, and $\|h_x\|\leq CW_x$,
where $W_x$ is the norm on its $2\radius$-ball, with boundary source
$h_x$ in the outer annulus as in
\eqref{eqn:local-cutoff-source-bounds}. On the right the annular
endpoints must be regular, so their layer blocks are uniformly
bounded; on the left require $\|\rmT_x\|\leq C_{b,\qr}$ when
applying the second factor. These assumptions verify the local inputs
\eqref{eqn:local-root-input-bounds}.
The errors in root recovery, including the removed root degree,
must be smaller than these stated margins. In the left clauses also
require $\tau_1<1<\tau_2^{\mathrm L}$, the two-factor tail inputs,
and $\delta_{\rm root}+|\bar\alpha_x-\alpha_x|\leq c_2s$,
where $c_2$ is sufficiently small for the clause's squared buffer.
Require $\eta_\gap/s$ and $\error_d/s$ to be sufficiently small
for the form and root-identification margins. In every clause,
\begin{align}
 \|\scA-\hatH_{\mathcal C}\|&\leq C\eta_\gap,
 &\|\Pi_{\mathcal C}^{(\vtau)}\scA^{(\vtau)}
        Q_{\mathcal C}^{(\vtau)}\|&\leq C\error_d.
 \label{eqn:shared-comparison-block-conclusion}
\end{align}
The complementary spectra satisfy:
\begin{enumerate}
 \item \emph{Right comparison.} Select
 $\mathcal C=\gV_{\mathrm R}^\star(\gap)$, or every right profile
 at the location scale, and let $\kappa>1/4+c_0$.
 With right targets $\mu\geq\redge+\kappa s+C\eta_\gap$,
 \begin{align}
  \|Q_{\mathcal C}^{(\vtau)}\scA^{(\vtau)}
           Q_{\mathcal C}^{(\vtau)}\|
       &\leq\redge+\kappa s+C\eta_\gap.
  \label{eqn:shared-right-comparison}
 \end{align}
 \item \emph{Left comparison on a positive window.} Fix
 $\sK\Subset(0,\ledge)$ of positive length and select
 $\mathcal C=\gV_{\mathrm L}^\star(\gap)$, with its positive
 lower profile cutoff. Use $\kappa=1/4$, $0<\gap<\ledge$.
 Then
 \begin{align}
  \operatorname{Spec}_+(Q_{\mathcal C}^{(\vtau)}\scA^{(\vtau)}
          Q_{\mathcal C}^{(\vtau)})
  \cap[3\inf\sK/4,\ledge-s/4-C\eta_\gap)=\varnothing.
  \label{eqn:shared-left-window-comparison}
 \end{align}
 Zero and smaller positive complementary eigenvalues are permitted.
 The squared buffer to be checked is
 $(\kappa-1/8)s[2\ledge-(\kappa+1/8)s]>0$.
 \item \emph{Global left comparison with a degree floor.} Select
 $\mathcal C=\gL_2$, retain all $m$ columns, and assume every
 normalized retained column degree is at least $a_0/2>0$.
 Use $\kappa=1$, $s=\gap=3\error_d^{1/4}/4<\ledge/2$.
 The tail comparison must hold on the interval through zero in
 \eqref{eqn:left-global-tail-domain}, with the positive pivot in
 \eqref{eqn:left-global-positive-schur-pivot} and the parameter
 recovery above $z_0$. For $Q=Q_1\oplus Q_2$ and
 $X_Q=Q_1\rmX^{(\vtau)}Q_2$, one has
 \begin{align}
  \|X_Qv\|&\geq(\ledge-s)\|v\|,
       \qquad v\in\operatorname{Ran}Q_2.
  \label{eqn:shared-left-global-comparison}
 \end{align}
 This clause includes the test at singular value zero.
\end{enumerate}
The constants have the dependence of the corresponding input clause:
$(B,\qr,\nu)$ on the right, $(b,\qr,\sK)$ for a left window,
and $(b,\qr,a_0)$ for the global left comparison.
\end{proposition}

\begin{proof}
Proposition~\ref{prop:semilocalized-block-approximation} proves
\eqref{eqn:shared-comparison-block-conclusion}.
For a target complementary eigenvector, root selection and the local
profile comparison bound the dangerous projection by
$C\varepsilon_xW_x/s$ on selected centres; it vanishes elsewhere. The assumed cutoff source bounds verify
Lemma~\ref{lem:semilocalized-weighted-mass}; its two budgets in
\eqref{eqn:shared-comparison-mass-budget} give weighted exceptional
mass at most $C\eta_\gap$. The regularized form tests yield
\[
 \mu\leq\redge+s/8+C\eta_\gap\quad\hbox{on the right},
 \qquad
 \mu^2\geq(\ledge-s/8)^2-C\eta_\gap\quad\hbox{on the left}.
\]
The right inequality excludes a target above
$\redge+\kappa s+C\eta_\gap$ after increasing the error constant;
bipartite symmetry gives the norm bound.
The positive-window squared buffer excludes its target interval.

For the global clause, take a unit eigenvector
$v\in\operatorname{Ran}Q_2$ of
$(X_Q^*X_Q)|_{\operatorname{Ran}Q_2}$ with eigenvalue
$\mu^2<(\ledge-s)^2$. Put $u=X_Qv/\mu$ for
$\mu>0$ and $u=0$ for $\mu=0$.
Then $w=(u,v)\in\operatorname{Ran}Q$, $\|u\|\leq1$,
$\|w\|\leq\sqrt2$, and the compression defect gives
\begin{align}
 \scA^{(\vtau)}w=\mu w+e,
 \qquad \|e\|\leq C_{b,\qr,a_0}\error_d.
 \label{eqn:appendix-left-column-eigenvector-defect}
\end{align}
At $\mu=0$ this follows from $X_Qv=0$ and requires no division.
The local defect bound \eqref{eqn:local-defect-bound}, homogeneous
in the vector, applies to $w$. The two adjacency factors give the
left squared cutoff source \eqref{eqn:left-cutoff-local-equation},
with bounded critical local norms, also at $\mu=0$. Orthogonality to the signed profiles uses the normalized even
profile \eqref{eqn:left-even-local-profile}.
The pivot excludes low root eigenvalues near zero, and recovery
above $z_0$ puts every other dangerous root in $\gL_2$.
Lemma~\ref{lem:semilocalized-weighted-mass}, applied to $w$, bounds
the weighted column mass by $C\eta_\gap$.
Lemma~\ref{lem:semilocalized-complement-form-tests}(ii), including
its test at $\mu=0$, gives the squared lower bound above.
The buffer \eqref{eqn:left-global-squared-buffer-margin}
exceeds $C\eta_\gap$ under the input budget, excluding
$\mu^2<(\ledge-s)^2$. This controls every complementary column
singular value, with all $m$ original column positions retained.
\end{proof}

\subsubsection{Window comparison}
\label{subsec:semilocalized-complement-proof}
\begin{proof}[Proof of Proposition~\ref{prop:semilocalized-complement-spectrum}]
Work on the common deterministic-cutoff event of
Subsection~\ref{subsec:semilocalized-pruning-proofs}, intersected with
the profile and regularization events and, on the left, the
giant-comparison event. Choose the right input exponents above $\nu$
and retain a common fixed left exponent, decreasing it if necessary.
Write $s=s(\gap)$ and $\eta_\gap=\error_d+\error_{\mathrm{prof}}(\gap)$.
Lemma~\ref{lem:semilocalized-pruning-input}, the profile estimates and
\eqref{eqn:pruning-generic-comparison-budgets} give the block,
boundary and column-loss budgets. The weighted residual estimate
\eqref{eqn:pruning-weighted-profile-acceptance} gives the sharp mass
budget since $s^{-2}\leq C\error_d^{-1}$.
The shared error-to-gap ratio is uniformly small on the right after
choosing $\gC_0$ and the window threshold, and tends to zero at
critical sparsity on the left. For $\gap\geq1$, use $s=1$ and the
fixed pruning construction.

On the right, Lemma~\ref{lem:uniform-right-envelope} with
$\kappa=1/2$ and $c_0=1/16$ gives the $5s/16$ profile buffer and
$3s/16-C\error_d$ separation in
\eqref{eqn:right-root-exception-is-active}, which also verifies activity.
On the left, the positive-window tail inputs and
\eqref{eqn:left-window-root-error-budget} verify recovery and
\eqref{eqn:left-root-exception-is-active}. Choose the window
threshold large enough to retain the upper profile cutoff
$\ledge-2\error_d^{1/4}$.
The cutoff source bounds \eqref{eqn:local-cutoff-source-bounds}
and the bounded layer blocks verify the remaining local assumptions.
Clauses (i) and (ii) of
Proposition~\ref{prop:semilocalized-shared-comparison} give the window
spectra and the common block bound, including the
$C\error_d$ compression defect.
\end{proof}

\subsubsection{Location-scale pruning and profile inputs}
\label{subsec:semilocalized-shrinking-gap-proofs}

\begin{lemma}[Pruning and profile inputs at the location scale]
\label{lem:shrinking-gap-profile-inputs}
Use the pruning scale and error in \eqref{eqn:calibrated-comparison-scales} with
$\mathcal C_{\mathrm L}=\gL_2$. Maxima over an empty set are
understood as zero. The location-scale identities and budgets are
\eqref{eqn:location-pruning-scale-identities}--\eqref{eqn:location-exponential-column-budgets}.
\begin{enumerate}
    \item \emph{Right edge.}
    For every $\nu>0$, choose the constant $\gC_0$ in
    \eqref{eqn:subcritical_regime} sufficiently large depending on
    $(B,\qr,\nu)$. With probability at least
    $1-C_{B,\qr,\nu}N^{-\nu}$, the pruning construction is valid at
    gap $\frac34\error_d^{1/4}$, its active right centres are precisely
    $\mathcal C_{\mathrm R}$, and the radial profiles are orthonormal.
    For $\mathcal C=\mathcal C_{\mathrm R}$, they satisfy
    \begin{align}
        \|\scA-\scA^{(\vtau)}\|
        &\leq C_{B,\qr,\nu}\left(\error_d+\frac{4\error_d^{3/4}}{3\log d}\right),
        \label{eqn:shrinking-gap-pruning-norm}\\
        \max_{x\in\mathcal C,\,\sigma=\pm}
        \|\bigl(\scA^{(\vtau)}-\sigma\Lambda_x\bigr)
          \hatvtau_\sigma(x)\|
        &\leq C_{B,\qr,\nu}\error_d,
        \label{eqn:shrinking-gap-profile-residual}\\
        \frac1{\error_d}\max_{x\in\mathcal C}(1+\alpha_x)
        \max_{\sigma=\pm}
        \|\bigl(\scA^{(\vtau)}-\sigma\Lambda_x\bigr)
          \hatvtau_\sigma(x)\|^2
        &\leq C_{B,\qr,\nu}\left(\error_d+\frac{4\error_d^{3/4}}{3\log d}\right).
        \label{eqn:shrinking-gap-weighted-profile-error}
    \end{align}
    Increasing $\gC_0$ makes $\left(\error_d+\frac{4\error_d^{3/4}}{3\log d}\right)/\left(\frac34\error_d^{1/4}\right)$ uniformly as small as needed.

    \item \emph{Critical left edge with a degree floor.}
    Suppose that $d=b\log N$, $\qr>1$, and fix $a_0>0$.
    There are constants $\nu_0,C>0$, depending on $(b,\qr,a_0)$,
    and an event of probability at least $1-CN^{-\nu_0}$ on which
    the following implication holds: if
    $\min_{x\in\gV_2}\alpha_x\geq a_0$, the radial
    construction at $\frac34\error_d^{1/4}$ is valid at every $x\in\mathcal C_{\mathrm L}$.
    The pruning norm satisfies
    \eqref{eqn:shrinking-gap-pruning-norm}, and the profiles are
    orthonormal and satisfy \eqref{eqn:shrinking-gap-profile-residual}
    and \eqref{eqn:shrinking-gap-weighted-profile-error} with
    $\mathcal C=\mathcal C_{\mathrm L}$ and constants depending on
    $(b,\qr,a_0)$. These estimates include restriction to the
    giant component and extension by zero. Every retained root degree
    satisfies
    $d^{-1}|\layer_1^{(\vtau)}(x)|\geq a_0/2$ for
    $x\in\Vtwolow{\tau_2^{\mathrm{L}}(\frac34\error_d^{1/4})}$. Under critical sparsity, $\error_d+\frac{4\error_d^{3/4}}{3\log d}=o\!\left(\frac34\error_d^{1/4}\right)$.
\end{enumerate}
\end{lemma}

\begin{proof}[Proof of Lemma~\ref{lem:shrinking-gap-profile-inputs}]
Write $\xi=\error_d$, $s=3\xi^{1/4}/4$ and
$\eta_{\rm loc}=\xi+4\xi^{3/4}/(3\log d)$.
Choose $\gC_0$ so that $s\leq\min\{1,4/\qr\}$ and $\xi\leq1$.
Equations \eqref{eqn:location-pruning-scale-identities}--
\eqref{eqn:location-exponential-column-budgets} verify admissibility
and the error-to-gap, boundary, and column-loss budgets at this gap.

For item~(i), apply Lemmas~\ref{lem:existence_pruned_graph},
\ref{lem:semilocalized-layer-concentration} and
\ref{lem:semilocalized-pruning-input} with $a=3/4$.
They give the pruning construction, retained-layer event and
\eqref{eqn:shrinking-gap-pruning-norm}, with arbitrary prescribed
polynomial failure. The admissibility just checked gives the
nonbacktracking and regularized-form inputs at this scale.
The pruning cutoffs retain every centre in $\gR_1\cup\gR_2$.
Their degrees satisfy $\rD_x\geq c_\qr d$, and their profile values
satisfy $\Lambda_x\geq\redge+\xi^{1/4}/8$.
Equation~\eqref{eqn:profile-energy-degree-separation} and its first-side
analogue give $1-|q_x|\geq c_\qr\xi^{1/8}$.
The resulting terminal correction is a polynomial in $\log N$ times
$e^{-c_\qr\xi^{1/8}\radius}$, absorbed by the shared radius scale.
These degree and decay bounds verify
Lemma~\ref{lem:semilocalized-selected-profile-estimates} for
$\mathcal C=\mathcal C_{\mathrm R}$, giving orthonormality and both
residual estimates.
Choose larger input exponents before intersecting the events.

For item~(ii), use $d=b\log N$, $B=2b$, and intersect the preceding
inputs with the component-size and giant-comparison events.
The common event has probability at least
$1-C_{b,\qr,a_0}N^{-\nu_0}$ for some fixed $\nu_0>0$.
Conditional on $\min_{x\in\gV_2}\alpha_x\geq a_0$, the removed-degree
bound and \eqref{eqn:location-pruning-scale-identities} give
\begin{align}
 \bar\alpha_x=d^{-1}|\layer_1^{(\vtau)}(x)|
 \geq a_0-C_{b,\qr}\frac{\log N}{d^2s^2}
 \geq a_0/2,\qquad x\in\gV_2.
 \label{eqn:appendix-left-retained-root-degree}
\end{align}
This includes every low exceptional root.
Every $x\in\gL_2$ passes the low-degree cutoff because $\qr s/8<1$.
Its profile value is at most $\ledge-\xi^{1/4}/8$, so the same
coefficient-decay estimate applies. The retained degree floor
verifies the lower-degree input of
Lemma~\ref{lem:semilocalized-selected-profile-estimates}, which gives
both profile bounds at all $x\in\gL_2$, with constants depending on
$(b,\qr,a_0)$.

Corollary~\ref{cor:bipartite-giant-component-reduction} puts every
column in the giant under this floor and pads only isolated rows.
Equation~\eqref{eqn:giant-component-centered-comparison} adds
$O_{b,\qr}(\xi)$ to the norm and residual estimates. Since selected
left degrees are bounded, this also preserves the weighted budget.
Thus item~(ii) holds as an implication on the common input event;
\eqref{eqn:location-error-gap-ratio} gives $\eta_{\rm loc}=o(s)$ at
critical sparsity.
\end{proof}

\Needspace{4\baselineskip}
\subsubsection{Location comparison}
\label{subsec:semilocalized-calibrated-proofs}

\begin{proof}[Proof of Corollary~\ref{cor:semilocalized-calibrated-comparison}]
Write $s=3\error_d^{1/4}/4$ and
$\eta_{\rm loc}=\error_d+4\error_d^{3/4}/(3\log d)$.

\emph{Right edge.}
Fix $\nu>0$ and choose $\gC_0$ sufficiently large depending on
$(B,\qr,\nu)$ for the following margins.
Lemma~\ref{lem:shrinking-gap-profile-inputs}(i) selects all
$\gR_1\cup\gR_2$ and gives the pruning and profile inputs with
probability at least $1-C_{B,\qr,\nu}N^{-\nu}$.
Its weighted estimate and
\eqref{eqn:location-exponential-column-budgets} verify the sharp mass
budgets because $s^{-2}\leq C\error_d^{-1}$.
The error-to-gap budget \eqref{eqn:location-error-gap-ratio} is
uniformly small after increasing $\gC_0$.
For any local eigenvalue within $s/16$ of a target above
$\redge+s$, Lemma~\ref{lem:uniform-right-envelope} with $\kappa=1$
gives $\Lambda_{\mathrm{br}}(\alpha_x)\geq\redge+13s/16$.
Its root is in $\gR_1\cup\gR_2$, and its ideal value is separated
from the regular local interval by $11s/16-C\error_d$.
The cutoff source algebra \eqref{eqn:right-cutoff-local-equation}
and bounded nonroot layer blocks apply directly at this radius,
independently of the window lower-gap constant.
Clause (i) of Proposition~\ref{prop:semilocalized-shared-comparison}
gives the block and complementary bounds and a
$C_{B,\qr,\nu}\error_d$ compression defect.
Increase $\gC_0$ so that the combined block-and-complement error is
at most $\error_d^{1/4}/4$, as required in the right location theorem.

\emph{Left edge.}
Assume $b>\qr^{-2}$.
Corollary~\ref{cor:bipartite-giant-component-reduction} gives
$a_0,\nu_0>0$ and an event of probability at least
$1-C_{b,\qr}N^{-\nu_0}$ on which, decreasing $a_0$ to lie in
$(0,\qr^2-1)$,
\begin{align}
 \min_{y\in\gV_2}\alpha_y\geq a_0.
 \label{eqn:left-global-minimum-column-degree}
\end{align}
Intersect this event with
Lemma~\ref{lem:shrinking-gap-profile-inputs}(ii), decreasing the
common fixed exponent if necessary. All $m$ columns remain in the
giant and only isolated rows are padded. The input lemma gives the
pruning, residual and centering budgets and the retained column floor
\eqref{eqn:appendix-left-retained-root-degree}.

Under the retained degree floor,
Lemma~\ref{lem:global-left-root-recovery} provides the tail comparison
through zero, the positive pivot on $[0,z_0]$, and root recovery for
$z\geq z_0$ with error $o(s)$ in the squared singular value.
The squared buffer \eqref{eqn:left-global-squared-buffer-margin}
absorbs the $O(\error_d)$ form correction.
For a target $0\leq\mu\leq\ledge-s$, choose
$c_0=c_0(b,\qr,a_0)>0$ so small that
$|z-\mu^2|<c_0s$ implies
$|\sqrt z-\mu|\leq3c_0\error_d^{1/4}/(4\sqrt{z_0})
\leq\error_d^{1/4}/8$.
Hence
\begin{align}
 \Lambda_\qr(\alpha_x)
 &\leq\mu+\error_d^{1/4}/8+o(\error_d^{1/4})
 <\ledge-\error_d^{1/4}/8,
 \label{eqn:appendix-left-dangerous-root-selected}
\end{align}
so the root is selected in $\gL_2$.
The root-zero floor and profile cutoff
$\ledge-4\error_d^{1/4}/32$ give separation from the regular
spectrum of at least $\error_d^{1/4}/64$ for large $N$.

The location profile and scale budgets give $\eta_{\rm loc}=o(s)$.
The two-factor source equation has bounded critical local norm, and
the zero-singular-value construction is already included in the
common proof \eqref{eqn:appendix-left-column-eigenvector-defect}.
Regularization gives column loss $O(\eta_{\rm loc})$ and a fixed
positive retained-row inverse margin.
Clause (iii) of Proposition~\ref{prop:semilocalized-shared-comparison}
therefore proves both left comparisons, on an event of probability
at least $1-C_{b,\qr}N^{-\nu_0}$, retaining all $m$ column positions.
\end{proof}

%% file: app-transition.tex

\section{Deferred Proofs in
Section~\ref{sec:transition-from-delocalization-to-localization}}
\label{sec:transition-deferred-proofs}

\subsection{General degree estimates}
\label{subsec:transition-general-degree-proofs}

\begin{proof}[Proof of Proposition~\ref{prop:uniform-degree-counts}]
Within a vertex class, the degrees depend on disjoint sets of independent
edges. Hence $X_i(J)$ is binomial with mean $\mu_i(J)$ and variance at
most $\mu_i(J)$. Lemma~\ref{lem:Bernstein} gives, for every $u>0$,
\begin{align}
    \P\left(
      |X_i(J)-\mu_i(J)|>2\sqrt{\mu_i(J)u}+2u
    \right)\leq2\mathrm e^{-u}.
\end{align}
There are $O(N^2)$ distinct intervals on each integer degree lattice,
including singletons and one-sided intervals. Taking
$u=(\nu+4)\log N$ and a union bound proves
\eqref{eqn:uniform-degree-count-concentration}. The upper bound
\eqref{eqn:uniform-degree-count-upper} follows from
$\sqrt{\mu\log N}\leq(\mu+\log N)/2$.
Dividing the deviation estimate by $\mu_i(J)$ gives the asserted
relative concentration when the smallest mean is $\gg\log N$.

For a deterministic $J$, Markov's inequality gives
$\P(X_i(J)\geq1)\leq\mu_i(J)$. If
$q_i(J)=\P(\rD^{(i)}/d\in J)$, the same binomial law gives
\begin{align}
    \P(X_i(J)=0)
    =(1-q_i(J))^{|\gV_i|}
    \leq\exp[-\mu_i(J)],
\end{align}
which proves \eqref{eqn:degree-window-presence-absence}.

Finally, for a binomial random variable $Z$ with mean $\mu>0$, the
moment-generating function satisfies
\begin{align}
    \E\mathrm e^{\theta Z}
    \leq\exp[\mu(\mathrm e^\theta-1)],
    \qquad\theta\in\R.
\end{align}
For $k>0$, use exponential Markov with
$\theta=\log(k/\mu)$, positive for the upper tail and negative for
the lower tail. This gives
\eqref{eqn:degree-chernoff-upper} and
\eqref{eqn:degree-chernoff-lower}. At $k=\mu$ the bound is trivial;
at $k=0$ the lower-tail bound follows by letting $\theta\to-\infty$.
No independence between vertex classes is used.
\end{proof}

\begin{proof}[Proof of Lemma~\ref{lem:sparse-degree-window-means}]
For $1\leq k\leq Ad$, Lemma~\ref{lem:bipartite-poisson-approx}
gives, uniformly in $i$ and $k$,
\begin{align}
    \log\P(\rD^{(i)}=k)
    =-dI_{s_i}(k/d)+O_{\qr,A,B}(\log d).
    \label{eqn:bounded-degree-point-rate}
\end{align}
For $k=0$, Lemma~\ref{lem:bipartite-poisson-approx} gives
$\log\P(\rD^{(i)}=0)=-d_i+O_\qr(d^2/N)$, so the same estimate
holds because $I_{s_i}(0)=s_i^2$ and
$d_i=s_i^2d+O_\qr(d/N)$. In a nonempty $J_d$, the largest mass
gives a lower bound on the window probability; summing at most $Ad+1$
masses gives the matching upper bound. This proves
\eqref{eqn:sparse-degree-window-mean}.

For the tails, compare the Chernoff rate based on the exact mean in
\eqref{eqn:exact-degree-tail-rate} with $I_{s_i}$. Uniformly for
integers $0\leq k\leq N$,
\begin{align}
    \Psi_i(k)
    &=dI_{s_i}(k/d)
      +k\log(s_i^2d/d_i)+d_i-s_i^2d
      \notag\\
    &=dI_{s_i}(k/d)+O_\qr((k+d)/N).
    \label{eqn:exact-normalized-degree-rate}
\end{align}
For $a>s_i^2$, take $k=\lceil da\rceil$ in
\eqref{eqn:degree-chernoff-upper}. For $0<a<s_i^2$, take
$k=\lfloor da\rfloor$ in \eqref{eqn:degree-chernoff-lower}.
The derivative $I_{s_i}'(a)=\log(a/s_i^2)$ has the appropriate sign
on each branch, so rounding increases the rate. The error in
\eqref{eqn:exact-normalized-degree-rate} is bounded, and multiplying
the tail probabilities by $|\gV_i|$ proves
\eqref{eqn:sparse-degree-tail-mean-bound}. The comparison with $d_i/d$
is uniform on compact subsets separated from $s_i^2$.
\end{proof}

\begin{proof}[Proof of Lemma~\ref{lem:right-uniform-degree-tail-count}]
For every $a\geq s_i^2+1$, the upper-tail argument above applies to
$k=\lceil da\rceil$ whenever $k\leq N$, giving
\begin{align}
    \mu_i([a,\infty))\leq C_\qr N\mathrm e^{-dI_{s_i}(a)}
\end{align}
uniformly in $a$. If $k>N$, the count is zero. Apply
\eqref{eqn:uniform-degree-count-upper} to these one-sided intervals
to obtain \eqref{eqn:right-uniform-degree-tail-count} on the single
simultaneous event of Proposition~\ref{prop:uniform-degree-counts}.
\end{proof}

\subsection{Critical resonant-degree counts}
\label{subsec:transition-critical-count-proofs}

\begin{proof}[Proof of Corollary~\ref{cor:critical-resonant-upper-counts}
and Proposition~\ref{prop:fixed-exterior-resonant-counts}]
\leavevmode\par
\paragraph{Degree windows and rates.}
At the right edge, write
$\mathcal U_{1,\lambda,\delta}=\gW_{\lambda,\delta}^{(1)}$ and
$\mathcal U_{2,\lambda,\delta}=\gW_{\lambda,\delta}^{(2,\mathrm R)}$.
The inverse profile maps give the degree intervals
\begin{align}
    J_{i,\lambda}^{\mathrm R}
    \coloneqq
    \left[\invR{s_i}(\lambda-\delta_{\mathrm R,N}(\lambda)),
          \invR{s_i}(\lambda+\delta_{\mathrm R,N}(\lambda))\right],
    \qquad i\in\{1,2\}.
    \label{eqn:right-uniform-degree-window}
\end{align}
The half-gap condition places these intervals on the upper-degree
branches, in a fixed bounded degree range. Since $I_{s_i}$ is increasing
there, the mean lemma gives an upper estimate with boundary rate
$I_{s_i}(\inf J_{i,\lambda}^{\mathrm R})$ and logarithmic error
$O_{\qr,L}(\log d)$. This estimate also holds for empty lattice windows.
For $\redge<t\leq L+1$,
\begin{align}
    \left|\frac{d}{dt}(I_{s_i}\circ\invR{s_i})(t)\right|
    \leq C_{\qr,L}\bigl(1+(t-\redge)^{-1/2}\bigr).
\end{align}
Thus \eqref{eqn:right-uniform-resonance-window} and
$\gapR-\delta_{\mathrm R,N}(\lambda)\geq\gapR/2$ imply
\begin{align}
    I_{s_i}\bigl(\invR{s_i}(\lambda-\delta_{\mathrm R,N}(\lambda))\bigr)
    =I_{s_i}\bigl(\invR{s_i}(\lambda)\bigr)+o(1)
    \label{eqn:right-uniform-rate-comparison}
\end{align}
uniformly on the prescribed spectral interval.

At the left edge, the degree interval is
\begin{align}
    J_\lambda^{\mathrm L}
    \coloneqq
    \left[\invL(\lambda-\delta_{\mathrm L,N}(\lambda)),
          \invL(\lambda+\delta_{\mathrm L,N}(\lambda))\right].
    \label{eqn:left-uniform-degree-window}
\end{align}
Its endpoints lie in a fixed compact subset of $(0,\qr^2-1)$.
The rate $I_\qr$ is decreasing there, so the mean lemma gives an upper
estimate with boundary rate $I_\qr(\sup J_\lambda^{\mathrm L})$ and
logarithmic error $O_{\qr,K_0}(\log d)$, again including empty windows.
The bounded derivative of $\rateL(\cdot;\qr)$ near $K_0$ gives
\begin{align}
    \rateL(\lambda+\delta_{\mathrm L,N}(\lambda);\qr)
    =\rateL(\lambda;\qr)+o(1)
    \label{eqn:left-uniform-rate-comparison}
\end{align}
uniformly on $K_0$.

Define the untruncated branch exponents on their respective intervals by
\begin{align}
    r_i^{\mathrm R}(\lambda)
    &\coloneqq1-bI_{s_i}\bigl(\invR{s_i}(\lambda)\bigr),
    &r^{\mathrm L}(\lambda)&\coloneqq1-b\rateL(\lambda;\qr).
    \label{eqn:resonant-untruncated-branch-exponents}
\end{align}
Since $\log|\gV_i|=\log N+O_\qr(1)$ and $d=b\log N$, the mean estimates
become
\begin{align}
    \E|\mathcal U_{i,\lambda,\delta_{\mathrm R,N}(\lambda)}|
    &\leq N^{r_i^{\mathrm R}(\lambda)+o(1)},
    \qquad i\in\{1,2\},
    \label{eqn:right-resonant-window-mean}\\
    \E|\gW_{\lambda,\delta_{\mathrm L,N}(\lambda)}^{(2,\mathrm L)}|
    &\leq N^{r^{\mathrm L}(\lambda)+o(1)}.
    \label{eqn:left-resonant-window-mean}
\end{align}
The exponent errors are deterministic and uniform. The untruncated
exponents describe means that may tend to zero; by definition,
$r_i^{\mathrm R}(\lambda)\leq\rho_{\mathrm R}(\lambda,b)$ and
$r^{\mathrm L}(\lambda)\leq\rho_{\mathrm L}(\lambda,b)$.

\paragraph{Uniform upper counts.}
On the simultaneous event of Proposition~\ref{prop:uniform-degree-counts},
the right branch counts sum to at most
\begin{align}
    C_\nu\left(
      \sum_{i=1}^2N^{r_i^{\mathrm R}(\lambda)+o(1)}+\log N
    \right)
    \leq N^{\rho_{\mathrm R}(\lambda,b)+\varepsilon}
\end{align}
for all sufficiently large $N$, uniformly in $\lambda$. The same
estimate with one branch gives the left bound. The last inequalities
absorb the constants and $\log N$ because $\rho_{\mathrm E}\geq0$.
This proves Corollary~\ref{cor:critical-resonant-upper-counts} with every
prescribed failure exponent $\nu>0$.

\paragraph{Counts at a fixed spectral parameter.}
Fix $\mathrm E$ and $\lambda$ as in
Proposition~\ref{prop:fixed-exterior-resonant-counts}, and abbreviate
$\rho=\rho_{\mathrm E}(\lambda,b)>0$.
For each width $\delta_{\sigma,N}$, use the corresponding inverse
windows above. At this fixed $\lambda$, the derivatives of the inverse
maps are bounded above and below by positive constants. Hence each
window has length at least $c_\lambda\delta_{\sigma,N}$, and
\eqref{eqn:fixed-exterior-resonant-widths} ensures that its degree
lattice contains a point. Lemma~\ref{lem:sparse-degree-window-means}
now gives matching mean estimates. Writing
$X_{i,\sigma}^{\mathrm R}=|\mathcal U_{i,\lambda,\delta_{\sigma,N}}|$
and $X_\sigma^{\mathrm L}=|\gW_{\lambda,\delta_{\sigma,N}}^{(2,\mathrm L)}|$
in the respective cases, we obtain
\begin{align}
    \E X_{i,\sigma}^{\mathrm R}
    &=N^{r_i^{\mathrm R}(\lambda)+e_{i,\sigma,N}},
    \qquad i\in\{1,2\},\quad\sigma\in\{-1,+1\},
    \label{eqn:right-fixed-window-count-means}\\
    \E X_\sigma^{\mathrm L}
    &=N^{r^{\mathrm L}(\lambda)+e_{\sigma,N}},
    \qquad\sigma\in\{-1,+1\},
    \label{eqn:left-fixed-window-count-means}
\end{align}
with deterministic exponent errors tending to zero.
Let $a_N$ be their maximum absolute value. Choose a deterministic
$\eta_N\downarrow0$ such that $a_N=o(\eta_N)$,
$\eta_N\log N\to\infty$, and $\eta_N<\rho/2$ for large $N$.

The largest right branch exponent, or the single left branch exponent,
equals $\rho$. For either width, the corresponding mean is at least
$N^{\rho-a_N}$ and therefore $\gg\log N$. On the event of
Proposition~\ref{prop:uniform-degree-counts}, this branch count is
$(1+o(1))$ times its mean, supplying the lower bound
$N^{\rho-\eta_N}$. Every branch mean is at most $N^{\rho+a_N}$;
\eqref{eqn:uniform-degree-count-upper} bounds their sum by
$C_\nu(N^{\rho+a_N}+\log N)\leq N^{\rho+\eta_N}$ for large $N$.
The sequence $\eta_N$ is independent of $\nu$, since each fixed constant
$C_\nu$ is absorbed by $\eta_N\log N\to\infty$. This proves
\eqref{eqn:fixed-exterior-resonant-counts} for both widths on the same
simultaneous event.
\end{proof}

\subsection{Extremal degree input}
\label{subsec:transition-extremal-degree-proof}

\begin{proof}[Proof of Corollary~\ref{cor:critical-extreme-degrees}]
We use the root characterizations from
Subsection~\ref{subsec:transition-degree-roots}.

Fix a sufficiently small $\varepsilon>0$ so that
$a_{\mathrm R}(b;s_i)-\varepsilon>s_i^2$ for both classes.
The tail-mean bound \eqref{eqn:sparse-degree-tail-mean-bound} and
\eqref{eqn:degree-window-presence-absence} give
\begin{align}
    \P\left(
      \max_{x\in\gV_i}\alpha_x
      \geq a_{\mathrm R}(b;s_i)+\varepsilon
    \right)
    \leq
    N^{1-bI_{s_i}(a_{\mathrm R}(b;s_i)+\varepsilon)+o(1)}
    =o(1),
\end{align}
Conversely, choose an integer degree $k_i$ with
\begin{align}
    \frac{k_i}{d}
    =a_{\mathrm R}(b;s_i)-\varepsilon+o(1).
\end{align}
The singleton case of Lemma~\ref{lem:sparse-degree-window-means}
gives expected count
\begin{align}
    N^{1-bI_{s_i}(a_{\mathrm R}(b;s_i)-\varepsilon)+o(1)}
    \geq N^{c_\varepsilon}.
\end{align}
By \eqref{eqn:degree-window-presence-absence}, the probability that
degree $k_i$ does not occur is at most $\exp(-N^{c_\varepsilon})$.
These bounds for small $\varepsilon$ imply convergence for every
fixed deviation, proving
\eqref{eqn:extreme-degree-V1} and
\eqref{eqn:extreme-degree-V2}.

For the lower extreme, first suppose $b>\qr^{-2}$.  For every
sufficiently small fixed $\varepsilon>0$,
the same tail-mean and presence estimates, together with the
monotonicity of $I_\qr$ on $(0,\qr^2)$, give
\begin{align}
    \P\left(
      \min_{x\in\gV_2}\alpha_x
      \leq a_{\mathrm L}(b;\qr)-\varepsilon
    \right)
    =o(1).
\end{align}
An integer degree $\ell$ with
$\ell/d=a_{\mathrm L}(b;\qr)+\varepsilon+o(1)$ has expected count
at least $N^{c_\varepsilon}$ by the singleton mean estimate.
Equation~\eqref{eqn:degree-window-presence-absence} shows that it
occurs with probability at least $1-\exp(-N^{c_\varepsilon})$.
Hence the minimum converges in probability to $a_{\mathrm L}(b;\qr)$.

If $b<\qr^{-2}$, \eqref{eqn:bipartite-isolated-V2-count} in
Lemma~\ref{lem:bipartite-isolated-vertices} gives an isolated vertex in
$\gV_2$ with overwhelming probability. Thus the minimum normalized degree
equals zero, as does $a_{\mathrm L}(b;\qr)$, on that event.
At $b=\qr^{-2}$, for every fixed
$0<\varepsilon<\qr^2$, choose an integer $\ell$ with
$\ell/d=\varepsilon+o(1)$.  Since
\begin{align}
    1-bI_{\qr}(\varepsilon)>0,
\end{align}
the singleton mean estimate gives polynomially diverging expectation.
Equation~\eqref{eqn:degree-window-presence-absence} then shows that
this degree occurs with high probability. Hence
the minimum normalized degree is at most $\varepsilon+o(1)$ with high
probability.  Letting $\varepsilon\downarrow0$ proves convergence to
$a_{\mathrm L}(b;\qr)=0$ at the endpoint and completes
\eqref{eqn:extreme-degree-V2-lower}.
\end{proof}

\subsection{Spectral comparison}
\label{subsec:transition-spectral-comparison-proof}

\begin{proof}[Proof of Lemma~\ref{lem:spectral-counting}]
By Lemma~\ref{lem:weyl}, corresponding ordered eigenvalues of
$\rmS$ and $\rmT$ differ by at most $\varepsilon$. The inclusions
between the buffered intervals give
\eqref{eqn:transition-counting-sandwich-general}, with multiplicities
counted by their eigenvalue indices.

For the projection estimate, $P^\perp$ commutes with $\rmT$, so
$\rmS u=\lambda u$ implies
\begin{align}
    (\rmT-\lambda)P^\perp u=P^\perp(\rmT-\rmS)u.
\end{align}
The inverse of $\rmT-\lambda$ on $\operatorname{Ran}P^\perp$ has
norm at most $g^{-1}$, giving $\|P^\perp u\|\leq\varepsilon/g$.
Orthogonality yields the mass bound
\eqref{eqn:general-comparison-profile-mass}.
\end{proof}

\subsection{Profile-to-spectrum transfer}
\label{subsec:transition-spectral-transfer-proof}

\begin{proof}[Proof of Proposition~\ref{prop:transition-counting-transfer}]
At either fixed exterior spectral parameter, the edge gap $\gap>0$
is independent of $N$, so
\begin{align}
    \error_{\mathrm{prof}}(\gap)
    =\frac{\error_d}{\gap\log d}
    =o(\error_d).
\end{align}
Thus the gap-dependent errors in the block and complementary-spectrum
estimates are absorbed into the stated $C\error_d$ buffers for all
sufficiently large $N$.

For all sufficiently large $N$, the buffered widths satisfy the
window restrictions in
Subsubsection~\ref{subsec:semilocalized-eigenvector-semilocalization}.
The resonant-centre inclusion stated there therefore identifies all
centres counted below as active.

For the right edge, set $\gap=\gapR$ and construct
$\hatH_{\mathrm{R}}(\gap)$ as in \eqref{eqn:right-hatH}.  Since $\gap$ is fixed,
invoke the right-edge block and complementary-spectrum estimates with
target failure exponent $\nu+1$. Their intersection has probability at
least $1-\gC_{B,\qr,\nu}N^{-\nu}$, and
Proposition~\ref{prop:semilocalized-complement-spectrum} gives
\begin{align}
    \|\scA-\hatH_{\mathrm{R}}(\gap)\|
    \leq\varepsilon_{\mathrm{R},N}.
\end{align}
The complementary-spectrum bound \eqref{eqn:right-complement-spectrum}
puts the complementary block of $\hatH_{\mathrm{R}}(\gap)$ below
\begin{align}
    \redge+\frac{\gap\wedge1}{2}+o(1)
    <\lambda-c_{\lambda,\qr}.
    \label{eqn:transition-right-complement-gap}
\end{align}
Thus, for all sufficiently large $N$, the only eigenvalues of
$\hatH_{\mathrm{R}}(\gap)$ in $[\lambda-2\delta,\lambda+2\delta]$ are the positive
profile eigenvalues $\Lambda_x$.  Equations
\eqref{eqn:right-hatH} and \eqref{eqn:profile-energy} give
\begin{align}
    \mathrm{N}_{\lambda,t}(\hatH_{\mathrm{R}}(\gap))
    =|\gW_{\lambda,t}^{(1)}|
     +|\gW_{\lambda,t}^{(2,\mathrm{R})}|,
    \qquad
    t\in[
      \delta-\varepsilon_{\mathrm{R},N},
      \delta+\varepsilon_{\mathrm{R},N}
    ].
    \label{eqn:transition-profile-count-right}
\end{align}
Applying
\eqref{eqn:transition-counting-sandwich-general} with
$\rmS=\scA$ and $\rmT=\hatH_{\mathrm{R}}(\gap)$ proves
\eqref{eqn:transition-counting-sandwich-right}.

For the left edge, apply the constructions with
$\sK=K$, set $\gap=\gapL$, and construct
$\hatH_{\mathrm{L}}(\gap)$. The left block-approximation and complementary-spectrum
events hold jointly with probability at least
$1-\gC_{b,\qr,K}N^{-\nu_0}$ after decreasing $\nu_0$ if necessary. On
this event, the left block approximation gives
\begin{align}
    \|\scA-\hatH_{\mathrm{L}}(\gap)\|
    \leq\varepsilon_{\mathrm{L},N}.
\end{align}
The complementary-spectrum estimate
\eqref{eqn:left-complement-spectrum} separates the positive
complementary spectrum from $\lambda$ by a constant.  Zero, the
negative spectrum, and the remaining positive spectrum below
$3\inf K/4$ are also separated from $\lambda$.  Hence
\begin{align}
    \mathrm{N}_{\lambda,t}(\hatH_{\mathrm{L}}(\gap))
    =|\gW_{\lambda,t}^{(2,\mathrm{L})}|,
    \qquad
    t\in[
      \delta-\varepsilon_{\mathrm{L},N},
      \delta+\varepsilon_{\mathrm{L},N}
    ].
    \label{eqn:transition-profile-count-left}
\end{align}
Another application of
\eqref{eqn:transition-counting-sandwich-general} proves
\eqref{eqn:transition-counting-sandwich-left}.
\end{proof}


\section{Technical Lemmas}

\begin{lemma}[Bernstein's inequality, {\cite[Theorem $2.8.4$]{vershynin2018high}}]\label{lem:Bernstein}
    Let $\rX_1,\dots, \rX_n$ be independent mean zero random variables such that $|\rX_i|\leq K$ for all $i$. Let $\sigma^2 = \sum_{i=1}^{n}\E \rX_i^2$. Then for every $t \geq 0$,
    \begin{align}
        \P \Bigg( \Big|\sum_{i=1}^{n} \rX_i \Big| \geq t \Bigg) \leq 2 \exp \Bigg( - \frac{t^2/2}{\sigma^2 + Kt/3} \Bigg)\,.
    \end{align}
\end{lemma}

\begin{lemma}[Bennett's inequality, {\cite[Theorem $2.9.2$]{vershynin2018high} }]\label{lem:Bennett}
    Let $\rX_1,\dots, \rX_n$ be independent random variables. Assume that $|\rX_i - \E \rX_i| \leq K$ almost surely for every $i$. Then for any $t>0$, we have
    \begin{align}
        \P \Bigg( \sum_{i=1}^{n} (\rX_i - \E \rX_i) \geq t \Bigg) \leq \exp \Bigg( - \frac{\sigma^2}{K^2} \cdot \benrate \bigg( \frac{Kt}{\sigma^2} \bigg)\Bigg)\,, \notag
    \end{align}
    where $\sigma^2 = \sum_{i=1}^{n}\Var(\rX_i)$ is the variance of the
    sum and $\benrate(t) \coloneqq(1+t)\log(1+t)-t$.
    Furthermore, define $\rY_i=-\rX_i$ and apply the inequality above;
    for any $t>0$, we then have
    \begin{align}
        \P \Bigg( \sum_{i=1}^{n} (\rX_i - \E \rX_i) \leq -t \Bigg) \leq \exp \Bigg( - \frac{\sigma^2}{K^2} \cdot \benrate \bigg( \frac{Kt}{\sigma^2} \bigg)\Bigg)\,.
    \end{align}
\end{lemma}

\begin{lemma}[Weyl's inequality, \cite{Weyl1912DasAV}]
    \label{lem:weyl}
    Let $\rmA, \rmE \in \R^{m \times n}$. Then
    \begin{align}
        |\sigma_i(\rmA+\rmE)-\sigma_i(\rmA)|
        \leq \|\rmE\|,
        \qquad 1\leq i\leq\min\{m,n\}.
    \end{align}
    If $m=n$ and $\rmA,\rmE\in\R^{n\times n}$ are symmetric, then
    \begin{align}
        |\lambda_i(\rmA+\rmE)-\lambda_i(\rmA)|
        \leq\|\rmE\|,
        \qquad 1\leq i\leq n.
    \end{align}
\end{lemma}

\begin{lemma}[Schur test in matrix version, \cite{schur1911bemerkungen}]\label{lem:schur_test}
    Let $\{a_j\}_{j=1}^{n}$ and $\{b_l\}_{l=1}^{m}$ be two sequences of positive real numbers, and let $\lambda, \mu$ be positive real numbers. For matrix $\rmX \in \R^{n \times m}$, we have $\|\rmX\| \leq \sqrt{\lambda \mu}$ if
    \begin{align}
        \sum_{j=1}^{n} |\ermX_{jl}| \cdot a_{j} \leq \lambda \, b_{l}, \,\, \forall l\in [m], \quad \textnormal{ and }\quad \sum_{l=1}^{m} |\ermX_{jl}| \cdot b_{l} \leq \mu\, a_{j}, \,\, \forall j \in [n].
    \end{align}
    As a consequence, we take $a_{j} = b_{l} = 1$ and choose $\lambda$ (resp.\ $\mu$) as the maximum column (resp.\ row) sum, then
    \begin{align}
        \|\rmX\| \leq \bigg( \max_{j\in [n]} \sum_{l=1}^{m} |\ermX_{jl}| \bigg)^{1/2} \cdot \bigg( \max_{l\in [m]}\sum_{j=1}^{n} |\ermX_{jl}| \bigg)^{1/2}.
    \end{align}
\end{lemma}

\begin{lemma}[Schur's complement formula]\label{lem:shur-complement}
    The inverse of the block matrix $\mM$ is given by
    \begin{align}
        \mM^{-1} =&\, \begin{bmatrix}
            \mA &\, \mB \\
            \mC &\, \mD
        \end{bmatrix}^{-1}\\
        =&\, \begin{bmatrix}
            \mA^{-1} + \mA^{-1}\mB (\mM/\mA)^{-1}\mC \mA^{-1} &\, -\mA^{-1}\mB(\mM/\mA)^{-1} \\
            - (\mM/\mA)^{-1}\mC \mA^{-1} &\, (\mM/\mA)^{-1}
        \end{bmatrix}\quad (\textnormal{if $\mA$ is invertible})\\
        =&\, \begin{bmatrix}
            (\mM/\mD)^{-1} &\, - (\mM/\mD)^{-1}\mB \mD^{-1} \\
            -\mD^{-1}\mC(\mM/\mD)^{-1} &\, \mD^{-1} + \mD^{-1}\mC (\mM/\mD)^{-1}\mB \mD^{-1}
        \end{bmatrix}\quad (\textnormal{if $\mD$ is invertible})
    \end{align}
    where we denote
    \begin{align}
        \mM/\mA = \mD - \mC \mA^{-1}\mB, \quad \mM/\mD = \mA - \mB \mD^{-1}\mC.
    \end{align}
\end{lemma}
\begin{proof}[Proof of Lemma~\ref{lem:shur-complement}]
    The results can be verified via Gaussian elimination.
\end{proof}

For a rectangular matrix $M$ and an interval $I\subset(0,\infty)$, put
\[
 Q_I(M)=\mathbf1_{I^2}(M^{\sT}M),\qquad
 \widehat Q_I(M)=\mathbf1_{I^2}(MM^{\sT}),\qquad
 I^2=\{s^2:s\in I\}.
\]
For $\mX$ and the centered giant block, we suppress the matrix argument.
Giant-block projectors are extended by zero to the full row and column
coordinate spaces.

\begin{lemma}[Schur reduction and spectral projection]
    \label{lem:low-column-Schur-spectrum}
    \begin{enumerate}
    \item[(i)] \emph{Schur reduction.}
    Let $D=[C\ B]$, where $C$ has $\ell$ columns, and suppose that
    $B^{\sT}B\geq g^2d\,\id$ for some $g,d>0$. Put
    \[
     G=B^{\sT}B,\qquad Y=B^\dagger C,\qquad
     P=\operatorname{proj}_{\ker B^{\sT}},\qquad K=C^{\sT}PC.
    \]
    Suppose $s\coloneqq\|K\|<g^2d$, and write
    $0\leq\theta_1\leq\cdots\leq\theta_\ell$ for the eigenvalues of $K$.
    There are exactly $\ell$ eigenvalues
    $0\leq\lambda_1\leq\cdots\leq\lambda_\ell<g^2d$ of $D^{\sT}D$,
    and every remaining eigenvalue is at least $g^2d$. Moreover,
    \begin{align}
     \frac{\theta_i}{1+\tau}\leq\lambda_i\leq\theta_i,
     \qquad
     \tau=\frac{\|Y\|^2}{1-s/(g^2d)},\qquad 1\leq i\leq\ell.
     \label{eqn:core-Schur-eigenvalue-bounds}
    \end{align}
    In particular, $\ker D=\{(x,-Yx):x\in\ker K\}$.
    For $0\leq\lambda<g^2d$, the Schur complement of
    $G-\lambda\id$ in $D^{\sT}D-\lambda\id$ is
    \begin{align}
     K-\lambda\left[\id+
     Y^{\sT}(\id-\lambda G^{-1})^{-1}Y\right].
     \label{eqn:core-Schur-complement}
    \end{align}
    \item[(ii)] \emph{Projection of approximate eigenvectors.}
    Let $H$ be real symmetric, let $T$ have $t$ columns, and let $\Lambda$
    be a real symmetric $t\times t$ matrix. Suppose
    $\|HT-T\Lambda\|\leq r$. If $I=[a,b]$ and
    $\operatorname{Spec}\Lambda\subset[a+\eta,b-\eta]$ for $\eta>0$,
    then
    \[
     \|(\id-\mathbf1_I(H))T\|\leq\sqrt2\,r/\eta.
    \]
    If $H$ is the linearization of a rectangular matrix, $I\subset(0,\infty)$,
    and $T=2^{-1/2}(U,V)$ with $U^{\sT}U=V^{\sT}V=\id_t$, then
    \[
     \|(\id-\widehat Q_I)U\|+
     \|(\id-Q_I)V\|\leq4r/\eta.
    \]
    Each summand is at most $2r/\eta$. If $2r/\eta<1$, then
    $\operatorname{rank}Q_I\geq t$. Moreover, for orthonormal vectors
    $f_1,\ldots,f_k$ and every orthogonal projection $Q$,
    \begin{align}
     \operatorname{rank}Q\geq\sum_{i=1}^k\langle f_i,Qf_i\rangle
     =\sum_{i=1}^k\bigl(1-\|(\id-Q)f_i\|^2\bigr).
     \label{eqn:trace-rank-count}
    \end{align}
    \end{enumerate}
    \end{lemma}
    \begin{proof}[Proof of Lemma~\ref{lem:low-column-Schur-spectrum}]
    \emph{(i)} Since $C=PC+BY$, expansion of the Gram matrix and elimination of its
    $B$ block give \eqref{eqn:core-Schur-complement}. Congruence preserves
    inertia, so below $g^2d$ the negative eigenvalue count is the negative
    eigenvalue count of that Schur complement. On the subspace
    $\{(x,-Yx):x\in\R^\ell\}$ the Gram Rayleigh quotient is
    \[
     \frac{x^{\sT}Kx}{\|x\|^2+\|Yx\|^2}.
    \]
    Min--max therefore gives $\lambda_i\leq\theta_i\leq s$ for
    $i\leq\ell$. Principal-submatrix interlacing gives
    $\lambda_{\ell+1}(D^{\sT}D)\geq\lambda_1(G)\geq g^2d$ whenever
    this index exists. For $0\leq\lambda\leq s$,
    \[
     \id\leq\id+Y^{\sT}(\id-\lambda G^{-1})^{-1}Y
     \leq(1+\tau)\id.
    \]
    The inertia identity, applied at parameters on either side of each
    ordered eigenvalue, gives the lower bound in
    \eqref{eqn:core-Schur-eigenvalue-bounds}. The kernel identity also
    follows from $D(x,z)=PCx+B(Yx+z)$ and orthogonality of these two terms.

    \emph{(ii)} For the spectral subspace of $H$ below $a$, the equation
    $HT-T\Lambda=R$ is a separated Sylvester equation. Its solution is
    the integral of $e^{t(H-a)}R e^{-t(\Lambda-a)}$, up to sign, so the
    operator norm is at most $r/\eta$. The subspace above $b$ has the
    same bound. These two output subspaces are orthogonal, giving
    $\sqrt2r/\eta$. The projector onto the linearization eigenvalues
    in $I\cup(-I)$ is $\operatorname{diag}(\widehat Q_I,Q_I)$ and has
    no larger error. Taking its two blocks proves the stated estimates.
    Finally, $Q_IV$ is injective if $\|(\id-Q_I)V\|<1$, because $V$
    is an isometry. This argument has no factor depending on $t$.
    For \eqref{eqn:trace-rank-count}, complete $f_1,\ldots,f_k$ to an
    orthonormal basis. Since $Q$ is positive semidefinite,
    $\operatorname{rank}Q=\Tr Q$ is the sum of $\langle f,Qf\rangle\geq0$
    over this basis, which is at least the sum over $f_1,\ldots,f_k$. For a
    unit vector $f$, $\langle f,Qf\rangle=\|Qf\|^2=1-\|(\id-Q)f\|^2$.
    \end{proof}